\documentclass{agtart_a}
\pdfoutput=1

\usepackage{amscd}
\usepackage{pinlabel}
\usepackage[all]{xy}

\title{Bottom tangles and universal invariants}

\author{Kazuo Habiro}
\givenname{Kazuo}
\surname{Habiro}
\address{Research Institute for Mathematical Sciences\\
Kyoto University\\\newline
Kyoto 606--8502\\
Japan}
\email{habiro@kurims.kyoto-u.ac.jp}
\urladdr{}

\volumenumber{6}
\issuenumber{}
\publicationyear{2006}
\papernumber{41}
\startpage{1113}
\endpage{1214}

\doi{}
\MR{}
\Zbl{}

\keyword{knots}
\keyword{links}
\keyword{tangles}
\keyword{braided categories}
\keyword{ribbon Hopf algebras}
\keyword{braided Hopf algebras}
\keyword{universal link invariants}
\keyword{transmutation}
\keyword{local moves}
\keyword{Hennings invariants}
\keyword{bottom tangles}
\keyword{claspers}
\subject{primary}{msc2000}{57M27}
\subject{secondary}{msc2000}{57M25}
\subject{secondary}{msc2000}{18D10}

\received{20 December 2005}
\revised{}
\accepted{4 May 2006}
\published{7 September 2006}
\publishedonline{7 September 2006}
\proposed{}
\seconded{}
\corresponding{}
\editor{}
\version{}

\arxivreference{math.GT/0505219}

\makeatletter
\def\cnewtheorem#1[#2]#3{\newtheorem{#1}{#3}[section]
\expandafter\let\csname c@#1\endcsname\c@theorem}
\makeatother

\let\xysavmatrix\xymatrix
\def\xymatrix{\disablesubscriptcorrection\xysavmatrix}
\AtBeginDocument{\let\bar\wbar\let\tilde\wtilde\let\hat\what}

\def\nocolon{}

\numberwithin{equation}{section}

\theoremstyle{plain}
\newtheorem{Theorem}{Theorem}
\makeautorefname{Theorem}{Theorem}
\newtheorem{theorem}{Theorem}[section]
\cnewtheorem{lemma}[theorem]{Lemma}
\cnewtheorem{corollary}[theorem]{Corollary}
\cnewtheorem{proposition}[theorem]{Proposition}
\cnewtheorem{conjecture}[theorem]{Conjecture}
\theoremstyle{definition}
\cnewtheorem{definition}[theorem]{Definition}
\cnewtheorem{remark}[theorem]{Remark}
\cnewtheorem{example}[theorem]{Example}
\newtheorem*{acknowledgments}{Acknowledgments}

\newcommand\FIG[3]{\begin{figure}[t]
    \begin{center}
      \includegraphics[#3]{\figdir/#1}
    \end{center}
    \caption{#2}
    \label{fig:#1}
\end{figure}}

\newcommand\FIGn[3]{\begin{figure}[t]
    \begin{center}
      \includegraphics[#3]{\figdir/#1}
    \end{center}
    \nocolon\caption{#2}
    \label{fig:#1}
\end{figure}}

\newcommand\incl[2]{{\includegraphics[height=#1]{\figdir/#2}}}
\newcommand\evdn{\raisebox{-0.2em}{\incl{1em}{}}}
\newcommand\cvdn{\raisebox{-0.2em}{\incl{1em}{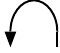}}}
\newcommand\evup{\raisebox{-0.2em}{\incl{1em}{}}}
\newcommand\cvup{\raisebox{-0.2em}{\incl{1em}{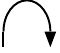}}}
\newcommand\col {\colon\thinspace}
\newcommand\ver {\ |\ }

\newcommand\simh{\sim_h}
\newcommand\trr{\triangleright}
\newcommand\ul{\underline}
\newcommand\one {{\mathbf 1}}
\newcommand\modA {{\mathcal A}}
\newcommand\sfA{{\mathsf A}}
\newcommand\modb {{\mathsf b}}
\newcommand\modB {{\mathsf B}}
\newcommand\B{{\mathcal B}}
\newcommand\cB{\check{\modB }}
\newcommand\modC {{\mathcal C}}
\newcommand\modH {{\langle \HH\rangle }}
\newcommand\modh {{\mathsf h}}
\newcommand\modI {{\mathcal I}}
\newcommand\modJ {{\mathsf J}}
\newcommand\modk {{\mathbf k}}
\newcommand\modL {{\mathcal L}}
\newcommand\LL{{\mathsf L}}
\newcommand\HH{{\mathsf H}}
\newcommand\modM {{\mathcal M}}
\newcommand\modr {{\mathbf r}}
\newcommand\modS {{\mathcal S}}
\newcommand\modT {{\mathcal T}}
\newcommand\T{{\mathsf T}}
\newcommand\modV {{\mathsf V}}
\newcommand\V{{\mathcal V}}
\newcommand\modY {{\mathcal Y}}
\newcommand\modZ {{\mathbb Z}}
\newcommand\g{{\mathfrak g}}

\def\bS{\mskip1mu\underline{\mskip-1mu{S}\mskip-2mu}\mskip2mu}
\def\bD{\mskip1.5mu\underline{\mskip-1.5mu{\Delta}\mskip-1.5mu}\mskip1.5mu}
\def\bH{\mskip1.5mu\underline{\mskip-1.5mu{H}\mskip-5mu}\mskip5mu}

\newcommand\simeqto{\stackrel{\simeq}{\to}}

\makeop{Ab}
\makeop{Ob}
\makeop{Mor}
\makeop{ad}
\makeop{coad}
\makeop{ev}
\makeop{coev}
\makeop{id}
\makeop{Span}
\makeop{End}
\makeop{source}
\makeop{target}
\newcommand\opint{\operatorname{int}}
\newcommand\clo{\operatorname{cl}}
\makeop{Hom}
\makeop{tr}
\newcommand\Mod{{\mathsf{Mod}}}
\newcommand\Sets{{\operatorname{\sf Sets}}}
\newcommand\BT{\operatorname{\sf BT}}
\newcommand\ABT{\operatorname{\sf ABT}}
\newcommand\SL{\operatorname{\sf SL}}

\newcommand\coadb{\underline{\coad}}
\newcommand\adb{\underline{\ad}}

\begin{document}

\begin{asciiabstract}
A bottom tangle is a tangle in a cube consisting only of arc components,
each of which has the two endpoints on the bottom line of the cube,
placed next to each other.  We introduce a subcategory B of the category
of framed, oriented tangles, which acts on the set of bottom tangles.
We give a finite set of generators of B, which provides an especially
convenient way to generate all the bottom tangles, and hence all the
framed, oriented links, via closure.  We also define a kind of "braided
Hopf algebra action" on the set of bottom tangles.

Using the universal invariant of bottom tangles associated to each ribbon
Hopf algebra H, we define a braided functor J from B to the category Mod_H
of left H-modules.  The functor J, together with the set of generators of
B, provides an algebraic method to study the range of quantum invariants
of links.  The braided Hopf algebra action on bottom tangles is mapped
by J to the standard braided Hopf algebra structure for H in Mod_H.

Several notions in knot theory, such as genus, unknotting number,
ribbon knots, boundary links, local moves, etc are given algebraic
interpretations in the setting involving the category B.  The functor
J provides a convenient way to study the relationships between these
notions and quantum invariants.
\end{asciiabstract}

\begin{abstract}
  A \emph{bottom tangle} is a tangle in a cube consisting only of arc
  components, each of which has the two endpoints on the bottom line
  of the cube, placed next to each other.  We introduce a subcategory
  $\mathsf{B}$ of the category of framed, oriented tangles, which acts
  on the set of bottom tangles.  We give a finite set of generators of
  $\mathsf{B}$, which provides an especially convenient way to generate
  all the bottom tangles, and hence all the framed, oriented links,
  via closure.  We also define a kind of ``braided Hopf algebra action''
  on the set of bottom tangles.

  Using the universal invariant of bottom tangles associated to each
  ribbon Hopf algebra $H$, we define a braided functor $\mathsf{J}$ from
  $\mathsf{B}$ to the category $\mathsf{Mod}_H$ of left $H$--modules.
  The functor $\mathsf{J}$, together with the set of generators
  of $\mathsf{B}$, provides an algebraic method to study the range of
  quantum invariants of links.  The braided Hopf algebra action on bottom
  tangles is mapped by~$\mathsf{J}$ to the standard braided Hopf algebra
  structure for $H$ in $\mathsf{Mod}_H$.

  Several notions in knot theory, such as genus, unknotting number,
  ribbon knots, boundary links, local moves, etc are given algebraic
  interpretations in the setting involving the category $\mathsf{B}$.
  The functor $\mathsf{J}$ provides a convenient way to study the
  relationships between these notions and quantum invariants.
\end{abstract}

\maketitle

\section{Introduction}
\label{sec1}

The notion of category of tangles (see Yetter \cite{Yetter:88} and Turaev
\cite{Turaev:89}) plays a
crucial role in the study of the quantum link invariants.  One can
define most quantum link invariants as braided functors from the
category of (possibly colored) framed, oriented tangles to other
braided categories defined algebraically.  An
important class of such functorial tangle invariants is introduced by
Reshetikhin and Turaev~\cite{Reshetikhin-Turaev:90}: Given a ribbon
Hopf algebra $H$ over a field $k$, there is a canonically defined
functor $F\col \modT _H\rightarrow \Mod^f_H$ of the category $\modT _H$ of framed, oriented
tangles colored by finite-dimensional representations of $H$ into the
category $\Mod^f_H$ of finite-dimensional left $H$--modules.  The
Jones polynomial \cite{Jones} and many other polynomial link
invariants (see Freyd--Yetter--Hoste--Lickorish--Millett--Ocneanu
\cite{FYHLMO}, Przytycki--Traczyk \cite{PT}, Brandt--Lickorish--Millett
\cite{BLM}, Ho \cite{Ho} and Kauffman \cite{Kauffman:poly}) can be understood in
this setting, where $H$ is a quantized enveloping algebra of simple
Lie algebra.

One of the fundamental problems in the study of quantum link
invariants is {\em to determine the range of a given invariant over a
given class of links}.  So far, the situation is far from
satisfactory.  For example, the range of the Jones polynomial for
knots is not completely understood yet.

The purpose of the present paper is to provide a useful algebraic
setting for the study of the range of quantum invariants of links and
tangles.  The main ingredients of this setting are
\begin{itemize}
\item a special kind of tangles of arcs, which we call {\em bottom
  tangles},
\item a braided subcategory $\modB $ of the category $\modT $ of (uncolored)
  framed, oriented tangles, which ``acts'' on the set of bottom
  tangles by composition, and provides a convenient way to generate
  all the bottom tangles,
\item for each ribbon Hopf algebra $H$ over a commutative, unital ring
  $\modk $, a braided functor $\modJ \col \modB \rightarrow \Mod_H$ from $\modB $ to the category
  of left $H$--modules.
\end{itemize}

\subsection{Bottom tangles}
\label{sec:bottom-tangles}
When one studies links in the $3$--sphere, it is often useful to
represent a link as the {\em closure} of a tangle consisting only of
arc components.  In such an approach, one first study tangles, and
then obtains results for links from those for tangles, via the closure
operation.  The advantage of considering tangles of arcs is that one
can paste tangles together to obtain another tangle, and such pasting
operations produce useful algebraic structures.  For example, the set
of $n$--component string links, up to ambient isotopy fixing endpoints,
forms a monoid with multiplication induced by vertical pasting.

Bottom tangles, which we study in the present paper, are another kind
of tangles of arcs.  An $n$--component bottom tangle
$T=T_1\cup \cdots\cup T_n$ is a framed tangle consisting of $n$ arcs
$T_1,\ldots,T_n$ in a cube such that all the endpoints of~$T$ are on a
line at the bottom square of the cube, and for each $i=1,\ldots,n$ the
component $T_i$ runs from the $2i$th endpoint on the bottom to the
$(2i-1)$st endpoint on the bottom, where the endpoints are counted
from the left.  The component $T_i$ is called the $i$th component of
$T$.  For example, see \fullref{fig:bottom-tangle} (a).
\labellist\small
\pinlabel {$T_1$} [br] at 20 115
\pinlabel {$T_2$} [bl] at 55 117
\pinlabel {$T_3$} [b] at 104 100
\pinlabel {$L_1$} [br] at 182 120
\pinlabel {$L_2$} [bl] at 215 121
\pinlabel {$L_3$} [b] at 265 105
\pinlabel {(a)} [b] at 60 0
\pinlabel {(b)} [b] at 225 0
\endlabellist
\FIG{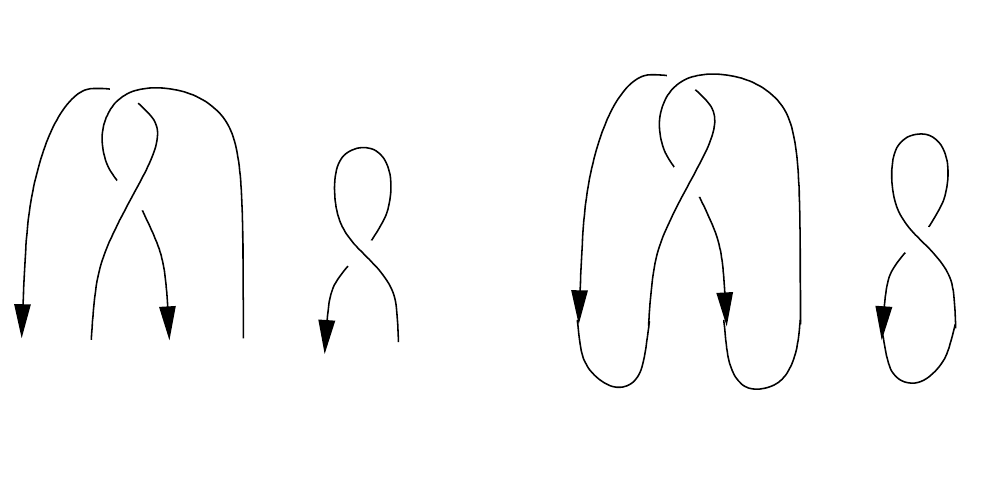}{(a) A $3$--component bottom tangle
$T=T_1\cup T_2\cup T_3$.  (b) The closure $\clo(T)=L_1\cup L_2\cup L_3$ of
$T$.}{height=32mm}

For $n\ge 0$, let $\BT_n$ denote the set of the ambient isotopy classes,
relative to endpoints, of $n$--component bottom tangles.  (As usual,
we often confuse a tangle and its ambient isotopy class.)

There is a natural closure operation which transforms an $n$--component
bottom tangle $T$ into an $n$--component framed, oriented, ordered link
$\clo(T)$, see \fullref{fig:bottom-tangle}~(b).  This operation
induces a function
\begin{equation*}
  \clo_n=\clo\col \BT_n \rightarrow  \LL_n,
\end{equation*}
where $\LL_n$ denotes the set of the ambient isotopy classes of
$n$--component, framed, oriented, ordered links in the $3$--sphere.  It is clear
that $\clo_n$ is surjective, ie, for any link $L$ there is a bottom
tangle $T$ such that $\clo(T)=L$.  Consequently, one can use bottom
tangles to represent links.  In many situations, one can divide the
study of links into the study of bottom tangles and the study of the
effect of closure operation.

\begin{remark}
  \label{r14}
  The notion of bottom tangle has appeared in many places in knot
  theory, both explicitly and implicitly, and is essentially
  equivalent to the notion of based links, as is the case with string
  links.  We establish a specific one-to-one correspondence between
  bottom tangles and string links in
  \fullref{sec:string-links-bottom}.
\end{remark}

\subsection{An approach to quantum link invariants using universal
  invariants of bottom tangles}
\label{sec:an-approach-quantum}

\subsubsection{Universal link invariants associated to ribbon Hopf algebras}
\label{sec:color-link-invar}

For each ribbon Hopf algebra $H$, there is an invariant of links and
tangles, which is called the {\em universal invariant} associated to
$H$, introduced by Lawrence \cite{Lawrence:89,Lawrence:90} in the case
of links and quantized enveloping algebras.  Around the same time,
Hennings \cite{Hennings:96} formulated link invariants associated to
quasitriangular Hopf algebras which do not involve representations but
involve trace functions on the algebras.  Reshetikhin \cite[Section
4]{Reshetikhin:89} and Lee \cite{Lee:90} considered universal
invariants of $(1,1)$--tangles (ie, tangles with one endpoint on the
top and one on the bottom) with values in the center of a quantum
group, which can be thought of as the $(1,1)$--tangle version of the
universal link invariant.  Universal invariants are further
generalized to more general oriented tangles by Lee
\cite{Lee:92a,Lee:92b,Lee:96} and Ohtsuki \cite{Ohtsuki:93}.  Kauffman
\cite{Kauffman:93} and Kauffman and Radford \cite{Kauffman-Radford:01}
defined functorial versions of universal tangle invariant for a
generalization of ribbon Hopf algebra which is called ``(oriented)
quantum algebra''.

The universal link invariants are defined at the quantum group level
and they do not require any representations.  The universal link
invariant have a universality property that colored link invariants
can be obtained from the universal link invariants by taking traces in
the representations attached to components.  Thus, in order to study
the range of the representation-colored link invariants, it suffices,
in theory, to study the range of the universal invariant.

\subsubsection{Universal invariant of bottom tangles and their closures}
\label{sec:universal-invariant}
Here we briefly describe the relationship between the colored link
invariants and the universal invariants, using bottom tangles.

Let $H$ be a ribbon Hopf algebra over a commutative ring $\modk $ with
unit.  For an $n$--component bottom tangle $T\in \BT_n$, the universal
invariant $J_T=J_T^H$ of $T$ associated to $H$ takes value in the
$n$--fold tensor product $H^{\otimes n}$ of $H$.  Roughly speaking,
$J_T\in H^{\otimes n}$ is obtained as follows.  Choose a diagram $D$ of $T$.
At each crossing of $D$, place a copy of universal $R$--matrix
$R\in H^{\otimes 2}$, which is modified in a certain rule using the antipode
$S\col H\rightarrow H$.  Also, place some grouplike elements on the local maxima
and minima.  Finally, read off the elements placed on each component
of $H$.  An example is given in \fullref{fig:example-univ-inv} in
\fullref{sec:defin-univ-invar}.  For a more precise and more
general definition, see \fullref{sec:defin-univ-invar}.

For each $n\ge 0$, the universal invariant defines a function
\begin{equation*}
  J\col \BT_n \rightarrow  H^{\otimes n},\quad  T\mapsto J_T.
\end{equation*}
In this section we do not give the definition of the universal
invariant of links, but it can be obtained from the universal
invariant of bottom tangles, as we explain below.  Set
\begin{equation}
  \label{e15}
  N= \Span_{\modk} \{x y-yS^2(x)\ver x,y\in H\}\subset H.
\end{equation}
The projection
\begin{equation*}
  \tr_q\col H\rightarrow H/N
\end{equation*}
is called the {\em universal quantum trace}, since if $\modk $ is a field
and $V$ is a finite-dimensional left $H$--module, then the quantum trace
$\tr_q^V\col H\rightarrow \modk $ in $V$ factors through $\tr_q$.

The universal link invariant $J_L\in (H/N)^{\otimes n}$ for an $n$--component
framed link $L\in \LL_n$, which we define in \fullref{sec:defin-univ-invar}, satisfies
\begin{equation*}
  J_L = \tr_q^{\otimes n}(J_T),
\end{equation*}
where $T\in \BT_n$ satisfies $\clo(T)=L$.

\subsubsection{Reduction to the colored link invariant}
\label{sec:recov-repr-color}
Let $\modk $ be a field, and $V_1,\ldots,V_n$ be finite-dimensional
left $H$--modules.  Then the quantum invariant $J_{L;V_1,\ldots,V_n}$ of an
$n$--component colored link $(L;V_1,\ldots,V_n)$ can be obtained from
$J_L$ by
\begin{equation*}
  J_{L;V_1,\ldots,V_n} = \bigl(\bar{\tr}_q^{V_1}\otimes \cdots\otimes
\bar{\tr}_q^{V_n}\bigr)(J_L).
\end{equation*}
where $\bar{\tr}_q^{V_i}\col H/N\rightarrow \modk $ is induced by the quantum trace
$\tr_q^{V_i}\col H\rightarrow \modk $.  Hence if $\clo(T)=L$, $T\in \BT_n$, we have
\begin{equation*}
  J_{L;V_1,\ldots,V_n} = \bigl(\tr_q^{V_1}\otimes \cdots\otimes
\tr_q^{V_n}\bigr)(J_T).
\end{equation*}
These facts can be summarized into a commutative diagram:
\begin{equation}
  \label{eq:18}
  \xymatrix{
    \BT_n\ar[r]^{J\quad}\ar[d]_{\clo_n} &{\quad H^{\otimes n}\quad}
    \ar[drr]^{\bigotimes_{i=1}^n\tr_q^{V_i}}\ar[d]_{\tr_q^{\otimes n}}&&\\
    \;\LL_n\;\ar[r]_{J\quad} & (H/N)^{\otimes n}
    \ar[rr]_{\bigotimes_{i=1}^n\bar{\tr}_q^{V_i}} && \modk .
  }
\end{equation}
Given finite-dimensional left $H$--modules $V_1,\ldots,V_n$, we are
interested in the range of $J_{L;V_1,\ldots,V_n}\in \modk $ for $L\in \LL_n$.
Since $\clo_n$ is surjective, it follows from \eqref{eq:18} that
\begin{equation}
  \label{e17}
  \begin{split}
    \{J_{L;V_1,\ldots,V_n}\ver L\in \LL_n\}
    &=\bigl(\bar{\tr}_q^{V_1}\otimes \cdots\otimes \bar{\tr}_q^{V_n}\bigr)(J(\LL_n))\\
    &=\bigl(\tr_q^{V_1}\otimes \cdots\otimes \tr_q^{V_n}\bigr)(J(\BT_n)).
  \end{split}
\end{equation}
Hence, to determine the range of the representation-colored link
invariants, it suffices to determine the images $J(\BT_n)\subset H^{\otimes n}$
for $n\ge 0$ and to study the maps $\tr_q^{V_i}$.

\subsection{The category $\modB $ acting on bottom tangles}
\label{sec:category--acting}
To study bottom tangles, and, in particular, to determine the images
$J(\BT_n)\subset H^{\otimes n}$, it is useful to introduce a braided subcategory
$\modB $ of the category $\modT $ of (uncolored) framed, oriented tangles
which {\em acts on bottom tangles by composition}.

Here we give a brief definition of $\modB $, assuming familiarity with the
braided category structure of $\modT $.  For the details, see \fullref{sec3}.

The objects of $\modB $ are the expressions $\modb ^{\otimes n}$ for $n\ge 0$.
(Later, $\modb $ is identified with an object $\modb =\downarrow \otimes \uparrow $ in $\modT $.)  For
$m,n\ge 0$, a morphism $T$ from $\modb ^{\otimes m}$ to $\modb ^{\otimes n}$ is the ambient
isotopy class relative to endpoints of a framed, oriented tangle $T$
satisfying the following.
\begin{enumerate}
\item There are $2m$ (resp. $2n$) endpoints on the top (resp. bottom),
  where the orientations are as $\downarrow \uparrow \cdots\downarrow \uparrow $.
\item For any $m$--component bottom tangle $T'$, the composition $TT'$
  (obtained by stacking $T'$ on the top of $T$) is an $n$--component
  bottom tangle.
\end{enumerate}
It follows that $T$ consists of $m+n$ arc components and no circle
components.  For example, see \fullref{fig:example-bottom} (a).
\labellist\small
\pinlabel {(a)} [b] at 60 0
\pinlabel {(b)} [b] at 205 0
\pinlabel {(c)} [b] at 335 0
\pinlabel {$T'$} at 392 128
\pinlabel {$T$} at 392 70
\endlabellist
\FIG{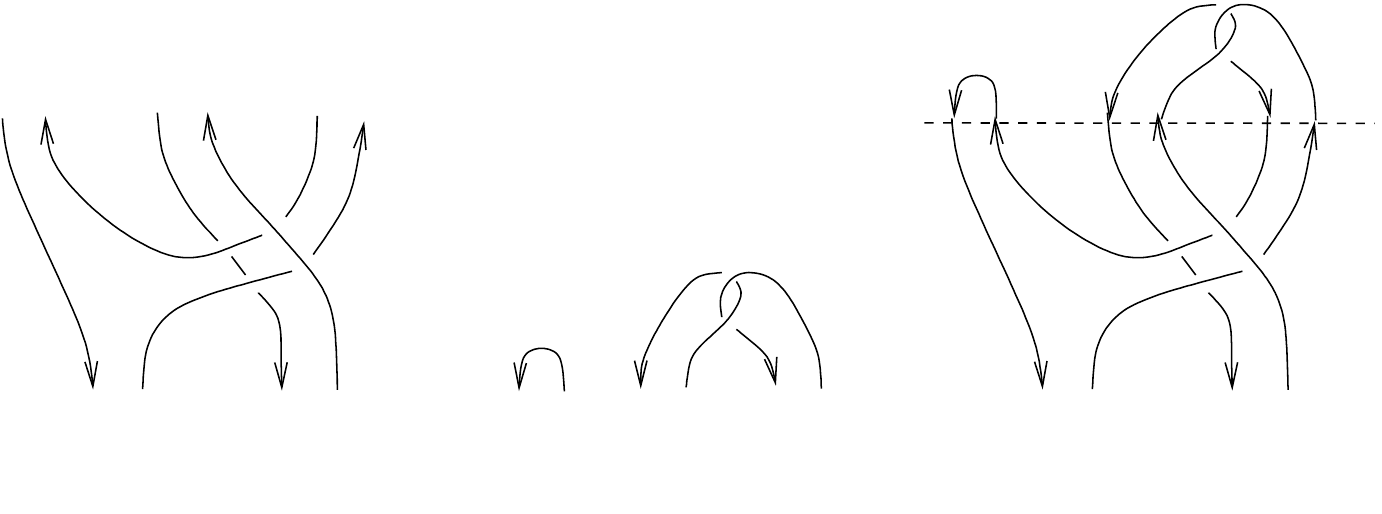}{(a) A morphism $T\in \modB (3,2)$.  (b) A bottom tangle
$T'\in \BT_3$.  (c) The composition $TT'\in \BT_2$.}{height=32mm} The set
$\modB (\modb ^{\otimes m},\modb ^{\otimes n})$ of morphisms from $\modb ^{\otimes m}$ to $\modb ^{\otimes n}$ in
$\modB $ is often denoted simply by $\modB (m,n)$.  The composition of two
morphisms in $\modB $ is pasting of two tangles vertically, and the
identity morphism $1_{\modb ^{\otimes m}}=\downarrow \uparrow \cdots\downarrow \uparrow $ is a tangle consisting
of $2m$ vertical arcs.  The monoidal structure is given by pasting two
tangles side by side.  The braiding is defined in the usual way; ie,
the braiding for the generating object $\modb \in \Ob(\modB )$ and itself is
given by
$$\psi _{\modb ,\modb }=\raisebox{-0.8em}{\incl{2em}{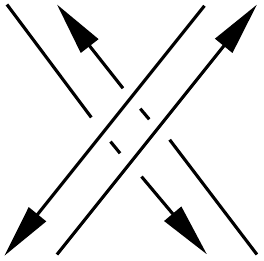}}.$$  For each $n\ge 0$, we can
identify $\BT_n$ with $\modB (0,n)$.  The category $\modB $ acts on
$\BT=\coprod_{n\ge 0}\BT_n$ via the functions
\begin{equation*}
  \modB (m,n)\times \BT_m \rightarrow  \BT_n,\quad (T,T')\mapsto TT'.
\end{equation*}
In the canonical way, one may regard this action as a functor
\begin{equation*}
  \modB (\one ,-)\col \modB \rightarrow \Sets
\end{equation*}
from $\modB $ to the category $\Sets$ of sets, which maps an object
$\modb ^{\otimes n}$ into $\BT_n$.

\subsection{The braided functor $\modJ \col \modB \rightarrow \Mod_H$}
\label{sec:functor-mod_h}
Let $\Mod_H$ denote the category of left $H$--modules, with the
standard braided category structure.  Unless otherwise stated, we
regard $H$ as a left $H$--module with the left adjoint action.

We study a braided functor
\begin{equation*}
  \modJ \col \modB \rightarrow \Mod_H,
\end{equation*}
which is roughly described as follows.  For the details, see \fullref{sec8}.  For objects, we set $\modJ (\modb ^{\otimes n})=H^{\otimes n}$, where
$H^{\otimes n}$ is given the standard tensor product left $H$--module
structure.  For $T\in \modB (m,n)$, the left $H$--module homomorphism
$\modJ (T)\col H^{\otimes m}\rightarrow H^{\otimes n}$ maps $x=\sum x_1\otimes \cdots\otimes x_m\in H^{\otimes m}$ to the
element $\modJ (T)(x)$ in $H^{\otimes n}$ obtained as follows.  Set
\begin{equation*}
  \eta _{\modb} =\raisebox{-0.5em}{\incl{1.5em}{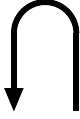}}\in \modB (0,1),
\end{equation*}
and set for $m\ge 0$
\begin{equation*}
  \eta _m = \eta _{\modb} ^{\otimes
m}=\raisebox{-0.2em}{\incl{1em}{}}\cdots\raisebox{-0.2em}{\incl{1em}{bot0}}\in \modB (0,m).
\end{equation*}
Consider a diagram for the composition $T\eta _m$, and put the element
$x_i$ on the $i$th component (from the left) of $\eta _m$ for
$i=1,\ldots,m$, and put the copies of universal $R$--matrix and the
grouplike element on the strings of $T$ as in the definition of $J_T$.
Then the element $\modJ (T)(x)\in H^{\otimes n}$ is read off from the diagram.
(See \fullref{fig:FT} in \fullref{sec:braided-functor-h} to get
a hint for the definition.)  We see in \fullref{sec:braided-functor-h} that $\modJ (T)$ is a left $H$--module
homomorphism, and $\modJ $ is a well-defined braided functor.

\subsection{Generators of the braided category $\modB $}
\label{sec:gener--appl}

Recall that a (strict) braided category $M$ is said to be generated by
a set $X\subset \Ob(M)$ of objects and a set $Y\subset \Mor(M)$ of morphisms if
every object of~$M$ is a tensor product of finitely many copies of
objects from $X$, and if every morphism of $M$ is an iterated tensor
product and composition of finitely many copies of morphisms from $Y$
and the identity morphisms of the objects from~$X$.  In the category
$\modT $, ``tensor product and composition'' is ``horizontal and vertical
pasting of tangles''.

\begin{Theorem}[\fullref{thm:2}]
  \label{thm:39}
  As a braided subcategory of $\modT $, $\modB $ is generated by the object
  $\modb $ and the morphisms
  \begin{equation*}
    \def\sss{1.5em}
    \eta _{\modb} =\raisebox{-0.5em}{\incl{\sss}{bot0}},\quad
    \def\sss{2em}
    \mu _{\modb} =\raisebox{-0.6em}{\incl{\sss}{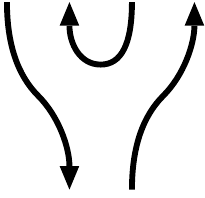}},\quad
    v_+=\raisebox{-0.6em}{\incl{\sss}{}},\quad
    v_-=\raisebox{-0.6em}{\incl{\sss}{}},\quad
    c_+=\raisebox{-0.6em}{\incl{\sss}{}},\quad
    c_-=\raisebox{-0.6em}{\incl{\sss}{}}.
  \end{equation*}
\end{Theorem}

Consequently, any bottom tangle can be obtained by horizontal and
vertical pasting from finitely many copies of the above-listed
tangles, the braidings $\psi _{\modb ,\modb },\psi _{\modb ,\modb }^{-1}$ and the identity
$1_{\modb} =\downarrow \;\uparrow $.  \fullref{thm:39} implies that, as a category, $\modB $
is generated by the morphisms
\begin{equation*}
  f_{(i,j)} = \modb ^{\otimes i}\otimes f\otimes \modb ^{\otimes j}
\end{equation*}
for $i,j\ge 0$ and $f\in \bigl\{\eta _{\modb} ,\mu _{\modb} ,v_{\pm}
,c_{\pm}
,\psi _{\modb ,\modb }^{\pm 1}\bigr\}$.
\fullref{thm:39} provides a convenient way to generate all the
bottom tangles, and hence all the links via closure operation.

We can use \fullref{thm:39} to determine the range
$J(\BT)=\{J_T\ver T\in \BT\}$ of the universal invariant of bottom tangles
as follows.  \fullref{thm:39} and functoriality of~$\modJ $ implies
that any morphism $f$ in $\modB $ is the composition of finitely many copies
of the left $H$--module homomorphisms
\begin{equation}
  \label{eq:26}
  \modJ (f_{(i,j)})= 1_H^{\otimes i}\otimes \modJ (f)\otimes 1_H^{\otimes j},
\end{equation}
for $i,j\ge 0$ and $f\in \{\eta _{\modb} ,\mu _{\modb} ,v_{\pm} ,c_{\pm} ,\psi _{\modb ,\modb }^{\pm 1}\}$.
Hence we have the following characterization of the range of $J$ for
bottom tangles.

\begin{Theorem}
  \label{r6}
  The set $J(\BT)$ is equal to the smallest subset of
  $\coprod_{n\ge 0}H^{\otimes n}$ containing $1\in \modk =H^{\otimes 0}$ and stable under
  the functions $\modJ (f_{(i,j)})$ for $i,j\ge 0$ and
  $f\in \bigl\{\eta _{\modb} ,\mu _{\modb} ,v_{\pm} ,c_{\pm} ,\psi _{\modb ,\modb
}^{\pm 1}\bigr\}$.
\end{Theorem}

See \fullref{sec:values-jt} for some variants of \fullref{r6},
which may be more useful in applications than \fullref{r6}.

\subsection{Hopf algebra action on bottom tangles}
\label{sec:hopf-algebra-action-1}

We define a kind of ``Hopf algebra action'' on the set $\BT$, which is
formulated as an ``exterior Hopf algebra'' in the category $\modB $.
Roughly speaking, this exterior Hopf algebra in $\modB $ is an
``extension'' of an algebra structure for the object $\modb $ in $\modB $ to a
Hopf algebra structure at the level of sets of morphisms in $\modB $.  This
``extension'' is formulated as a functor
\begin{equation*}
  F_{\modb} \col \modH \rightarrow \Sets,
\end{equation*}
from the free strict braided category $\modH $ generated by a Hopf algebra
$\HH$ to the category $\Sets$ of sets, where each object $\HH^{\otimes m}$,
$m\ge 0$, is mapped to the set $\modB (\one ,\modb ^{\otimes m})=\BT_m$.  This functor
essentially consists of functions
\begin{equation}
  \label{e11}
  \modH (\HH^{\otimes m},\HH^{\otimes n})\times \BT_m \rightarrow  \BT_n\quad \text{for $m,n\ge 0$},
\end{equation}
which can be regarded as a left action of the category $\modH $ on the
graded set~$\BT$.

The two functors $\modB (\one ,-)\col \modB \rightarrow \Sets$ and
$F_{\modb} \col \modH \rightarrow \Sets$ are related as follows.  Let $\langle \sfA\rangle $ denote the braided
category freely generated by an algebra $\sfA$.  Note that the
morphisms $\mu _{\modb} \col \modb ^{\otimes 2}\rightarrow \modb $ and
$\eta _{\modb} \col \one \rightarrow \modb $ in $\modB $ define an
algebra structure for the object $\modb $.  Consider the following
diagram
\begin{equation*}
  \begin{CD}
    \langle \sfA\rangle  @>i_{\modB ,\modb }>> \modB \\
    @Vi_{\modH ,\HH}VV @VV{\modB (\one ,-)}V\\
    \modH  @>>{F_{\modb} }> \Sets.
  \end{CD}
\end{equation*}
Here the two arrows $i_{\modB ,\modb }$ and $i_{\modH ,\HH}$ are the unique
braided functors that map the algebra structure for $\sfA$ into those
of $\modb $ and $\HH$, respectively.  Both $i_{\modB ,\modb }$ and $i_{\modH ,\HH}$
are faithful.  The above diagram turns out to be commutative.
In other words, the action of the algebra
structure in $\modB $ on~$\BT$ extends to an action by a Hopf algebra
structure on $\BT$.

\begin{remark}
  The action of $\modH $ on $\BT$ mentioned above extends to an action of
  a category $\B$ of ``bottom tangles in handlebodies'' which we
  shortly explain in \fullref{sec:category-0}.  Also, the above
  Hopf algebra action is closely related to the Hopf algebra structure
  in the category $\modC $ of cobordisms of connected, oriented surfaces
  with boundary parameterized by a circle
  (see Crane--Yetter \cite{Crane-Yetter:99} and Kerler \cite{Kerler:99}),
  and also to the Hopf algebra properties satisfied by claspers
  (see Habiro \cite{Habiro:claspers}).
\end{remark}

Now we explain the relationship between the Hopf algebra action on
$\BT$ and the braided functor $\modJ \col \modB \rightarrow \Mod_H$.  Let $H$ be a ribbon
Hopf algebra, hence in particular $H$ is quasitriangular.  Recall that
the {\em transmutation} $\bH$ of a quasitriangular Hopf algebra $H$ is
a braided Hopf algebra structure in $\Mod_H$ defined for the object
$H$ with the left adjoint action in the braided category $\Mod_H$,
with the same algebra structure and the same counit as $H$ but with a
twisted comultiplication $\bD\col H\rightarrow H\otimes H$ and a twisted antipode
$\bS\col H\rightarrow H$.  For details, see Majid
\cite{Majid:algebras,Majid}.  The
transmutation $\bH$ yields a braided functor $F_{\bH}\col \modH \rightarrow \Mod_H$,
which maps the Hopf algebra structure of $\HH$ into that of $\bH$.
Via $F_{\bH}$, the category $\modH $ acts on the $H^{\otimes n}$ as
\begin{equation*}
  \modH (\HH^{\otimes m},\HH^{\otimes n}) \times  H^{\otimes m} \rightarrow  H^{\otimes n},\quad
  (f,x)\mapsto F_{\bH}(f)(x).
\end{equation*}
We have the following commutative diagram
\begin{equation*}
  \begin{CD}
    \modH (\HH^{\otimes m},\HH^{\otimes n})\times  \BT_m @>>> \BT_n \\
    @V{\id\times J}VV @VVJV \\
    \modH (\HH^{\otimes m},\HH^{\otimes n}) \times  H^{\otimes m} @>>> H^{\otimes n}.
  \end{CD}
\end{equation*}
This means that the Hopf algebra action on the bottom tangles is
mapped by the universal invariant maps $J\col \BT_n\rightarrow H^{\otimes n}$ into the
Hopf algebra action on the $H^{\otimes n}$, defined by transmutation.  This
fact may be considered as a {\em topological interpretation of
transmutation}.

\subsection{Local moves}
\label{sec:local-moves}

The setting of bottom tangles and the category $\modB $ is also useful in
studying local moves on bottom tangles, and hence on links.  Recall
that a type of local move can be defined by specifying a pair of
tangles $(t,t')$ in a $3$--ball with the same set of endpoints and with
the same orientations and framings at the endpoints.  Two tangles $u$
and $u'$ in another $3$--ball $B$ are said to be related by a
$(t,t')$--move if $u'$ is, up to ambient isotopy which fixes endpoints,
obtained from $u$ by replacing a ``subtangle'' $t$ in $u$ with a copy
of $t'$.  Here, by a subtangle of a tangle $u$ in $B$ we mean a tangle
$u\cap B'$ contained in a $3$--ball $B'\subset \opint B$.

In the following, we restrict our attention to the case where both the
tangles $t$ and $t'$ consist only of arcs.  In this case, we can
choose two bottom tangles $T,T'\in \BT_n$ such that the notion of
$(t,t')$--move is the same as that of $(T,T')$--move.

The following result implies that the notion of $(T,T')$--move for
bottom tangles, where $T,T'\in \BT_n$, can be formulated in a totally
algebraic way within the category $\modB $.

\begin{Theorem}[Consequence of \fullref{lem:2} and
    \fullref{thm:5}]
  \label{thm:44}
  Let $T,T'\in \BT_m$ and $U,U'\in \BT_n$ with $m,n\ge 0$.  Then $U$ and
  $U'$ are related by a $(T,T')$--move if and only if there is
  $W\in \modB (m,n)$ such that $U=WT$ and $U'=WT'$.
\end{Theorem}

The setting of bottom tangles and the category $\modB $ is useful in
studying local moves on bottom tangles.  In \fullref{sec5}, we
formulate in algebraic terms the following typical questions in the
theory of local moves, under some mild conditions.
\begin{itemize}
\item When are two bottom tangles related by a sequence of a given set
  of local moves?  (\fullref{lem:7}.)
\item When are two bottom tangles equivalent under the equivalence
  relation generated by a given set of local moves?  (\fullref{lem:4}.)
\item When is a bottom tangle related to the trivial bottom tangle
  $\eta _n$ by just one local move of a given type?  (\fullref{thm:4}.)
\item When is a bottom tangle equivalent to $\eta _n$ under the
  equivalence relation generated by a given set of local moves?
  (\fullref{thm:19}.)
\end{itemize}
It is easy to modify the answers to the above questions for bottom
tangles into those for links, via closure.  Some of the algebraic
formulations of the above questions are combined with the functor $\modJ $
in \fullref{sec:values-jt}.

A remarkable application of the algebraic formulation of local moves
is the following.  A {\em delta move} (see Murakami and Nakanishi
\cite{Murakami-Nakanishi}) or a
{\em Borromean transformation} (see Matveev \cite{Matveev}) can be defined as a
$(\eta _3,B)$--move, where $B\in \BT_3$ is the {\em Borromean tangle}
depicted in \fullref{fig:borromean}.
\FIG{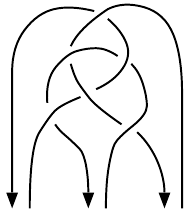}{The Borromean
tangle $B$}{height=20mm} An $n$--component link has a zero linking
matrix if and only if it is related to the $n$--component unlink by a
sequence of delta moves \cite{Murakami-Nakanishi}.  Using the obvious
variation of this fact for bottom tangles, we obtain the following.

\begin{Theorem}[Part of \fullref{thm:24}]
  \label{thm:26}
  A bottom tangle with zero linking matrix is obtained by pasting
  finitely many copies of
  $1_{\modb} ,\psi _{\modb ,\modb },\psi _{\modb ,\modb }^{-1},\mu
_{\modb} ,\eta _{\modb} ,\gamma _+,\gamma _-,B$, where
  \begin{gather*}
    \gamma _+ =\raisebox{-5mm}{\incl{15mm}{}},\quad
    \gamma _- =\raisebox{-5mm}{\incl{15mm}{}}\in \modB (1,2).
  \end{gather*}
\end{Theorem}

For applications of \fullref{thm:26} to the universal invariant of
bottom tangles with zero linking matrix, see \fullref{thm:29},
which states that the universal invariant associated to a ribbon Hopf
algebra $H$ of $n$--component bottom tangles with zero linking matrix
is contained in a certain subset of $H^{\otimes n}$.  The case of the
quantized enveloping algebra $U_h(\g)$ of simple Lie algebra $\g$ will
be used in future publications \cite{Habiro:in-preparation} (for
$\g=sl_2$) and \cite{Habiro-Le:in-preparation} (for general $\g$) to
show that there is an invariant of integral homology spheres with
values in the completion
$\varprojlim_{n}\modZ [q]/((1-q)(1-q^2)\cdots(1-q^n))$ studied in
\cite{Habiro:cyclotomic}, which unifies the quantum $\g$ invariants at
all roots of unity, as announced in \cite{Habiro:rims2001},
\cite[Conjecture 7.29]{Ohtsuki:problem}.

\subsection{Other applications}
\label{sec:other-applications-1}
The setting of bottom tangles, the category $\modB $ and the functor $\modJ $
can be applied to the following notions in knot theory: unknotting
number (\fullref{sec:unknotting-number}), Seifert surface, knot
genus and boundary links (\fullref{sec:commutators}), unoriented
spanning surface, crosscap number and $\modZ _2$--boundary links (\fullref{sec:unorient-spann}), clasper moves (\fullref{sec:cn-moves}),
Goussarov--Vassiliev finite type invariants (\fullref{sec:vass-gouss-filtr-2}), twist moves and Fox's congruence
(\fullref{sec:twist-moves}), ribbon knots (\fullref{sec:ribbon-discs}), and the Hennings $3$--manifold invariants
(\fullref{sec11}).

\subsection{Organization of the paper}
\label{sec:organization-paper}

Here we briefly explain the organization of the rest of the paper.
\fullref{sec2} provides some preliminary facts and notations about
braided categories and Hopf algebras.  In \fullref{sec3}, we
recall the definition of the category $\modT $ of framed oriented
tangles, and define the subcategory $\modB $ of $\modT $.  In \fullref{sec4}, we study a subcategory $\modB _0$ of $\modB $, which is
used in \fullref{sec5}, where we study local moves on tangles and
give a set of generators of $\modB $.  \fullref{sec6} deals with
the Hopf algebra action on bottom tangles.  In \fullref{sec7}, we
recall the definition of the universal tangle invariant and provide
some necessary facts.  In \fullref{sec8}, we define and study the
braided functor $\modJ \col \modB \rightarrow \Mod_H$.  In
\fullref{sec9}, we give several applications of the results in the
previous sections to the values of the universal invariant.  In
\fullref{sec:knots-tensor-product}, we modify the functor $\modJ $
into a braided functor $\tilde{\modJ }\col \modB \rightarrow
\Mod_H$, which is useful in the study of the universal invariants of
bottom knots.  In \fullref{sec10}, we define a refined version of
the universal link invariant by giving a necessary and sufficient
condition that two bottom tangles have the same closures.  In \fullref{sec11}, we reformulate the Hennings invariant of $3$--manifolds
in our setting, using an algebraic version of Kirby calculus.  In
\fullref{sec:string-links-bottom}, we relate the structure of the
sets of bottom tangles to the sets of string links.  In \fullref{sec12}, we give some remarks.

\begin{remark}
  \label{r12}
  Around the same time as the first preprint version of the present
  paper is completed, a paper by Brugui\`eres and Virelizier \cite{BV}
  appeared on the arXiv.  Part of \cite{BV} is closely related to part
  of the present paper.  The present paper is independent of \cite{BV}
  (except this \fullref{r12}).

  The notion of bottom tangle is equivalent to that of ``ribbon
  handle'' in \cite{BV}.  \fullref{thm:39} is related to
  \cite[Theorems 1.1 and 3.1]{BV}.  \fullref{sec11} is related to
  \cite[Section 5]{BV}.  \fullref{sec:string-links-bottom} is
  closely related to \cite{BV}.
\end{remark}

\begin{acknowledgments}
  This work was partially supported by the Japan Society for the
  Promotion of Science, Grant-in-Aid for Young Scientists (B),
  16740033.

  The author would like to thank Thang Le, Gregor Masbaum,
  Jean-Baptiste Meilhan, and Tomotada Ohtsuki for helpful discussions
  and comments.
\end{acknowledgments}

\section{Braided categories and Hopf algebras}
\label{sec2}
In this section, we fix some notations concerning monoidal and braided
categories, and braided Hopf algebras.

If $C$ is a category, then the set (or class) of objects in $C$ is
denoted by $\Ob(C)$, and the set (or class) of morphisms in $C$ is
denoted by $\Mor(C)$.  For $a,b\in \Ob(C)$, the set of morphisms from
$a$ to $b$ is denoted by $C(a,b)$.  For $a\in \Ob(C)$, the identity
morphism $1_a\in C(a,a)$ of $a$ is sometimes denoted by $a$.

\subsection{Monoidal and braided categories}
\label{sec:mono-braid-categ}

We use the following notation for monoidal and braided categories.
See Mac\,Lane \cite{MacLane} for the definitions of monoidal and braided
categories.  If $\modM $ is a monoidal category (also called tensor
category), then the tensor functor is denoted by $\otimes _{\modM} $ and
the unit object by $\one _{\modM} $.  We omit the subscript $\modM $
and write $\otimes _{\modM} =\otimes $ and $\one _{\modM} =\one $ if
there is no fear of confusion.  If $\modM $ is a braided category,
then the braiding of $a,b\in \Ob(\modM )$ is denoted by
\begin{equation*}
  \psi _{a,b}^{\modM} =\psi _{a,b}\col a\otimes b\rightarrow b\otimes a
\end{equation*}
for $a,b\in \Ob(\modM )$.

Definition of monoidal category involves also the associativity and
the unitality constraints, which are functorial isomorphisms
\begin{gather*}
  a\otimes (b\otimes c) \simeqto (a\otimes b)\otimes c,\quad
  \one \otimes a\simeqto a,\quad a\otimes \one  \simeqto a.
\end{gather*}
 A monoidal category is called {\em strict} if these constraints are
identity morphisms.  In what follows, a strict monoidal category is
simply called a ``monoidal category''.  Similarly, a strict braided
category is called a ``braided category''.  Also, we sometimes need
non-strict monoidal or braided categories, such as the category of left
modules over a Hopf algebra.  When this is the case, we usually
suppress the associativity and the unitality constraints, and we argue
as if they are strict monoidal or braided categories.  This should not
cause confusion.

A {\em strict monoidal functor} from a (strict) monoidal category $\modM $
to a (strict) monoidal category $\modM '$ is a functor $F\col \modM \rightarrow \modM '$ such
that
\begin{align*}
  F\otimes  &= \otimes (F\times F) \col  \modM \times \modM  \rightarrow  \modM ',\\
  F(\one _{\modM} ) &= 1_{\modM '},
\end{align*}
ie, a strict monoidal functor is a functor which (strictly)
preserves $\otimes $ and $\one $.  Unless otherwise stated, by a ``monoidal
functor'' we mean a strict monoidal functor.  A {\em braided functor}
from a braided category $\modM $ to a braided category $\modM '$ is a monoidal
functor $F\col \modM \rightarrow \modM '$ such that $F(\psi ^{\modM} _{a,b})=\psi ^{\modM '}_{F(a),F(b)}$
for all $a,b\in \Ob(\modM )$.

A monoidal category $M$ is said to be {\em generated} by a set
$X\subset \Ob(M)$ of objects and a set $Y\subset \Mor(M)$ of morphisms if every
object of $M$ is a tensor product of finitely many copies of objects
from $X$, and every morphism of $M$ is an iterated tensor product and
composition of finitely many copies of morphisms from $Y$ and the
identity morphisms of the objects from $X$.  A braided category~$M$ is
said to be generated by $X\subset \Ob(M)$ and $Y\subset \Mor(M)$ if $M$ is
generated as a monoidal category by $X$ and
$Y\cup \bigl\{\psi _{x,y}^{\pm 1}\ver x,y\in X\bigr\}$.

\subsection{Algebras and Hopf algebras in braided categories}
\label{sec:braid-hopf-algebr}
Here we fix the notations for algebras and Hopf algebras in a monoidal
or braided category.  We refer the reader to
Majid \cite{Majid:algebras,Majid} for the details.

An {\em algebra} (also called {\em monoid}) in a monoidal category
$\modM $ is an object $A$ equipped with morphisms $\mu \col A\otimes A\rightarrow A$
({\em multiplication}) and $\eta \col \one \rightarrow A$ ({\em unit}) satisfying
\begin{gather*}
  \mu (\mu \otimes A)=\mu (A\otimes \mu ),\qquad \mu (\eta \otimes A)=1_A=\mu (A\otimes \eta ).
\end{gather*}
A {\em coalgebra} $C$ in $\modM $ is an object $C$ equipped with morphisms
$\Delta \col A\rightarrow A\otimes A$ ({\em comultiplication}) and $\epsilon \col A\rightarrow \one $ ({\em counit})
satisfying
\begin{gather*}
  (\Delta \otimes A)\Delta =(A\otimes \Delta )\Delta ,\qquad (\epsilon \otimes A)\Delta =1_A=(A\otimes \epsilon )\Delta .
\end{gather*}
A {\em Hopf algebra} in a braided category $\modM $ is an object $H$ in $\modM $
equipped with an algebra structure $\mu ,\eta $, a coalgebra structure
$\Delta ,\epsilon $ and a morphism $S\col A\rightarrow A$ ({\em antipode}) satisfying
\begin{gather*}
  \epsilon \eta  = 1_{\one} ,\qquad
  \Delta \eta  = \eta \otimes \eta ,\qquad  \epsilon \mu  = \epsilon \otimes \epsilon ,\\
  \Delta \mu  = (\mu \otimes \mu )(A\otimes \psi _{A,A}\otimes A)(\Delta \otimes \Delta ),\\
  \mu (A\otimes S)\Delta =\mu (S\otimes A)\Delta =\eta \epsilon .
\end{gather*}
Later, we sometimes use the notations $\mu _A=\mu $, $\eta _A=\eta $,
$\Delta _A=\Delta $, $\epsilon _A=\epsilon $, $S_A=S$, to distinguish structure morphisms from
different Hopf algebras.

A Hopf algebra (in the usual sense) over a commutative, unital ring
$\modk $ can be regarded as a Hopf algebra in the symmetric monoidal
category of $\modk $--modules.

If $A$ is an algebra in $\modM $, then let
$\smash{\mu_A^{[n]}}=\smash{\mu^{[n]}}\col A^{\otimes n}\rightarrow A$ ($n\ge 0$) denote the {\em $n$--input
multiplication} defined by $\mu ^{[0]}=\eta $, $\mu ^{[1]}=1_A$, and
\begin{equation*}
  \mu ^{[n]} = \mu (\mu \otimes 1)\cdots(\mu \otimes 1^{\otimes (n-2)})
\end{equation*}
for $n\ge 2$.  Similarly, if $C$ is an coalgebra, then let
$\Delta _C^{[n]}=\Delta ^{[n]}\col C\rightarrow C^{\otimes n}$ denote the {\em $n$--output
comultiplication} defined by $\Delta ^{[0]}=\epsilon $, $\Delta ^{[1]}=1_C$ and
\begin{equation*}
  \Delta ^{[n]} = (\Delta \otimes 1^{\otimes (n-2)})\cdots(\Delta \otimes 1)\Delta
\end{equation*}
for $n\ge 2$.

\section{The category $\modT $ of tangles and the subcategory $\modB $}
\label{sec3}

In this section, we first recall the definition of the category $\modT $
of framed, oriented tangles.  Then we give a precise definition of the
braided subcategory $\modB $ of $\modT $.

In the rest of the paper, by an ``isotopy'' between two tangles we
mean ``ambient isotopy fixing endpoints''.  Thus, two tangles are said
to be ``isotopic'' if they are ambient isotopic relative to endpoints.

\subsection{The category $\modT $ of tangles}
\label{sec:category-tangles}

Here we recall the definition of the braided category $\modT $ of framed,
oriented tangles, and fix the notations.  For details, see
Yetter \cite{Yetter:88,Yetter}, Turaev \cite{Turaev:89,Turaev}, Freyd and Yetter
\cite{Freyd-Yetter}, Shum \cite{Shum} and Kassel \cite{Kassel}.

The objects in the category $\modT $ are the tensor words of symbols $\downarrow $
and $\uparrow $, ie, the expressions $x_1\otimes \cdots\otimes x_n$ with
$x_1,\ldots,x_n\in \{\downarrow ,\uparrow \}$, $n\ge 0$.  The tensor word of length $0$ is
denoted by $\one =\one _{\modT} $.  The morphisms $T\col w\rightarrow w'$ between
$w,w'\in \Ob(\modT )$ are the isotopy classes of framed, oriented tangles in
a cube $[0,1]^3$ such that the endpoints at the top are described by
$w$ and those at the bottom by~$w'$, see \fullref{fig:tangle} for
example.  \FIG{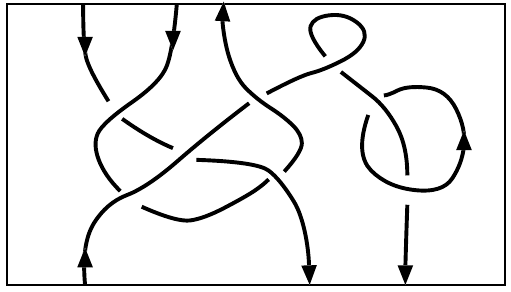}{A tangle
$T\col \downarrow \otimes \downarrow \otimes \uparrow \rightarrow \uparrow \otimes \downarrow \otimes \downarrow $}{height=20mm} We use the blackboard
framing convention in the figures.  In what follows, by ``tangles'' we
mean framed, oriented tangles unless otherwise stated.  As usual, we
systematically confuse a morphism in $\modT $ with a tangle representing
it.

The composition $gf$ of a composable pair $(f,g)$ of morphisms in $\modT $
is obtained by placing $g$ below $f$, and the tensor product $f\otimes g$
of two morphisms $f$ and~$g$ is obtained by placing $g$ on the right
of $f$.  Graphically,
\begin{equation*}
\labellist\small
\pinlabel {$f$} at 100 746
\pinlabel {$g$} at 100 711
\endlabellist
  gf = \raisebox{-1.3em}{\incl{3.6em}{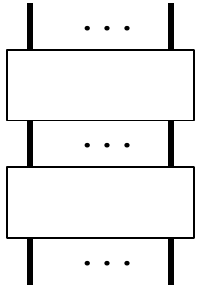}},\quad
\labellist\small
\pinlabel {$f$} at 100 701
\pinlabel {$g$} at 160 701
\endlabellist
  f\otimes g = \raisebox{-0.9em}{\incl{2.5em}{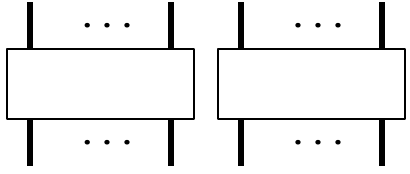}}.
\end{equation*}
The braiding $\psi _{w,w'}\col w\otimes w'\rightarrow w'\otimes w$ for $w,w'\in \Ob(\modT )$ is
the positive braiding of parallel families of strings.  For
$w\in \Ob(\modT )$, the dual $w^*\in \Ob(\modT )$ of $w$ is defined by $\one ^*=\one $,
$\downarrow ^*=\uparrow $, $\uparrow ^*=\downarrow $, and
\begin{gather*}
  (x_1\otimes \cdots\otimes x_n)^*=x_n^*\otimes \cdots\otimes x_1^*\quad (x_1,\ldots,x_n\in \{\downarrow ,\uparrow \}, n\ge 2).
\end{gather*}
For $w\in \Ob(\modT )$, let
\begin{gather*}
  \ev_w \col w^*\otimes w\rightarrow \one ,\quad \coev_w\col \one \rightarrow w\otimes w^*
\end{gather*}
denote the duality morphisms.
  For each object $w$ in $\modT $, let
$t_w\col w\rightarrow w$ denote the positive full twist defined by
\begin{equation*}
\labellist\small
\pinlabel {$w$} [b] at 3 68
\pinlabel {$w$} [t] at 3 12
\endlabellist
  t_w = (w\otimes \ev_{w^*})(\psi _{w,w}\otimes w^*)(w\otimes
\coev_w)=\raisebox{-6mm}{\incl{15mm}{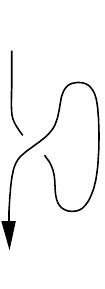}}.
\end{equation*}
It is well known that $\modT $ is generated as a monoidal category by
the objects $\downarrow ,\uparrow $ and the morphisms
\begin{gather*}
  \psi _{\downarrow ,\downarrow }=\raisebox{-0.4em}{\incl{1.3em}{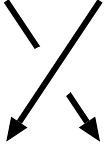}},\quad
  \psi _{\downarrow ,\downarrow }^{-1}=\raisebox{-0.4em}{\incl{1.3em}{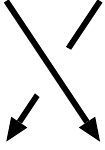}},\quad
  \ev_{\downarrow} =\raisebox{-0.3em}{\incl{1em}{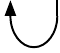}},\quad
  \ev_{\uparrow} =\raisebox{-0.3em}{\incl{1em}{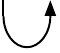}},\quad
  \coev_{\downarrow} =\raisebox{-0.3em}{\incl{1em}{}},\quad
  \coev_{\uparrow} =\raisebox{-0.3em}{\incl{1em}{evup}}.
\end{gather*}

\subsection{The braided subcategory $\modB $ of $\modT $}
\label{sec:subcategories}

Two tangles $T,T'\in \modT (w,w')$, $w,w'\in \Ob(\modT )$, are said to be {\em
  homotopic} (to each other) if there is a homotopy between $T$ and
  $T'$ which fixes the endpoints, where the framings are ignored.  We
  write $T\simh T'$ if $T$ and $T'$ are homotopic.  Note that two
  tangles are homotopic if and only if they are related by a finite
  sequence of isotopies, crossing changes, and framing changes.

Now we define a braided subcategory $\modB $ of $\modT $.  Set
\begin{equation*}
  \Ob(\modB )=\{\modb ^{\otimes m}\ver m\ge 0\}\subset \Ob(\modT ),
\end{equation*}
where $\modb =\downarrow \otimes \uparrow \in \Ob(\modT )$.  Set
\begin{gather*}
  \eta _{\modb} =\coev_{\downarrow} =\cvdn\in \modT (\one ,\modb ),\\
  \eta _n = \eta _{\modb} ^{\otimes n}\in \modT (\one ,\modb ^{\otimes n})\quad \text{for $n\ge 0$}.
\end{gather*}
For $m,n\ge 0$, set
\begin{equation*}
  \modB (\modb ^{\otimes m},\modb ^{\otimes n})
  =\{T\in \modT (\modb ^{\otimes m},\modb ^{\otimes n})\ver T\eta _m\simh\eta _n\}.
\end{equation*}
If $T\in \modB (\modb ^{\otimes l},\modb ^{\otimes m})$ and $T'\in \modB (\modb ^{\otimes m},\modb ^{\otimes n})$, then we
have $T'T\eta _l\simh T'\eta _m\simh \eta _n$, hence
$T'T\in \modB (\modb ^{\otimes l},\modb ^{\otimes n})$.  Clearly, we have
$1_{\modb ^{\otimes n}}\in \modB (\modb ^{\otimes n},\modb ^{\otimes n})$.  Hence $\Ob(\modB )$ and the
$\modB (\modb ^{\otimes m},\modb ^{\otimes n})$ form a subcategory of $\modT $.

If $T\in \modB (\modb ^{\otimes m},\modb ^{\otimes n})$ and $T'\in \modB (\modb ^{\otimes m'},\modb ^{\otimes n'})$, then
we have
\begin{gather*}
  (T\otimes T')\eta _{m+m'}=(T\eta _m)\otimes (T'\eta _{m'})\simh\eta _n\otimes \eta _{n'}=\eta _{n+n'}.
\end{gather*}
We also have $\one \in \Ob(\modB )$.  Hence $\modB $ is a monoidal subcategory of
$\modT$.

We have $\psi _{\modb ,\modb }^{\pm 1}\eta _2=\eta _2$, hence
$\psi _{\modb ,\modb }^{\pm 1}\in \modB (\modb ^{\otimes 2},\modb ^{\otimes 2})$.  Since the object $\modb $
generates the monoid $\Ob(\modB )$, it follows that
$\psi _{\modb ^{\otimes m},\modb ^{\otimes n}}\in \modB (\modb ^{\otimes (m+n)},\modb ^{\otimes (m+n)})$.  Hence $\modB $ is a
braided subcategory of $\modT $.

For simplicity of notation, we set
\begin{equation*}
  \modB (m,n) = \modB (\modb ^{\otimes m},\modb ^{\otimes n})
\end{equation*}
for $m,n\ge 0$.  We use the similar notation for subcategories of $\modB $
defined later.

For $n\ge 0$, we have
\begin{equation*}
  \modB (0,n)=\{T\in \modT (\one ,\modb ^{\otimes n})\ver T\simh\eta _n\}.
\end{equation*}
Hence we can naturally identify $\modB (0,n)$ with the set $\BT_n$ of
isotopy classes of $n$--component bottom tangles.

\section{The subcategory $\modB _0$ of $\modB $}
\label{sec4}

In this section, we introduce a braided subcategory $\modB _0$ of $\modB $,
and give a set of generators of $\modB _0$.  The category $\modB _0$ is used
in \fullref{sec5}.

\subsection{Definition of $\modB _0$}
\label{sec:definition-0}

Let $\modB _0$ denote the subcategory of $\modB $ with $\Ob(\modB _0)=\Ob(\modB )$
and
\begin{equation*}
  \modB _0(m,n) =
  \{T\in \modB (m,n)\ver T\eta _m=\eta _n\}
\end{equation*}
for $m,n\ge 0$.  It is straightforward to check that the category $\modB _0$
is well-defined and it is a braided subcategory of $\modB $.  Note that we
have $\modB _0(0,n)=\{\eta _n\}$.

Set
\begin{align*}
  \gamma _+&=\bigl(\downarrow \otimes \psi _{\modb ,\uparrow }\psi
_{\uparrow ,\modb }\bigr)(\coev_{\downarrow} \otimes \modb )
  =\raisebox{-0.8em}{\incl{2.5em}{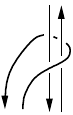}}\in \modB _0(1,2),\\
  \gamma _-&=\bigl(\downarrow \otimes \psi _{\uparrow ,\modb }^{-1}\psi
_{\modb ,\uparrow }^{-1}\bigr)(\coev_{\downarrow} \otimes \modb )
  =\raisebox{-0.8em}{\incl{2.5em}{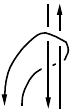}}\in \modB _0(1,2),\\
  t_{+-} &= t_{\downarrow} \otimes t_{\uparrow} ^{-1}\in \modB _0(1,1).
\end{align*}
Note that $t_{+-}$ is an isomorphism.

The purpose of this section is to prove the following theorem, which
is used in \fullref{sec5}.

\begin{theorem}
  \label{thm:9}
  As a braided subcategory of $\modB $, $\modB _0$ is generated by the
  object~$\modb $ and the morphisms
  $\mu _{\modb} ,\eta _{\modb} ,\gamma _+,\gamma _-,t_{+-},t_{+-}^{-1}$.
\end{theorem}

The proof of \fullref{thm:9} is given in
\fullref{sec:proof-theor-refthm:9}, after giving a lemma on string
links in \fullref{sec:pres-string-links}.

\subsection{Clasper presentations for string links}
\label{sec:pres-string-links}
To prove the case of ``doubled string links'' of \fullref{thm:9}
(see \fullref{sec:case-doubled-string}), we need a lemma which
presents an $n$--component string link as the result of surgery on
$1_{\downarrow ^{\otimes n}}$ along some claspers
(see Goussarov \cite{Goussarov:Y-graphs} and Habiro
\cite{Habiro:claspers}).  In this and the next
subsections (but not elsewhere in this paper), a ``clasper'' means a
``strict tree clasper of degree $1$'' in the sense of
\cite{Habiro:claspers}, ie, a clasper consisting of two disc-leaves
and one edge which looks as depicted in \fullref{fig:clasper2} (a).
\labellist\small
\pinlabel {(a)} [b] at 15 0
\pinlabel {(b)} [b] at 75 0
\pinlabel {(c)} [b] at 115 0
\pinlabel {(d)} [b] at 175 0
\pinlabel {$C$} [b] at 14 22
\pinlabel {$C'$} [b] at 115 23
\tiny
\pinlabel {\parbox{25pt}{surgery along $C$}} [b] at 43 21
\pinlabel {\parbox{26pt}{surgery along $C'$}} [b] at 143 20
\endlabellist
\FIG{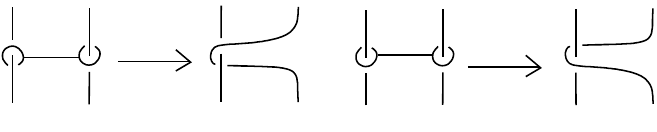}{Here each string may be replaced with parallel
strings.}{height=20mm} One can perform surgery on a clasper as
depicted in \fullref{fig:clasper2} (b), see \cite[Remark
2.4]{Habiro:claspers}.  We also use the fact that the result of
surgery on another clasper $C'$ depicted in \fullref{fig:clasper2}
(c) is as depicted in \fullref{fig:clasper2} (d).

For $n\ge 0$, let $\SL_n$ denote the submonoid of
$\modT (\downarrow ^{\otimes n},\downarrow ^{\otimes n})$ consisting of the isotopy classes of the
$n$--component framed string links.  Thus we have
\begin{equation*}
  \SL_n = \bigl\{T\in \modT (\downarrow ^{\otimes n},\downarrow ^{\otimes
n})\ver T\simh\downarrow ^{\otimes n}\bigr\}.
\end{equation*}

\begin{lemma}
  \label{lem:12}
  If $T\in \SL_n$, then there are mutually disjoint claspers
  $C_1,\ldots,C_r$ ($r\ge 0$) for $1_{\downarrow ^{\otimes n}}$ satisfying the following
  properties.
  \begin{enumerate}
  \item The tangle $T$ is obtained from $1_{\downarrow ^{\otimes n}}$ by surgery
  along $C_1,\ldots,C_r$ and framing change.
  \item $1_{\downarrow ^{\otimes n}}$ and $C_1,\ldots,C_r$ is obtained by pasting
  horizontally and vertically finitely many copies of
  the following:
  \end{enumerate}
  \begin{equation*}
    \def\size{18mm}
    \incl{\size}{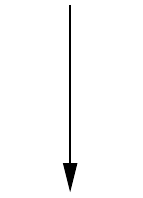}\raisebox{10mm}{,\quad}
    \labellist\small \pinlabel {edge} [t] at 19 13 \endlabellist
    \incl{\size}{}\raisebox{10mm}{,\quad}
    \labellist\small \pinlabel {edge} [t] at 7 14 \endlabellist
    \incl{\size}{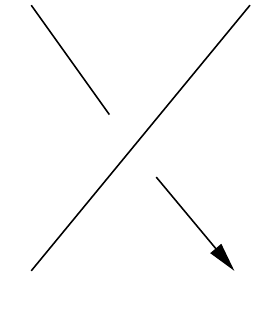}\raisebox{10mm}{,\quad}
    \labellist\small \pinlabel {edge} [t] at 65 13 \endlabellist
    \incl{\size}{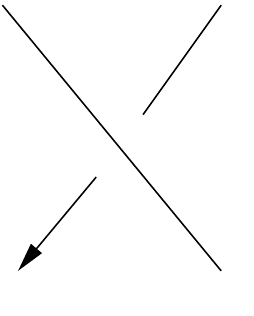}\raisebox{10mm}{,\quad}
    \labellist\small \pinlabel {edge} [t] at 23 13 \endlabellist
    \incl{\size}{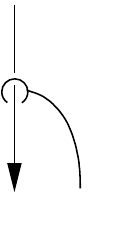}\raisebox{10mm}{,\quad}
    \labellist\small \pinlabel {edge} [t] at 23 13 \endlabellist
    \incl{\size}{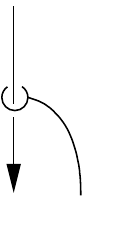}\raisebox{10mm}{,\quad}
    \labellist\small \pinlabel {edge} [b] at 15 55 \endlabellist
    \incl{\size}{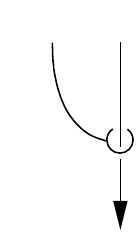}\raisebox{10mm}{.}
  \end{equation*}
\end{lemma}

\begin{proof}
  In this proof, we can ignore the framings.

  As is well known, we can express $T$ as a ``partially closed braid''
  in the sense that there is an integer $p\ge 1$ and a pure braid
  $\beta \in \modT (\downarrow ^{\otimes np},\downarrow ^{\otimes np})$ of $np$ strings such that
  \begin{equation*}
    T =
    \bigl(\downarrow ^{\otimes n}\otimes \ev_{\uparrow ^{\otimes n(p-1)}}\bigr)
    \bigl(\psi _{\downarrow ^{\otimes n(p-1)},\downarrow ^{\otimes
n}}\beta \otimes \uparrow ^{\otimes n(p-1)}\bigr)
    \bigl(\downarrow ^{\otimes n}\otimes \coev_{\downarrow ^{\otimes
n(p-1)}}\bigr),
  \end{equation*}
  see \fullref{fig:pf-lem12-1-flip} (a).
\labellist\small
\pinlabel {\tiny$\beta$} at 30 128
\pinlabel {\tiny$\beta'$} at 239 116
\pinlabel {(a)} [b] at 80 0
\pinlabel {(b)} [b] at 230 0
\pinlabel {(c)} [b] at 370 0
\pinlabel {(d)} [b] at 495 0
\endlabellist
\FIGn{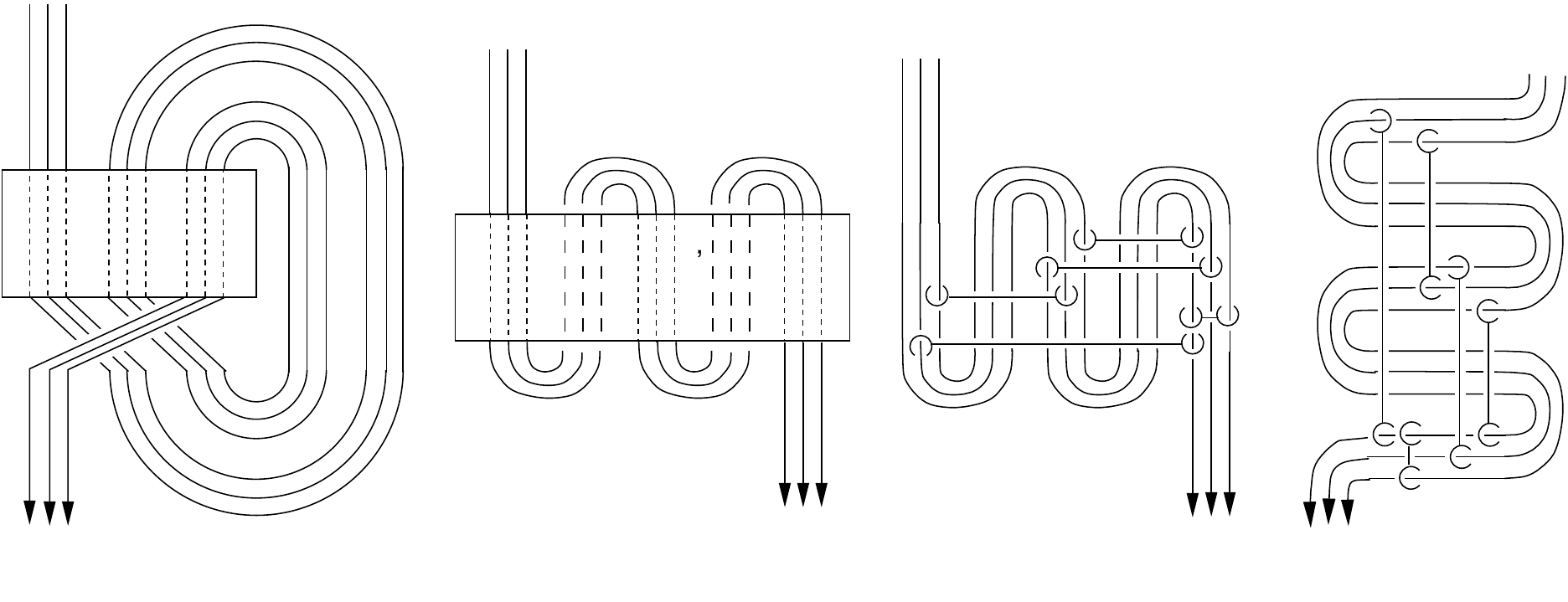}{}{height=35mm}
  By isotopy, $T$ can be
  expressed as in \fullref{fig:pf-lem12-1-flip} (b), where
  $\beta '=((t_{\downarrow ^{\otimes n}})^{\otimes (p-1)}\otimes \downarrow ^{\otimes n})\beta $ is a pure braid, and
  where the upward parts of the strings run under, and are not
  involved in, the pure braid $\beta '$.  We express $\beta '$ as the product
  of copies of generators $A_{i,j}$ ($1\le i<j\le np$) of the $np$--string
  pure braid group and their inverses.  Here $A_{i,j}$ is as depicted
  in \fullref{fig:Aij-new} (a).
\labellist\small
\pinlabel {(a)} [b] at 85 0
\pinlabel {(b)} [b] at 275 0
\pinlabel {(c)} [b] at 480 0
\tiny
\pinlabel {$A_{i,j} =$} [r] at 27 66
\pinlabel {$1$} [b] at 32 89
\pinlabel {$i$} [b] at 70 89
\pinlabel {$j$} [b] at 117 89
\pinlabel {$np$} [b] at 154 89
\pinlabel {$\cdots$} at 44 64
\pinlabel {$\cdots$} at 95 76
\pinlabel {$\cdots$} at 95 42
\pinlabel {$\cdots$} at 145 64
\pinlabel {$A_{i,j} =$} [r] at 217 66
\pinlabel {$1$} [b] at 221 89
\pinlabel {$i$} [b] at 259 89
\pinlabel {$j$} [b] at 307 89
\pinlabel {$np$} [b] at 343 89
\pinlabel {$\cdots$} at 233 64
\pinlabel {$\cdots$} at 285 76
\pinlabel {$\cdots$} at 285 42
\pinlabel {$\cdots$} at 335 64
\pinlabel {$A_{i,j}^{-1} =$} [r] at 417 66
\pinlabel {$1$} [b] at 419 89
\pinlabel {$i$} [b] at 455 89
\pinlabel {$j$} [b] at 503 89
\pinlabel {$np$} [b] at 541 89
\pinlabel {$\cdots$} at 432 64
\pinlabel {$\cdots$} at 483 76
\pinlabel {$\cdots$} at 483 42
\pinlabel {$\cdots$} at 533 64
\endlabellist
\FIGn{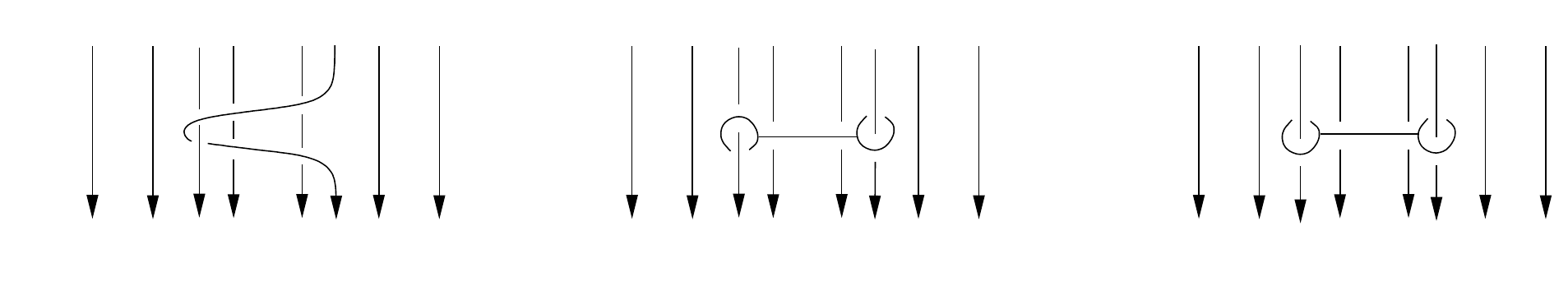}{}{height=23mm}
  (See Birman \cite{Birman} for the generators of the pure braid group.)
  Using claspers, we can express $A_{i,j}^{\pm 1}$ as depicted in
  \fullref{fig:Aij-new} (b), (c).  Let $T_0$ denote the string link
  obtained from the tangle depicted in \fullref{fig:pf-lem12-1-flip} (b) by replacing the pure braid $\beta '$ with
  $1_{\downarrow ^{\otimes np}}$.  There are claspers $C'_1,\ldots,C'_r$ ($r\ge 0$) for
  $T_0$ corresponding to the generators and inverses involved in $\beta '$
  such that surgery on $T_0$ along $C'_1,\ldots,C'_r$ yields $T$, see eg
  \fullref{fig:pf-lem12-1-flip} (c).  We can isotop
  $T_0,C'_1,\ldots,C'_r$ to the identity braid $1_{\downarrow ^{\otimes n}}$ and claspers
  $C_1,\ldots,C_r$ satisfying the desired properties, as is easily seen
  from \fullref{fig:pf-lem12-1-flip} (d).
\end{proof}

The rest of this subsection is not necessary in the rest of the paper,
but seems worth mentioned.  Let $\modS $ denote the monoidal subcategory
of $\modT $ generated by the objects $\downarrow $, $\modb $ and the following
morphisms
\begin{gather*}
  t_{\downarrow} ,  t_{\downarrow} ^{-1}\col \downarrow \rightarrow \downarrow ,\quad
  \psi _{\downarrow ,\modb }\col \downarrow \otimes \modb \rightarrow \modb \otimes \downarrow ,\quad
  \psi _{\downarrow ,\modb }^{-1}\col \modb \otimes \downarrow \rightarrow \downarrow \otimes \modb ,\\
  \delta _+= \raisebox{-0.8em}{\incl{3em}{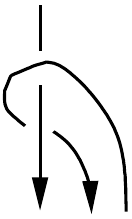}},\quad
  \delta _-= \raisebox{-0.8em}{\incl{3em}{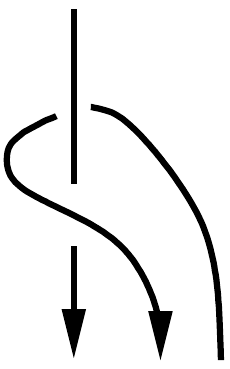}}\col \downarrow \rightarrow \downarrow \otimes \modb ,\quad
  \alpha =\raisebox{-1.0em}{\incl{3em}{}}\col \modb \otimes \downarrow \rightarrow \downarrow .
\end{gather*}

\begin{proposition}
  \label{thm:18}
  For $n\ge 0$, we have $\SL_n=\modS (\downarrow ^{\otimes n},\downarrow ^{\otimes n})$.
\end{proposition}

\begin{proof}
  The inclusion $\SL_n\subset \modS (\downarrow ^{\otimes n},\downarrow ^{\otimes n})$ easily follows from
  \fullref{lem:12}.  We prove the other inclusion
  $\modS (\downarrow ^{\otimes n},\downarrow ^{\otimes n})\subset \SL_n$.  Let $\modS _0$ denote the monoidal
  subcategory of $\modT $ generated by the objects $\downarrow $, $\modb $ and the
  morphisms $\psi _{\downarrow ,\modb }$, $\psi _{\downarrow ,\modb
}^{-1}$, $\downarrow \otimes \eta _{\modb} $, $\alpha $.  We
  can prove that $\modS _0(\downarrow ^{\otimes n},\downarrow ^{\otimes n})=\{1_{\downarrow ^{\otimes n}}\}$.  Since
  any morphism in~$\modS $ is homotopic to a morphism in $\modS _0$, the
  assertion follows.
\end{proof}

\fullref{thm:18} may be useful in studying quantum invariants
of string links.  For another approach to string links, see \fullref{sec:string-links-bottom}.

\subsection[Proof of \ref{thm:9}]{Proof of \fullref{thm:9}}
\label{sec:proof-theor-refthm:9}

In this subsection we prove \fullref{thm:9}.  Let $\modB _0'$
denote the braided subcategory of $\modB $ generated by the object $\modb $
and the morphisms
\begin{equation*}
  \eta _{\modb} ,\mu _{\modb} ,\gamma _+,\gamma _-, t_{+-},t_{+-}^{-1},
t_{\modb} ,t_{\modb} ^{-1}, \gamma '_+,\gamma '_-,
\end{equation*}
where
\begin{gather*}
  \gamma '_+=\raisebox{-0.8em}{\incl{3em}{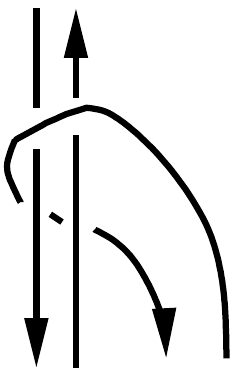}},\quad
  \gamma '_-=\raisebox{-0.8em}{\incl{3em}{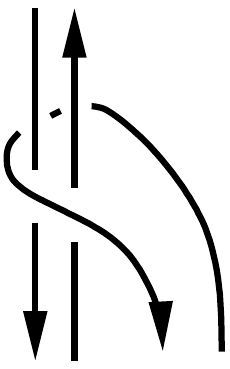}}\col \modb \rightarrow \modb ^{\otimes 2}.
\end{gather*}
Since these morphisms are in $\modB _0$, it follows that $\modB '_0$ is a
subcategory of $\modB _0$.  Since
\begin{gather*}
  t_{\modb} ^{\pm 1}=\mu _{\modb}  \gamma _{\mp} t_{+-}^{\mp1},\quad
  \gamma '_{\pm}  =\psi _{\modb ,\modb }^{\pm 1}\gamma _{\mp},
\end{gather*}
it follows that $\modB '_0$ is generated as a braided subcategory of $\modB $
by the object $\modb $ and the morphisms $\mu _{\modb} ,\eta _{\modb} ,\gamma _+,\gamma _-,
t_{+-},t_{+-}^{-1}$.  Hence it suffices to prove that any morphism in
$\modB _0$ is in $\modB _0'$.

\subsubsection{The case of doubled string links}
\label{sec:case-doubled-string}
We here prove that if $T\in \modB _0(n,n)$ is obtained from a framed
string link $T'\in \SL_n$ by doubling each component, then we have
$T\in \modB _0'(n,n)$.

By \fullref{lem:12}, there are mutually disjoint claspers
$C_1,\ldots,C_r$ ($r\ge 0$) for $1_{\modb ^{\otimes n}}$ and integers
$l_1,\ldots,l_n\in \modZ $ satisfying the following properties.
\begin{enumerate}
\item $\tilde T=T\bigl(t_{\modb} ^{l_1}\otimes \cdots\otimes t_{\modb}
^{l_n}\bigr)$ is obtained from
  $1_{\downarrow ^{\otimes n}}$ by surgery along $C_1,\ldots,C_r$.
\item $1_{\downarrow ^{\otimes n}}$ and $C_1,\ldots,C_r$ is obtained by pasting
  horizontally and vertically finitely many copies of
  the following:
  \begin{equation*}
    \def\size{4em}
    \raisebox{-1.5em}{\incl{\size}{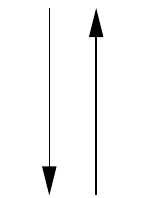}},\quad
    \raisebox{-2.0em}{\incl{\size}{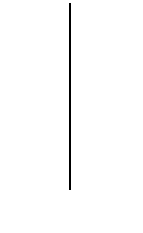}}, \quad
    \labellist\tiny \pinlabel {edge} [t] at 9 13 \endlabellist
    \raisebox{-2.0em}{\incl{\size}{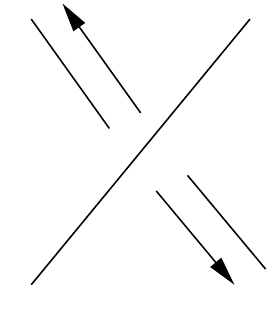}}, \quad
    \labellist\tiny \pinlabel {edge} [t] at 68 13 \endlabellist
    \raisebox{-2.0em}{\incl{\size}{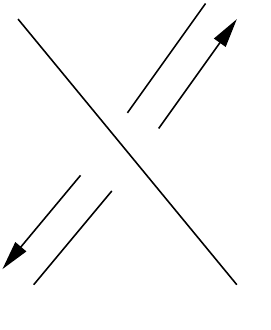}}, \quad
    \labellist\tiny \pinlabel {edge} [t] at 30 13 \endlabellist
    \raisebox{-2.0em}{\incl{\size}{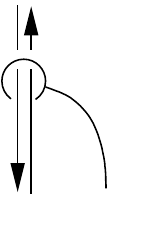}}, \quad
    \labellist\tiny \pinlabel {edge} [t] at 30 13 \endlabellist
    \raisebox{-2.0em}{\incl{\size}{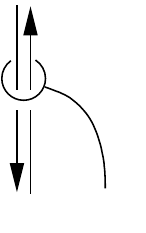}},\quad
    \labellist\tiny \pinlabel {edge} [b] at 17 54 \endlabellist
    \raisebox{-2.0em}{\incl{\size}{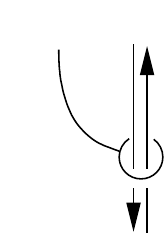}}.
  \end{equation*}
\end{enumerate}
Surgery on each $C_i$ moves the band intersecting the lower leaf of
$C_i$ and let it clasp with the band intersecting the upper leaf of
$C_i$, and we can isotop the result of surgery to the tangle
representing a morphism in $\modB _0'$ as depicted in \fullref{fig:pf-dsl}.
\labellist\tiny
\pinlabel {surgery} [b] at 78 82
\pinlabel {isotopy} [b] at 186 80
\pinlabel {surgery} [b] at 388 82
\pinlabel {isotopy} [b] at 496 80
\endlabellist
\FIGn{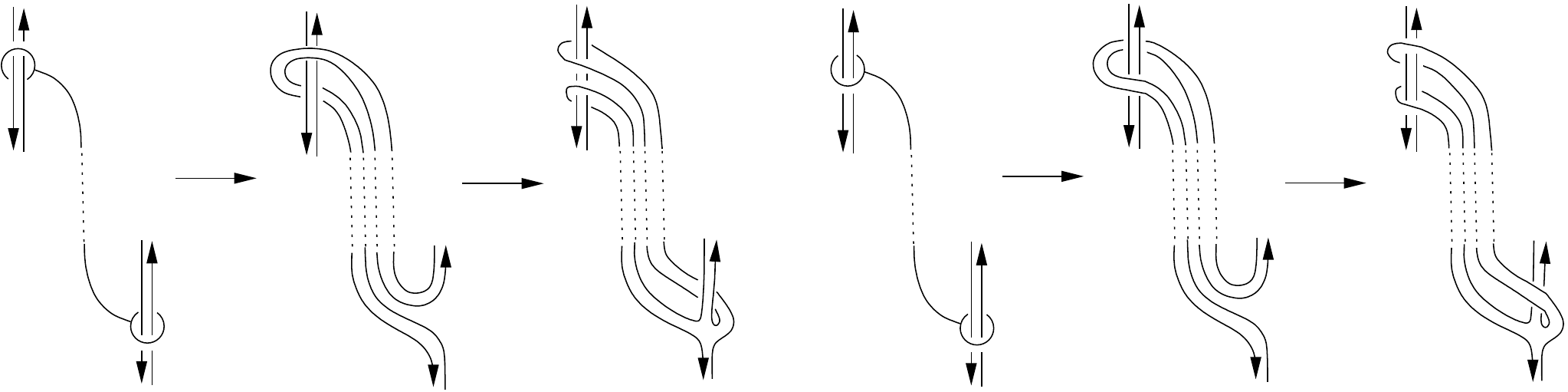}{}{height=32mm}
Hence it follows that
$\tilde T$ is in $\modB _0'$.  Since
$t_{\modb} ^{l_1}\otimes \cdots\otimes t_{\modb} ^{l_n}\in \modB _0'$, we have $T\in \modB _0'$.

\subsubsection{The general case}
\label{sec:general-case-1}

Suppose that a tangle $T$ in $[0,1]^3$ represents a morphism
$T\in \modB _0(m,n)$.  We also assume that the endpoints of $T$ are
contained in the two intervals $\bigl\{\frac12\bigl\}\times [0,1]\times \{\xi \}$, $\xi =0,1$.

For $i=1,\ldots,m$, let $c_i$ denote the interval in
$\bigl\{\frac12\bigr\}\times [0,1]\times \{1\}$ bounded by the $(2i{-}1)$st and the $2i$th
upper endpoints of $T$.  Set $c=c_1\cup \cdots \cup c_m$.  Similarly, for
$j=1,\ldots,n$, let $d_j$ denote the interval in
$\bigl\{\frac12\bigr\}\times [0,1]\times \{0\}$ bounded by the $(2j{-}1)$st and the $2j$th
lower endpoints of $T$.  Set $d=d_1\cup \cdots \cup d_n$.  Note that $T\cup c$
consists of $n$ mutually disjoint arcs $e_1,\ldots,e_n$, such that
$\partial e_j=\partial d_j$ for $j=1,\ldots,n$.  Set $e=e_1\cup \cdots \cup e_n$.

Consider $T\eta _m$, which is regarded as a tangle in $[0,1]^2\times [0,2]$,
where the lower cube $[0,1]^2\times [0,1]$ contains $T$ and the upper cube
$[0,1]^2\times [1,2]$ contains $\eta _m$.  Note that $e\subset [0,1]^2\times [0,2]$ can
be regarded as a tangle, and is equivalent to $T\eta _m$, and hence, by
the assumption, equivalent to $\eta _n$ (after identifying $[0,1]^3$ and
$[0,1]^2\times [0,2]$ in a natural way).  Hence for $j=1,\ldots,n$,
$e_j\cup d_j$ bounds a disc $D_j$ in $[0,1]^2\times [0,2]$, where
$D_1,\ldots,D_n$ are mutually disjoint.  Here each~$D_i$ is chosen so
that the framing of $e_i\cup d_i$ induced by $D_i$ is the same as the
framing of $e_i\cup d_i$ induced by that of $T$.  (Here we use the
convention that the framing of oriented tangle component is given by
the blackboard framing convention, see \fullref{fig:framing}.)  Set
$D=D_1\cup \cdots \cup D_n$.
\labellist\small
\pinlabel {$c_i$} [b] at 76 242
\pinlabel {$c_{i'}$} [b] at 142 242
\pinlabel {$c_{i''}$} [b] at 238 242
\pinlabel {$d_j$} [t] at 140 38
\tiny
\pinlabel {$\cdots$} at 30 212
\pinlabel {$\cdots$} at 110 212
\pinlabel {$\cdots$} at 190 212
\pinlabel {$\cdots$} at 280 212
\pinlabel {$\cdots$} at 80 62
\pinlabel {$\cdots$} at 200 62
\endlabellist
\FIG{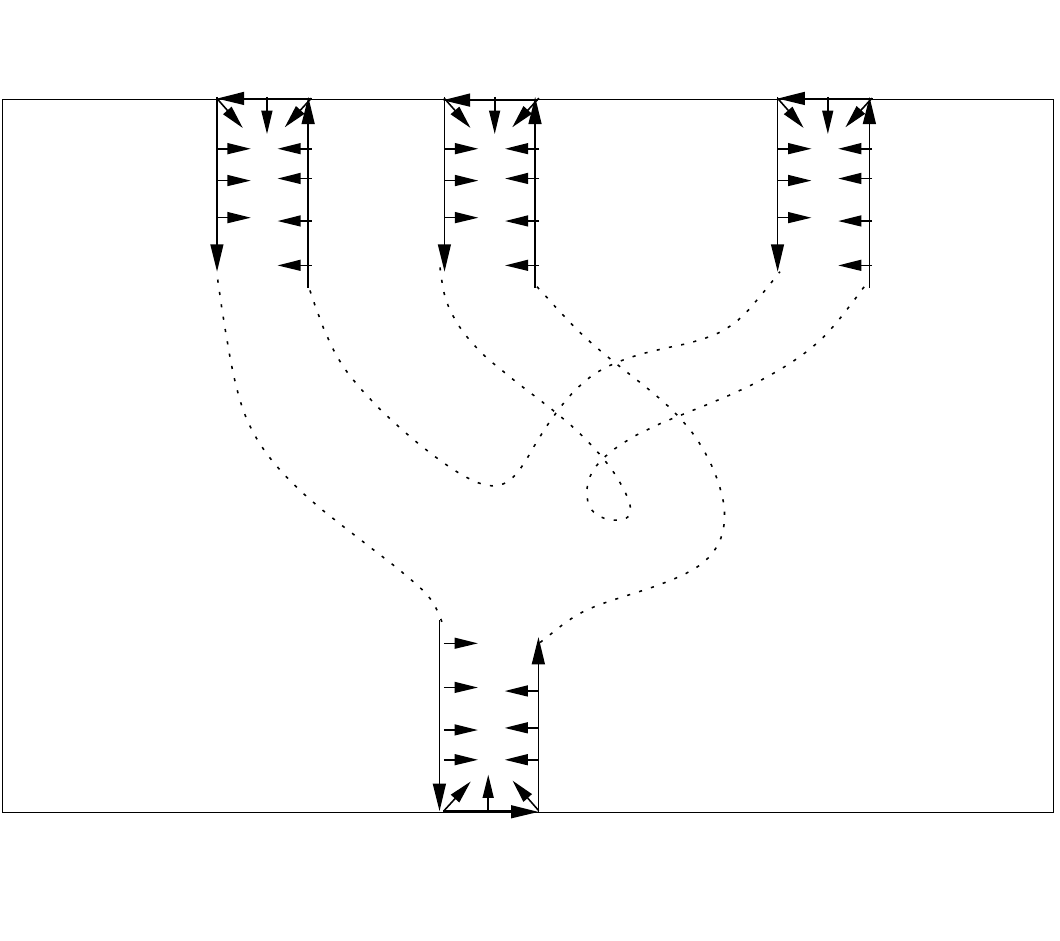}{Here the small arrows determines the framing of $e\cup d$
near $c\cup d$.}{height=40mm}

Let $\pi \col [0,1]^2\times \{1\}\rightarrow [0,1]^2$ denote the projection.  Using a
small isotopy if necessary, we may assume that for small $\epsilon >0$
we have the following.
\begin{itemize}
\item $N=\pi (c)\times [1-\epsilon ,1]$ is a regular neighborhood of $c$ in $D$.
\item  $e\setminus N\subset [0,1]^2\times [0,1-\epsilon )$.
\end{itemize}

For $i=1,\ldots,m$, let $U_i$ denote a small regular neighborhood of $c_i$
in $[0,1]^2\times \{1\}$.  Using an isotopy of $[0,1]^2\times [0,2]$ fixing
$[0,1]^3$, we can assume that for each $i=1,\ldots,m$, we have
\begin{equation*}
  \bigl(\pi (U_i)\times (1,2]\bigr)\cap D = \pi (U_i)\times \{p_{i,1},\ldots,p_{i,l_i}\},
\end{equation*}
where $1<p_{i,1}<\cdots<p_{i,l_i}<2$, $l_i\ge 0$.  Define a piecewise-linear
homeomorphism $f\col [0,2]\rightarrow [0,1]$ by
\begin{equation*}
  f(t) =
  \begin{cases}
    t, &\text{if $0\le t\le 1-\epsilon $},\\
    \frac{\epsilon t+(1-\epsilon )}{1+\epsilon }, &\text{if $1-\epsilon \le t\le 2$}.
  \end{cases}
\end{equation*}
Define $\tilde{f}\col [0,1]^2\times [0,2]\rightarrow [0,1]^3$ by
$\tilde{f}(x,y,t)=(x,y,f(t))$ for $x,y\in [0,1]$, $t\in [0,2]$.  (Thus
$\tilde{f}$ fixes $[0,1]^2\times [0,1-\epsilon ]$, and maps $[0,1]^2\times [1-\epsilon ,2]$
onto $[0,1]^2\times [1-\epsilon ,1]$ linearly.)  Set
\begin{equation*}
  D' = \tilde{f}(D\setminus N)\cup N.
\end{equation*}
Note that $D'$ is a union of $n$ immersed discs whose only
singularities are ribbon singularities
\begin{equation*}
  \pi (c_i)\times \{f(p_{i,k})\}\quad \text{for $1\le i\le m$, $1\le k\le l_i$}.
\end{equation*}
For $j=1,\ldots,n$, let $D'_j\subset D'$ denote the unique immersed disc
containing $d_j$.  Note that $\partial D'_j=e_j\cup d_j$.

We prove the assertion by induction on the number $l=l_1+\cdots+l_m$ of
ribbon singularities.  There are two cases.

\textbf{Case 1}\qua $l=0$, ie, there are no ribbon singularities in $D'$.
  In this case, we claim that $T$ is equivalent to
\begin{equation*}
  (\mu _{\modb} ^{[m_1]}\otimes \cdots\otimes \mu _{\modb} ^{[m_n]})\sigma \beta ,
\end{equation*}
where $m_j$ is the number of arcs $c_i$ contained in $D'_j$ for
$j=1,\ldots,n$, $\sigma $ is a doubled braid, and $\beta $ is a doubled string
link.  This claim can be proved as follows.  Choose points
$x_i\in \opint c_i$ for $i=1,\ldots,m$, and $y_j\in \opint d_j$ for
$j=1,\ldots,n$.  For each $i=1,\ldots,m$, let $b_i$ denote a proper arc in
$D'$ which connects $x_i$ and $y_{j(i)}$, where $j(i)$ is such that
$x_i$ and $y_j$ are in the same component of~$D'$.  Set
$b=b_1\cup \cdots \cup b_m$.  We may assume that for $i\neq i'$ the intersection
$b_i\cap b_{i'}$ is either empty if $y(i)\neq y(i')$, or the common
endpoint $y_j$ if $y(i)=y(i')$.  By small isotopy, we may assume that
for sufficiently small $\epsilon '>0$ the intersection
$([0,1]^2\times [0,\epsilon '])\cap b$ consists of ``stars'' rooted at $y_1,\ldots,y_n$.
Here each star at~$y_j$ consists of $m_j$ line segments, and each
$b'_i=([0,1]^2\times [\epsilon ',1])\cap b_i$ is an arc.  Let $D''$ be a sufficiently
small regular neighborhood of $c\cup b\cup d$ in $D'$ so that
$D''\cap ([0,1]^2\times [\epsilon ',1])$ consists of mutually disjoint bicollar neighborhoods
of $b'_1,\ldots,b'_m$.  There is an isotopy of $[0,1]^3$ which fixes
$c\cup b\cup d$ and deforms $D'$ to $D''$.  Set $T'=\partial D''\setminus (\opint c\cup \opint
d)$, which is a tangle isotopic to $T$.  By an isotopy fixing
$[0,1]^2\times \{\epsilon '\}$ as a set, we may assume that
$([0,1]^2\times [0,\epsilon '])\cap T'$ is a tangle of the form
$\mu _{\modb} ^{[m_1]}\otimes \cdots\otimes \mu _{\modb} ^{[m_n]}$, and $([0,1]^2\times [\epsilon ',1])\cap T'$ is a
a tangle of the form $\sigma \beta $, as desired.  Hence we have the claim.

It follows from \fullref{sec:case-doubled-string} that $\beta $ is a
morphism in $\modB '_0$.  Obviously, $\sigma $ and
$\mu _{\modb} ^{[m_1]}\otimes \cdots\otimes \mu _{\modb} ^{[m_n]}$ are morphisms in $\modB '_0$.  Hence
$T$ is in $\modB '_0$.

\textbf{Case 2}\qua  $l\ge 1$, ie, there is at least one ribbon
singularity in $D'$.
Suppose $i\in \{1,\ldots,m\}$ with $l_i\ge 1$.  Using
isotopy of $[0,1]^3$ fixing $(\pi (c)\times [1-\epsilon ,1])\cup d$, we see that $T$
is equivalent to
\begin{equation}
  \label{eq:1}
  T'(\modb ^{\otimes (i-1)}\otimes \gamma _{\pm} \otimes \modb ^{\otimes (n-i)}),
\end{equation}
where $T'\in \modT (\modb ^{\otimes (m+1)},\modb ^{\otimes n})$.  Here the ribbon singularity
$\pi (c_i)\times \{f(p_{i,l_i})\}$ of~$D$ is isotoped to the obvious ribbon
singularity involved in the copy of $\gamma _{\pm} $ in \eqref{eq:1}.  Since we
have
\begin{gather*}
  T'\eta _{m+1}=T'(\modb ^{\otimes (i-1)}\otimes \gamma _{\pm} \otimes \modb ^{\otimes (n-i)})\eta _m
  =T\eta _m=\eta _n,
\end{gather*}
it follows that $T'\in \modB _0(m+1,n)$.  Note that $T'$ bounds ribbon discs
with less singularities than $T$ by $1$.  By the induction assumption,
it follows that $T'$ is in $\modB _0'$.  Hence we have $T\in \modB _0'(m,n)$.

This completes the proof of \fullref{thm:9}.

\section{Local moves}
\label{sec5}

In this section, we explain how the category $\modB $ can be used in the
study of local moves on links and tangles.

\subsection{Monoidal relations and monoidal congruences}
\label{sec:monoidal-relations}
In this subsection, we recall the notions of monoidal relations and
monoidal congruences in monoidal categories.

Two morphisms in a monoidal category $\modM $ are said to be {\em
compatible} if they have the same source and the same target.

A {\em monoidal relation} in a monoidal category $\modM $ is a binary
relation $R\subset \Mor(\modM )\times \Mor(\modM )$ on $\Mor(\modM )$ satisfying the
following conditions.
\begin{enumerate}
\item If $(f,f')\in R$, then $f$ and $f'$ are compatible.
\item If $(f,f')\in R$, $a\in \Ob(\modM )$, then $(a\otimes f,a\otimes f')\in R$ and
  $(f\otimes a,f'\otimes a)\in R$.
\item If $(f,f')\in R$, $g\in \Mor(\modM )$ and $\target(f)=\source(g)$
  (resp. $\source(f)=\target(g)$), then $(gf,gf')\in R$
  (resp. $(fg,f'g)\in R$).
\end{enumerate}

For any relation $X\subset \Mor(\modM )\times \Mor(\modM )$ satisfying the condition (1)
above, there is the smallest monoidal relation $R_X$ containing $X$,
which is called the monoidal relation in $\modM $ generated by $X$.  If
$X=\{(f,f')\}$, then $R_X$ is also said to be generated by the pair
$(f,f')$.

Suppose $f,f'\in \modM (a,b)$ and $g,g'\in \modM (c,d)$ with $a,b,c,d\in \Ob(\modM )$.
Then $g$ and $g'$ are related by the monoidal relation generated by
$(f,f')$, if and only if there are $z,z'\in \Ob(\modM )$ and morphisms
$h_1\in \modM (c,z\otimes a\otimes z')$, $h_2\in \modM (z\otimes b\otimes z',d)$ such that
\begin{equation}
  \label{eq:4}
  g = h_2(z\otimes f\otimes z')h_1,\quad g' = h_2(z\otimes f'\otimes z')h_1.
\end{equation}

\begin{lemma}
  \label{lem:6}
  Let $\modM $ be a braided category and let $(f,f')\in \modM (\one ,b)$ with
  $b\in \Ob(\modM )$.  Then $g,g'\in \modM (c,d)$ ($c,d\in \Ob(\modT )$) are related by
  the monoidal relation generated by $(f,f')$ if and only if there is
  $h\in \modM (c\otimes b,d)$ such that
  \begin{equation*}
    g=h(c\otimes f),\quad g'=h(c\otimes f').
  \end{equation*}
\end{lemma}

\begin{proof}
  The ``if'' part is obvious.  We prove the ``only if'' part.
  By assumption, there are $z,z'\in \Ob(\modM )$ and $h_1\in \modM (c,z\otimes z')$,
  $h_2\in \modM (z\otimes b\otimes z',d)$ satisfying \eqref{eq:4}.
  Set $h=h_2(z\otimes \psi _{z',b})(h_1\otimes b)$.  Then we have the assertion.
\end{proof}

A {\em monoidal congruence}, also called {\em four-sided congruence},
in a monoidal category $\modM $ is a monoidal relation in $\modM $ which is an
equivalence relation.  If $\sim$ is a monoidal congruence in a
monoidal (resp. braided) category $\modM $, then the quotient category
$\modM /\sim$ is equipped with a monoidal (resp. braided) category
structure induced by that of $\modM $.

\begin{example}
  \label{r18}
  The notion of homotopy (see \fullref{sec:subcategories}) for
  morphisms in $\modT $ is a monoidal congruence in $\modT $ generated by
  $\{(\psi _{\downarrow ,\downarrow },\psi _{\downarrow ,\downarrow
}^{-1}),(1_{\downarrow} ,t_{\downarrow} )\}$.
\end{example}

\subsection{Topological and algebraic definitions of local
  moves}
\label{sec:topol-defin-algebr}
In this subsection, we first recall a formulation of local moves on
links and tangles, and then we reformulate it in the setting of the
category $\modT $.

Informally, a ``local move'' is an operation on a link (or a tangle)
which replaces a tangle in a link (or a tangle) contained in a
$3$--ball $B$ with another tangle.  The following is a precise
definition of the notion of local moves.

\begin{definition}
  \label{def:3}
  Two tangles $t$ and $t'$ in a $3$--ball $B$ are said to be {\em
  compatible} if we have $\partial t=\partial t'$ and $t$ and $t'$ have the same
  framings and the same orientations at the endpoints.  Let $(t,t')$
  be a compatible pair of tangles in $B$.  For two compatible tangles
  $u$ and~$u'$ in another $3$--ball $D$, we say that $u$ and~$u'$ are
  {\em $(t,t')$--related}, or $u'$ is obtained from $u$ by a {\em
  $(t,t')$--move}, if there is an orientation-preserving embedding
  $f\col B\hookrightarrow \opint D$ and a tangle $u''$ in $D$ isotopic to
  $u'$ such that
  \begin{equation}
    \label{eq:2}
    f(t)=u\cap f(B),\quad
    f(t')=u''\cap f(B),\quad
    u\setminus \opint f(B)=u''\setminus \opint f(B),
  \end{equation}
  where the orientations and framings are the same in the two
  $1$--submanifolds in each side of these three identities.
\end{definition}

Now we give an algebraic formulation of local moves.

\begin{definition}
  \label{def:4}
  Let $(T,T')$ be a compatible pair of morphisms in $\modT $.  Let
  $R_{(T,T')}$ be the monoidal relation generated by the pair
  $(T,T')$.  For two morphisms $U$ and $U'$ in $\modT $, we say that $U$
  and $U'$ are {\em $(T,T')$--related}, or $U'$ is obtained from $U$ by
  a {\em $(T,T')$--move}, if $(U,U')\in R_{(T,T')}$.
\end{definition}

The following shows that we can reduce the study of local moves
defined by compatible pairs of tangles in a $3$--ball $B$ to the study
of local moves defined by compatible pairs of morphisms in $\modT $.

\begin{proposition}
  \label{lem:2}
 Let $(t,t')$ be a compatible pair of tangles in a $3$--ball $B$, and
  let $C$ be a compatibility class of tangles in another $3$--ball $D$.  Then
  there is a compatible pair $(T,T')$ of morphisms in $\modT $ and an
  orientation-preserving homeomorphism $g\col D\cong[0,1]^3$ such that
  two tangles $u,u'\in C$ are $(t,t')$--related if and only if the
  morphisms in $\modT $ represented by $g(u)$ and $g(u')$ are
  $(T,T')$--related as morphisms in $\modT $.
\end{proposition}

\begin{proof}
  We choose an orientation-preserving homeomorphism $h\col B\cong[0,1]^3$
  such that we have $h(\partial t)=h(\partial t')\subset
\bigl\{\frac12\bigr\}\times (0,1)\times \{0,1\}$ so
  that the tangles $h(\partial t),h(\partial t')\subset [0,1]^3$ represents (compatible)
  morphisms in $\modT $.  Similarly, we choose an orientation-preserving
  homeomorphism $g\col D\cong[0,1]^3$ such that for any $u\in C$ we have
  $g(\partial u)\subset \bigl\{\frac12\bigr\}\times (0,1)\times \{0,1\}$ so that for $u\in C$ the tangle
  $g(u)\subset [0,1]^3$ represents a morphism in $\modT $.  Set
  $T=h(t),T'=h(t')\subset [0,1]^3$, which are regarded as morphisms in $\modT $.

  It is easy to verify the ``if'' part.  We prove the ``only if'' part
  below.  Suppose $u,u'\in C$ are $(t,t')$--related.  By the definition,
  there is an orientation-preserving embedding $f\col B\hookrightarrow D$
  and a tangle $u''\in C$ isotopic to $u'$ satisfying \eqref{eq:2}.  We
  can assume without loss of generality that $u''=u'$.  Set
  \begin{equation*}
    f'=gfh^{-1}\col [0,1]^3\hookrightarrow[0,1]^3.
  \end{equation*}
  Then $g(u)$ and $g(u')$, as tangles in $[0,1]^3$, are related by a
  $(T,T')$--move via~$f'$.  By applying an appropriate
  self-homeomorphism of $[0,1]^3$ fixing boundary to both $g(u)$ and
  $g(u')$, we have in $\modT $
  \begin{gather*}
    g(u) = W_2(z\otimes T\otimes z')W_1,\quad
    g(u') = W_2(z\otimes T'\otimes z')W_1,
  \end{gather*}
  where $z,z'\in \Ob(\modT )$ and
  \begin{gather*}
    W_1\in \modT (\source(g(u)),z\otimes \source(T)\otimes z'),\\
    W_2\in \modT (z\otimes \target(T)\otimes z',\target(g(u))).  
  \end{gather*}
  This shows that $g(u)$
  and $g(u')$ are $(T,T')$--related as morphisms in $\modT $.
\end{proof}

\begin{definition}
  Two compatible pairs $(T_1,T_1')$ and $(T_2,T_2')$ of morphisms
  in~$\modT $ are said to be {\em equivalent} to each other if there is an
  orientation-preserving self-homeomorphism $\phi $ (not necessarily
  fixing the boundary) of the cube $[0,1]^3$ such that $\phi (T_1)=T_2$
  and $\phi (T_1')=T_2'$.
\end{definition}

Note that if $(T_1,T_1')$ and $(T_2,T_2')$ are equivalent pairs of
mutually compatible morphisms in $\modT $, then the notions of
$(T_1,T_1')$--move and $(T_2,T_2')$--move are the same, ie, two
tangles $U$ and $U'$ are $(T_1,T_1')$--related if and only if $U$ and
$U'$ are $(T_2,T_2')$--related.

\subsection{Local moves defined by pairs of bottom tangles}
\label{sec:local-moves-defined}
In the following we restrict our attention to local moves defined by
pairs of mutually homotopic tangles $T$ and $T'$ consisting only of
arc components.  Let us call such a local move an {\em arc local
move}.  Arc local moves fit nicely into the setting of the category
$\modB $.

The following implies that, to study the arc local moves, it suffices
to study the local moves defined by pairs of bottom tangles.

\begin{proposition}
  \label{thm:3}
  Let $(T,T')$ be a pair of mutually homotopic morphisms in~$\modT $, each
  consisting of $n$ arc components and no circle components.  Then
  there is a pair $(T_1,T_1')$ of mutually homotopic $n$--component
  bottom tangles which is equivalent to $(T,T')$.
\end{proposition}

\begin{proof}
  Set $a=\source(T)$ and $b=\target(T)$.  There is a (not unique)
  framed braid $\beta \in \modT (b\otimes a^*,\modb ^{\otimes n})$ such that
\begin{equation*}
  T_1=\beta (T\otimes a^*)\coev_a\quad \text{and}\quad T'_1=\beta (T'\otimes a^*)\coev_a
\end{equation*}
are bottom
  tangles.  Clearly, the two pairs $(T,T')$ and $(T_1,T'_1)$ are
  equivalent, and $T_1$ and $T_1'$ are homotopic to each other.
  Hence we have the assertion.
\end{proof}

We are in particular interested in arc local moves on bottom tangles.
The following theorem implies that the study of arc local moves on
tangles in $\modB $ is reduced to the study of monoidal relations in $\modB $
generated by pairs of bottom tangles.

\begin{theorem}
  \label{thm:5}
  For $T,T'\in \BT_n$ and $U,U'\in \modB (k,l)$, the following conditions
  are equivalent.
  \begin{enumerate}
  \item $U$ and $U'$ are $(T,T')$--related.
  \item $U$ and $U'$ are ``$(T,T')$--related in $\modB $'', ie, related
    by the monoidal relation in $\modB $ generated by $(T,T')$.
  \item There is a morphism $W\in \modB (k+n,l)$ such that
    \begin{equation}
      \label{eq:3}
      U=W(\modb ^{\otimes k}\otimes T),\quad U'=W(\modb ^{\otimes k}\otimes T').
    \end{equation}
  \end{enumerate}
\end{theorem}

\begin{proof}
  By \fullref{lem:6}, (2) and (3) are equivalent.  Obviously, (2)
  implies~(1).  We show that (1) implies (3).  By assumption, $U$ and
  $U'$ are $(T,T')$--related.  By \fullref{lem:6}, there is
  $W\in \modT (\modb ^{\otimes (k+n)},\modb ^{\otimes l})$ satisfying \eqref{eq:3}.  We have
  \begin{equation*}
    W\eta _{k+n}
    =W(\modb ^{\otimes k}\otimes \eta _n)\eta _k
    \simh W(\modb ^{\otimes k}\otimes T)\eta _k
    = U\eta _k
    \simh \eta _l,
  \end{equation*}
  where we used the fact that $T$ is a bottom tangle and $U$ is a
  morphism in $\modB $.  Hence we have $W\in \modB (k+n,l)$.
\end{proof}

\subsection{Admissible local moves}
\label{sec:admiss-local-moves}

An $n$--component tangle $t$ in a $3$--ball $B$ is said to be {\em
admissible} if the pair $(B,t)$ is homeomorphic to the pair
$([0,1]^3,\eta _n)$.  (In the literature, such a tangle is sometimes
called ``trivial tangle'', but here we do not use this terminology,
since it may give an impression that a tangle is equivalent to a
``standard'' tangle such as $\eta _n$.)

A compatible pair $(t,t')$ of tangles in a $3$--ball $B$ is called {\em
  admissible} if both $t$ and $t'$ are admissible.  A local move
  defined by admissible pair is called an admissible local move.  In
  this subsection, we translate some well-known properties of
  admissible local moves into our category-theoretical setting.

It follows from the previous subsections that, to study admissible
local moves on morphisms in $\modB $, it suffices to study the monoidal
relations in $\modB $ generated by pairs of admissible bottom tangles with
the same number of components.

For $n\ge 0$, let $\ABT_n$ denote the subset of $\BT_n$ consisting of
admissible bottom tangles.  We set
\begin{equation*}
  \ABT =\bigcup_{n\ge 0}\ABT_n\subset \BT.
\end{equation*}

\begin{lemma}
  \label{lem:1}
  If $T,T'\in \ABT_n$, then there is $V\in \ABT_n$ such that the two pairs
  $(T,T')$ and $(\eta _n,V)$ are equivalent.
\end{lemma}

\begin{proof}
  Since $T\in \ABT_n$, there is a framed pure braid
  $\beta \in \modT (\modb ^{\otimes n},\modb ^{\otimes n})$ of $2n$--strings such that $\beta T=\eta _n$.
  Setting $V=\beta T'$, we can easily verify the assertion.
\end{proof}

\fullref{lem:1} above implies that, to study admissible local moves
on morphisms in $\modB $, it suffices to study the admissible local moves
defined by pairs $(\eta _n,T)$ for $T\in \ABT_n$.  Hence it is useful to
make the following definition.  If two tangles $U$ and $U'$ in $\modB $
are $(\eta _n,T)$--related, then we simply say that $U$ and $U'$ are {\em
$T$--related}, or $U'$ is obtained from $U$ by a {\em $T$--move}.

Now we consider sequences of admissible local moves.

\begin{proposition}
  \label{thm:6}
  Let $T_1,\ldots,T_r\in \ABT$, $r\ge 0$, and let $U,U'\in \modB (k,l)$ be two
  morphisms in $\modB $.  Then the following conditions
  are equivalent.
  \begin{enumerate}
  \item There is a sequence $U_0=U,U_1,\ldots,U_r=U'$ of morphisms in $\modB $
    from $U$ to $U'$ such that, for $i=1,\ldots,r$, the tangles $U_{i-1}$
    and $U_i$ are $T_i$--related.
  \item $U$ and $U'$ are $(T_1\otimes \cdots\otimes T_r)$--related.
  \end{enumerate}
\end{proposition}

\begin{proof}
  Obviously, (2) implies (1).  We show that (1) implies (2).  It is
  well known (see, for example, \cite[Lemma 3.21]{Habiro:claspers})
  that if there is a sequence from a tangle $U$ to another tangle $U'$
  of admissible local moves, then $T'$ can be obtained from $T$ by
  simultaneous application of admissible local moves of the same types
  as those appearing in the sequence.  Hence, after suitable isotopy
  of $[0,1]^3$ fixing the boundary, there are mutually disjoint small cubes
  $C_1,\ldots,C_r$ in $[0,1]^3$ such that
  \begin{itemize}
  \item for $i=1,\ldots,r$, the tangle $C_i\cap U$ in $C_i$ is equivalent to
  $\eta _{n_i}$, where we set $n_i=|T_i|$,
  \item the tangle obtained from $U$ by replacing the copy of
  $\eta _{n_i}$ in $C_i$ with a copy of $T_i$ for all $i=1,\ldots,r$ is
  equivalent to $U'$.
  \end{itemize}
  Using an isotopy of $[0,1]^3$ fixing the boundary which move the
  cubes $C_1,\ldots,C_r$ to the upper right part of $[0,1]^3$, we can
  express $U$ and $U'$ as
  \begin{equation*}
    U = W(\modb ^{\otimes k}\otimes \eta _{n_1+\cdots+n_r}), \quad
    U' =W(\modb ^{\otimes k}\otimes T_1\otimes \cdots\otimes T_r),
  \end{equation*}
  where $W\in \modT (\modb ^{\otimes (k+n_1+\cdots+n_r)},\modb ^{\otimes l})$.  One can easily
  verify $W\eta _{k+n_1+\cdots+n_r}\simh \eta _l$, hence
  $W\in \modB (k+n_1+\cdots+n_r,l)$.  Hence we have the assertion.
\end{proof}

In the study of local moves, it is often useful to consider the
relations on tangles defined by several types of moves.  Let $M\subset \ABT$
be a subset.  For two tangles $U$ and $U'$ in $\modB $, we say that $U$
and $U'$ are {\em $M$--related}, or $U'$ is obtained from~$U$ by an
{\em $M$--move}, if there is  $T\in M$ such that $U$ and $U'$
are $T$--related.

For $M\subset \ABT$, let $M^*$ denote the subset of $\ABT$ of the form
$T_1\otimes \cdots\otimes T_r$ with $T_i\in M$ for $i=1,\ldots,r$, $r\ge 0$.  Note that
$M^*\subset \ABT$.  The following immediately follows from \fullref{thm:6}.

\begin{proposition}
  \label{lem:7}
  Let $M\subset \ABT$.  Then $U,U'\in \modB (k,l)$ are related by a finite
  sequence of $M$--moves if and only if $U$ and $U'$ are $M^*$--related.
\end{proposition}

For $M\subset \ABT$, the {\em $M$--equivalence} is the equivalence
relation on tangles generated by the $M$--moves.  Note that, for
morphisms in $\modB $, the $M$--equivalence is the same as the
monoidal congruence in $\modB $ generated by the set
$\{(\eta _{|T|},T)\ver T\in M\}$.

A subset $M\subset \ABT$ is said to be {\em inversion-closed} if for each
$T\in M$, there is a sequence of $M$--moves from $T$ to $\eta _{|T|}$.
In this case, two tangles $U$ and $U'$ are $M$--equivalent if and
only if there is a sequence of $M$--moves from $U$ to $U'$.

Given any subset $M\subset \ABT$, one can construct an inversion-closed
subset $M'\subset \ABT$ such that any two morphisms $U$ and $U'$ in $\modB $ are
$M$--equivalent if and only if there is a sequence from $U$ to $U'$ of
$M'$--moves.  For example, $M'=M\cup \{\bar{T}\ver T\in M\}$ satisfies this
condition, where $\bar{T}=\beta ^{-1}\eta _n$ with $n=|T|$ and $\beta \in \modB (n,n)$
a (not unique) framed pure braid such that $T=\beta \eta _n$.  (Note that the
pair $(\bar{T},\eta _n)$ is equivalent to $(\eta _{\modb} ,T)$.)

By \fullref{lem:7}, we have the following.

\begin{proposition}
  \label{lem:4}
  Let $M\subset \ABT$ be inversion-closed.  Then $U,U'\in \modB (k,l)$ are
  $M$--equivalent if and only if they are $M^*$--related.
\end{proposition}

\subsection{Tangles obtained from $\eta _n$ by an admissible local move}
\label{sec:tangl-obta-from}

Let $\cB_0$ denote the braided subcategory of $\modB _0$ generated by the
object $\modb $ and the morphisms $\mu _{\modb} ,\eta _{\modb} ,\gamma _+,\gamma _-$, and let $\T$
denote the monoidal subcategory of $\modB _0$ generated by the object~$\modb $
and the morphisms $t_{+-}$ and $t_{+-}^{-1}$.  It is easy to see that
for $n\ge 0$ we have
\begin{equation*}
  \T(n,n) = \{t_{+-}^{k_1}\otimes \cdots\otimes t_{+-}^{k_n}\ver k_1,\ldots,k_n\in \modZ \},
\end{equation*}
and $\T(m,n)$ is empty if $m\neq n$.

\begin{lemma}
  \label{lem:8}
  For any morphism $T\in \modB _0(m,n)$ with $m,n\ge 0$, there are
  $T'\in \cB_0(m,n)$ and $T''\in \T(m,m)$ such that $T=T'T''$.  (The
  decomposition $T=T'T''$ is not unique.)
\end{lemma}

\begin{proof}
  Using \fullref{thm:9} and the identities
  \begin{gather*}
    t_{+-}^k\mu _{\modb} =\mu _{\modb} (t_{+-}^k\otimes t_{+-}^k),\quad
    t_{+-}^k\eta _{\modb} =\eta _{\modb} ,\quad
    (t_{+-}^k\otimes t_{+-}^l)\gamma _{\pm} =\gamma _{\pm} t_{+-}^l,
  \end{gather*}
  for $k,l\in \modZ $, we can easily prove the assertion.
\end{proof}

\begin{theorem}
  \label{thm:4}
  Let $T\in \ABT_m$ and $U\in \BT_n$.  Then $U$ is obtained from $\eta _n$ by
  one $T$--move if and only if there is $W\in \cB_0(m,n)$ such that
  $U=WT$.
\end{theorem}

\begin{proof}
  The ``if'' part is obvious.  We prove the ``only if'' part.  Suppose
  that $\eta _n$ and $U\in \BT_n$ are $T$--related.  By \fullref{thm:5},
  there is $W'\in \modB (m,n)$ such that $\eta _n=W'\eta _m$ and $U=W'T$.  The
  first identity means that $W'\in \modB _0(m,n)$.  By \fullref{lem:8}, we
  have $W'=WV$, where $W\in \cB_0(m,n)$ and $V\in \T(m,m)$.  It is easy to
  see that $VT=T$.  Hence we have $U=W'T=WVT=WT$.
\end{proof}

For $M\subset \ABT$, a bottom tangle $T\in \BT_n$ is said to be {\em
$M$--trivial} if $T$ is $M$--equivalent to $\eta _n$.  The following
immediately follows from \fullref{lem:4} and \fullref{thm:4}.

\begin{corollary}
  \label{thm:19}
  Let $M\subset \ABT$ be inversion-closed, and let $U\in \BT_n$.  Then $U$ is
  $M$--trivial if and only if there are $T\in M^*$ and $W\in \cB_0(|T|,n)$
  such that $U=WT$.
\end{corollary}

\subsection{Generators of $\modB $}
\label{sec:generators--1}

Here we use the results in the previous subsections to obtain a simple
set of generators of $\modB $.  Define morphisms $v_{\pm} \in \BT_1$ and
$c_{\pm} \in \BT_2$ by
\begin{equation*}
  v_{\pm} =(t_{\downarrow} ^{\mp1}\otimes \uparrow )\eta _{\modb} ,\quad
  c_{\pm} =(\downarrow \otimes (\psi _{\downarrow ,\uparrow }\psi
_{\uparrow ,\downarrow })^{\pm 1}\otimes \uparrow )(\eta _{\modb} \otimes
\eta _{\modb} ).
\end{equation*}
Graphically, we have
\begin{equation*}
  v_+ = \raisebox{-1.0em}{\incl{2.5em}{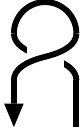}},\qquad
  v_- = \raisebox{-1.0em}{\incl{2.5em}{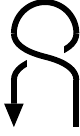}},\qquad
  c_+ = \raisebox{-1.0em}{\incl{2.5em}{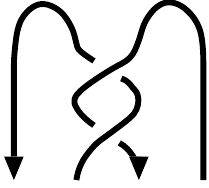}},\qquad
  c_- = \raisebox{-1.0em}{\incl{2.5em}{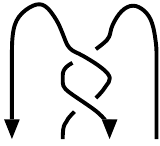}}.
\end{equation*}
Note that a $v_{\pm} $--move is change of framing by $1$, and a $c_{\pm} $--move
is a crossing change.  Hence two morphisms in $\modB $ are homotopic if
and only if they are $\{v_+,v_-,c_+,c_-\}$--equivalent.

\begin{theorem}
  \label{thm:2}
  As a braided subcategory of $\modT $, $\modB $ is generated by the
  object~$\modb $ and the morphisms $\mu _{\modb} ,\eta _{\modb} ,v_+,v_-,c_+,c_-$.
\end{theorem}

\begin{proof}
  Let $\modB '$ denote the braided subcategory of $\modT $ generated by the
  object~$\modb $ and the morphisms $\mu _{\modb} ,\eta _{\modb} ,v_+,v_-,c_+,c_-$.  It
  suffices to show that any tangle $U'\in \modB (m,n)$ is a morphism in $\modB '$.

  Choose a tangle $U\in \modB _0(m,n)$ which is homotopic to $U'$, ie,
  $\{v_+,v_-,c_+,c_-\}$--equivalent to $U'$.  Since
  $\{v_+,v_-,c_+,c_-\}$ is an inversion-closed subset of $\ABT$,
  \fullref{thm:5} and \fullref{lem:4} imply that
  there are $T\in \{v_+,v_-,c_+,c_-\}^*$ and $W\in \modB (m+|T|,n)$ such that
  \begin{eqnarray}
    \label{eq:80}
    U &=& W(\modb ^{\otimes m}\otimes \eta _{|T|}),\\
    \label{eq:81}
    U'&=& W(\modb ^{\otimes m}\otimes T).
  \end{eqnarray}
  Since $U$ is a morphism in $\modB _0$, \eqref{eq:80} implies that $W$ is
  a morphism in $\modB _0$.  The generators of $\modB _0$ given in
  \fullref{thm:9} are in $\modB '$, since we have
  \begin{align*}
    \gamma _{\pm}  &=
    (\mu _{\modb} \otimes \modb )(\modb \otimes \psi _{\modb ,\modb }^{\pm
1})(\modb \otimes \mu _{\modb} ^{[3]}\otimes \modb )(c_{\pm} \otimes \modb
\otimes c_{\mp}),\\
    t_{+-}^{\pm 1} &= \mu _{\modb} ^{[3]}(v_{\mp}\otimes \modb \otimes
v_{\pm} ).
  \end{align*}
  Hence $W$ is in~$\modB '$.  Since $\modb ^{\otimes p}\otimes T$ is in $\modB '$,
  it follows from \eqref{eq:81} that $U'$ is in~$\modB '$.  This completes
  the proof.
\end{proof}

\begin{remark}
  \label{rem:1}
  The set of generators of $\modB $ given in \fullref{thm:2} is not
  minimal.  One can show, for example, that $\modB $ is minimally
  generated as a braided subcategory of $\modT $ by the object $\modb $ and
  the morphisms $\mu _{\modb} $, $v_+$, $c_+$ and $c_-$.
\end{remark}

\fullref{thm:2} implies that each bottom tangle can be obtained as
a result of horizontal and vertical pasting of finitely many copies of
the tangles
$1_{\modb} $, $\psi _{\modb ,\modb }^{\pm 1}$, $\mu _{\modb} $, $\eta
_{\modb} $, $v_+$, $v_-$, $c_+$, $c_-$.

In the following we give several corollaries to \fullref{thm:2}.

The following notation is useful in the rest of the paper.
For $f\in \modB (m,n)$ and $i,j\ge 0$, set
\begin{equation*}
  f_{(i,j)} = \modb ^{\otimes i}\otimes f\otimes \modb ^{\otimes j} \in \modB (i+m+j,i+n+j).
\end{equation*}
The following corollary to \fullref{thm:2} is sometimes useful.

\begin{corollary}
  \label{thm:36}
  As a subcategory of $\modT $, $\modB $ is generated by the objects
  $\modb ^{\otimes n}$, $n\ge 0$, and the morphisms
  \begin{eqnarray*}
    (\psi _{\modb ,\modb })_{(i,j)},(\psi _{\modb ,\modb }^{-1})_{(i,j)}
\quad &\text{for } & i,j\ge 0,\\
    f_{(i,0)} \quad &\text{for }& f\in \{\mu _{\modb} ,\eta _{\modb}
,v_{\pm} ,c_{\pm} \}, i\ge 0.
  \end{eqnarray*}
\end{corollary}

\begin{proof}
  By \fullref{thm:2}, $\modB $ is generated as a subcategory of $\modT $
  by the morphisms $f_{(i,j)}$ with
  $f\in \{\psi _{\modb ,\modb }^{\pm 1},\mu _{\modb} ,\eta _{\modb}
,v_{\pm} ,c_{\pm} \}$, $i,j\ge 0$.  For
  $f\neq \psi _{\modb ,\modb }^{\pm 1}$, we can express $f_{(i,j)}$ as a conjugate
  of $f_{(i+j,0)}$ by a doubled braid.  (Here a {\em doubled braid}
  means a morphism in the braided subcategory of $\modB $ generated by the
  object $\modb $.)  This implies the assertion.
\end{proof}

\begin{corollary}
  \label{thm:41}
  (1) Each $T\in \BT_n$ can be expressed as
  \begin{equation}
    \label{eq:21}
    T = \bigl(t_{+0}^{p_1}\otimes \cdots\otimes t_{+0}^{p_n}\bigr)
    \bigl(\mu _{\modb} ^{[j_1]}\otimes \cdots\otimes \mu _{\modb}
^{[j_n]}\bigr) \beta \bigl(c_+^{\otimes l_+}\otimes c_-^{\otimes l_-}\bigr)
  \end{equation}
  with $t_{+0}=t_{\downarrow} \otimes \uparrow (=\mu _{\modb} (v_-\otimes \modb ))\in \modB (1,1)$, $p_1,\ldots,p_n\in \modZ $,
  $j_1,\ldots,j_n\ge 0$, $l_+,l_-\ge 0$, $2(l_++l_-)=j_1+\cdots+j_n$, and
  $\beta \in \modB (2(l_++l_-),2(l_++l_-))$ a doubled braid.  (For example, see
  \fullref{fig:example-thm41}.)  \FIG{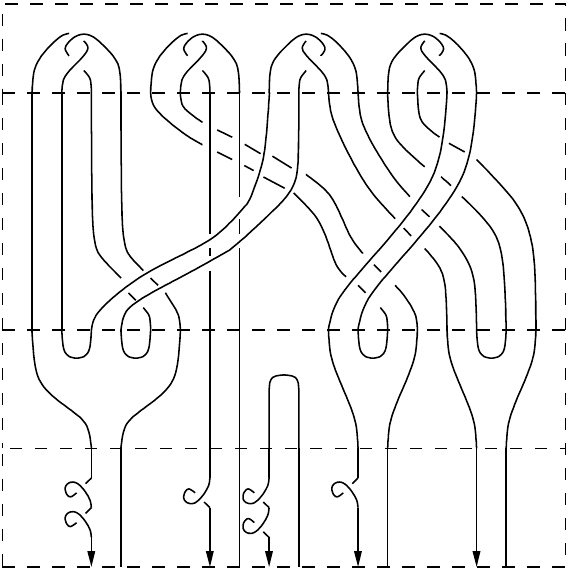}{An example of
  $T$ with $n=5$, $(p_1,\ldots,p_5)=(-2,1,2,-1,0)$,
  $(j_1,\ldots,j_5)=(3,1,0,2,2)$, and $l_+=l_-=2$}{height=35mm}

  (2) Each $T\in \BT_n$ can be expressed up to framing change as
  \begin{equation}
    \label{eq:13}
    T = \bigl(\mu _{\modb} ^{[j_1]}\otimes \cdots\otimes \mu _{\modb}
^{[j_n]}\bigr) \beta \bigl(c_+^{\otimes l_+}\otimes c_-^{\otimes l_-}\bigr)
  \end{equation}
  where $j_1,\ldots,j_n\ge 0$, $l_+,l_-\ge 0$ with
  $2(l_++l_-)=j_1+\cdots+j_n$, and where $\beta \in \modB (2(l_++l_-),2(l_++l_-))$ is
  a doubled braid.
\end{corollary}

\begin{proof}
Note that composing $\modb ^{\otimes (i-1)}\otimes t_{+0}^{p_i}\otimes \modb ^{\otimes (n-i)}$ from the
  left to $T\in \BT_n$ just changes the framing of the $i$th component
  of $T$ by $p_i$.  Hence we have only to prove (2).  In the following
  we ignore the framing.  It suffices to prove that if $T$ is as in
  \eqref{eq:13}, and $U=f_{(r,s)}$ with
  $f\in \{\psi _{\modb ,\modb }^{\pm 1},\mu _{\modb} ,\eta _{\modb}
,c_{\pm} \}$ and $r,s\ge 0$ such that
  $UT$ is well defined, then $UT$ has a decomposition similar to
  \eqref{eq:13}.  In the following we use the
  notation
  \begin{equation*}
    \mu _{\modb} ^{[a_1,a_2,\ldots,a_k]} = \mu _{\modb} ^{[a_1]}\otimes
\mu _{\modb} ^{[a_2]}\otimes \cdots\otimes \mu _{\modb} ^{[a_k]}
  \end{equation*}
  for $a_1,\ldots,a_k\ge 0$.

  The case $f=\psi _{\modb ,\modb }$ follows from
  \begin{multline*}
      (\psi _{\modb ,\modb })_{(r,s)}\mu _{\modb} ^{[j_1,\ldots,j_n]}
       =\\
       \quad \mu _{\modb} ^{[j_1,\ldots,j_r,j_{r+2},j_{r+1},j_{r+3},\ldots,j_n]}
      (\psi _{\modb ^{\otimes j_{r+1}},\modb ^{\otimes j_{r+2}}})_{(j_1+\cdots+j_r,j_{r+3}+\cdots+j_n)}.
  \end{multline*}
  The case $f=\psi _{\modb ,\modb }^{-1}$ is similar.  The cases
  $f=\mu _{\modb} ,\eta _{\modb} $ follow from
  \begin{align*}
    \mu _{(r,s)}\mu _{\modb} ^{[j_1,\ldots,j_n]}
    &=\mu _{\modb} ^{[j_1,\ldots,j_r,j_{r+1}+j_{r+2},j_{r+3},\ldots,j_n]},\\
    \eta _{(r,s)}\mu _{\modb} ^{[j_1,\ldots,j_n]}
    &=\mu _{\modb} ^{[j_1,\ldots,j_r,0,j_{r+1},\ldots,j_n]}.
  \end{align*}
  For $f=c_{\pm} $, we have
  \begin{multline*}
      f_{(r,s)}\mu _{\modb} ^{[j_1,\ldots,j_n]}\beta (c_+^{\otimes l_+}\otimes c_-^{\otimes l_-})=\\
      \mu _{\modb} ^{[j_1,\ldots,j_r,1,1,j_{r+1},\ldots,j_n]}
      (1^{\otimes (j_1+\cdots+j_r)}\otimes \psi _{\modb ^{\otimes
(j_{r+1}+\cdots+j_n)},\modb ^{\otimes 2}}) \\
      (\beta \otimes \modb ^{\otimes 2})(c_+^{\otimes l_+}\otimes c_-^{\otimes l_-}\otimes f).
  \end{multline*}
  If $f=c_-$, then we are done.  The other case $f=c_+$ follows from
$$c_+^{\otimes l_+}\otimes c_-^{\otimes l_-}\otimes c_+
    =(\modb ^{\otimes l_+}\otimes \psi _{\modb ^{\otimes 2},\modb
^{\otimes 2l_-}})(c_+^{\otimes (l_++1)}\otimes c_-^{\otimes l_-}).\proved$$
\end{proof}

\begin{remark}
  \label{r2}
  In \fullref{thm:41} (1), we may assume that
  $p_1,\ldots,p_n\in \{0,1\}$.  This follows from the identity
  $t_{+0}^{\pm 2}=\mu _{\modb} ^{[3]}(c_{\mp}\otimes \modb )$.  In particular, it follows
  that if each component of $T\in \BT_n$ is of even framing, then $T$
  can be expressed as in \fullref{thm:41} (1) with
  $p_1=\cdots=p_n=0$.  This fact is used in \fullref{sec:quasitriangular}.
\end{remark}

Let $\modA $ denote the braided subcategory of $\modB $ generated by the
object $\modb $ and the morphisms $\mu _{\modb} $ and $\eta _{\modb} $.  ($\modA $ is naturally
isomorphic to the braided category $\langle \sfA\rangle $ freely generated by an
algebra $\sfA$, defined later in \fullref{sec:exter-braid-hopf},
but we do not need this fact.)  Clearly, $\modA $ is a subcategory of
$\cB_0$ (and hence of $\modB _0$).  We need the following corollary later.

\begin{corollary}
  \label{thm:14}
  Any $T\in \BT_n$ can be expressed as a composition
  $T=T'T''$ with $T'\in \modA (m,n)$ and $T''\in \{v_{\pm} ,c_{\pm} \}^*\cap \BT_m$,
  $m\ge 0$.
\end{corollary}

\begin{proof}
  This easily follows from \fullref{thm:2}, similarly to
  \fullref{thm:41}.
\end{proof}

\section{Hopf algebra action on bottom tangles}
\label{sec6}

\subsection{The braided category $\modH $ freely generated by a Hopf
  algebra $\HH$}
\label{sec:redef-noti-hopf}

Let $\modH $ denote the braided category freely generated by a Hopf
algebra $\HH$.  In other words, $\modH $ is a braided category with a Hopf
algebra $\HH$ such that if $\modM $ is a braided category and $H$ is a
Hopf algebra in $\modM $, then there is a unique braided functor
$F_H\col \modH \rightarrow \modM $ that maps the Hopf algebra structure of $\HH$ into that
of $H$.  Such $\modH $ is unique up to isomorphism.

A more concrete definition of $\modH $ (up to isomorphism) is sketched as
follows.  Set $\Ob(\modH )=\{\HH^{\otimes n}\ver n\ge 0\}$.  Consider the expressions
obtained by compositions and tensor products from copies of the
morphisms
\begin{gather*}
  1_{\one} \col \one \rightarrow \one ,\quad
  1_{\HH}\col \HH\rightarrow \HH,\quad
  \psi _{\HH,\HH}^{\pm 1}\col \HH^{\otimes 2}\rightarrow \HH^{\otimes2},\quad
  \mu _{\HH}\col \HH^{\otimes 2}\rightarrow \HH,\quad
  \eta _{\HH}\col \one \rightarrow \HH,\\
  \Delta _{\HH}\col \HH\rightarrow \HH^{\otimes 2},\quad
  \epsilon _{\HH}\col \HH\rightarrow \one ,\quad
  S_{\HH}\col \HH\rightarrow \HH,
\end{gather*}
where we understand $\one =\HH^{\otimes 0}$ and $\HH=\HH^{\otimes 1}$, and define an
equivalence relation on such expressions generated by the axioms of
braided category and Hopf algebra.  Then the morphisms in $\modH $ are the
equivalence classes of such expressions.

Note that if $F\col \modH \rightarrow \modM $ is a braided functor of $\modH $ into a braided
category $\modM $, then $F(\HH)\in \Ob(\modM )$ is equipped with a Hopf algebra
structure
\begin{equation*}
  (F(\mu _{\HH}),F(\eta _{\HH}),F(\Delta _{\HH}),F(\epsilon
_{\HH}),F(S_{\HH})).
\end{equation*}
Conversely, if $H$ is a Hopf algebra in $\modM $, then there is a unique
braided functor $F\col \modH \rightarrow \modM $ such that $F$ maps the Hopf algebra
structure of $\HH$ into that of~$H$.  Hence there is a canonical
one-to-one correspondence between the Hopf algebras in a braided
category $\modM $ and the braided functors from $\modH $ to $\modM $.

For $f\in \modH (\HH^{\otimes m},\HH^{\otimes n})$ ($m,n\ge 0$) and $i,j\ge 0$, set
\begin{equation*}
  f_{(i,j)} =
  \HH^{\otimes i}\otimes f\otimes \HH^{\otimes j}\in \modH (\HH^{\otimes (m+i+j)},\HH^{\otimes (n+i+j)}).
\end{equation*}
Note that $\modH $ is generated as a category by the objects $\HH^{\otimes i}$,
$i\ge 0$, and the morphisms $f_{(i,j)}$ with
$f\in \{\psi _{\HH,\HH},\psi _{\HH,\HH}^{-1},\mu ,\eta ,\Delta ,\epsilon ,S\}$ and $i,j\ge 0$.
In the following, we write
$(\psi _{\HH,\HH}^{\pm 1})_{(i,j)}=\psi _{(i,j)}^{\pm 1}$.

\begin{lemma}
  \label{lem:22}
  As a category, $\modH $ has a presentation with the generators
  $f_{(i,j)}$ for
  $f\in \{\psi _{\HH,\HH},\psi _{\HH,\HH}^{-1},\mu ,\eta ,\Delta ,\epsilon ,S\}$ and
  $i,j\ge 0$, and the relations
  \begin{gather}
    \label{eq:62}
    f_{(i,j+q'+k)}g_{(i+p+j,k)}=g_{(i+p'+j,k)}f_{(i,j+q+k)},\\
    \label{eq:63}
    \psi _{(i,j)}\psi _{(i,j)}^{-1}=\psi _{(i,j)}^{-1}\psi _{(i,j)}=1_{\HH^{\otimes (i+j+2)}},\\
    \label{eq:64}
    \left\{\begin{aligned}
    (\psi _{p',1})_{(i,j)}f_{(i,j+1)}\negthinspace&=\negthinspace
      f_{(i+1,j)}(\psi _{p,1})_{(i,j)},\\
    (\psi _{1,p'})_{(i,j)}f_{(i+1,j)}\negthinspace&=\negthinspace
      f_{(i,j+1)}(\psi _{1,p})_{(i,j)},
    \end{aligned}\right.\\
    \left\{\begin{aligned}
    \mu _{(i,j)}\eta _{(i,j+1)}&=\mu _{(i,j)}\eta _{(i+1,j)}
      =1_{\HH^{\otimes (i+j+1)}},\\
    \mu _{(i,j)}\mu _{(i,j+1)}&=\mu _{(i,j)}\mu _{(i+1,j)},
    \end{aligned}\right.\\
    \left\{\begin{aligned}
    \epsilon _{(i,j+1)}\Delta _{(i,j)}&=\epsilon _{(i+1,j)}\Delta _{(i,j)}
      =1_{\HH^{\otimes (i+j+1)}},\\
    \Delta _{(i,j+1)}\Delta _{(i,j)}&=\Delta _{(i+1,j)}\Delta _{(i,j)},
    \end{aligned}\right.\\
    \left\{\begin{aligned}
    \epsilon _{(i,j)}\eta _{(i,j)}&=1_{\HH^{\otimes (i+j)}},\\
    \epsilon _{(i,j)}\mu _{(i,j)}&=\epsilon _{(i,j)}\epsilon _{(i+1,j)},\\
    \Delta _{(i,j)}\eta _{(i,j)}&=\eta _{(i+1,j)}\eta _{(i,j)},
    \end{aligned}\right.\\
    \label{eq:11}
    \Delta _{(i,j)}\mu _{(i,j)}=
    \mu _{(i+1,j)}\mu _{(i,j+2)}\psi _{(i+1,j+1)}\Delta _{(i,j+2)}\Delta _{(i+1,j)},\\
    \label{eq:12}
    \mu _{(i,j)}S_{(i,j+1)}\Delta _{(i,j)}=\mu _{(i,j)}S_{(i+1,j)}\Delta _{(i,j)}
  =  \eta _{(i,j)}\epsilon _{(i,j)},
  \end{gather}
  for $i,j,k\ge 0$ and
  $f,g\in \{\psi _{\HH,\HH},\psi _{\HH,\HH}^{-1},\mu ,\eta ,\Delta ,\epsilon ,S\}$ with
  $f\col \HH^{\otimes p}\rightarrow \HH^{\otimes p'}$, $g\col \HH^{\otimes q}\rightarrow \HH^{\otimes q'}$.  Here,
  $(\psi _{p,1})_{(i,j)}$ and $(\psi _{1,p})_{(i,j)}$ for $p=0,1,2$ and
  $i,j\ge 0$ are defined by
  \begin{gather*}
    (\psi _{0,1})_{(i,j)}=(\psi _{1,0})_{(i,j)}=1_{\HH^{\otimes(i+j+1)}},\quad
    (\psi _{1,1})_{(i,j)}=\psi _{(i,j)},\\
    (\psi _{2,1})_{(i,j)}=\psi _{(i,j+1)}\psi _{(i+1,j)},\quad
    (\psi _{1,2})_{(i,j)}=\psi _{(i+1,j)}\psi _{(i,j+1)}.
  \end{gather*}
\end{lemma}

\begin{proof}
  We only give a sketch proof, since a detailed proof is long though
  straightforward.  The relations given in the lemma are the ones
  derived from the axioms for braided category and Hopf algebra, hence
  valid in $\modH $.  We have to show, conversely, that all the relations
  in $\modH $ can be derived from the relations given in the lemma.  It
  suffices to show that the category $\modH '$ with the presentation given
  in the lemma is a braided category with a Hopf algebra $\HH$.  The
  relation \eqref{eq:62} implies that $\modH '$ is a monoidal category,
  since we can define the monoidal structure for $\modH '$ by
  \begin{equation*}
    f_{(i,j)}\otimes f'_{(i',j')}
    = f_{(i,j+i'+n'+j')}f'_{(i+m+j+i',j')}
  \end{equation*}
  for $f,f'\in \{\psi _{\HH,\HH},\psi _{\HH,\HH}^{-1},\mu ,\eta ,\Delta ,\epsilon ,S\}$,
  $f\col \HH^{\otimes m}\rightarrow \HH^{\otimes n}$, $f'\col \HH^{\otimes m'}\rightarrow \HH^{\otimes n'}$.  The relations
  \eqref{eq:63} and \eqref{eq:64} imply that $\modH '$ is a braided category,
  and the other relations imply that $\HH$ is a Hopf algebra in $\modH '$.
\end{proof}

\subsection{External Hopf algebras in braided categories}
\label{sec:exter-braid-hopf}

Let $\langle \sfA\rangle $ denote the braided category freely generated by an
algebra $\sfA=(\sfA,\mu _{\sfA},\eta _{\sfA})$.  For a braided category $\modM $
and an algebra $A$ in $\modM $, let
\begin{equation*}
  i_{\modM ,A}\col \langle \sfA\rangle \rightarrow \modM
\end{equation*}
denote the unique braided functor that maps the algebra structure of
$\sfA$ into the algebra structure of $A$.

\begin{definition}
  An {\em external Hopf algebra} $(H,F)$ in a braided category $\modM $ is
  a pair of an algebra $H=(H,\mu _H,\eta _H)$ in $\modM $ and a functor
  $F\col \modH \rightarrow \Sets$ into the category $\Sets$ of sets and functions such
  that we have a commutative square
  \begin{equation*}
    \begin{CD}
      \langle \sfA\rangle  @>i_{\modM ,H}>> \modM \\
      @Vi_{\modH ,\HH}VV @VV\modM (\one ,-)V\\
      \modH  @>>F> \Sets.
    \end{CD}
  \end{equation*}
\end{definition}

\begin{remark}
  \label{rem:2}
  The functor $i_{\modH ,\HH}\col \langle \sfA\rangle \rightarrow \modH $ is faithful.  Hence we can
  identify $\langle \sfA\rangle $ with the braided subcategory of $\modH $ generated
  by $\HH$, $\mu _{\HH}$ and $\eta _{\HH}$.  However, we do not need this fact in
  the rest of the paper.
\end{remark}

\begin{remark}
  \label{rem:3}
  To each Hopf algebra $H$ in a braided category $\modM $, we can
  associate an external Hopf algebra in $\modM $ as follows.  Let
  $F_H\col \modH \rightarrow \modM $ be the braided functor which maps the Hopf
  algebra $\HH$ into $H$.  Then the pair
  $((H,\mu _H,\eta _H),\modM (\one ,F_H(-)))$ is an external Hopf algebra in $\modM $.
  Therefore, the notion of external Hopf algebra in $\modM $ can be
  regarded as a generalization of the notion of Hopf algebra in $\modM $.
\end{remark}

\subsection{The external Hopf algebra structure in $\modB $}
\label{sec:an-external braided}

In this subsection, we define an external Hopf algebra
$(\modb ,F_{\modb} )$ in $\modB $.  First note that $(\modb ,\mu _{\modb}
,\eta _{\modb} )$ is an algebra
in $\modB $.  (This algebra structure of $\modb $ cannot be extended to a Hopf
algebra structure in $\modB $ in the usual sense, since there is no
morphism $\epsilon \col \modb \rightarrow \one $ in $\modB $.)

For $i,j\ge 0$ and $T\in \BT_{i+j+1}$, set
\begin{gather*}
  \begin{split}
  \check\Delta _{(i,j)}(T)
  &= (\modb ^{\otimes i}\otimes (\downarrow \otimes \psi _{\modb ,\uparrow })\otimes \modb ^{\otimes j})T',\\
  \check\epsilon _{(i,j)}(T) &= T'',\\
  \check S_{(i,j)}(T) &=
  (\modb ^{\otimes i}\otimes \psi _{\uparrow ,\downarrow }(\uparrow
\otimes t_{\downarrow} )\otimes \modb ^{\otimes j}) T''',
  \end{split}
\end{gather*}
where $T'\in \modT (\one ,\modb ^{\otimes i}\otimes \downarrow \otimes \downarrow \otimes \uparrow \otimes \uparrow \otimes \modb ^{\otimes j})$ is obtained
from~$T$ by duplicating the $(i+1)$st component of $T$,
$T''\in \modT (\one ,\modb ^{\otimes (i+j)})$ is obtained from $T$ by removing the
$(i+1)$st component of $T$, and
$T'''\in \modT (\one ,\modb ^{\otimes i}\otimes \uparrow \otimes \downarrow \otimes \modb ^{\otimes j})$ is obtained from $T$ by
reversing orientation of the $(i+1)$st component of $T$, see \fullref{fig:parallel-and-twist}.
\labellist\tiny
\pinlabel {$1$} [t] at 26 736
\pinlabel {$\cdots$} [t] at 45 732
\pinlabel {$i{+}1$} [t] at 68 736
\pinlabel {$\cdots$} [t] at 87 732
\pinlabel {$i{+}j{+}1$} [t] at 115 736
\small
\pinlabel {$T$} [b] at 70 675
\pinlabel {$\check{\Delta}_{(i,j)}(T)$} [b] at 220 675
\pinlabel {$T'$} at 260 790
\pinlabel {$T''$} at 405 790
\pinlabel {$\check{\epsilon}_{(i,j)}(T)$} [b] at 360 675
\pinlabel {$\check{S}_{(i,j)}(T)$} [b] at 515 675
\pinlabel {$T'''$} at 555 790
\endlabellist
\FIGn{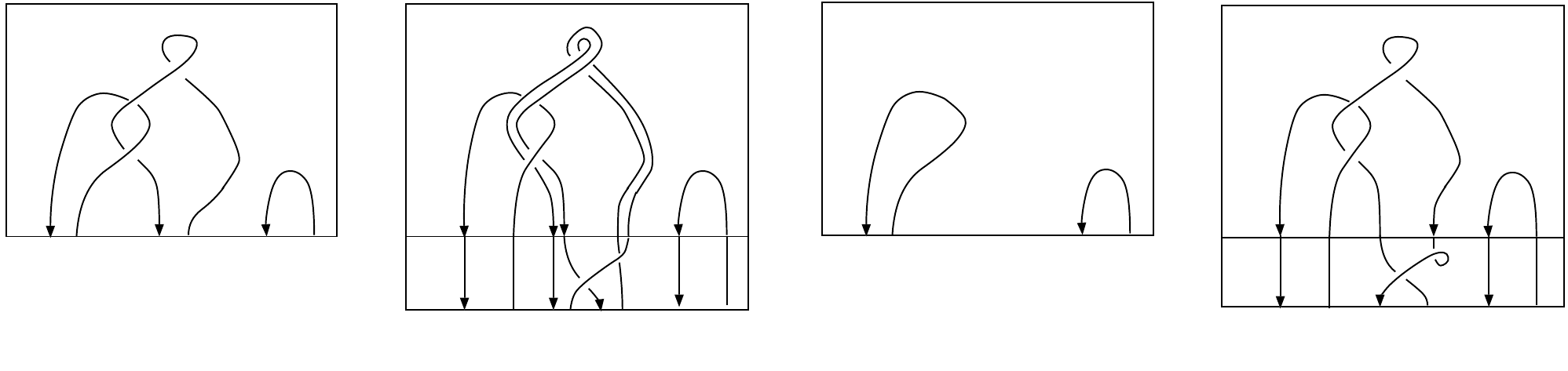}{}{height=30mm}
We have
\begin{equation*}
  \check\Delta _{(i,j)}(T)\in \BT_{i+j+2},\quad
  \check\epsilon _{(i,j)}(T)\in \BT_{i+j},\quad
  \check S_{(i,j)}(T)\in \BT_{i+j+1}.
\end{equation*}
Hence there are functions
\begin{gather*}
  \begin{split}
  \check\Delta _{(i,j)}&\col \BT_{i+j+1}\rightarrow \BT_{i+j+2},\\
  \check\epsilon _{(i,j)}&\col \BT_{i+j+1}\rightarrow \BT_{i+j},\\
  \check S_{(i,j)}&\col \BT_{i+j+1}\rightarrow \BT_{i+j+1}.
  \end{split}
\end{gather*}

\begin{theorem}
  \label{thm:60}
  There is a unique external Hopf algebra
  $((\modb ,\mu _{\modb} ,\eta _{\modb} ),F_{\modb} )$ in $\modB $ with
$F_{\modb} \col \modH \rightarrow \Sets$ satisfying
  \begin{equation}
    \label{eq:65}
    \begin{split}
      F_{\modb} (\Delta _{(i,j)}) = \check\Delta _{(i,j)},\quad
      F_{\modb} (\epsilon _{(i,j)}) = \check\epsilon _{(i,j)},\quad
      F_{\modb} (S_{(i,j)}) = \check S_{(i,j)}
    \end{split}
  \end{equation}
  for $i,j\ge 0$.
\end{theorem}

\begin{proof}
  We claim that there is a functor $F_{\modb} \col \modH \rightarrow \Sets$ satisfying
  \eqref{eq:65} and
  \begin{align*}
    F_{\modb} (\psi _{(i,j)}^{\pm 1}) &= (\psi _{\modb ,\modb }^{\pm 1})_{(i,j)}(-),\\
    F_{\modb} (\mu _{(i,j)}) &= (\mu _{\modb} )_{(i,j)}(-),\\
    F_{\modb} (\eta _{(i,j)}) &= (\eta _{\modb} )_{(i,j)}(-),
  \end{align*}
  for $i,j\ge 0$.  If this claim is true, then one can easily check that
  $F_{\modb} $ is unique and $(\modb ,F_{\modb} )$ is an external Hopf algebra
  in $\modB $.

  To prove the above claim, it suffices to check that each relation in
  \fullref{lem:22} are mapped into a relation in $\modB $.  For example,
  the relation \eqref{eq:11} is mapped into
  \begin{equation}
    \label{eq:66}
    \begin{split}
      \check\Delta _{(i,j)}((\mu _{\modb} )_{(i,j)}g)
      =&(\mu _{\modb} )_{(i+1,j)}
      (\mu _{\modb} )_{(i,j+2)}
      (\psi _{\modb ,\modb })_{(i+1,j+1)}
      \check\Delta _{(i,j+2)}\check\Delta _{(i+1,j)}(g)\\
      \Bigl(=&
      ((\mu _{\modb} \otimes \mu _{\modb} )(\modb \otimes \psi _{\modb ,\modb }\otimes \modb ))_{(i,j)}
      \check\Delta _{(i,j+2)}\check\Delta _{(i+1,j)}(g)\Bigr)\\
    \end{split}
  \end{equation}
  for $g\in \BT_{i+j+2}$.  This can be proved in a graphical way.  If
\labellist\tiny
\pinlabel {$1$} [b] at 42 718
\pinlabel {$i{+}1$} [b] at 78 718
\pinlabel {$i{+}2$} [b] at 109 718
\pinlabel {$i{+}j{+}2$} [b] at 152 718
\pinlabel {$\cdots$} at 63 747
\pinlabel {$\cdots$} at 128 747
\endlabellist
$$g=\raisebox{-6mm}{\incl{10mm}{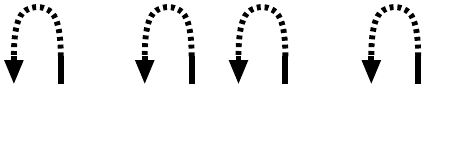}},$$
  then the left and the right hand sides
  of \eqref{eq:66} are as depicted in \fullref{fig:delta-mu} (a)
  and (b), respectively.
  \begin{figure}[t]
    \begin{equation*}
      \begin{matrix}
\labellist\tiny
\pinlabel {$1$} [b] at 163 710
\pinlabel {$i{+}1$} [b] at 215 710
\pinlabel {$i{+}2$} [b] at 246 710
\pinlabel {$i{+}j{+}2$} [b] at 294 710
\pinlabel {$\cdots$} at 190 760
\pinlabel {$\cdots$} at 265 760
\endlabellist
	\incl{20mm}{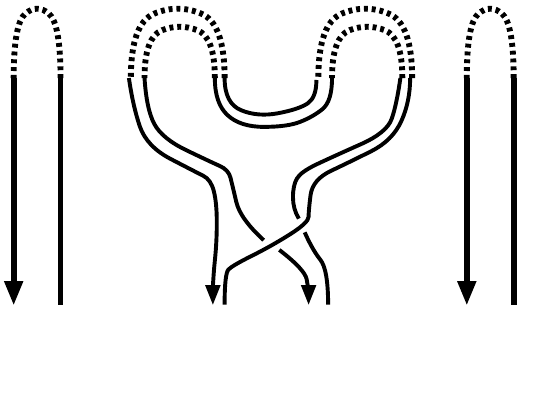}&\hspace{10mm}
\labellist\tiny
\pinlabel {$1$} [b] at 333 710
\pinlabel {$i{+}1$} [b] at 375 710
\pinlabel {$i{+}2$} [b] at 430 710
\pinlabel {$i{+}j{+}2$} [b] at 473 710
\pinlabel {$\cdots$} at 352 760
\pinlabel {$\cdots$} at 454 760
\endlabellist
	&\incl{20mm}{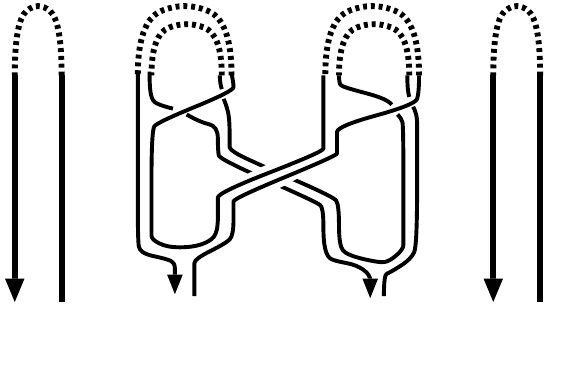}\\
	\text{(a)}& &\text{(b)}
      \end{matrix}
  \end{equation*}
    \nocolon\caption{}
    \label{fig:delta-mu}
  \end{figure}
  As another example, \eqref{eq:12} is mapped into the equivalence of
  \fullref{fig:antipode} (a), (b), and (c).
  \begin{figure}[t]
    \begin{equation*}
      \begin{matrix}
	\incl{25mm}{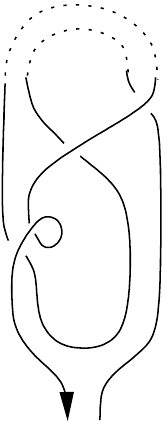}&\quad &
	\incl{25mm}{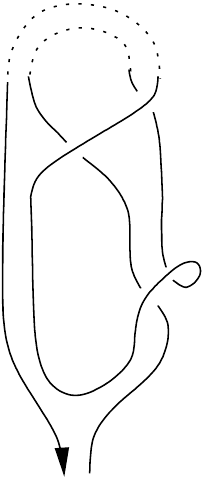}& &
	\incl{25mm}{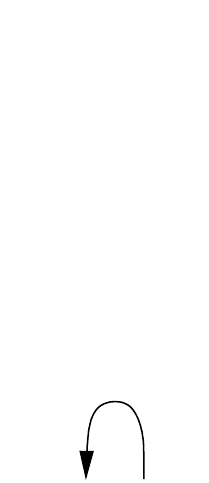}\\
	\text{(a)}& &
	\text{(b)}& &
	\text{(c)}
      \end{matrix}
    \end{equation*}
    \caption{Here, only the $i$th component is depicted in each figure.}
    \label{fig:antipode}
  \end{figure}
  It is straightforward to check the other relations, and we leave it
  to the reader.  The detail is to a certain extent similar to the proof of
  the existence of the Hopf algebra in the category of cobordisms of
  surfaces with connected boundary (see Crane and Yetter
\cite{Crane-Yetter:99} and Kerler \cite{Kerler:99})
  and also to the proof of the Hopf algebra relations satisfied by
  claspers \cite{Habiro:claspers}.  See also \fullref{sec:category-0}.
\end{proof}

\fullref{thm:60} implies that there is a Hopf algebra action on
the bottom tangles as we mentioned in \fullref{sec:hopf-algebra-action-1}.

\subsection{Adjoint actions and coactions}
\label{sec:some-identities}
In this subsection, we consider the image by $F_{\modb} $ of the left
adjoint action and the left adjoint coaction of $\HH$, which we
need later.

Let $\ad_{\HH}\col \HH^{\otimes 2}\rightarrow \HH$ denote the {\em left adjoint action}
for $\HH$, which is defined by
\begin{equation*}
  \ad_{\HH} = \mu _{\HH}^{[3]}(\HH\otimes \psi _{\HH,\HH})(\HH\otimes
S_{\HH}\otimes \HH)(\Delta _{\HH}\otimes \HH).
\end{equation*}
For $i,j\ge 0$ and $T\in \BT_{i+j+2}$, we have
\begin{equation*}
  \begin{split}
    &F_{\modb} ((\ad_{\HH})_{(i,j)})(T)\\
    =&F_{\modb} ((\mu _{\HH}^{[3]}(\HH\otimes \psi _{\HH,\HH})(\HH\otimes
S_{\HH}\otimes \HH)(\Delta _{\HH}\otimes \HH))_{(i,j)})(T)\\
    =&(\mu _{\modb} ^{[3]})_{(i,j)}(\modb \otimes \psi _{\modb ,\modb })_{(i,j)}
    F_{\modb} ((\HH\otimes S_{\HH}\otimes \HH)_{(i,j)})
    (F_{\modb} ((\Delta _{\HH}\otimes \HH)_{(i,j)})(T)).
  \end{split}
\end{equation*}
Hence $F_{\modb} ((\ad_{\HH})_{(i,j)})$ maps $T\in \BT_{i+j+1}$ to a bottom
tangle as depicted in \fullref{fig:ad} (a).
\labellist\tiny
\endlabellist
\FIG{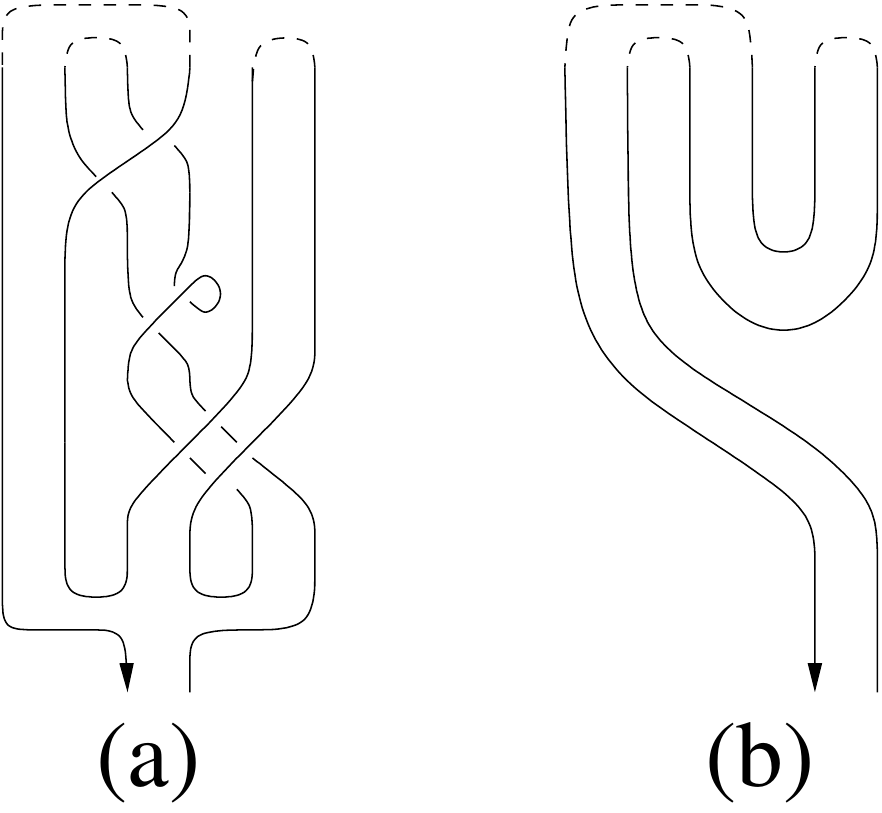}{(a) The tangle $F_{\modb} ((\ad_{\HH})_{(i,j)})(T)\in \BT_{i+j+1}$,
with only the $(i+1)$st component depicted.  The dotted lines denote
(parallel copies
of) the $(i+1)$st and $(i+2)$nd components of $T$.  (b) A tangle isotopic to (a).}{height=30mm}
By isotopy, we obtain a simpler tangle as in
\fullref{fig:ad} (b).  (Note that the closures $\clo(T)$ and
$\clo(F_{\modb} ((\ad_{\HH})_{(i,j)})(T))$ are isotopic.  This fact is used in
\fullref{sec10}.)

Let $\coad_{\HH}\col \HH\rightarrow \HH^{\otimes 2}$ denote the {\em left adjoint coaction}
defined by
\begin{equation*}
  \coad_{\HH}=(\mu _{\HH}\otimes \HH)(\HH\otimes \psi
_{\HH,\HH})(\HH\otimes \HH\otimes S_{\HH})\Delta _{\HH}^{[3]}.
\end{equation*}
For $i,j\ge 0$ and $T\in \BT_{i+j+1}$, we have
\begin{equation*}
  \begin{split}
    &F_{\modb} ((\coad_{\HH})_{(i,j)})(T)\\
    =&F_{\modb} (((\mu _{\HH}\otimes \HH)(\HH\otimes \psi
_{\HH,\HH})(\HH\otimes \HH\otimes S_{\HH})\Delta _{\HH}^{[3]})_{(i,j)})(T)\\
    =&(\mu _{\modb} \otimes \modb )_{(i,j)}(\modb \otimes \psi _{\modb ,\modb })_{(i,j)}
    F_{\modb} ((\HH\otimes \HH\otimes S_{\HH})_{(i,j)})(F_{\modb} ((\Delta
_{\HH}^{[3]})_{(i,j)})(T)).
  \end{split}
\end{equation*}
Hence $F_{\modb} ((\coad_{\HH})_{(i,j)})$ maps $T\in \BT_{i+j+1}$ to a
$(i+j+2)$--component bottom tangle as depicted in \fullref{fig:coad} (a).
\FIG{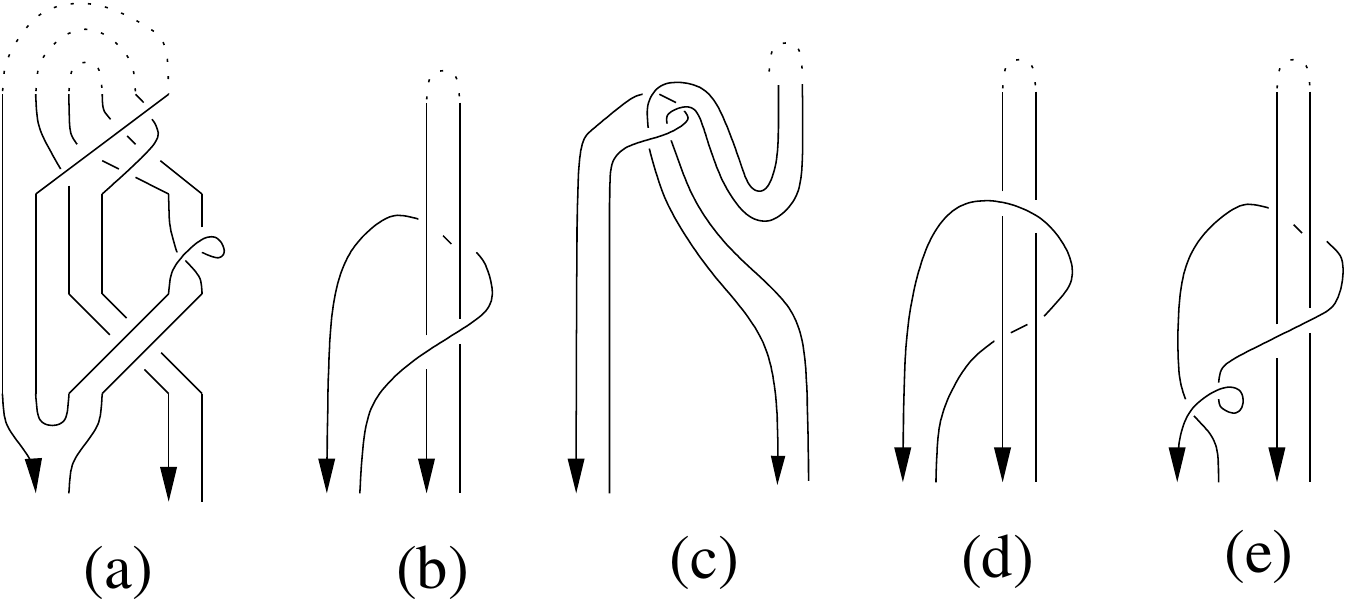}{(a) The tangle $F_{\modb}
((\coad_{\HH})_{(i,j)})(T)\in \BT_{i+j+2}$, where only the $(i+1)$st and
$(i+2)$nd components are depicted.  (b) A tangle isotopic to (a),
which is $(\gamma _+)_{(i,j)}T$.  (c) Another tangle isotopic to (a),
which is $F_{\modb} ((\HH\otimes \ad_{\HH})_{(i,j)})((c_+\otimes \modb
)_{(i,j)}T)$.  (d) The tangle $(\gamma _-)_{(i,j)}T$ with only
$(i+1)$st component depicted.  (e) A tangle isotopic to
(d).}{height=35mm} Since it is isotopic to \fullref{fig:coad} (b)
and (c), we have
\begin{equation}
  \label{eq:7}
  F_{\modb} ((\coad_{\HH})_{(i,j)})(T)
  =(\gamma _+)_{(i,j)}T
  =F_{\modb} ((\HH\otimes \ad_{\HH})_{(i,j)})((c_+\otimes \modb )_{(i,j)}T).
\end{equation}
Since \fullref{fig:coad} (d) and (e) are isotopic, we have
\begin{equation}
  \label{eq:17}
  (\gamma _-)_{(i,j)}T = F_{\modb} ((S_{\HH}\otimes \HH)_{(i,j)})((\gamma _+)_{(i,j)}T).
\end{equation}

\section{Universal tangle invariant associated to a ribbon Hopf algebra}
\label{sec7}

In this section, we give a definition of a universal tangle invariant
associated to a ribbon Hopf algebra.

\subsection{Ribbon Hopf algebras}
\label{sec:ribbon-hopf-algebra-2}
In this subsection, we  recall the definition of ribbon Hopf
algebra (see Reshetikhin and Turaev \cite{Reshetikhin-Turaev:90}).

Let $H=(H,\mu ,\eta ,\Delta ,\epsilon ,S)$ be a Hopf algebra over a commutative, unital
ring $\modk $.

A {\em universal $R$--matrix} for $H$ is an invertible element
$R\in H^{\otimes 2}$ satisfying the following properties:
\begin{gather}
  \label{eq:46}
  R\Delta (x)R^{-1} = P_{H,H}\Delta (x) \quad \text{for all $x\in H$},\\
  \label{eq:47}
  (\Delta \otimes 1)(R) = R_{13}R_{23},\quad
  (1\otimes \Delta )(R) = R_{13}R_{12}.
\end{gather}
Here $P_{H,H}\col H^{\otimes 2}\rightarrow H^{\otimes 2}$, $x\otimes y\mapsto y\otimes x$, is the
$\modk $--module homomorphism which permutes the tensor factors, and
\begin{gather*}
  R_{12}= R\otimes 1\in H^{\otimes 3},\qua
  R_{13}= (1\otimes P_{H,H})(R_{12})\in H^{\otimes 3},\qua
  R_{23}= 1\otimes R\in H^{\otimes 3}.
\end{gather*}
A Hopf algebra equipped with a universal $R$--matrix is called a {\em
  quasitriangular Hopf algebra}.

In what follows, we freely use the notations
$$R=\sum\alpha \otimes \beta \quad\text{and}\quad
R^{-1}=\sum\bar{\alpha }\otimes \bar{\beta }\quad(=(S\otimes 1)(R)).$$
We also use the notations
\begin{gather*}
  R=  \sum \alpha _i\otimes \beta _i,\quad
  R^{-1}= \sum\bar{\alpha }_i\otimes \bar{\beta }_i,
\end{gather*}
where $i$ is any index, used to distinguish several copies of $R^{\pm 1}$.

A {\em ribbon element} for $(H,R)$ is a central element $\modr \in H$ such
that
\begin{gather*}
  \modr ^2 = uS(u),\quad S(\modr )=\modr ,\quad \epsilon (\modr )=1,\quad
  \Delta (\modr ) = (R_{21}R)^{-1}(\modr \otimes \modr ),
\end{gather*}
where $R_{21}=P_{H,H}(R)=\sum \beta \otimes \alpha $ and $u=\sum S(\beta )\alpha $.  Since $u$ is
invertible, so is $\modr $.  The triple $(H,R,\modr )$ is called a {\em
  ribbon Hopf algebra}.

The element $\kappa =u\modr ^{-1}$ is grouplike, ie, $\Delta (\kappa )=\kappa \otimes \kappa $,
$\epsilon (\kappa )=1$.  We also have
\begin{equation*}
  \kappa x\kappa ^{-1} = S^2(x) \quad \text{for all $x\in H$}.
\end{equation*}
In what follows, we often use the Sweedler notation for
comultiplication.  For $x\in H$, we write
\begin{align*}
  \Delta (x) &= \sum x_{(1)}\otimes x_{(2)},\\
  \Delta ^{[n]}(x) &= \sum x_{(1)}\otimes \cdots\otimes x_{(n)}\quad \text{for $n\ge 1$}.
\end{align*}

\subsection{Adjoint action and universal quantum trace}
\label{sec:adjo-acti-univ}

We regard $H$ as a left $H$--module with the (left) adjoint action
\begin{equation*}
  \ad=\trr \col H^{\otimes 2}\rightarrow H,\quad x\otimes y\mapsto x\trr y,
\end{equation*}
defined by
\begin{equation*}
  x\trr y= \sum x_{(1)}yS(x_{(2)})\quad \text{for $x,y\in H$}.
\end{equation*}
Recall that $\ad$ is a left $H$--module homomorphism.

The function
\begin{equation*}
  H^{\otimes 2}\rightarrow H,\quad x\otimes y\mapsto x\trr y -\epsilon (x)y,
\end{equation*}
is a left $H$--module homomorphism.  Hence the image
\begin{equation*}
  N = \Span_{\modk} \{x\trr y-\epsilon (x)y\ver x,y\in H\}\subset H
\end{equation*}
is a left $H$--submodule of $H$.  Note that $H/N$, the {\em module of
coinvariants}, inherits from $H$ a trivial left $H$--module structure.

The definition of $N$ above is compatible with \eqref{e15} as follows.

\begin{lemma}
  \label{lem:5}
  We have
  \begin{equation*}
    N = \Span_{\modk} \{xy-yS^2(x)\ver x,y\in H\}.
  \end{equation*}
\end{lemma}

\begin{proof}
  The assertion follows from
  \begin{align*}
    x\trr y-\epsilon (x)y &=
      \sum\Bigl(x_{(1)}yS(x_{(2)})-yS(x_{(2)})S^2(x_{(1)})\Bigr),\\
      xy-y S^2(x) &=
      \sum\Bigl(x_{(1)}\trr yS^2(x_{(2)})-\epsilon (x_{(1)})yS^2(x_{(2)})\Bigr),
  \end{align*}
  for $x,y\in H$.
\end{proof}

As in \fullref{sec:universal-invariant}, let $\tr_q\col H\rightarrow H/N$
  denote the projection, and call it the {\em universal quantum trace}
  for $H$.  If $\modk $ is a field and $V$ is a finite-dimensional left
  $H$--module, then the {\em quantum trace} in $V$
\begin{equation*}
  \tr_q^V\col H\rightarrow \modk
\end{equation*}
factors through $\tr_q$.  Here $\tr_q^V$ is defined by
\begin{equation*}
  \tr_q^V(x)=\tr^V(\rho (\kappa x))\quad \text{for $x\in H$},
\end{equation*}
where $\rho \col H\rightarrow \End_{\modk} (V)$ denotes the left action of $H$ on $V$, and
\begin{equation*}
  \tr^V\col \End_{\modk} (V)\rightarrow \modk
\end{equation*}
denotes the trace in $V$.

\subsection{Definition of the universal invariant}
\label{sec:defin-univ-invar}

In this subsection, we recall the definition of the universal
invariant of tangles associated to a ribbon Hopf algebra $H$.  The
definition below is close to Ohtsuki's one \cite{Ohtsuki:93}, but
we use different conventions and we make some modifications.  In
particular,  for closed components, we use the universal quantum trace
instead of the universal trace.

Let $T=T_1\cup \cdots \cup T_l\cup L_1\cup \cdots \cup L_m$ with $l,m\ge 0$ be a (framed,
oriented) tangle in a cube, consisting of $l$ arc components
$T_1,\ldots,T_l$ and $m$ {\em ordered} circle components $L_1,\ldots,L_m$.
First, assume that $T$ is given by pasting copies of the tangles
$$\downarrow,\quad \uparrow,\quad \psi_{\downarrow,\downarrow}^{\pm 1},
\quad \ev_{\downarrow},\quad \ev_{\uparrow},\quad \coev_{\downarrow},\quad
\coev_{\uparrow}.$$
(Later we consider a more general case.)  We formally put
elements of $H$ on the strings of $T$ according to the rule depicted
in \fullref{fig:placing3}.
\labellist\small
\pinlabel {$\alpha$} [tr] at 48 688
\pinlabel {$\beta$} [tl] at 68 688
\pinlabel {$\wbar{\alpha}$} [br] at 120 669
\pinlabel {$\wbar{\beta}$} [bl] at 140 669
\pinlabel {$1$} [tl] at 202 678
\pinlabel {$1$} [br] at 242 680
\pinlabel {$\kappa$} [tr] at 303 678
\pinlabel {$\kappa^{-1}$} [bl] at 393 680
\endlabellist
\FIG{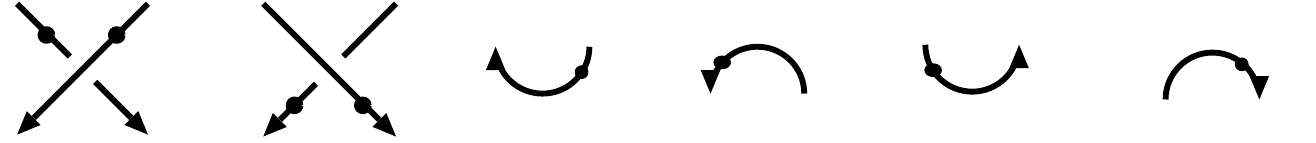}{How to place elements in
$H$ on the strings in the tangle $T$}{height=12mm}
\labellist\tiny
\pinlabel {$S(\alpha)$} [tr] at 106 592
\pinlabel {$\beta$} [tl] at 126 592
\pinlabel {$\alpha$} [tr] at 166 592
\pinlabel {$S(\beta)$} [tl] at 186 592
\pinlabel {$S(\alpha)$} [tr] at 239 592
\pinlabel {$S(\beta)$} [tl] at 260 592
\pinlabel {$S(\wbar{\alpha})$} [br] at 316 572
\pinlabel {$\wbar{\beta}$} [bl] at 336 572
\pinlabel {$\wbar{\alpha}$} [br] at 376 572
\pinlabel {$S(\wbar{\beta})$} [bl] at 396 572
\pinlabel {$S(\wbar{\alpha})$} [br] at 450 572
\pinlabel {$S(\wbar{\beta})$} [bl] at 469 572
\endlabellist
\FIG{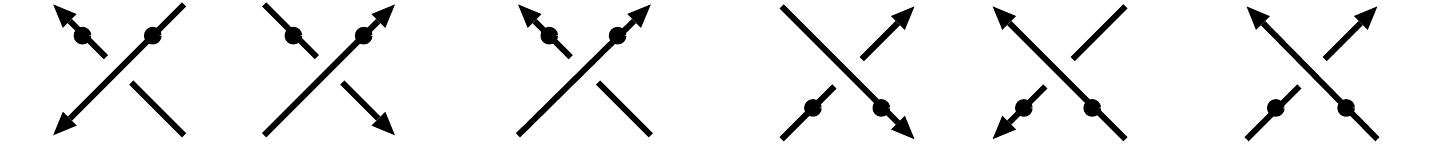}{The cases of crossings with upward strings}{height=12mm}
We define
\begin{equation}
  \label{eq:48}
  J_T = \sum
  J_{(T_1)}\otimes \cdots\otimes J_{(T_l)}\otimes J_{(L_1)}\otimes \cdots\otimes J_{(L_m)}
  \in H^{\otimes l}\otimes (H/N)^{\otimes m}
\end{equation}
as follows.  For each $i=1,\ldots,l$, we formally set $J_{(T_i)}$ to
be the product of the elements put on the component $T_i$ obtained by
reading the elements using the order determined by the opposite
orientation of $T_i$ and writing them down from left to right.  For
each $j=1,\ldots,s$, we define $J_{(L_j)}$ by first obtaining a word
$w$ by reading the elements put on $L_j$ starting from any point on
$L_j$, and setting formally $J_{(L_j)}=\tr_q(\kappa ^{-1}w)$.  (Here, it
should be noted that each of the $J_{(T_i)}$ and the $J_{(L_j)}$ has
only notational meaning and does not define an element of $H$ or $H/N$
by itself.)  For example, see \fullref{fig:example-univ-inv}.
\labellist\tiny
\pinlabel {$\wbar{\alpha}_c$} [br] at 170 393
\pinlabel {$\wbar{\alpha}_d$} [r] at 168 364
\pinlabel {$\wbar{\alpha}_e$} [br] at 164 322
\pinlabel {$\wbar{\beta}_c$} [l] at 183 394
\pinlabel {$\wbar{\beta}_d$} [bl] at 185 366
\pinlabel {$\wbar{\beta}_e$} [bl] at 186 333
\pinlabel {$T_1$} [r] at  161 310
\pinlabel {$\kappa$} [tr] at  202 315
\pinlabel {$S(\alpha_a)$} [r] at  246 389
\pinlabel {$\alpha_b$} [r] at 247 371
\pinlabel {$\beta_a$} [l] at  262 388
\pinlabel {$S(\beta_b)$} [l] at  263 370
\pinlabel {$\kappa$} [tr] at  272 333
\pinlabel {$1$} [br] at 274 411
\pinlabel {$L_1$} [l] at  295 351
\pinlabel {$1$} [b] at 194 416
\endlabellist
\FIG{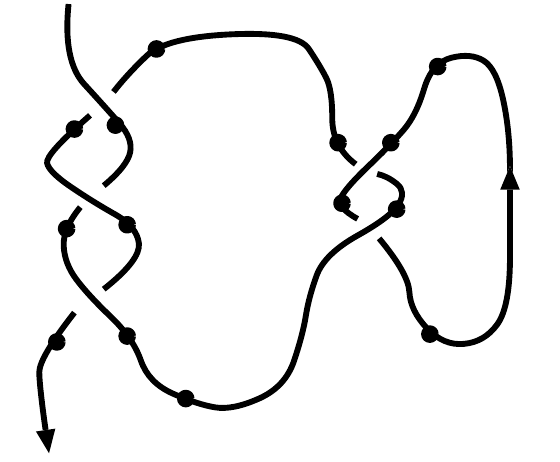}{For the tangle $T=T_1\cup L_1$, we have
$J_T=\sum\bar{\alpha}_e\bar{\beta}_d\bar{\alpha}_c1S(\alpha_a)S(\beta_b)
  \kappa\bar{\beta}_e\bar{\alpha}_d\bar{\beta}_c
  \otimes\tr_q(\kappa^{-1}\kappa\alpha_b\beta_a1)$,
where
$R=\sum\alpha_a\otimes\beta_a$, and $R^{-1}=\sum\bar{\alpha }_c\otimes \bar{\beta }_c$, etc.}{height=35mm}

Now we check that $J_T$ does not depend on where we start reading the
elements on the closed components.  Let $L_j$ be a closed component of
$T$ and let $x_1,\ldots,x_r$ be the elements read off from $L_j$.  Then we
have formally $J_{(L_j)}=\tr_q(\kappa ^{-1}x_1x_2\cdots x_r)$.  If we start from
$x_2$, then the right hand side becomes
\begin{gather*}
  \tr_q(\kappa ^{-1}x_2\cdots x_rx_1)
  =\tr_q(S^{-2}(x_1)\kappa ^{-1}x_2\cdots x_r)
  =\tr_q(\kappa ^{-1}x_1x_2\cdots x_r).
\end{gather*}
It follows that $J_{(L_i)}$ does not depend on where we
start reading the elements.

Now we can follow Ohtsuki's arguments \cite{Ohtsuki:93} to check that
$J_T$ does not depend on how we decompose $T$ into copies of
$\downarrow $, $\uparrow $, $\psi _{\downarrow ,\downarrow }^{\pm 1}$,
$\ev_{\downarrow} $, $\ev_{\uparrow} $, $\coev_{\downarrow} $ and
$\coev_{\uparrow} $, and that $J_T$ defines an isotopy invariant of framed,
oriented, ordered tangles.

It is convenient to generalize the above definition to the case where
$T$ is given as pasting of copies of the tangles
$\downarrow $, $\uparrow $,
$\psi _{a,b}^{\pm 1}$ ($a,b\in \{\downarrow ,\uparrow \}$),
$\ev_{\downarrow} $, $\ev_{\uparrow} $, $\coev_{\downarrow} $
and $\coev_{\uparrow} $.  In this case, we put elements of $H$ on the
components of $T$ as depicted in Figures \ref{fig:placing3} and
\ref{fig:placing4}.  Then $J_T$ is defined in the same way as above.
We can check that $J_T$ is well defined as follows.  For each tangle
diagram $T$ in \fullref{fig:placing4}, we choose a tangle diagram
$T'$ isotopic to $T$ obtained by pasting copies of $\downarrow $, $\uparrow $,
$\psi _{\downarrow ,\downarrow }^{\pm 1}$, $\ev_{\downarrow} $,
$\ev_{\uparrow} $, $\coev_{\downarrow} $ and $\coev_{\uparrow} $.
Then we can verify that $J_{T'}$ in the first definition, is equal to
$J_T$ in the second definition given by \fullref{fig:placing4}.
For example, consider the second tangle in \fullref{fig:placing4}.
Then $T$ and $T'$ are as depicted in \fullref{fig:placing5}.
\labellist\tiny
\pinlabel {$\alpha$} [tr] at 47 473
\pinlabel {$S(\beta)$} [tl] at 67 473
\pinlabel {$\wbar{\alpha}$} [br] at 147 454
\pinlabel {$\wbar{\beta}$} [bl] at 162 454
\pinlabel {$\kappa$} [t] at 170 445
\pinlabel {$\kappa^{-1}$} [b] at 138 482
\small
\pinlabel {$T$} [b] at 55 415
\pinlabel {$T'$} [b] at 150 415
\endlabellist
\FIGn{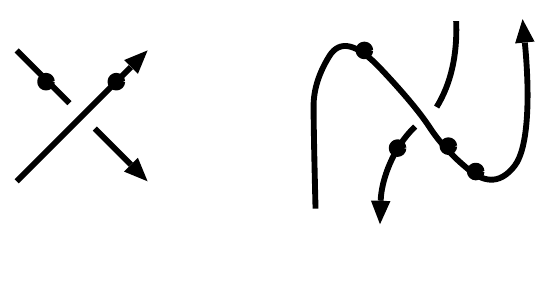}{}{height=25mm}
We have
\begin{equation*}
  J_{T'}
  =\sum\bar{\alpha }\otimes \kappa \bar{\beta }\kappa ^{-1}
  =\sum \alpha \otimes \kappa S^{-1}(\beta )\kappa ^{-1}
  =\sum \alpha \otimes S(\beta )
  =J_T,
\end{equation*}
where we used $\sum\bar{\alpha }\otimes \bar{\beta }=\sum\alpha \otimes S^{-1}(\beta )$.
The other cases can be similarly proved.

\begin{remark}
  \label{thm:27}
  To study invariants of tangles, it is sometimes useful to define a
  {\em functorial} invariant.  One can modify Kauffman and Radford's
  functorial universal regular isotopy invariant
  \cite{Kauffman-Radford:01} to define a functorial universal
  invariant defined on $\modT $, ie, a braided functor
  $F\col \modT \rightarrow \operatorname{Cat}(H)$ of $\modT $ into a category
  $\operatorname{Cat}(H)$ defined as in \cite{Kauffman-Radford:01}.
  However, we do not do so here, since in the present paper we are
  interested in {\em ordered} links.  Note that the categories $\modT $
  and $\operatorname{Cat}(H)$ do not care about the order of the
  circle components.  One can still define a functorial universal
  invariant which distinguishes circle components by using the
  category of colored tangles, but it would cause unnecessary
  complication and we do not take this approach here.
\end{remark}

\subsection{Effect of closure operation}
\label{sec:univ-invar-tangl}

In Ohtsuki's definition \cite{Ohtsuki:93} of his version of the
universal invariant, the universal trace $H\rightarrow H/I$, with
$I=\Span_{\modk} \{xy-yx\ver x,y\in H\}$, is used.
For our purposes, the universal quantum trace is more natural and
more useful than the universal trace.
Note that $I$ is not a left $H$--submodule of
$H$ in general.  The following proposition shows another reason why
the universal quantum trace is more convenient.

\begin{proposition}
  \label{thm:28}
  If $T\in \BT_n$, then we have
  \begin{equation*}
    J_{\clo(T)} = \tr_q^{\otimes n}(J_T).
  \end{equation*}
\end{proposition}

\begin{proof}
  Set $L=\clo(T)=L_1\cup \cdots \cup L_n$.
  We write
  \begin{align*}
    J_T &= \sum J_{(T_1)}\otimes \cdots\otimes J_{(T_n)}\in H^{\otimes n},\\
    J_{L} &= \sum J_{(L_1)}\otimes \cdots\otimes J_{(L_n)}\in (H/N)^{\otimes n},
  \end{align*}
  For $i=1,\ldots,n$, the part $J_{(L_i)}$ is computed as follows.  The
  diagram of $L_i$ is divided into the diagram of $T_i$ and the
  diagram of $\ev_{\uparrow} $.  See \fullref{fig:closure}.
\labellist\small
\pinlabel {$J_{(T_i)}$} [br] at 31 59
\pinlabel {$T_i$} at 95 60
\pinlabel {$\ev_{\uparrow}$} at 95 20
\pinlabel {$\kappa$} [tr] at 32 12
\endlabellist
  \FIGn{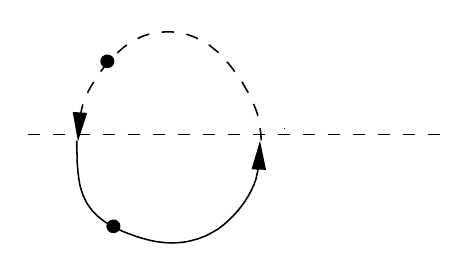}{}{height=20mm} Since we have $J_{\ev_{\uparrow} }=\kappa $, it
  follows that
  \begin{equation*}
    J_{(L_i)}
    = \tr_q(\kappa ^{-1}J_{\ev_{\uparrow} }J_{(T_i)})
    = \tr_q(\kappa ^{-1}\kappa J_{(T_i)})
    =\tr_q(J_{(T_i)}).
  \end{equation*}
  This implies the assertion.
\end{proof}

In \fullref{sec10}, we give a definition of a more
refined version of the universal invariants of links.

\subsection{Duplication, removal, and orientation-reversal}
\label{sec:dupl-remov-orient}

Let $T$ be a tangle and let $T_i$ be an arc component of $T$.  Define
a $\modk $--module homomorphism $\tilde{S}_{T_i}\col H\rightarrow H$ by
\begin{equation*}
  \tilde{S}_{T_i}(x) = \kappa ^{-r(T_i)}S(x)\kappa ^{s(T_i)}\quad \text{for $x\in H$},
\end{equation*}
where $r(T_i)=0$ if $T_i$ starts at the top and $r(T_i)=1$ otherwise,
and $s(T_i)=0$ if $T_i$ ends at the bottom, and $s(T_i)=1$ otherwise.
For example, we have
\begin{align*}
  \tilde{S}_{\downarrow} (x)&=S(x), &
  \tilde{S}_{\uparrow} (x)&=\kappa ^{-1}S(x)\kappa =S^{-1}(x),\\
  \tilde{S}_{\incl{.5em}{evdown}}(x)
  &=\tilde{S}_{\incl{.5em}{evup}}(x)=S(x)\kappa , &
  \tilde{S}_{\incl{.5em}{coevdown}}(x)
  &=\tilde{S}_{\incl{.5em}{}}(x)
  =\kappa ^{-1}S(x)
\end{align*}
for $x\in H$.
We need the following result, which is almost standard.

\begin{lemma}
  \label{lem:11}
  Let $T=T_1\cup \cdots \cup T_n$ be a tangle with $n$ arcs with $n\ge 1$.  For
  $i=1,\ldots,n$, let $\Delta _i(T)$ (resp.  $\epsilon _i(T)$, $S_i(T)$) denote the
  tangle obtained from~$T$ by duplicating (resp. removing,
  orientation-reversing) the $i$th component $T_i$.  Then we have
  \begin{align}
    \label{eq:36}
    J_{\Delta _i(T)} &= (1^{\otimes (i-1)}\otimes \Delta \otimes 1^{\otimes (n-i)})(J_T),\\
    \label{eq:38}
    J_{\epsilon _i(T)} &= (1^{\otimes (i-1)}\otimes \epsilon \otimes 1^{\otimes (n-i)})(J_T),\\
    \label{eq:41}
    J_{S_i(T)} &= (1^{\otimes (i-1)}\otimes \tilde{S}_{T_i}\otimes 1^{\otimes (n-i)})(J_T).
  \end{align}
\end{lemma}

\begin{proof}
  The cases of $\Delta _i(T)$ and $\epsilon _i(T)$ are standard.  We prove
  the case of $S_i(T)$, which may probably be well known to the
  experts but does not seem to have appeared in a way as general as
  here.

  We can easily check \eqref{eq:41} for
  $T=\psi _{\modb ,\modb },\psi _{\modb ,\modb }^{-1},\evdn,\evup,\cvdn,\cvup$.
  For the general case, we express $T$ as an iterated composition and
  tensor product of finitely many copies of the morphisms $\downarrow $, $\uparrow $,
  $\psi _{\downarrow ,\downarrow }^{\pm 1}$, $\evdn$, $\evup$, $\cvdn$, $\cvup$.  We may
  assume that $T_i$ involves at least one crossing or critical point,
  since otherwise the assertion is obvious.
  We decompose the component $T_i$ into finitely many intervals
  $T_{i,1},\ldots,T_{i,p}$ with $p\ge 1$, where
  \begin{itemize}
  \item if one goes along $T_i$ in the opposite direction to the
    orientation, then one encounter the intervals in the order
    $T_{i,1},\ldots,T_{i,p}$, and
  \item for each $j=1,\ldots,p$, there is just one crossing or critical
    points in $T_{i,j}$.
  \end{itemize}
  For each $j=1,\ldots,p$, let $x_j=J_{(T_{i,j})}$ denote the formal
  element put on the interval $T_{i,j}$ in the definition of $J_T$.
  Then we have $J_{(T_i)}=x_1x_2\cdots x_p$.  Let $-T_{i,j}$ denote the
  orientation reversal of $T_{i,j}$.  Then it follows from the cases
  of $T=\psi _{\modb ,\modb },\psi _{\modb ,\modb }^{-1},\evdn,\evup,\cvdn,\cvup$ that the
  formal element put on $-T_{i,j}$ in the definition of $J_{S_i(T)}$
  is $\tilde{S}_{T_{i,j}}(x_j)$.  Hence we have $J_{(-T_i)}=x'_p\cdots x'_1$,
  where
  \begin{equation*}
    x'_j= J_{(-T_{i,j})}
    =\tilde{S}_{T_{i,j}}(x_j)
    =\kappa ^{-r(T_{i,j})}S(x_j)\kappa ^{s(T_{i,j})}
  \end{equation*} for
  $j=1,\ldots,p$.  We have $s(T_{i,j})=r(T_{i,j-1})$ for $j=2,\ldots,p$.  We
  also have $s(T_{i,1})=s(T_i)$ and $r(T_{i,p})=r(T_i)$.  Hence it
  follows that
  \begin{align*}
      J_{(-T_i)}
      &=x'_px'_{p-1}\ldots x'_1\\
      &=\bigl(\kappa ^{-r(T_{i,p})}S(x_p)\kappa ^{s(T_{i,p})}\bigr)
      \bigl(\kappa ^{-r(T_{i,p-1})}S(x_{p-1})\kappa ^{s(T_{i,p-1})}\bigr)\\
      &\qquad\qquad
      \ldots \bigl(\kappa ^{-r(T_{i,1})}S(x_1)\kappa ^{s(T_{i,1})}\bigr)\\
      &=\kappa ^{-r(T_i)}S(x_p)
      S(x_{p-1})
      \cdots S(x_1)\kappa ^{s(T_i)}\\
      &=\kappa ^{-r(T_i)}S(x_1\ldots x_{p-1}x_p)
      \kappa ^{s(T_i)}\\
      &=\kappa ^{-r(T_i)}S(J_{(T_i)})
      \kappa ^{s(T_i)}\\
      &=\tilde{S}_{T_i}(J_{(T_i)}).
  \end{align*}
  Now the assertion immediately follows.
\end{proof}

\section{The braided functor $\modJ \col \modB \rightarrow \Mod_H$}
\label{sec8}
In this section, we fix a ribbon Hopf algebra $H$ over a commutative,
unital ring~$\modk $.

\subsection{The category $\Mod_H$ of left $H$--modules}
\label{sec:category-modh-left}
In this subsection, we recall some algebraic facts about the category
$\Mod_H$ of left $H$--modules.  For details, see
Majid \cite{Majid:algebras,Majid}.

Let $\Mod_H$ denote the category of left $H$--modules and left
$H$--module homomorphisms.  The category $\Mod_H$ is equipped with a
(non-strict) monoidal category structure with the tensor functor
$\otimes \col \Mod_H\times \Mod_H\rightarrow \Mod_H$ given by tensor product over $\modk $ with
the usual left $H$--module structure defined using comultiplication.
The unit object is $\modk $ with the trivial left $H$--module structure.
The braiding $\psi _{V,W}$ and its inverse of two objects $V$ and $W$ are
given by
\begin{gather}
  \label{eq:19}
  \psi _{V,W}(v\otimes w) = \sum \beta w \otimes  \alpha v,\qquad
  \psi _{V,W}^{-1}(v\otimes w) = \sum \bar{\alpha }w \otimes  \bar{\beta }v,
\end{gather}
for $v\in V$, $w\in W$.

We regard $H$ as a left $H$--module using the adjoint action
$\ad=\trr\col H^{\otimes 2}\rightarrow H$.  By \eqref{eq:19}, the braiding
$\psi _{H,H}\col H^{\otimes 2}\rightarrow H^{\otimes 2}$ and its inverse
are given by
\begin{gather*}
  \psi _{H,H}(x\otimes y) = \sum (\beta \trr y) \otimes  (\alpha \trr x),\qquad
  \psi _{H,H}^{-1}(x\otimes y) = \sum (\bar{\alpha }\trr y) \otimes  (\bar{\beta }\trr x)
\end{gather*}
for $x,y\in H$.

The {\em transmutation} \cite{Majid:algebras,Majid} of a
quasitriangular Hopf algebra $H$ is a Hopf algebra
$\bH=(H,\mu ,\eta ,\bD,\epsilon ,\bS)$ in the braided category $\Mod_H$, which is
obtained by modifying the Hopf algebra structure of $H$ as follows.
The algebra structure morphisms $\mu $ and $\eta $, and the counit $\epsilon $ of
$\bH$ are the same as those of $H$.  The comultiplication
$\bD\col H\rightarrow H^{\otimes 2}$ and the antipode $\bS\col H\rightarrow H$ are defined by
\begin{align}
  \label{eq:8}\bD(x) &= \sum x_{(1)}S(\beta )\otimes (\alpha \trr x_{(2)}),\\
  \label{eq:9}\bS(x) &= \sum \beta S(\alpha \trr x),
\end{align}
for $x\in H$.  The morphisms $\mu ,\eta ,\bD,\epsilon ,\bS$ are all left $H$--module
homomorphisms, and $\bH$ is a Hopf algebra in the braided category
$\Mod_H$.

Define $c^H_{\pm} \in H^{\otimes 2}$ by
\begin{equation}
  \label{eq:31}
  c^H_{\pm} =(S\otimes 1)((R_{21}R)^{\pm 1}).
\end{equation}
By abuse of notation, we denote by $c^H_{\pm} $ the $\modk $--module
homomorphism $\modk \rightarrow H^{\otimes 2}$ which maps $1$ to
$c^H_{\pm} $.

Using $\Delta (x)R_{21}R=R_{21}R\Delta (x)$, $x\in H$, one can verify
\begin{equation}
  \label{e12}
  c^H_{\pm} \in \Mod_H(\modk ,H^{\otimes 2}).
\end{equation}

\subsection{Definition of $\modJ \col \modB \rightarrow \Mod_H$}
\label{sec:braided-functor-h}
In this subsection, we define a braided functor
\begin{equation}
  \label{eq:16}
  \modJ \col \modB \rightarrow \Mod_H,
\end{equation}
which maps $\modb \in \Ob(\modB )$ to $H\in \Ob(\Mod_H)$.

For $T\in \modB (m,n)$ with $m,n\ge 0$, we define a $\modk $--module homomorphism
\begin{equation*}
  \modJ (T)\col H^{\otimes m}\rightarrow H^{\otimes n}
\end{equation*}
as follows.  Consider a tangle diagram of $T\eta _m$, see
\fullref{fig:FT}.
\labellist\tiny
\pinlabel {$x_1$} [r] at 198 610
\pinlabel {$x_2$} [r] at 243 609
\pinlabel {$1$} [tl] at 233 568
\pinlabel {$1$} [br] at 294 571
\pinlabel {$1$} [br] at 279 496
\pinlabel {$S(\alpha\!)$} [bl] at 306 571
\pinlabel {$\alpha'$} [r] at 315 535
\pinlabel {$\beta$} [tl] at 334 569
\pinlabel {$S(\!\beta'\!)$} [tl] at 327 529
\pinlabel {$1$} [br] at 388 491
\pinlabel {$x_3$} [r] at 396 610
\small
\pinlabel {$\eta_3$} [l] at 438 610
%\pinlabel {$S_0^{\otimes 3}$} [l] at 438 610
\pinlabel {$T$} [l] at 438 510
\endlabellist
\FIG{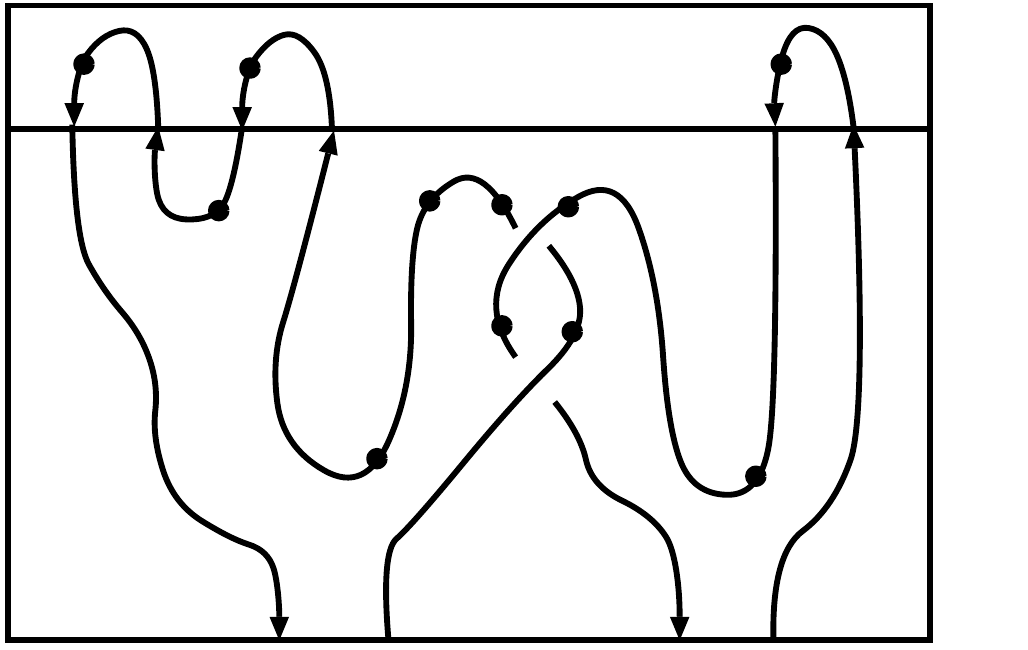}{For the tangle
$T=(\mu _{\modb} \otimes \mu _{\modb} )(\mu _{\modb} \otimes c_+\otimes \modb )\col \modb ^{\otimes 3}\rightarrow \modb ^{\otimes 2}$ depicted, and
$\sum x_1\otimes x_2\otimes x_3\in H^{\otimes 3}$, we have $\modJ (T)(\sum x_1\otimes x_2\otimes x_3) =
\sum x_1x_2S(\alpha )S(\beta ')\otimes \alpha '\beta x_3$, where
$R=\sum\alpha \otimes \beta =\sum\alpha '\otimes \beta '$.}{height=35mm}
Given an element
$\sum x_1\otimes \cdots\otimes x_m\in H^{\otimes m}$, we put $x_i$ on the $i$th component in
$\eta _m$ for each $i=1,\ldots,m$.  Moreover, we put elements in $H$
to the components in $T$ as in the definition of $J_T$.  Then we obtain a
tangle diagram consisting of $n$ arcs, decorated with elements of $H$.
For $i=1,\ldots,n$, let $y_i$ denote the word obtained by reading the
elements on the $i$th component of $T\eta _m$.  Then we set
\begin{equation*}
  \modJ (T)\Bigl(\sum x_1\otimes \cdots\otimes x_m\Bigr) = \sum y_1\otimes \cdots\otimes y_n.
\end{equation*}
Clearly, $\modJ (T)$ is a $\modk $--module homomorphism, and does not depend on
the choice of the diagram of $T$.  Note that if $m=0$, then we have
\begin{equation*}
  \modJ (T)(1) = J_T.
\end{equation*}
It is also clear that $\modJ (TT')=\modJ (T)\modJ (T')$
for any two composable pair of morphisms $T$ and $T'$ in $\modB $, and
that $\modJ (1_{\modb} ^{\otimes n})=1_{H^{\otimes n}}$.  This means that the correspondence
$T\mapsto \modJ (T)$ defines a functor
\begin{equation*}
  \modJ _{\modk} \col \modB \rightarrow \Mod_{\modk} ,
\end{equation*}
where $\Mod_{\modk} $ denotes the category of $\modk $--modules and $\modk $--module
homomorphisms.  We give the category $\Mod_{\modk} $ the standard symmetric
monoidal category structure.  Then we can easily check that $\modJ
_{\modk} $ is
a monoidal functor.

To prove that $\modJ _{\modk} $ lifts along the forgetful functor
$\Mod_H\rightarrow \Mod_{\modk} $ to a monoidal functor \eqref{eq:16}, it suffices to
show that if $T$ is a morphism in $\modB $, then $\modJ (T)$ is a left
$H$--module homomorphism.  By \fullref{thm:2}, we have only to
check this property for
$T\in \{\psi _{\modb ,\modb }^{\pm 1},\mu _{\modb} ,\eta _{\modb} ,v_{\pm}
,c_{\pm} \}$.  This follows from
\fullref{thm:7} below, since
$\eta ,v^{\pm 1},c_{\pm} ^H,\mu ,\psi _{H,H}^{\pm 1}$ are left $H$--module
homomorphisms.  \fullref{thm:7} also shows that $\modJ $ is a
braided functor.

\begin{proposition}
  \label{thm:7}
  We have
  \begin{gather*}
    \modJ (\eta _{\modb} )= \eta ,\quad
    \modJ (v_{\pm} )= \modr ^{\pm 1},\quad
    \modJ (c_{\pm} )= c^H_{\pm} ,\quad
    \modJ (\mu _{\modb} ) = \mu ,\quad
    \modJ (\psi_{\modb,\modb}^{\pm 1}) = \psi _{H,H}^{\pm 1},\quad
  \end{gather*}
  for $x,y\in H$.  Here, by abuse of notation, we denote by $\modr ^{\pm 1}$
  the corresponding morphism in $\Mod_H(\modk ,H)$.
\end{proposition}

\begin{proof}
  For $T=\eta _{\modb} , v_{\pm} , c_{\pm} ,\mu _{\modb} $, the homomorphism $\modJ (T)$ are easily
  computed using \fullref{fig:Fgen1} (a)--(f).
\labellist
\tiny
\hair=1pt
\pinlabel {\small(a)} [b] at 65 400
\pinlabel {\small(b)} [b] at 125 400
\pinlabel {$S(\alpha)$} [r] at 109 473
\pinlabel {$\beta$} [l] at 138 473
\pinlabel {$\kappa^{-1}$} [bl] at 130 488
\pinlabel {\small(c)} [b] at 180 400
\pinlabel {$\wbar\alpha$} [r] at 167 458
\pinlabel {$S(\wbar\beta)$} [bl] at 191 458
\pinlabel {$\kappa^{-1}$} [bl] at 184 487
\pinlabel {\small(d)} [b] at 250 400
\pinlabel {$S(`\alpha`)$} [tr] at 245 471
\pinlabel {$\alpha'$} [tr] at 241 454
\pinlabel {$S(`\beta'`)$} [tl] at 257 452
\pinlabel {$\beta$} [tl] at 262 474
\pinlabel {\small(e)} [b] at 325 400
\pinlabel {$\wbar\alpha$} [r] at 315 460
\pinlabel {$S(`\wbar\alpha'`)$} [br] at 319 441
\pinlabel {$\wbar\beta'$} [l] at 337 442
\pinlabel {$S(\wbar\beta)$} [bl] at 333 461
\pinlabel {\small(f)} [b] at 395 400
\pinlabel {$x_1$} [br] at 372 489
\pinlabel {$x_2$} [br] at 402 489
\pinlabel {\small(g)} [b] at 510 400
\pinlabel {$x$} [r] at 481 499
\pinlabel {$y$} [r] at 528 499
\pinlabel {$(1{\otimes}S)(\Delta(\alpha))$} [tr] at 492 471
\pinlabel {$(1{\otimes}S)(\Delta(\beta))$} [tl] at 526 471
\endlabellist
  \FIGn{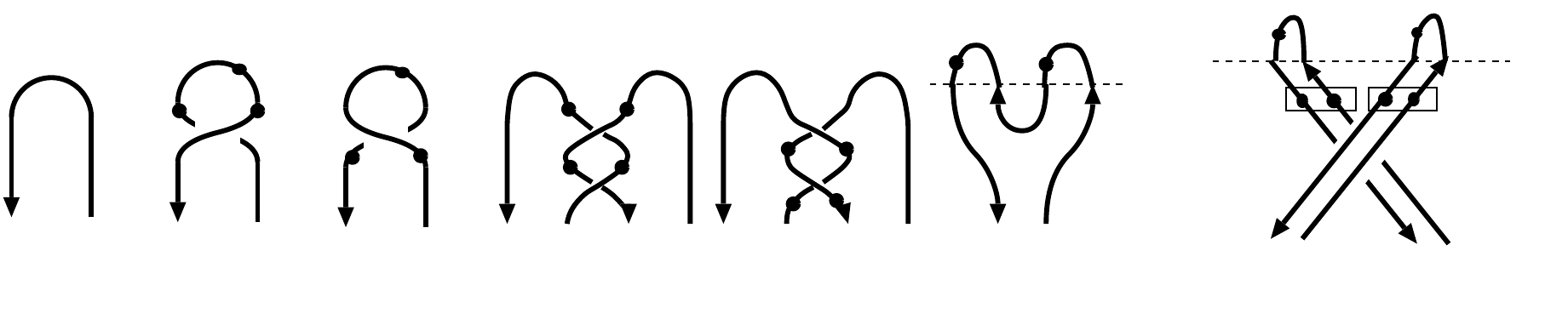}{}{height=26mm}
  The case $T=\psi _{H,H}$ is computed using
  \fullref{fig:Fgen1} (g), where
$$
\labellist\tiny
\pinlabel {$x{\otimes}y$} [r] at 142 519
\endlabellist
\incl{2.5em}{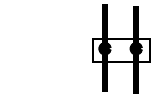}
\quad\raisebox{1.2em}{means}\quad
\labellist\tiny
\pinlabel {$x$} [r] at 247 522
\pinlabel {$y$} [r] at 260 522
\endlabellist
\incl{2.5em}{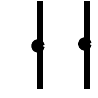}.$$
We have
\begin{align*}
\modJ (\psi _{\modb ,\modb })(x\otimes y)
&=\sum\beta _{(1)}yS(\beta _{(2)})\otimes \alpha _{(1)}xS(\alpha _{(2)})\\
&=\sum(\beta \trr y)\otimes (\alpha \trr x)=\psi _{H,H}(x\otimes y),
\end{align*}
where we write
$(\Delta \otimes \Delta )(R)=\sum\alpha _{(1)}\otimes \alpha _{(2)}\otimes
\beta_{(1)}\otimes \beta_{(2)}$.  We can
similarly check the case $T=\psi _{\modb ,\modb }^{-1}$.
\end{proof}

An easy consequence of the braided functor $\modJ $ is the following,
which is essentially well known.

\begin{proposition}[See Kerler {\cite[Corollary 12]{Kerler:97}}]
  \label{thm:1}
  If $T\in \BT_n$, then we have $J_T\in \Mod_H(\modk ,H^{\otimes n})$.  In
  particular, if $T\in \BT_1$, then $J_T\in H$ is central.
\end{proposition}

\subsection{The functor $\modJ $ as a morphism of external Hopf algebras}
\label{sec:transmutation}

Note that, by \fullref{rem:3}, the Hopf algebra structure of the
transmutation $\bH=(H,\mu ,\eta ,\bD,\epsilon ,\bS)$ of $H$ determines an external
Hopf algebra $((H,\mu ,\eta ),F_{\bH})$ in the canonical way.

\begin{theorem}
  \label{thm:63}
  The braided functor $\modJ \col \modB \rightarrow \Mod_H$ maps the external Hopf algebra
   $(\modb ,F_{\modb} )$ in $\modB $ into the external Hopf algebra
$(H,F_{\bH})$ in
   $\Mod_H$ in the following sense.
  \begin{enumerate}
  \item $\modJ $ maps the algebra $(\modb ,\mu_{\modb} ,\eta_{\modb})$ into the algebra
    $(H,\mu ,\eta )$.
  \item By defining
  $\modJ '_{\HH^{\otimes m}}=\modJ \col \BT_m\rightarrow \Mod_H(\modk ,H^{\otimes m})$ for
  $m\ge 0$, we obtain a natural transformation $\modJ '\col F_{\modb}
\Rightarrow F_{\bH}$.
  \end{enumerate}
\end{theorem}

\begin{proof}
  The condition (1) follows immediately from \fullref{thm:7}.

  The condition (2) is equivalent to that, for any morphism
  $f\col \HH^{\otimes m}\rightarrow \HH^{\otimes n}$ in~$\modH $, the diagram
  \begin{equation}
    \label{eq:15}
    \begin{CD}
      \BT_m @>\modJ >> \Mod_H(\modk ,H^{\otimes m}) \\
      @VF_{\modb} (f)VV @VVF_{\bH}(f)V\\
      \BT_n @>\modJ >> \Mod_H(\modk ,H^{\otimes n})
    \end{CD}
  \end{equation}
  commutes.  It suffices to prove \eqref{eq:15} for $f$ in a set of
  generators of $\modH $ as a category.  Hence we can assume $f=g_{(i,j)}$
  with
  $g\in \{\psi_{\HH,\HH},\psi_{\HH,\HH}^{-1},
    \mu_{\HH},\eta_{\HH},\Delta_{\HH},\epsilon_{\HH},S_{\HH}\}$
  and $i,j\ge 0$.  The condition (1) implies that we have \eqref{eq:15}
  if $g=\psi _{\HH,\HH}^{\pm 1}$, $\mu_{\HH}$ or~$\eta_{\HH}$.  To prove the cases
  $g=\Delta_{\HH},\epsilon_{\HH},S_{\HH}$, it suffices to prove
  \begin{gather}
    \label{eq:69}J_{\check\Delta _{(i,j)}(T)}=(1^{\otimes i}\otimes \bD\otimes 1^{\otimes j})(J_T),\\
    \label{eq:70}J_{\check\epsilon _{(i,j)}(T)}=(1^{\otimes i}\otimes \epsilon \otimes 1^{\otimes j})(J_T),\\
    \label{eq:71}J_{\check S_{(i,j)}(T)}=(1^{\otimes i}\otimes \bS\otimes 1^{\otimes j})(J_T)
  \end{gather}
  for $T\in \BT_{i+j+1}$, $i,j\ge 0$.  Note that \eqref{eq:70}
  follows from \fullref{lem:11}.

  We write
  \begin{equation*}
    J_T
    =\sum J_{(T_1)}\otimes \cdots\otimes J_{(T_{i+j+1})}
    = \sum y \otimes  x \otimes  y',
  \end{equation*}
  where
  \begin{align*}
  x&=J_{(T_{i+1})}\in H,\\
  y&=J_{(T_1)}\otimes \cdots\otimes J_{(T_i)}\in H^{\otimes i},\\
  y'&=J_{(T_{i+2})}\otimes \cdots\otimes J_{(T_{i+j+1})}\in H^{\otimes j}.
  \end{align*}
  By \fullref{fig:transm}, we have
  \begin{equation*}
    \begin{split}
      J_{\check\Delta _{(i,j)}(T)}
      &=\sum y\otimes x_{(1)}S(\beta )\otimes \alpha _{(1)}x_{(2)}S(\alpha _{(2)})\otimes y'\\
      &=\sum y\otimes x_{(1)}S(\beta )\otimes (\alpha \trr x_{(2)})\otimes y'\\
      &=\sum y\otimes \bD(x)\otimes y'\\
      &=(1^{\otimes i}\otimes \bD\otimes 1^{\otimes j})(J_T),
    \end{split}
  \end{equation*}
  and
  \begin{equation*}
    \begin{split}
      J_{\check S_{(i,j)}(T)}
      &=\sum y\otimes \beta \kappa \alpha _{(2)}\kappa ^{-1}S(x)S(\alpha _{(1)})\otimes y'\\
      &=\sum y\otimes \beta S^2(\alpha _{(2)})S(x)S(\alpha _{(1)})\otimes y'\\
      &=\sum y\otimes \beta S(\alpha _{(1)}xS(\alpha _{(2)}))\otimes y'\\
      &=\sum y\otimes \beta S(\alpha \trr x)\otimes y'\\
      &=\sum y\otimes \bS(x)\otimes y'\\
      &=(1^{\otimes i}\otimes \bS\otimes 1^{\otimes j})(J_T).
    \end{split}
  \end{equation*}
\labellist\tiny
\pinlabel {$x_{(1)}$} [r] at 84 646
\pinlabel {$x_{(2)}$} [l] at 104 646
\pinlabel {$\alpha_{(1)}$} [tr] at 115 577
\pinlabel {$S(`\alpha_{(2)}`)$} [bl] at 127 582
\pinlabel {$S(\beta)$} [tl] at 167 580
\pinlabel {$\kappa^{-1}`S(x)$} [r] at 311 646
\pinlabel {$S(`\alpha_{(1)}`)$} [r] at 317 574
\pinlabel {$\alpha_{(2)}$} [bl] at 335 572
\pinlabel {$\beta$} [tl] at 364 575
\pinlabel {$\kappa$} [t] at 376 538
\endlabellist
  \FIGn{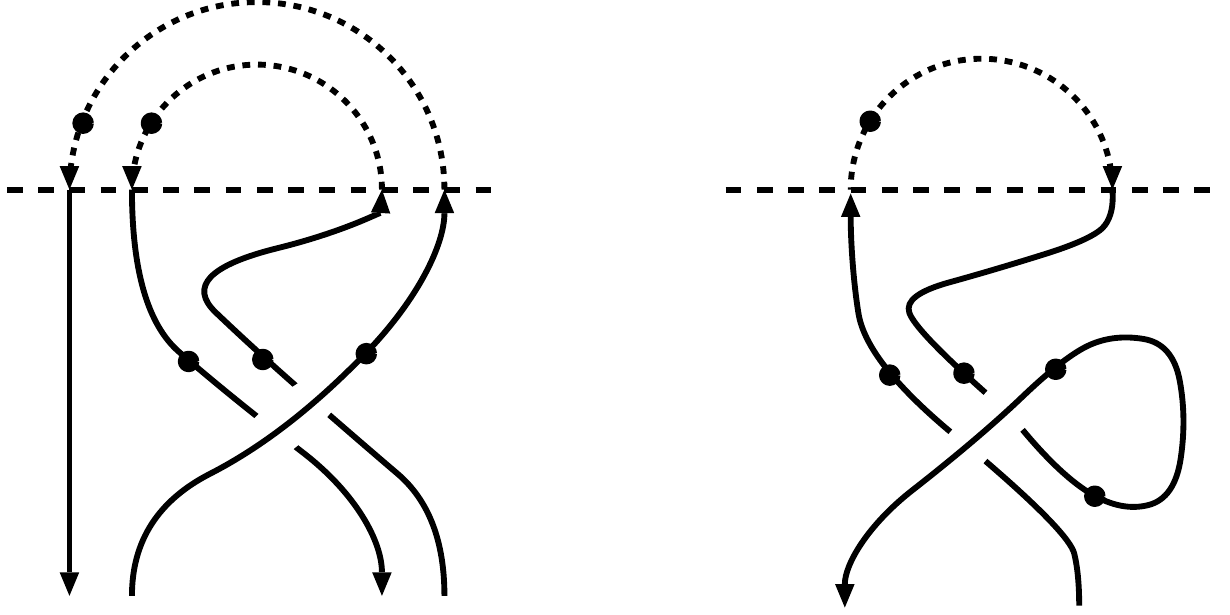}{}{height=35mm}
  Hence we have \eqref{eq:69} and \eqref{eq:71}.
\end{proof}

\fullref{thm:63} implies the relationship between the Hopf algebra
action on the bottom tangles and the functor $\modJ $ that we mentioned in
the latter half of \fullref{sec:hopf-algebra-action-1}.

\subsection{Topological proofs of algebraic identities}
\label{sec:topol-proofs-algebr}

\fullref{thm:63} means that, to a certain extent, the braided Hopf
algebra structure of $\bH$ is explained in terms of the external Hopf
algebra structure in $\modB $, which is defined topologically.  Thus,
\fullref{thm:63} can be regarded as a {\em topological
interpretation} of transmutation of a ribbon Hopf algebra.  We explain
below that \fullref{thm:63} can be used in proving various
identities for transmutation using isotopy of tangles.

For a $\modk $--module homomorphism $f\col H^{\otimes m}\rightarrow H^{\otimes n}$ and $i,j\ge 0$, we set
\begin{equation*}
  f_{(i,j)} = 1_H^{\otimes i}\otimes f\otimes 1_H^{\otimes j}\col H^{\otimes (m+i+j)}\rightarrow H^{\otimes (n+i+j)}.
\end{equation*}
For $f,g\in \Mod_H(H^{\otimes m},H^{\otimes n})$, we write $f\equiv g$ if we have
\begin{equation*}
  f_{(i,j)}(J_T)=g_{(i,j)}(J_T)
\end{equation*}
for all $i,j\ge 0$ and $T\in \BT_{i+j+m}$.  Note that if $m=0$, then
$f\equiv g$ and $f=g$ are equivalent.

\begin{remark}
  \label{thm:21}
  {\em All the formulas of the form ``$f\equiv g$'' that appear in
  what follows can be replaced with ``$f=g$''.}  One can prove this fact
  either by direct computation or using the functor $\modJ ^{\B}$ mentioned
  in \fullref{sec:category-0} below.  In the present paper, we
  content ourselves with the weaker form ``$f\equiv g$''.
\end{remark}

Let $\adb\in \Mod_H(H^{\otimes 2},H)$ denote the left adjoint action for $\bH$
defined by
\begin{equation*}
  \adb = \mu ^{[3]}(1\otimes \psi _{H,H})((1\otimes \bS)\bD\otimes 1).
\end{equation*}
It is well known that $\ad=\adb$.  As a first example of topological
proofs, we show that $\adb\equiv\ad$.  Since \fullref{fig:ad} (a)
and (b) are isotopic, we see that, for $i,j\ge 0$ and $T\in \BT_{i+j+2}$,
\begin{equation*}
  J_{F_{\modb} ((\ad_{\HH})_{(i,j)})(T)}
  =\adb_{(i,j)}(J_T)
\end{equation*}
\labellist\small
\pinlabel {$A$} [r] at 11 140
\pinlabel {\tiny$x_{i+2}$} [b] at 90 140
\pinlabel {$\kappa$} [t] at 67 80
\endlabellist
is calculated using \fullref{fig:ad-calc}. \FIG{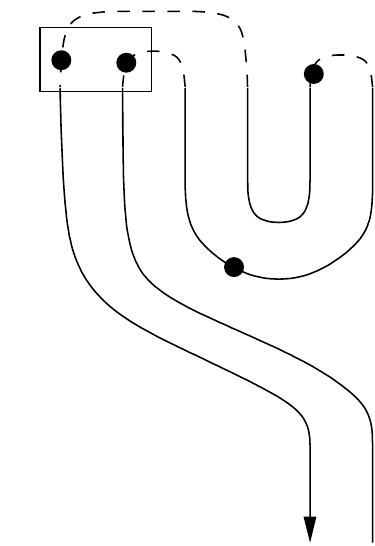}{The tangle
$F_{\modb} ((\ad_{\HH})_{(i,j)})(T)$ with only the $(i{+}1)$st component
depicted.  Here we write $J_T=\sum x_1\otimes \cdots\otimes
x_{i+j+2}$.  The box labeled ``A'' contains the tensor $\sum(x_{i+1})_{(1)}\otimes
\kappa ^{-1}S((x_{i+1})_{(2)})$.}{height=30mm} Hence we have
\begin{equation*}
  \begin{split}
    &\adb_{(i,j)}(J_T)\\
    &\quad=\sum
  x_1\otimes \cdots\otimes x_i\otimes (x_{i+1})_{(1)}x_{i+2}S((x_{i+1})_{(2)})
  \otimes x_{i+3}\otimes \cdots\otimes x_{i+j+2}\\
  &\quad =\sum
  x_1\otimes \cdots\otimes x_i\otimes (x_{i+1}\trr x_{i+2})\otimes x_{i+3}\otimes \cdots\otimes x_{i+j+2}\\
  &\quad =\ad_{(i,j)}(J_T).
\end{split}
\end{equation*}
Hence we have $\adb\equiv\ad$.

Let $\coadb\in \Mod_H(H,H^{\otimes 2})$ denote the left adjoint coaction for
$\bH$, defined by
\begin{equation*}
  \coadb = (\mu \otimes 1)(1\otimes \psi _{H,H}(1\otimes \bS))\bD^{[3]}.
\end{equation*}
By taking $J$ of \eqref{eq:7}, we have
\begin{equation*}
  \coadb_{(i,j)}(J_T)  =(\modJ (\gamma _+))_{(i,j)}(J_T)
  = (1_H\otimes \ad)_{(i,j)}(c^H_+\otimes 1_H)_{(i,j)}(J_T),
\end{equation*}
where we used $\ad\equiv\adb$.
Hence we have
\begin{equation}
  \label{eq:33}
  \modJ (\gamma _+) \equiv \coadb \equiv (1_H\otimes \ad)(c^H_+\otimes 1_H).
\end{equation}
By taking $J$ of \eqref{eq:17}, we have
\begin{equation*}
  \modJ (\gamma _-)_{(i,j)}(J_T) = (\bS\otimes 1_H)_{(i,j)}\modJ (\gamma _+)_{(i,j)}(J_T).
\end{equation*}
Hence we have
\begin{equation}
  \label{eq:32}
  \modJ (\gamma _-) \equiv (\bS\otimes 1_H)\coadb.
\end{equation}
We also note that the proof of external Hopf algebra axioms in $\modB $
yields a topological proof of the weak ``$\equiv$--version'' of the
identities in the axiom of Hopf algebra for $\bH$.  For example, one
can derive the formula
\begin{equation*}
  \bD\mu \equiv(\mu \otimes \mu )(1_H\otimes \psi _{H,H}\otimes 1_H)(\bD\otimes \bD)
\end{equation*}
from \fullref{fig:delta-mu}, and
\begin{equation*}
  \mu (1_H\otimes \bS)\bD\equiv\mu (\bS\otimes 1_H)\bD\equiv\eta \epsilon
\end{equation*}
from \fullref{fig:antipode}.

\section{Values of universal invariants of bottom tangles}
\label{sec9}
In this subsection, we study the set of values of universal invariants
of bottom tangles.  In \fullref{sec:values-jt}, we
give several general results, and in later subsections we give applications to
some specific cases.

We fix a ribbon Hopf algebra $H$ over a commutative, unital ring $\modk $.

\subsection{Values of $J_T$ }
\label{sec:values-jt}

We use the following notation.  Let $K\subset H^{\otimes m}$ and $L\subset H^{\otimes n}$ be
subsets.  Set
\begin{equation*}
  K\otimes L = \{x\otimes y\ver x\in K,y\in L\}\subset H^{\otimes (m+n)}.
\end{equation*}
If $x\in H$, then set
\begin{equation*}
  K\otimes x = K\otimes \{x\},\quad x\otimes K=\{x\}\otimes K.
\end{equation*}
The category $\modB $ acts on the left $H$--modules $H^{\otimes n}$, $n\ge 0$, by
the functions
\begin{equation*}
  j_{m,n}\col \modB (m,n)\times H^{\otimes m}\rightarrow H^{\otimes n},\quad
  (T,x)\mapsto Tx = \modJ (T)(x).
\end{equation*}
If $C$ is a subcategory of $\modB $, and if $K\subset \bigcup_{i\ge 0}H^{\otimes i}$, set
\begin{equation*}
    C\cdot K = \bigcup_{m,n\ge 0} j_{m,n}(C(m,n)\times (K\cap H^{\otimes m})).
\end{equation*}
Recall that $\modA $ denotes the braided subcategory of $\modB $ generated by
the object $\modb $ and the morphisms $\mu_{\modb} $ and $\eta_{\modb}$.

We have the following characterization of the possible values of the
universal invariants of bottom tangles.

\begin{theorem}
  \label{r5}
  The set $\{J_U\ver U\in \BT\}$ of the values of $J_U$ for all the bottom
  tangles $U\in \BT$ is given by
  \begin{equation}
    \label{eq:34}
    \{J_U\ver U\in \BT\}= \modA \cdot\{J_T\ver T\in \{v_{\pm} ,c_{\pm} \}^*\}.
  \end{equation}
\end{theorem}

\begin{proof}
  The result follows immediately from \fullref{thm:14} and
  functoriality of~$\modJ $.
\end{proof}

Using \fullref{r5}, we obtain the following, which will be useful
in studying the universal invariants of bottom tangles.

\begin{corollary}
  \label{thm:12}
  Let $K_i\subset H^{\otimes i}$ for $i\ge 0$, be subsets satisfying the following.
  \begin{enumerate}
  \item $1\in K_0$, $1,v^{\pm 1}\in K_1$, and $c^H_{\pm} \in K_2$.
  \item For $m,n\ge 0$, we have $K_m\otimes K_n\subset K_{m+n}$.
  \item For $p,q\ge 0$ we have
    \begin{align*}
      (\psi _{H,H}^{\pm 1})_{(p,q)}(K_{p+q+2})&\subset K_{p+q+2},\\
      \mu _{(p,q)}(K_{p+q+2})&\subset K_{p+q+1}.
    \end{align*}
  \end{enumerate}
  Then, for any $U\in \BT_n$, $n\ge 0$, we have $J_U\in K_n$.
\end{corollary}

\begin{proof}
  By (1) and (2), the $K_i$ contain $J_T$ for $T\in \{v_{\pm} ,c_{\pm} \}^*$.
  By (1), (2) and (3), the $K_i$ are invariant under the action of
  $\modA $.  Hence we have the assertion.
\end{proof}

Using \fullref{thm:4}, we see that, for $T\in \ABT$, the set of the
values of the universal invariant of bottom tangles obtained from
$\eta _n$ ($n\ge 0$) by one $T$--move is equal to $\cB_0\cdot\{J_T\}$.
Similarly, by \fullref{thm:19}, we see that, for $M\subset \ABT$
inversion-closed, the set of the values of the universal invariant of
the $M$--trivial bottom tangles is equal to
$\cB_0\cdot\{J_U\ver U\in M^*\}$.  From these observations we have the
following results, which will be useful in applications.

\begin{corollary}
  \label{thm:10}
  Let $T\in \ABT_m$, and let $K_i\subset H^{\otimes i}$ for $i\ge 0$ be subsets
  satisfying the following conditions.
  \begin{enumerate}
  \item $J_T\in K_m$.
  \item For $p,q\ge 0$ and
    $f\in \{\psi _{H,H}^{\pm 1},\eta ,\mu ,\coadb,(\bS\otimes 1)\coadb\}$ with
    $f\col \modb ^{\otimes i}\rightarrow \modb ^{\otimes j}$, we have
    \begin{equation*}
      f_{(p,q)}(K_{p+q+i})\subset K_{p+q+j}.
    \end{equation*}
  \end{enumerate}
  Then, for any $U\in \BT_n$ obtained from $\eta _n$ by one $T$--move, we
  have $J_U\in K_n$.
\end{corollary}

\begin{proof}
  We have to show that $\cB_0\cdot\{J_T\}\subset \bigcup_nK_n$.  Since
  $J_T\in K_m$, it suffices to show that $\bigcup_nK_n$ is stable under
  the action of generators of $\cB_0$.  Since $\cB_0$ is as a category
  generated by $f_{(p,q)}$ with $p,q\ge 0$ and
  $f\in \{\mu_{\modb},\eta_{\modb},\gamma_+,\gamma_-\}$, the condition (2) implies that
  $\bigcup_nK_n$ is stable under the action of $\cB_0$.  Here we use
  \eqref{eq:33} and \eqref{eq:32} with $\equiv$ replaced by $=$.  (See
  \fullref{thm:21}.)
\end{proof}

\begin{corollary}
  \label{thm:23}
  Let $M\subset \ABT$ be inversion-closed, and let $K_i\subset H^{\otimes i}$ for
  $i\ge 0$ be subsets satisfying the following conditions.
  \begin{enumerate}
  \item $1\in K_0$ and $1\in K_1$.
  \item For $V\in M$, we have $J_V\in K_{|V|}$.
  \item If $k,l\ge 0$, then we have $K_k\otimes K_l\subset K_{k+l}$.
  \item For $p,q\ge 0$ and
    $f\in \{\psi _{H,H}^{\pm 1},\mu ,\coadb,(\bS\otimes 1)\coadb\}$ with
    $f\col H^{\otimes i}\rightarrow H^{\otimes j}$, we have
    \begin{equation*}
      f_{(p,q)}(K_{p+q+i})\subset K_{p+q+j}.
    \end{equation*}
  \end{enumerate}
  Then, for any $M$--trivial $U\in \BT_n$, we have $J_U\in K_n$.
\end{corollary}

\begin{proof}
  It suffices to check the conditions in \fullref{thm:10}, where
  $T$ is an element of $M^*$.  The condition (1) in \fullref{thm:10}, ie, $J_T\in K_{|T|}$, follows from (1), (2) and (3).
  The condition (2) in \fullref{thm:10} with $f\neq\eta $ follows
  from~(4).  The condition (2) with $f=\eta $ in \fullref{thm:10}
  follows from (1), (3), and (4), since we have
$$\eta _{(p,q)}(K_{p+q}) =(1^{\otimes p}\otimes \psi _{H^{\otimes
q},H})(K_{p+q}\otimes 1).\proved$$
\end{proof}

The following will be useful in studying the set of $M$--equivalence
classes of bottom tangles for $M\subset \ABT$.

\begin{corollary}
  \label{thm:25}
  Let $M\subset \ABT$, not necessarily inversion-closed.  Let $K_i\subset H^{\otimes i}$
  for $i\ge 0$ be $\modZ $--submodules satisfying the following
  conditions.
  \begin{enumerate}
  \item For each $T\in M$ we have $J_T-1^{\otimes |T|}\in K_{|T|}$.
  \item For $i\ge 0$, we have
    \begin{equation*}
      K_i\otimes 1\subset K_{i+1},\quad
      K_i\otimes v^{\pm 1}\subset K_{i+1},\quad
      K_i\otimes c^H_{\pm} \subset K_{i+2}.
    \end{equation*}
  \item For $p,q\ge 0$, we have
    \begin{align}
      \label{eq:27}
      (\psi _{H,H}^{\pm 1})_{(p,q)}(K_{p+q+2})&\subset K_{p+q+2},\\
      \label{eq:28}
      \mu _{(p,q)}(K_{p+q+2})&\subset K_{p+q+1}.
    \end{align}
  \end{enumerate}
  Then, for any pair $U,U'\in \BT_n$ of $M$--equivalent
  bottom tangles, we have $J_{U'}-J_U\in K_n$.  Hence there is a
  well-defined function
  \begin{equation*}
    \BT_n/(\text{$M$--equivalence})\rightarrow H^{\otimes n}/K_n,\quad
    [U]\mapsto [J_U]
  \end{equation*}
  for each $n\ge 0$.
\end{corollary}

\begin{proof}
  Since $K_n\subset H^{\otimes n}$ are $\modZ $--submodules for $n\ge 0$, we may assume
  without loss of generality that $U$ and $U'$ are related by one
  $T$--move for $T\in M$.  Set $r=|T|$.  By \fullref{thm:5}, there is
  $W\in \modB (r,n)$ such that $U = W\eta _r$ and $U' = WT$.  Hence we have
  \begin{equation*}
    J_{U'}-J_U = \modJ (W)(J_T-1^{\otimes r}) \in \modJ (W)(K_r),
  \end{equation*}
  Therefore we have only to prove that $\modJ (W)(K_r)\subset K_n$.  By the
  assumptions, this holds for each generator $W$ of $\modB $ as a
  subcategory of $\modT $, described in \fullref{thm:36}.  Hence we
  have the assertion.
\end{proof}

\subsection{Unknotting number}
\label{sec:unknotting-number}
A {\em positive crossing change} is a local move on a tangle which
replaces a negative crossing with a positive crossing.  A {\em negative
crossing change} is the inverse operation.  In our terminology, a
positive (resp. negative) crossing change is equivalent to a
$c_-$--move (resp. $c_+$--move).

A {\em bottom knot} is a $1$--component bottom tangle.

\begin{corollary}
  \label{thm:15}
  Let $n_+,n_-\ge 0$, and let $K_i\subset H^{\otimes i}$, $i\ge 1$, be subsets
  satisfying the following conditions.
  \begin{enumerate}
  \item $(c^H_+)^{\otimes n_-}\otimes (c^H_-)^{\otimes n_+}\in K_{2(n_++n_-)}$.
  \item For $p,q\ge 0$ and
    $f\in \{\psi _{H,H}^{\pm 1},\eta ,\mu ,\coadb,(\bS\otimes 1)\coadb\}$ with
    $f\col \modb ^{\otimes i}\rightarrow \modb ^{\otimes j}$, we have
    \begin{equation*}
      f_{(p,q)}(K_{p+q+i})\subset K_{p+q+j}.
    \end{equation*}
  \end{enumerate}
  Then, if a bottom knot $T\in \BT_1$ of framing $0$
  is obtained from $\eta_{\modb} $ by $n_+$ positive crossing changes
  and $n_-$ negative crossing changes up to framing change, we have
  $J_T\in \modr ^{2(n_+-n_-)}K_1$.
\end{corollary}

\begin{proof}
  The result follows from \fullref{thm:10}, since $T$ is obtained
  from $\eta_{\modb} $ by a $(c_+^{\otimes n_-}\otimes c_-^{\otimes n_+})$--move and framing
  change by $-2(n_+-n_-)$.
\end{proof}

\fullref{thm:15} can be used to obtain an obstruction for a
bottom knot $T$ to be of unknotting number at most $n$, since a bottom
knot is of unknotting number~$n$ if and only if, for some $n_+,n_-\ge 0$
with $n_++n_-=n$, we have the situation in the statement of \fullref{thm:15}.  Also, we can use \fullref{thm:15} to obtain an
obstruction for a bottom tangle to be positively unknottable, ie,
obtained from~$\eta_{\modb}$ by finitely many positive crossing change up to
framing change.

In the literature, versions of unknotting numbers with respect to
various kinds of admissible local moves are studied, see, for example,
Murakami \cite{Murakami:sharp} or Murakami--Nakanishi
\cite{Murakami-Nakanishi}.  For $M\subset \ABT$, a
bottom knot $T\in \BT_1$ is said to be of ``$M$--unknotting number $n$''
if $T$ can be obtained from $\eta _1$ by $n$ applications of $M$--moves.
One can easily modify \fullref{thm:15} to give obstructions for
a bottom knot to be of $M$--unknotting number$\le n$.

\subsection{Commutators and Seifert surfaces}
\label{sec:commutators}

For any Hopf algebra $A$ in a braided category $\modM $, we define the
 {\em commutator morphism} \cite[Section 8.1]{Habiro:claspers}
 $Y_A\in \modM (A^{\otimes 2},A)$ by
\begin{gather}
  \label{eq:14}
  Y_A = \mu _A^{[4]}(A\otimes \psi _{A,A}\otimes A)
  (A\otimes S_A\otimes S_A\otimes A)(\Delta _A\otimes \Delta _A).
\end{gather}
Using adjoint action, we obtain a simpler formula, which is sometimes
more useful:
\begin{gather}
  \label{eq:6}
  Y_A = \mu _A(\ad_A\otimes A)(A\otimes S_A\otimes A)(A\otimes \Delta _A).
\end{gather}
The function
\begin{equation*}
  F_{\modb} ((Y_{\HH})_{(i,j)})\col \BT_{i+j+2}\rightarrow \BT_{i+j+1}
\end{equation*}
transforms a bottom tangle into another as illustrated in \fullref{fig:Y}.
\labellist\small
\pinlabel {(a)} [b] at 40 0
\pinlabel {\tiny $i{+}1$} [t] at 15 58
\pinlabel {\tiny $i{+}2$} [t] at 70 58
\pinlabel {(b)} [b] at 220 0
\pinlabel {\tiny $i{+}1$} [t] at 200 58
\pinlabel {$=$} at 302 140
\pinlabel {(c)} [b] at 385 0
\pinlabel {\tiny $i{+}1$} [t] at 415 58
\pinlabel {$=$} at 457 120
\pinlabel {(d)} [b] at 535 0
\pinlabel {\tiny $i{+}1$} [t] at 515 58
\pinlabel {(e)} [b] at 665 0
\pinlabel {\tiny $i{+}1$} [t] at 667 58
\pinlabel {\tiny surgery} [t] at 620 110
\endlabellist
\FIG{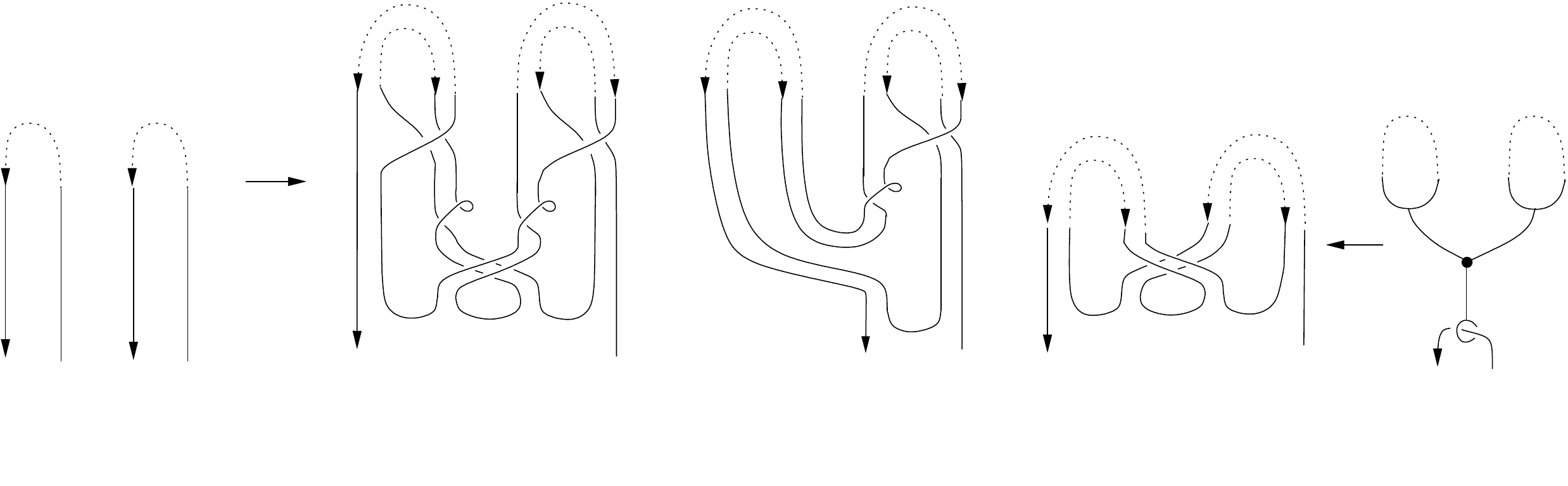}{(a) An $(i{+}j{+}2)$--component bottom tangle $T$.
(b) The $(i{+}j{+}1)$--component bottom tangle $T'=F((Y_{\HH})_{(i,j)})(T)$,
calculated using \eqref{eq:14}. (c) $T'$ calculated using
\eqref{eq:6}.  (d) Another picture of $T'$, in which the $(i+1)$st
component bounds a Seifert surface of genus $1$.  (d) A presentation
of $T'$ using a clasper.}{height=38mm}

For a ribbon Hopf algebra $H$, a (relatively) simple formula for the
left $H$--module homomorphism $Y_{\bH}$ is as follows.

\begin{proposition}
  \label{thm:13}
  For $\sum x\otimes y\in H^{\otimes 2}$, we have
  \begin{equation*}
    Y_{\bH}\Bigl(\sum x\otimes y\Bigr)=\sum \Bigl(x\trr \beta S((\alpha \trr y)_{(1)})\Bigr)(\alpha \trr y)_{(2)},
  \end{equation*}
  where $\Delta (\alpha \trr y)=\sum(\alpha \trr y)_{(1)}\otimes (\alpha \trr y)_{(2)}$.
\end{proposition}

\begin{proof}
  By computation, we have
  \begin{equation}
    \label{eq:10}
    (\bS\otimes 1)\bD(y) = \sum \beta S((\alpha \trr y)_{(1)})\otimes (\alpha \trr y)_{(2)}
  \end{equation}
  for $y\in H$.  Hence we have by \eqref{eq:6}
  \begin{equation*}
    \begin{split}
      Y_{\bH}\bigl(\sum x\otimes y\Bigr)
      =&\mu (\ad\otimes H)(1\otimes (\bS\otimes 1)\bD)\Bigl(\sum x\otimes y\Bigr)\\
      =& \sum (x\trr \beta S((\alpha \trr y)_{(1)}))(\alpha \trr y)_{(2)}.
    \end{split}
  \end{equation*}
This completes the proof.
\end{proof}

\begin{remark}
  \label{thm:16}
  \fullref{thm:13} holds also for the transmutation of a
  quasitriangular Hopf algebra $H$ which are not ribbon.
\end{remark}

A {\em Seifert surface} of a bottom knot $T$ in a cube
$[0,1]^3$ is a compact, connected, oriented surface $F$ in $[0,1]^3$
such that $\partial F=T\cup \gamma $ and $F\cap ([0,1]^2\times \{0\})=\gamma $, where
$\gamma \subset [0,1]^2\times \{0\}$ is the line segment with $\partial \gamma =\partial T$.  Note that a
Seifert surface of a bottom knot $T$ determines in the
canonical way a Seifert surface of the closure of $T$.

Recall that a link $L$ in $S^3$ is {\em boundary} if the components of
$L$ bounds mutually disjoint Seifert surfaces.  Similarly, a bottom
tangle $T\in \BT_n$ is said to be {\em boundary} if the components of
$T$ are of framing $0$ and bound mutually disjoint Seifert surfaces in
$[0,1]^3$.

\begin{theorem}
  \label{thm:17}
  Let $K_i\subset H^{\otimes i}$, $i\ge 0$, be as in \fullref{thm:12}.  Then, for
  any boundary bottom tangle $T=T_1\cup \cdots \cup T_n$ bounding mutually disjoint
  Seifert surfaces $F_1,\ldots,F_n$ of genus $g_1,\ldots,g_n$, we have
  \begin{equation*}
    J_T\in
    (\mu ^{[g_1]}\otimes \cdots\otimes \mu ^{[g_n]})Y_{\bH}^{\otimes (g_1+\cdots+g_n)}(K_{2(g_1+\cdots+g_n)}).
  \end{equation*}
  In particular, if a bottom knot bounds a Seifert surface of
  genus $g$, then we have
  \begin{equation*}
    J_T\in \mu ^{[g]}Y_{\bH}^{\otimes g}(K_{2g}).
  \end{equation*}
\end{theorem}

\begin{proof}
  Set $g=g_1+\cdots+g_n$.  By isotopy, we can arrange $F_1,\ldots,F_n$ as
  depicted in \fullref{fig:seifert}, where
  $D(T')\in \modT (\one ,\modb ^{\otimes 4g})$ is obtained from a bottom tangle
  $T'\in \BT_{2g}$ by doubling the components.
\labellist\small
\pinlabel {$T = $} [r] at 29 52
\pinlabel {$F_1$} [r] at 197 13
\pinlabel {$D(T')$} [t] at 218 80
\pinlabel {$F_n$} [r] at 395 13
\pinlabel {$\cdots$} at 118 28
\pinlabel {$\cdots$} at 118 58
\pinlabel {$\cdots$} at 218 8
\pinlabel {$\cdots$} at 218 28
\pinlabel {$\cdots$} at 218 58
\pinlabel {$\cdots$} at 316 28
\pinlabel {$\cdots$} at 316 58
\endlabellist
  \FIGn{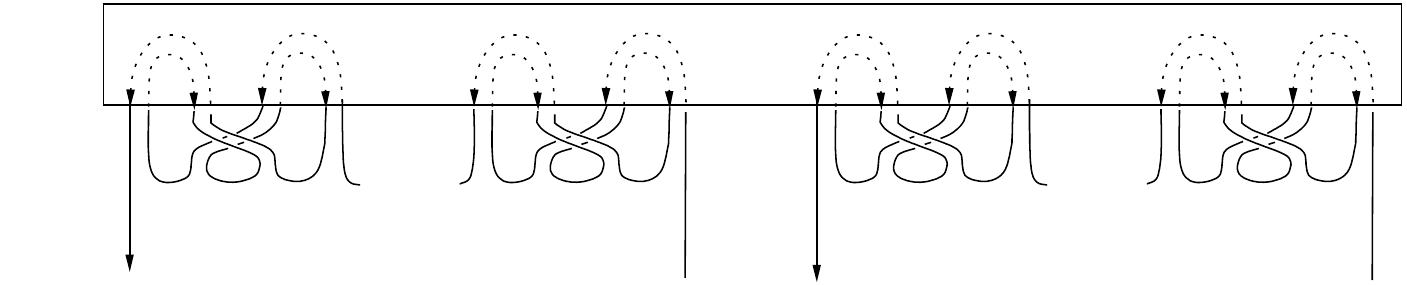}{}{height=23mm} (The surfaces bounded by the components
  of $T$ should be obvious from the figure.)  We have
  \begin{equation*}
    T = F_{\modb} ((\mu_{\HH}^{[g_1]}\otimes \cdots\otimes
\mu_{\HH}^{[g_n]})Y_{\HH}^{\otimes g})(T').
  \end{equation*}
  Hence
  \begin{equation*}
    \begin{split}
      J_T
      &= F_{\bH}((\mu_{\HH}^{[g_1]}\otimes \cdots\otimes
\mu_{\HH}^{[g_n]})Y_{\HH}^{\otimes g})(J_{T'})\\
      &= (\mu _H^{[g_1]}\otimes \cdots\otimes \mu
_H^{[g_n]})Y_{\bH}^{\otimes g}(J_{T'}).
    \end{split}
  \end{equation*}
  By \fullref{thm:12}, we have $J_{T'}\in K_{2g}$.  Hence we have
  the assertion.
\end{proof}

It is easy to verify that a link $L$ is boundary if and only if there
is a boundary bottom tangle $T$ such that the closure of $T$ is
equivalent to $L$.  (However, there are many non-boundary bottom
tangles whose closures are boundary.)  Hence we can use \fullref{thm:17} to obtain $J_L$ for boundary links $L$.  Also, the latter
part of \fullref{thm:17} can be used to obtain an obstruction for
a knot from being of genus$\le g$.

\subsection{Unoriented spanning surfaces}
\label{sec:unorient-spann}
Here we consider the ``unorientable version'' of the previous
subsection.

The {\em crosscap number} (see Clark \cite{Clark} and Murakami--Yasuhara
\cite{Murakami-Yasuhara}) of an unframed nontrivial knot $K$ is the
minimum number of the first Betti numbers of unorientable surfaces
bounded by $K$.  The crosscap number of an unknot is defined to be
$0$.

\begin{proposition}
  \label{r4}
  Let $T$ be a $0$--framed bottom knot of crosscap number $c\ge 0$ (ie,
  the closure of $T$ is of crosscap number $c$).  Then there is
  $T'\in \BT_c$ such that
  \begin{equation}
    \label{e8}
    J_T = \modr ^{4w(T')}(\mu \bD)^{\otimes c}(J_{T'}),
  \end{equation}
  where $w(T')\in \modZ $ is the writhe of the tangle $T'$.
\end{proposition}

\begin{proof}
  By assumption, the union of the bottom knot $T$ and the line segment
  bounded by the endpoints of $T$ bounds a connected, compact,
  unorientable surface $N$ of genus $c$ in the cube.  Here the framing
  of $T$ which is determined by $N$ may differ from the $0$--framing.
  (We ignore the framing until the end of this proof.)  As is well
  known, $N$ can be obtained from a disc $D$ by attaching $c$ bands
  $b_1,\ldots,b_c$ such that, for each $i=1,\ldots,c$, the union
  $D\cup b_i$ is a M\"obius band, and between the two components of
  $D\cap b_i$ there are no attaching region of the other band, see
  \fullref{fig:surface} (a).
\labellist\small
\pinlabel {(a)} [b] at 120 0
\pinlabel {(b)} [b] at 395 0
\pinlabel {$\cdots$} at 118 95
\pinlabel {$\cdots$} at 391 95
\tiny
\pinlabel {$b_1$} [r] at 44 86
\pinlabel {$b_c$} [r] at 156 86
\pinlabel {$D$} [r] at 225 55
\pinlabel {double of $T'$} [l] at 224 142
\pinlabel {double of $T'$} [l] at 482 142
\endlabellist
  \FIGn{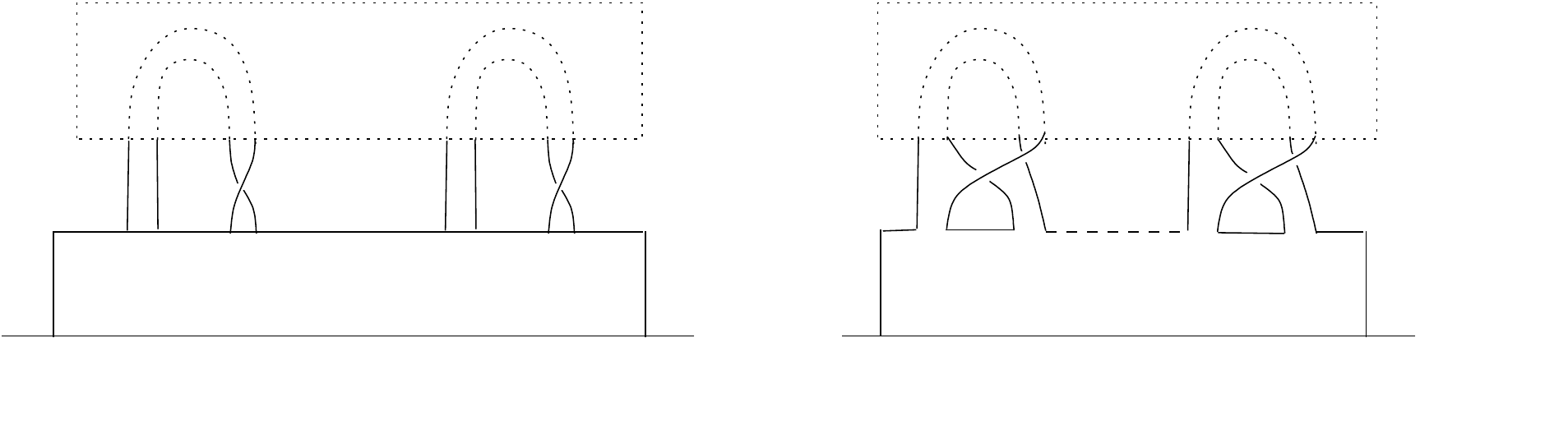}{}{height=30mm}
  Here the dotted
  part is obtained from a $c$--component bottom tangle $T'$ by
  replacing the components with bands, using the framings.  $T$ can be
  isotoped as in \fullref{fig:surface} (b).  Since the framing of
  the tangle depicted in \fullref{fig:surface} (b) is $4w(T')$, we
  have
  \begin{equation*}
    T =
    (t_{\downarrow} \otimes \uparrow )^{-4w(T')}
    F_{\modb} ((\mu_{\HH}\Delta_{\HH})^{\otimes c})(J_{T'}).
  \end{equation*}
  Hence we have the assertion.
\end{proof}

The unorientable version of boundary link is {\em $\modZ _2$--boundary
link} (see Hillman \cite{Hillman}).  A link $L$ in $S^3$ is called $\modZ _2$--boundary
if the components of $L$ bounds mutually disjoint possibly
unorientable surfaces.  Similarly, a bottom tangle $T\in \BT_n$ is said
to be {\em $\modZ _2$--boundary} if the components of $T$ are of framing
$0$ and bound mutually disjoint possibly unorientable surfaces in
$[0,1]^3$.  In the above definitions, ``possibly unorientable'' can
be replaced with ``unorientable''.  One can easily modify \fullref{thm:17} into the $\modZ _2$--boundary case.

\subsection{Borromean tangle and delta moves}
\label{sec:borr-tangle-delta}

We consider {\em delta moves} (see Murakami and Nakanishi
\cite{Murakami-Nakanishi}) or {\em
  Borromean transformation} (see Matveev \cite{Matveev}) on bottom tangles, which
  we mentioned in \fullref{sec:local-moves}.  In our setting, a
  delta move can be defined as a $B$--move, where $B\in \BT_3$ is the
  Borromean tangle defined in \fullref{sec:local-moves}.

The following is an easily verified variant of a theorem of Murakami
and Nakanishi \cite{Murakami-Nakanishi}, which makes delta moves
especially useful.

\begin{proposition}[Murakami--Nakanishi \cite{Murakami-Nakanishi}]
  \label{thm:11}
  Two $n$--component bottom tangles $T$ and $T'$ have the same linking
  matrix if and only if there is a sequence of finitely many delta
  moves (and isotopies) from $T$ to $T'$.  (Here the linking matrix of
  an $n$--component bottom tangle $T$ is defined to be the linking
  matrix of the closure of $T$.)
\end{proposition}

Using \fullref{thm:11}, we obtain the following results.

\begin{corollary}
  \label{thm:20}
  For two $n$--component bottom tangles $T$ and $T'$, the following
  conditions are equivalent.
  \begin{enumerate}
  \item $T$ and $T'$ have the same linking matrix.
  \item $T$ and $T'$ are delta move equivalent, ie, $B$--equivalent.
  \item For some $k\ge 0$ and $W\in \modB (3k,n)$, we have
    \begin{equation*}
      T=W\eta _{3k},\quad T'=W B^{\otimes k}.
    \end{equation*}
  \end{enumerate}
\end{corollary}

\begin{proof}
  The equivalence of (1) and (2) is just \fullref{thm:11}.  The
  equivalence of (2) and (3) follows from \fullref{lem:4},
  since the set $\{B\}$ is inversion-closed (see
  \cite{Murakami-Nakanishi}).
\end{proof}

\begin{corollary}
  \label{thm:24}
  The linking matrix of an $n$--component bottom tangle $T$ is zero if
  and only if there are $k\ge 0$ and $W\in \cB_0(3k,n)$ such that we have
  $T=WB^{\otimes k}$.

  In particular, a bottom tangle with zero linking matrix is obtained
  by pasting finitely many copies of
  $1_{\modb} ,\psi _{\modb ,\modb },\psi _{\modb ,\modb }^{-1},\mu_{\modb}
,\eta_{\modb} ,\gamma _+,\gamma _-,B$.
\end{corollary}

\begin{proof}
  The result follows from \fullref{thm:19}, since $T$ is of
  linking matrix $0$ if and only if $T$ is $B$--trivial.
\end{proof}

In some applications, the following form may be more useful.

\begin{corollary}
  \label{r8}
  The linking matrix of an $n$--component bottom tangle $T$ is zero if
  and only if there are $k\ge 0$ and $f\in \modH (3k,n)$ such that we have
  $T=F_{\modb} (f)(B^{\otimes k})$.
\end{corollary}

\begin{proof}
  The ``only if'' part follows easily from \fullref{thm:24}, using
  \begin{align*}
    (\gamma _+)_{(i,j)} U &= F_{\modb} ((\coad_{\HH})_{(i,j)})(U),\\
    (\gamma _-)_{(i,j)} U &= F_{\modb} (((S_{\HH}\otimes \HH)\coad_{\HH})_{(i,j)})(U),
  \end{align*}
  for $i,j\ge 0$, $U\in \BT_{i+j+1}$.

  The ``if'' part follows from the easily verified fact that the set
  of bottom tangles with zero linking matrices is closed under the
  Hopf algebra action.
\end{proof}

Now we apply the above results to the universal invariant.  First we
give a few formulas for $J_B\in H^{\otimes 3}$.  Using \fullref{fig:borromean}, we can easily see that
\begin{equation}
  \label{eq:24}
  J_B = \sum
  S^2(\alpha _5)     \beta _2            \alpha _6    S(\beta _1)\otimes
  \alpha _1          \beta _4            \alpha _2    S(\beta _3)\otimes
  \alpha _3          S^{-2}(\beta _6)    \alpha _4    S(\beta _5),
\end{equation}
where $R=\sum \alpha _i\otimes \beta _i$ for $i=1,\ldots,6$.  By \fullref{fig:borromean2}, we
have
\FIGn{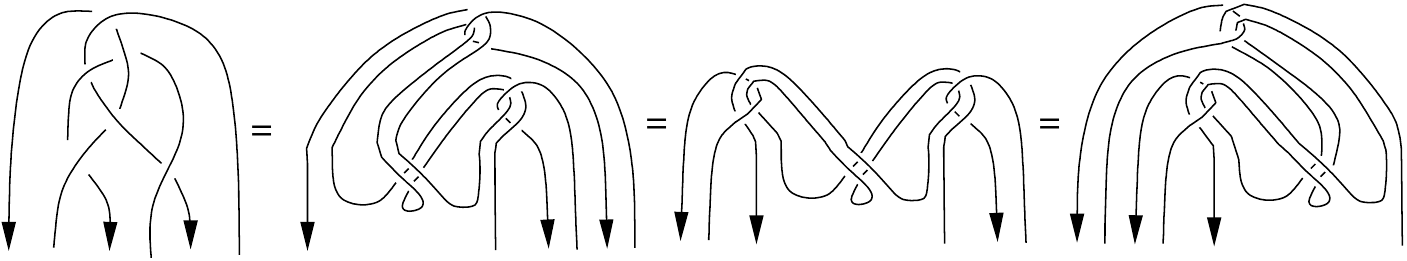}{}{height=20mm}
\begin{equation}
  \label{eq:25}
  B
  = F_{\modb} (Y_{\HH}\otimes \HH\otimes \HH)(c_{+,2})
  = F_{\modb} (\HH\otimes Y_{\HH}\otimes \HH)((c_+)^{\otimes 2})
  = F_{\modb} (\HH\otimes \HH\otimes Y_{\HH})(c_{+,2}),
\end{equation}
where $c_{+,2}=(\modb \otimes c_+\otimes \modb )c_+\in \BT_4$.  By \eqref{eq:25}, it follows
that
\begin{equation*}
  J_B
  = (Y_{\bH}\otimes 1_H^{\otimes 2})(c^H_{+,2})
  = (1_H\otimes Y_{\bH}\otimes 1_H)((c^H_+)^{\otimes 2})
  = (1_H^{\otimes 2}\otimes Y_{\bH})(c^H_{+,2}),
\end{equation*}
where $c^H_{+,2}=\sum (c^H_+)_{[1]}\otimes c^H_+\otimes (c^H_+)_{[2]}\in H^{\otimes 4}$ with
$c^H_+=\sum(c^H_+)_{[1]}\otimes (c^H_+)_{[2]}$.

One can apply \fullref{thm:25} to the case $M=\{B\}$ to obtain a
result about the difference of the universal invariants of bottom
tangles which have the same linking matrix.  We do not give the
explicit statement here.

For the bottom tangles with zero linking matrices, we can easily
derive the following result from Corollaries \ref{thm:24} and
\ref{r8}.

\begin{corollary}
  \label{thm:29}
  Set either
  \begin{equation*}
    X=\{\psi _{H,H}^{\pm 1},\;\mu ,\;\coadb,\;(\bS\otimes 1)\coadb\}
    \quad\text{or}\quad
    X=\{\psi _{H,H}^{\pm 1},\;\mu ,\;\bD,\;\bS\}.
  \end{equation*}
  Let $K_i\subset H^{\otimes i}$, $i\ge 0$, be
  subsets satisfying the following conditions.
  \begin{enumerate}
  \item $1\in K_0$, $1\in K_1$, and $J_B\in K_3$.
  \item If $k,l\ge 0$, then we have $K_k\otimes K_l\subset K_{k+l}$.
  \item For $p,q\ge 0$ and $f\in X$ with $f\col H^{\otimes i}\rightarrow H^{\otimes j}$, we have
    \begin{equation*}
      f_{(p,q)}(K_{p+q+i})\subset K_{p+q+j}.
    \end{equation*}
  \end{enumerate}
  Then, for any $U\in \BT_n$ with zero linking matrix, we have
  $J_U\in K_n$.
\end{corollary}

As mentioned in \fullref{sec:local-moves}, in future publications
we will apply \fullref{thm:29} to the case where $H$ is a
quantized enveloping algebra.

\subsection{Clasper moves}
\label{sec:cn-moves}

In this subsection, we apply the settings in this paper to the {\em
clasper moves} or {\em $C_n$--moves}
(see Goussarov \cite{Gusarov:variations} and Habiro
\cite{Habiro:claspers}) which are closely related to
the Goussarov--Vassiliev finite type link invariants
(see Vassiliev \cite{Vassiliev}, Goussarov
\cite{Gusarov:91,Gusarov:n-equivalence}, Birman \cite{Birman:93}, Birman
and Lin \cite{Birman-Lin} and Bar-Natan \cite{Bar-Natan}).

Recall that a simple $C_n$--moves in the sense of
\cite{Habiro:claspers} is a local move on a tangle~$T$ defined as
surgery on a strict tree clasper $C$ of degree $n$ (ie, with $n+1$
disc-leaves) such that each disc-leaf of $T$ intersects transversely
with $T$ by one point.  A simple $C_n$--move is a generalization of a
crossing change ($n=1$) and a delta move ($n=2$).  In this subsection,
for simplicity, we slightly modify the definition of a simple
$C_n$--move so that the sign of the intersection of $C$ (which is
defined as a surface homeomorphic to a disc) and the strings of $T$
are all positive or all negative.  It is known (see, for example, the
author's master's thesis \cite{Habiro:masterthesis}) that this does not
make any essential
difference if $n\ge 2$.  Ie, the relations on tangles defined by
the moves are the same.

We can use the results in the previous sections in the study of simple
$C_n$--moves, by redefining a simple $C_n$--move as an $M_n$--move, where
$M_n$ is an inversion-closed subset of $\ABT$ defined as follows.
Define $\modY _n\subset \modH (\HH^{\otimes n},\HH)$ for $n\ge 1$ inductively by
$\modY _1=\{1_{\HH}\}$ and
\begin{gather*}
  \modY _n = \{Y_{\HH}(f\otimes g)\ver f\in \modY _i,g\in \modY _j, i+j=n\}\quad \text{for $n\ge 2$}.
\end{gather*}
Thus $\modY _n$ is the set of iterated commutators of class $n$.  For
example, we have $\modY _2=\{Y_{\HH}\}$ and
$\modY _3=\{Y_{\HH}(Y_{\HH}\otimes \HH),Y_{\HH}(\HH\otimes Y_{\HH})\}$.  For $n\ge 1$, define
$M_n\subset \ABT_{n+1}$ by
\begin{equation*}
  M_n = \{F_{\modb} (f\otimes \HH^{\otimes n})(c_{+,n})\ver f\in \modY _n\},
\end{equation*}
where we set
\begin{equation*}
  c_{+,n}=
  (\modb ^{\otimes (n-1)}\otimes c_+\otimes \modb ^{\otimes (n-1)})\cdots(\modb \otimes c_+\otimes \modb )c_+
  \in \BT_{2n}
\end{equation*}
for $n\ge 1$.  (Here, the fact that each element of $M_n$ is admissible
follows from \cite[Lemma 3.20]{Habiro:claspers}.)  In particular, we
have $M_1=\{c_+\}$ and $M_2=\{B\}$.  For example, \fullref{fig:clasper} shows a clasper $C$ for $\eta _5$ such that surgery
along $C$ yields the tangle
\begin{equation*}
  F_{\modb} ((Y_{\HH}(Y_{\HH}\otimes \HH)(Y_{\HH}\otimes \HH^{\otimes 2}))\otimes \HH^{\otimes 4})(c_{+,4})\in M_4.
\end{equation*}
\labellist\small
\pinlabel {$C$} [r] at 340 693
\endlabellist
\FIG{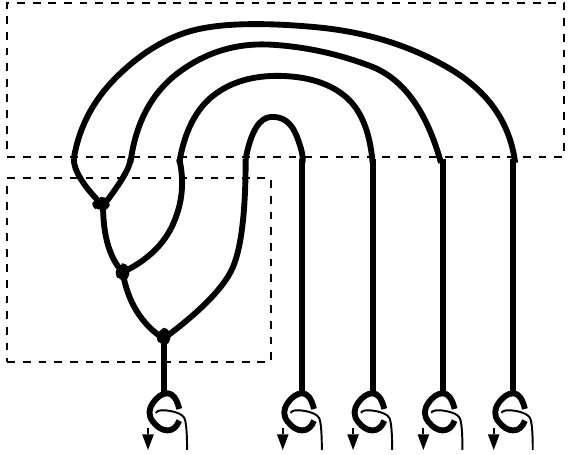}{The upper rectangle corresponds to $c_{+,4}\in \BT_8$.
  The lower rectangle corresponds to
  $Y_{\HH}(Y_{\HH}\otimes \HH)(Y_{\HH}\otimes \HH^{\otimes 2})\in \modY _4$.}{height=30mm}
We can also define the $M_n$ using the {\em cocommutator morphism}
\cite{Habiro:claspers}
\begin{equation*}
  Y_{\HH}^*\col \HH\rightarrow \HH\otimes \HH
\end{equation*}
defined by
\begin{equation}
  \label{e3}
  \begin{split}
    Y_{\HH}^*
    &= (\mu_{\HH}\otimes \mu_{\HH})(\HH\otimes S_{\HH}\otimes S_{\HH}\otimes \HH)
    (\HH\otimes \psi _{\HH,\HH}\otimes \HH)\Delta_{\HH}^{[4]}\\
    &=(\HH\otimes \mu_{\HH})(\coad_{\HH}\otimes \HH)(S_{\HH}\otimes
\HH)\Delta_{\HH}.
  \end{split}
\end{equation}
Note that the notion of cocommutator is dual to the notion of
commutator.  For $i,j\ge 0$, the function
\begin{equation*}
  F_{\modb} ((Y_{\HH}^*)_{(i,j)})\col \BT_{i+j+1}\rightarrow \BT_{i+j+2}
\end{equation*}
transforms a bottom tangle into another as illustrated in
\fullref{fig:Ystar}.
\labellist\tiny
\pinlabel {\small(a)} [b] at 17 0
\pinlabel {$i{+}1$} [t] at 17 48
\pinlabel {\small(b)} [b] at 125 0
\pinlabel {$i{+}1$} [t] at 105 48
\pinlabel {$i{+}2$} [t] at 147 48
\pinlabel {\small(c)} [b] at 230 0
\pinlabel {$i{+}1$} [t] at 200 48
\pinlabel {$i{+}2$} [t] at 248 48
\pinlabel {\small(d)} [b] at 353 0
\pinlabel {$i{+}1$} [t] at 332 48
\pinlabel {$i{+}2$} [t] at 372 48
\pinlabel {\small(e)} [b] at 460 0
\pinlabel {$i{+}1$} [t] at 435 48
\pinlabel {$i{+}2$} [t] at 482 48
\endlabellist
\FIG{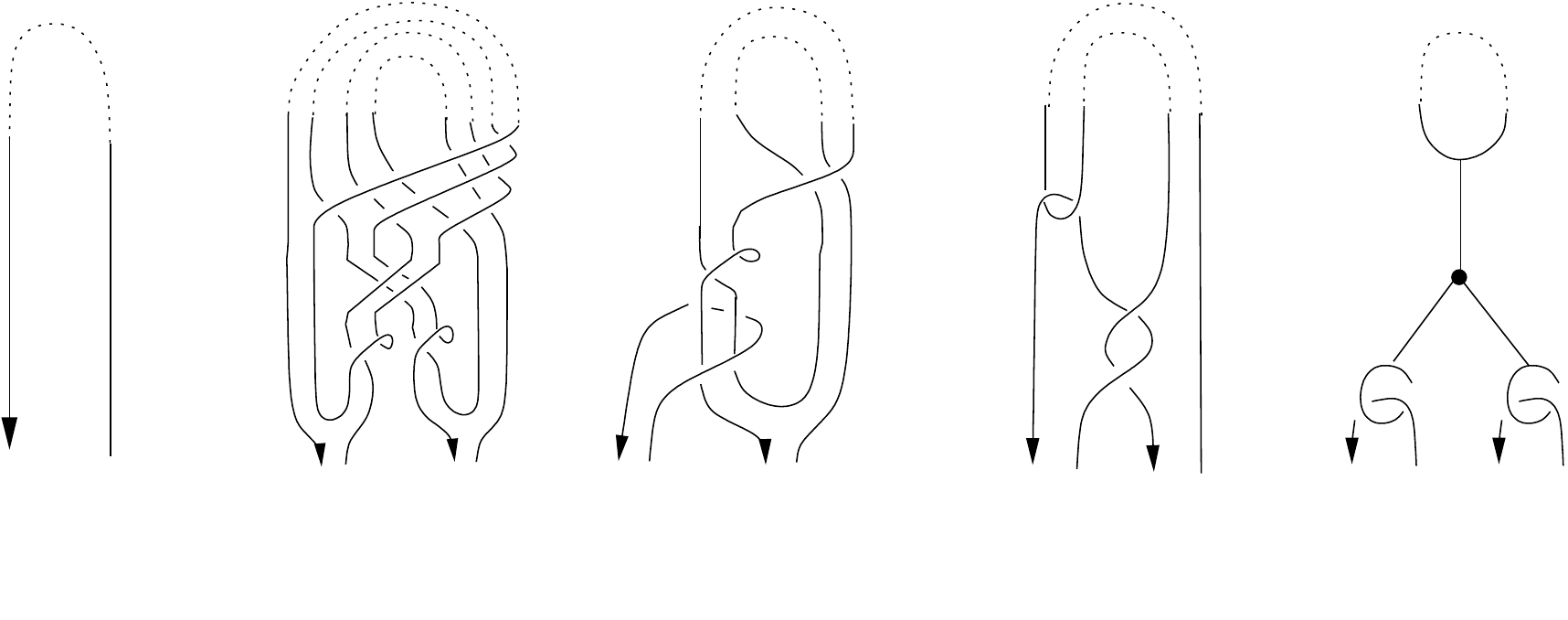}{(a) An $(i+j+1)$--component bottom tangle
$T$.  (b) The $(i+j+2)$--component bottom tangle
$T'=F((Y^*_{\HH})_{(i,j)})(T)$, calculated using \eqref{e3} (upper). (c)
$T'$ calculated using \eqref{e3} (lower).  (d) Another picture of
$T'$.  (d) A presentation of $T'$ using a clasper.}{height=37mm} For
$n\ge 1$, define $\modY ^*_n\subset \modH (\HH,\HH^{\otimes n})$ inductively by
$\modY ^*_1=\{1_{\HH}\}$ and
\begin{equation*}
  \modY ^*_n = \{(f\otimes g)Y^*_{\HH}\ver f\in \modY ^*_i,g\in \modY ^*_j, i+j=n\}\quad
  \text{for $n\ge 2$}.
\end{equation*}
Then we have for $n\ge 1$,
\begin{equation*}
  M_n = \{F_{\modb} (f\otimes g)(c_+)\ver f\in \modY ^*_i,g\in \modY ^*_j, i+j=n+1\},
\end{equation*}
which follows by induction using \eqref{eq:25} and
\begin{equation}
  \label{e5}
  B = F_{\modb} (Y_{\HH}^*\otimes \HH)(c_+) = F_{\modb} (\HH\otimes
Y_{\HH}^*)(c_+).
\end{equation}
(The above definition of $M_n$ using $Y^*_{\HH}$ is similar to the
definition of local moves in \cite{Habiro:masterthesis}, where we
defined a family of local moves without using claspers.  See also
Taniyama and Yasuhara \cite{Taniyama-Yasuhara} for a similar definition.)

One can show that the notion of simple $C_n$--move and that of
$M_n$--move are the same.

\begin{remark}
\label{r1}
A general $C_n$--move, which may not be simple, is obtained by allowing
removal, orientation reversal and parallelization of strings in the
tangles which define the move.  Hence it can be redefined as an
$M'_n$--move, where the set $M'_n\subset \ABT$ is defined by
\begin{equation*}
  \begin{split}
    M'_n =&
    \Bigl\{F_{\modb} \Bigl(\bigotimes_{i=1}^{n+1}\Delta_{\HH}^{[c_i;d_i]}\Bigr)(f)
      \;\Bigl|\; c_1,\ldots,c_{n+1},d_1,\ldots,d_{n+1}\ge 0, f\in M_n\Bigr\},
  \end{split}
\end{equation*}
where we set
\begin{equation*}
  \Delta_{\HH}^{[c;d]}=(\HH^{\otimes c}\otimes S_{\HH}^{\otimes d})\Delta ^{[c+d]}\in \modH (\HH,\HH^{\otimes (c+d)})
\end{equation*}
for $c,d\ge 0$.  As special cases, the following local moves in the
literature can be redefined algebraically:
\begin{enumerate}
\item A {\em pass-move} (see Kauffman \cite{Kauffman:pass} and \fullref{fig:moves}
   (a)), which characterizes the Arf invariant of knots, is the same
   as a $F_{\modb} ((\Delta_{\HH}^{[1;1]})^{\otimes 2})(c_+)$--move.
\item A {\em $\sharp$--move} (see Murakami \cite{Murakami:sharp} and  \fullref{fig:moves} (b)) is the same as a
  $F_{\modb} (\Delta_{\HH}^{\otimes 2})(c_+)$--move.  (This is a framed version. In
  applications to unframed or $0$--framed knots, one should take
  framings into account.)
\item A {\em $D(\Delta )$--move} (see Nakanishi \cite{Nakanishi} and
\fullref{fig:moves} (c)), which preserves the stable equivalence class of
 the Goeritz matrix (see Goeritz \cite{Goeritz} and Gordon--Litherland
\cite{Gordon-Litherland}) of (possibly
 unorientable) spanning surfaces of knots, can be redefined by
 setting
  \begin{equation*}
    D(\Delta )=
    \{F_{\modb} ((\Delta_{\HH})^{[i;2-i]}\otimes
(\Delta_{\HH})^{[j;2-j]}\otimes (\Delta_{\HH})^{[k;2-k]})(B)
      \ver 0\le i,j,k\le 2\}.
  \end{equation*}
  It is known that the $D(\Delta )$--equivalence is the same as an oriented
  version of it, which can be defined as the
  $F_{\modb} (\Delta_{\HH}^{\otimes 3})(B)$--equivalence, ie, the case $i=j=k=2$.
\item A {\em doubled-delta move} (see Naik--Stanford \cite{Naik-Stanford}
and \fullref{fig:moves} (d)), which characterizes the $S$--equivalence class
  of Seifert matrices of knots, can be defined as a
  $F_{\modb} ((\Delta_{\HH}^{[1;1]})^{\otimes 3})(B)$--move, which is the case $i=j=k=1$
  of $D(\Delta )$--move.
\end{enumerate}
\labellist\tiny
\pinlabel {\small(a)} [b] at 100 100
\pinlabel {pass-move} [b] at 100 154
\pinlabel {\small(b)} [b] at 335 100
\pinlabel {$\sharp$--move} [b] at 335 154
\pinlabel {\small(c)} [b] at 100 0
\pinlabel {$D(\Delta)$--move} [b] at 103 43
\pinlabel {\small(d)} [b] at 335 0
\pinlabel {\parbox{40pt}{doubled-delta move}} [b] at 340 43
\endlabellist
\FIG{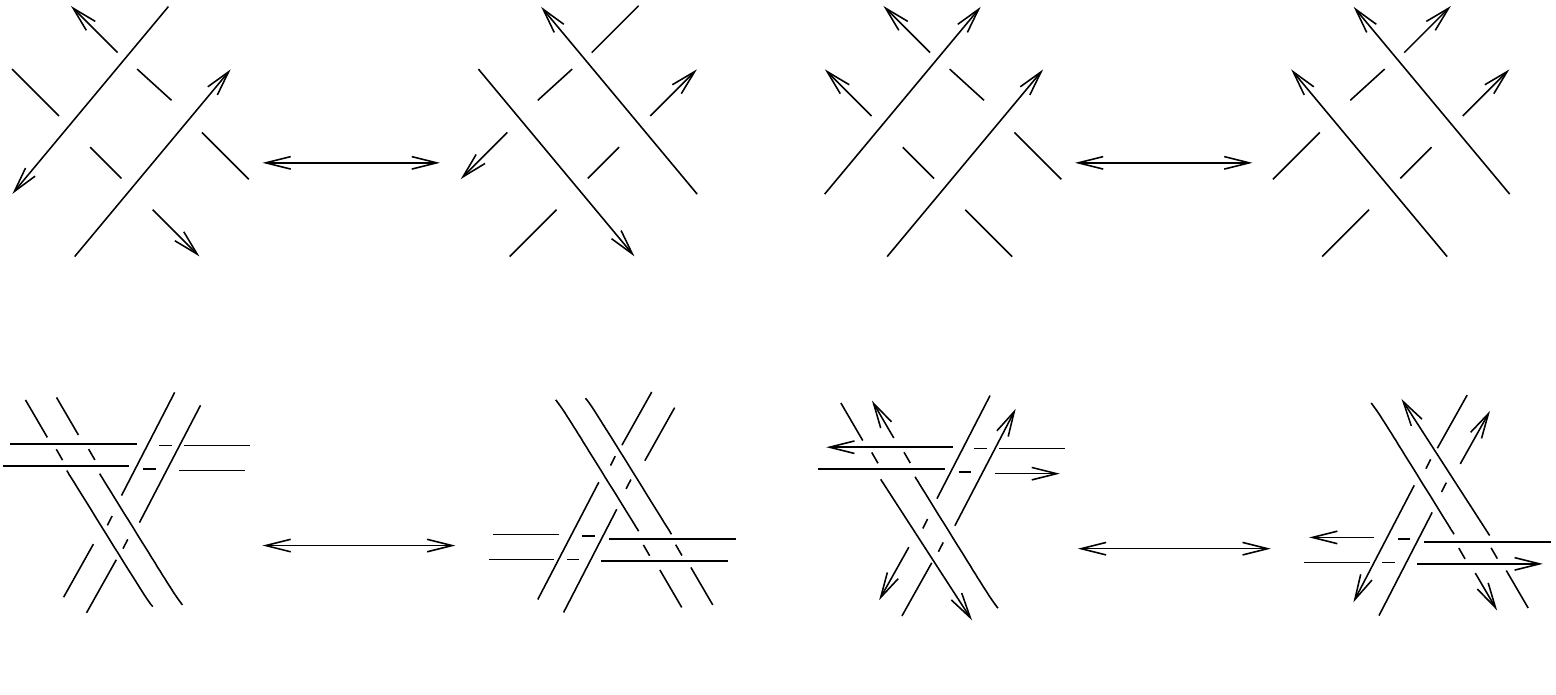}{(a) A pass-move.  (b) A $\sharp$--move.  (c) A $D(\Delta )$--move.
(Here the orientations of strings are arbitrary.)  (d) A doubled-delta
move.}{height=50mm}
\end{remark}

We postpone to future publications a more systematic study of the
clasper moves in a category-theoretical setting, which was announced
in \cite{Habiro:claspers}.  For this purpose the category $\B$ (see
\fullref{sec:category-0} below) is more useful than $\modB $.

\subsection{Goussarov--Vassiliev filtrations on tangles}
\label{sec:vass-gouss-filtr-2}
In this subsection, we give an algebraic formulation of
Goussarov--Vassiliev invariants using the setting of the category $\modB $.

\subsubsection{Four-sided ideals in a monoidal $\Ab$--category}
\label{sec:four-sided-ideals}

Here we recall the notion of {\em four-sided ideal} in a monoidal
$\Ab$--category, which can be regarded as the linearized version of
the notion of four-sided congruence in a monoidal category.

Let $C$ be a (strict) {\em monoidal $\Ab$--category}, ie, a monoidal
category $C$ such that for each pair $X,Y\in \Ob(C)$ the set $C(X,Y)$ is
equipped with a structure of a $\modZ $--module, and the composition and
the tensor product are bilinear.

A {\em four-sided ideal} $I=(I(X,Y))_{X,Y\in \Ob(C)}$ in a monoidal
$\Ab$--category $C$ is a family of $\modZ $--submodules $I(X,Y)$ of $C(X,Y)$
for $X,Y\in \Ob(C)$ such that
\begin{enumerate}
\item if $f\in I(X,Y)$ and $g\in C(Y,Z)$ (resp. $g\in C(Z,X)$),
  then we have $gf\in I(X,Z)$ (resp.  $fg\in I(Z,Y)$),
\item if $f\in I(X,Y)$ and $g\in C(X',Y')$, then we have
  $f\otimes g\in I(X\otimes X',Y\otimes Y')$ and $g\otimes f\in I(X'\otimes X,Y'\otimes Y)$.
\end{enumerate}
By abuse of notation, we denote by $I$, the union
$\bigcup_{X,Y\in \Ob(C)}I(X,Y)$.

Let $S\subset \Mor(C)$ be a set of morphisms in $C$.  Then there is the
smallest four-sided ideal $I_S$ in $C$ such that $S\subset I_S$.  The
four-sided ideal $I_S$ is said to be {\em generated} by $S$.  For
$X,X'\in \Ob(C)$, $I_S(X,X')$ is $\modZ $--spanned by the elements
\begin{equation}
  \label{e6}
  f'(g\otimes s\otimes g')f,
\end{equation}
where $s\in S$, and $f,f',g,g'\in \Mor(C)$ are such that the expression
\eqref{e6} gives a well-defined morphisms in $C(X,X')$.

For two four-sided ideals $I$ and $I'$ in a monoidal $\Ab$--category
$C$, the {\em product} $I'I$ of $I'$ and $I$ is defined to be the
smallest four-sided ideal in $C$ such that if $(g,f)\in I'\times I$ is a
composable pair, then $gf\in I I'$.  It follows that $f\in I$ and
$g\in I'$ implies $f\otimes g,g\otimes f\in I'I$.  For $X,Y\in \Ob(C)$, then we have
\begin{equation*}
  I'I(X,Y) = \sum_{Z\in \Ob(C)}I'(Z,Y)I(X,Z).
\end{equation*}
For $n\ge 0$, let $I^n$ denote the $n$th power of $I$, which is defined
by $I^0=\Mor(C)$, $I^1=I$, and $I^n=I^{n-1}I$ for $n\ge 2$.

\begin{lemma}
  \label{r17}
  Let $C$ be a braided $\Ab$--category and let $I$ be a four-sided
  ideal in $C$ generated by $S\subset \prod_{X\in \Ob(C)}C(\one ,X)$.  Then
  $I^n(X,Y)$ ($X,Y\in \Ob(C)$) is $\modZ $--spanned by the elements of the
  form $f(X\otimes s_1\otimes \cdots\otimes s_n)$, where $s_1,\ldots,s_n\in S$ and
  $f\in C(X\otimes \target(s_1\otimes \cdots\otimes s_n),Y)$.
\end{lemma}

\begin{proof}
  The proof is sketched as follows.  Each element of $I^n$ is a
  $\modZ $--linear combination of morphisms, each obtained as an iterated
  composition and tensor product of finitely many morphisms of $C$
  involving $n$ copies of elements $s_1,\ldots,s_n$ of $S$.  By the
  assumption, one can arrange (using braidings) the copies
  $s_1,\ldots,s_n$ involved in each term of an element of $I^n$ to be
  placed side by side as in $s_1\otimes \cdots\otimes s_n$ in the upper right
  corner, ie, we obtain a term of the form $f(X\otimes s_1\otimes \cdots\otimes s_n)$,
  as desired.
\end{proof}

\subsubsection{Goussarov--Vassiliev filtration for $\modZ \modT $}
\label{sec:vass-gouss-filtr}

Here we recall a formulation of Gouss\-ar\-ov--Vassiliev filtration using
the category~$\modT $ of framed, oriented tangles, which is given by Kassel and
Turaev \cite{Kassel-Turaev}.

Let $\modZ \modT $ denote the category of {\em $\modZ $--linear tangles}.  Ie, we
have $\Ob(\modZ \modT )=\Ob(\modT )$, and for $w,w'\in \Ob(\modT )$, the set
$\modZ \modT (w,w')$ is the free $\modZ $--module generated by the set $\modT (w,w')$.
$\modZ \modT $ is a braided $\Ab$--category.

Let $\modI $ denote the four-sided ideal in $\modZ \modT $ generated by the
morphism
\begin{equation*}
  \psi ^{\times} =\psi _{\downarrow ,\downarrow }-\psi _{\downarrow ,\downarrow }^{-1}\in \modZ \modT (\downarrow ^{\otimes 2},\downarrow ^{\otimes 2}).
\end{equation*}
For $n\ge 0$, let $\modI ^n$ denote the $n$th power of $\modI $.  $\modI ^n$ is
equal to the four-sided ideal in $\modZ \modT $ generated by the morphism
$(\psi ^{\times} )^{\otimes n}$.  For $w,w'\in \Ob(\modT )$, the filtration
\begin{equation}
  \label{e7}
  \modZ \modT (w,w')=\modI ^0(w,w')\supset \modI ^1(w,w')\supset \modI ^2(w,w')\supset \cdots
\end{equation}
is known \cite{Kassel-Turaev} to be the same as the
Goussarov--Vassiliev filtration for $\modZ \modT (w,w')$.

\begin{remark}
  \label{r20}
  There is an alternative, perhaps more natural, definition of the
  Gouss\-ar\-ov--Vassiliev filtration for framed tangles, which involves
  the difference $t_{\downarrow} -1_{\downarrow} $ of framing change as
well as $\psi ^{\times} $.  In
  the present paper, we do not consider this version for simplicity.
\end{remark}

\subsubsection{Goussarov--Vassiliev filtration for $\modZ \modB $}
\label{sec:vass-gouss-filtr-1}

Now we consider the case of tangles in $\modB $.  The definition of the
category $\modZ \modB $ of $\modZ $--linear tangles in $\modB $ is obvious.  For
$i,j\ge 0$, the Goussarov--Vassiliev filtration for the tangles in
$\modB (i,j)$ is given by the $\modZ $--submodules
\begin{equation*}
  (\modZ \modB \cap \modI ^n)(i,j) :=\modZ \modB (i,j)\cap \modI ^n(\modb ^{\otimes i},\modb ^{\otimes j})
\end{equation*}
for $n\ge 0$.  Clearly, this defines a four-sided ideal $\modZ \modB \cap \modI ^n$ in
$\modZ \modB $.

Set
\begin{equation*}
  c^{\times} =\eta _2-c_+\in \modZ \modB (0,2),
\end{equation*}
and let $\modI_{\modB} $ denote the four-sided ideal in $\modZ \modB $ generated by
$c^{\times} $.

The following result gives a definition of the Goussarov--Vassiliev
filtration for tangles in $\modB $, and in particular for bottom tangles,
{\em defined algebraically in $\modZ \modB $}.  Thus the setting in the
present paper is expected to be useful in the study of
Goussarov--Vassiliev finite type invariants.

\begin{theorem}
  \label{r13}
  For each $n\ge 0$, we have
  \begin{equation}
    \label{e14}
    \modZ \modB \cap \modI ^n=\modI_{\modB} ^n.
  \end{equation}
  For bottom tangles, we also have
  \begin{equation}
    \label{e18}
    \modI_{\modB} ^n(0,m) = \modZ \modB (2n,m)(c^{\times} )^{\otimes n}
    \Bigl(=\{f(c^{\times} )^{\otimes n}\ver f\in \modZ \modB (2n,m)\}\Bigr).
  \end{equation}
\end{theorem}

\begin{proof}
  We have
  \begin{align}
    \label{e10}
    c^{\times} &=(\downarrow \otimes \psi _{\downarrow ,\uparrow }\otimes
\uparrow )(\psi ^{\times} \otimes \uparrow \otimes \uparrow )
    \coev_{\downarrow \otimes \downarrow } \in \modI (\one ,\modb ^{\otimes 2}),\\
    \label{e13}
    \psi ^{\times} &=(\downarrow \otimes \downarrow \otimes \ev_{\downarrow \otimes \downarrow })(\downarrow \otimes \psi _{\downarrow ,\uparrow }^{-1}\otimes \uparrow \otimes \downarrow \otimes \downarrow )
    (c^{\times} \otimes \downarrow \otimes \downarrow ).
  \end{align}
  By \eqref{e10}, we have $\modI_{\modB} \subset \modI $, and hence
$\modI_{\modB} ^n\subset \modI ^n$ for
  $n\ge 0$.  Since $\modI_{\modB} ^n\subset \modZ \modB $, we have
$\modI_{\modB} ^n\subset \modZ \modB \cap \modI ^n$.

  We show the other inclusion.  Suppose that $f\in (\modZ \modB \cap \modI ^n)(l,m)$.
  By \eqref{e10} and \eqref{e13}, $\modI $ is generated by $c^{\times} $ as a
  four-sided ideal in $\modZ \modT $.  By \fullref{r17}, we have
  \begin{equation*}
    f= g'(\modb ^{\otimes l}\otimes (c^{\times} )^{\otimes n}),
  \end{equation*}
  where $g'\in \modZ \modT (\modb ^{\otimes (l+2n)},\modb ^{\otimes m})$.  We can write
  \begin{equation*}
    g'=\hspace{-10pt}\sum_{h\in \modT (\modb ^{\otimes (l+2n)},\modb
^{\otimes m})}\hspace{-10pt}p_hh,\quad p_h\in \modZ .
  \end{equation*}
  Set
  \begin{equation*}
    g=\hspace{-5pt}\sum_{h\in \modB (l+2n,m)}\hspace{-5pt}p_hh\in \modZ \modB (l+2n,m)
  \end{equation*}
  We have $(g-g')(\modb ^{\otimes l}\otimes (c^{\times} )^{\otimes n})=0$, since if
  $h\in \modT (\modb ^{\otimes (l+2n)},\modb ^{\otimes m})\setminus \modB (l+2n,m)$, then $h\eta _{l+2n}$ is
  not homotopic to $\eta _m$.  Hence
  \begin{equation*}
    f= g(\modb ^{\otimes l}\otimes (c^{\times} )^{\otimes n})\in
\modI_{\modB} ^n.
  \end{equation*}
  Hence we have $\modZ \modB \cap \modI ^n\subset \modI _{\modB} ^n$.

  The identity \eqref{e18} follows from the above argument with $l=0$.
  This completes the proof.
\end{proof}

\begin{remark}
  \label{r21}
  It is easy to generalize this subsection to the case of {\em skein
  modules} (see Przytycki \cite{Przytycki:skein}) involving bottom
  tangles.  Let $\modk $
  be a commutative, unital ring, and consider the $\modk $--linear braided
  categories $\modk \modT $ and $\modk \modB $.  A {\em skein element} is just a
  morphism $f\in \modk \modT (w,w')$, $w,w'\in \Ob(\modk \modT )=\Ob(\modT )$.  For a set
  $S\subset \Mor(\modk \modT )$ of skein elements, let $I_S$ denote the four-sided
  ideal in $\modk \modT $ generated by $S$.  Then the quotient ($\modk $--linear,
  braided) category $\modk \modT /I_S$ is known as the {\em skein category}
  defined by $S$ as the set of skein relations.

  Suppose $S\subset \bigcup_{n\ge 0}\modk \BT_n\subset \Mor(\modk \modB )$.  Thus $S$ is a set
  of skein elements involving only bottom tangles.  Let $I^{\modB} _S$
  denote the four-sided ideal in $\modk \modB $ generated by $S$.  Then we
  have the following generalization of \fullref{r13}:
  \begin{align*}
    \modk \modB \cap I_S&=I^{\modB} _S,\\
    I^{\modB} _S(0,n) &=\sum_{l\ge 0}\modk \modB (l,n)(S\cap \modk \BT_l).
  \end{align*}
  Thus, analogously to the case of local moves, it follows that skein
  theory defined by skein elements of compatible tangles consisting of
  arcs can be formulated within the setting of $\modk \modB $ using skein
  elements of bottom tangles.
\end{remark}

\subsection{Twist moves}
\label{sec:twist-moves}

A {\em twist move} is a local move on a tangle which performs a power
of full twist on a parallel family of strings.  A type of a twist move
is determined by a triple of integers $(n,i,j)$ with $i,j\ge 0$, where
the move performs $n$ full twists on a parallel family of $i$ downward
strings and $j$ upward strings.  Let us call it a $t^n_{i,j}$--move.
In our notation, a $t^n_{i,j}$--move is the same as a
$(\downarrow ^{\otimes i}\otimes \uparrow ^{\otimes j},t_{\downarrow ^{\otimes i}\otimes \uparrow ^{\otimes j}}^n)$--move.  For example,
see \fullref{fig:twist} (a).
\labellist\small
\pinlabel {(a)} [b] at 120 0
\pinlabel {$t^2_{2,1}$--move} [b] at 120 147
\pinlabel {(b)} [b] at 445 0
\endlabellist
\FIG{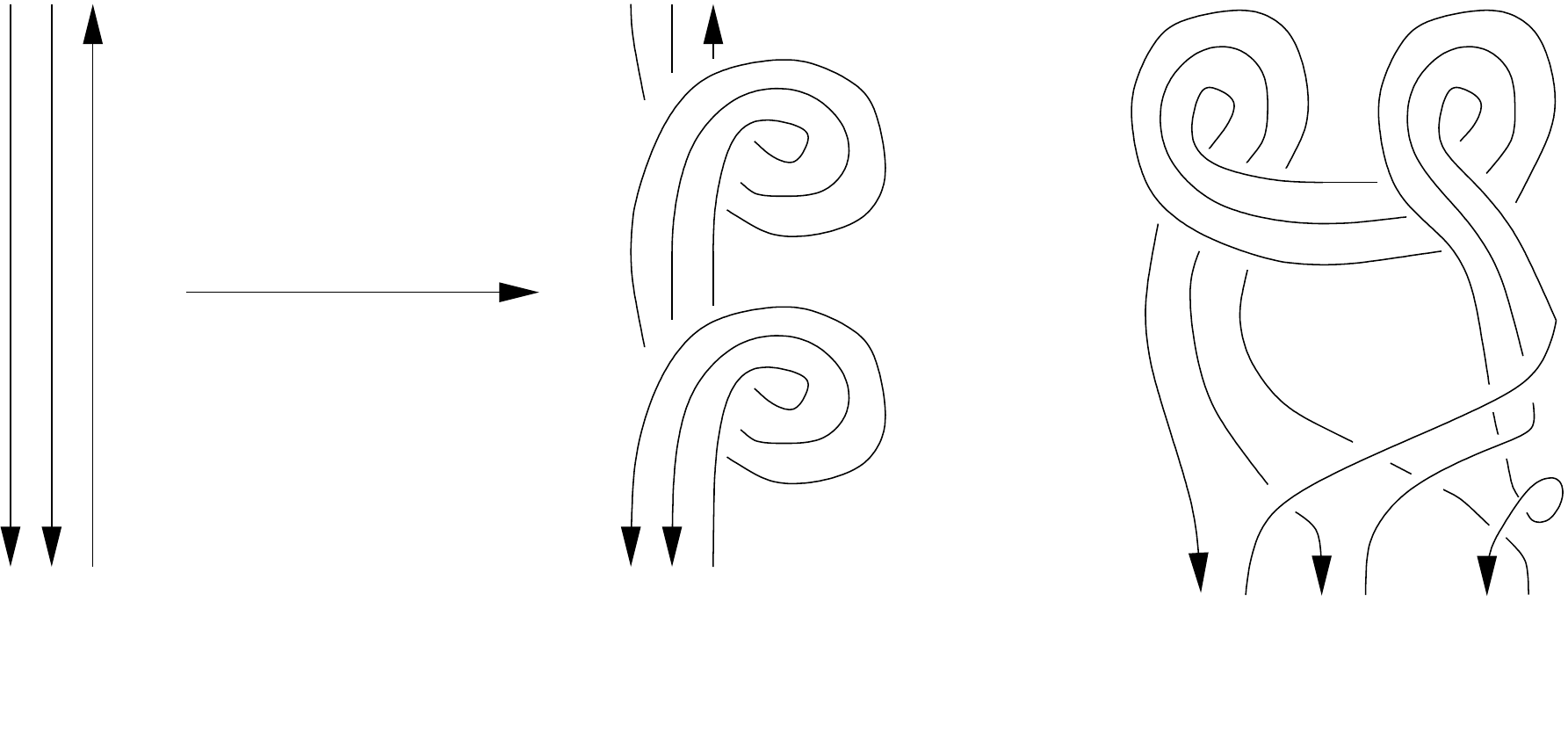}{(a) A $t^2_{2,1}$--move.  (b) The tangle $t^2_{2,1}\in
\BT_3$.}{height=28mm}

Note that a $t^n_{1,0}$--move is just an $n$--full twist of a string,
and is the same as $v_n$--move, where $v_n\in \BT_1$ is the $n$th
convolution power of $v_-$ defined by
\begin{equation*}
  v_n =
  \begin{cases}
    \mu _{\modb} ^{[n]}v_-^{\otimes n}&\text{if $n\ge 0$},\\
    \mu _{\modb} ^{[-n]}v_+^{\otimes (-n)}&\text{if $n\le 0$}.
  \end{cases}
\end{equation*}
Using an idea similar to the one in \fullref{r1}, we see that an
$t^n_{i,j}$--move is the same as $F_{\modb} (\Delta _{\HH}^{[i;j]})(v_n)$--move.
($\Delta _{\HH}^{[i;j]}$ is defined in \fullref{r1}.)  By abuse of
notation, set
\begin{equation*}
  t^n_{i,j}=F_{\modb} (\Delta _{\HH}^{[i;j]})(v_n)\in \BT_{i+j},
\end{equation*}
which should not cause confusion.    Note that $t^n_{i,j}$ is admissible.

Note that a $t^n_{i,j}$--move changes the writhe of a tangle by
$n(i-j)^2$.  (Here the {\em writhe} of a tangle is the number of
positive crossing minus the number of negative crossings.)  In the
literature, twist moves are often considered in the unframed context.
The modification to the unframed case is easy.  For example, two
$0$--framed knots are related by unframed $t^n_{i,j}$--move if they are
related by a sequence of a framed $t^n_{i,j}$--move and a framed
$t_{1,0}^{-n(i-j)^2}$--move.  The latter move multiplies the universal
invariant associated to a ribbon Hopf algebra by the factor of the
power $\modr ^{n(i-j)^2}$ of the ribbon element $\modr $.

Twist moves have long been studied in knot theory.  For a recent
survey, see Przytycki \cite{Przytycki:02}.  Here we give a few examples from the
literature with translations into our setting.  For simplicity, we
only give suitable framed versions of the notions in the literature.

For integers $n,k\ge 0$, a framed version of Fox's notion of {\em
congruence modulo $(n,k)$ } (see Fox \cite{Fox}, Nakanishi--Suzuki
\cite{Nakanishi-Suzuki} and Nakanishi \cite{Nakanishi:90}) can be defined as the
$FC_{n,k}$--equivalence, where we set
\begin{equation*}
  FC_{n,k} = \{t^n_{i,j}\ver i-j\equiv 0\pmod k\}\subset \ABT.
\end{equation*}
For integer $n$, a framed version of {\em $t_{2n}$--move} (see
Przytycki \cite{Przytycki:1}) can be defined as $t_{2,0}^n$--move, and a framed
version of {\em $\bar t_{2n}$--move} can be defined as
$t_{1,1}^n$--move.  Nakanishi's $4$--move conjecture
\cite{Nakanishi-Suzuki,Nakanishi:90}, which is still open, can be
restated that any knot is $\{t^2_{2,0},t^2_{1,1}\}$--equivalent to an
unknot.

We expect that the above ``algebraic redefinitions'' of twist moves
and equivalence relations are useful in the study of these notions in
terms of quantum invariants, by applying the results in \fullref{sec:values-jt}.

\section{The functor $\tilde{\modJ }\col \modB \rightarrow \Mod_H$ and universal
  invariants of bottom knots}
\label{sec:knots-tensor-product}

The following idea may be useful in studying the universal invariants
of bottom knots.

\subsection{The functor $\tilde{\modJ }\col \modB \rightarrow \Mod_H$}

Let $H$ be a ribbon Hopf algebra over a commutative, unital ring $\modk $,
and let $Z(H)$ denote the center of $H$.  Let $\tilde{H}$ denote $H$
regarded as a $Z(H)$--algebra.  For $n\ge 0$, let $\tilde{H}^{\otimes n}$
denote the $n$--fold iterated tensor product of $\tilde{H}$, ie,
$n$--fold tensor product of $H$ over $Z(H)$, regarded as a
$Z(H)$--algebra.  In particular, we have $\tilde{H}^{\otimes 0}=Z(H)$.  Let
\begin{equation*}
  \iota _n\col H^{\otimes n}\rightarrow  \tilde{H}^{\otimes n}
\end{equation*}
denote the natural map, which is surjective if $n\ge 1$.

The functor $\modJ \col \modB \rightarrow \Mod_H$ induces another functor
$\tilde{\modJ }\col \modB \rightarrow \Mod_H$ as follows.  For $n\ge 0$, set
$\tilde{\modJ }(\modb ^{\otimes n})=\tilde{H}^{\otimes n}$, which is given the left
$H$--module structure induced by that of $H^{\otimes n}$.  (This left
$H$--module structure of $\tilde{H}^{\otimes n}$ does {\em not} restrict to
the $Z(H)$--module structure of $\tilde{H}^{\otimes n}$.)

For each $f\in \modB (m,n)$, the left $H$--module homomorphism
$\tilde{\modJ }(f)\col \tilde{H}^{\otimes m}\rightarrow \tilde{H}^{\otimes n}$ is induced by
$\modJ (f)\col H^{\otimes m}\rightarrow H^{\otimes n}$ as follows.  If $m>0$, then $\tilde{\modJ }(f)$
is defined to be the unique map such that the following diagram
commutes
\begin{equation}
  \label{e4}
  \begin{CD}
    H^{\otimes m} @>\modJ (f)>> H^{\otimes n}\\
    @V\iota _mVV @VV\iota _nV\\
    \tilde{H}^{\otimes m} @>>\tilde{\modJ }(f)> \tilde{H}^{\otimes n}.
  \end{CD}
\end{equation}
If $m=0$, then set
\begin{equation}
  \label{e2}
  \tilde{\modJ }(f)(z)=z\iota _n(\modJ (f)(1))\quad \text{for $z\in Z(H)$}.
\end{equation}
Note that commutativity of the diagram \eqref{e4} holds also for
$m=0$.  It is straightforward to check that the above defines a
well-defined functor $\tilde{\modJ }$.

For $m,n\ge 0$, let
\begin{equation*}
  \xi _{m,n}\col \tilde{H}^{\otimes m}\otimes _{\modk} \tilde{H}^{\otimes n}\rightarrow \tilde{H}^{\otimes (m+n)}
\end{equation*}
denote the natural map.  The $\xi _{m,n}$ form a natural transformation
\begin{equation*}
  \xi \col \tilde{\modJ }(-)\otimes \tilde{\modJ }(-)\rightarrow \tilde{\modJ }(-\otimes -)
\end{equation*}
of functors from $\modB \times \modB $ to $\Mod_H$.

It is straightforward to see that the triple $(\tilde{\modJ },\xi ,\iota _0)$ is
an {\em ordinary} braided functor, ie a ordinary monoidal functor
which preserves braiding.  (By ``ordinary monoidal functor'', we mean
a ``monoidal functor'' in the ordinary sense, see Mac\,Lane \cite[Chapter XI,
Section 2]{MacLane}.)

The $\iota _n$ form a monoidal natural transformation $\iota \col \modJ \Rightarrow \tilde{\modJ }$
(in the ordinary sense) of ordinary monoidal functors from $\modB $ to
$\Mod_H$.

\subsection{Universal invariant of bottom knots}
\label{sec:univ-invar-bott}

Since $\iota _1\col H=H^{\otimes 1}\rightarrow \tilde{H}^{\otimes 1}=H$ is the identity, the
functor $\tilde{\modJ }$ can be used in computing $J_T=\modJ (T)(1)$ for a
bottom knot $T\in \BT_1$.  For example, we have the following version of
\fullref{thm:12} for $n=1$.

\begin{proposition}
  \label{r7}
  Let $K_i\subset \tilde{H}^{\otimes i}$ for $i\ge 0$, be $Z(H)$--submodules
  satisfying the following.
  \begin{enumerate}
  \item $1\in K_0$, $1,v^{\pm 1}\in K_1$, and $\iota _2(c^H_{\pm} )\in K_2$.
  \item For $m,n\ge 0$, we have $K_m\otimes _{Z(H)}K_n\subset K_{m+n}$.
  \item For $p,q\ge 0$ we have
    \begin{align*}
      (\tilde{\psi }_{H,H}^{\pm 1})_{(p,q)}(K_{p+q+2})&\subset K_{p+q+2},\\
      \tilde{\mu }_{(p,q)}(K_{p+q+2})&\subset K_{p+q+1},
    \end{align*}
    where
    \begin{equation*}
      (\tilde{\psi }_{H,H})_{(p,q)}
      =\tilde{\modJ }((\psi _{\modb ,\modb })_{(p,q)})
      \col \tilde{H}^{\otimes (p+q+2)}\rightarrow \tilde{H}^{\otimes (p+q+2)}
    \end{equation*}
    is induced by $(\psi _{H,H})_{(p,q)}\col H^{\otimes (p+q+2)}\rightarrow H^{\otimes (p+q+2)}$,
    and
    \begin{equation*}
      \tilde{\mu }_{(p,q)}=\tilde{\modJ }((\mu _{\modb} )_{(p,q)})
      \col \tilde{H}^{\otimes (p+q+2)}\rightarrow \tilde{H}^{\otimes (p+q+1)}
    \end{equation*} is induced by
    $\mu _{(p,q)}\col H^{\otimes (p+q+2)}\rightarrow H^{\otimes (p+q+1)}$.
  \end{enumerate}
  Then, for any bottom knot $U\in \BT_1$, we have $J_U\in K_1$.
\end{proposition}

\begin{proof}
  This is easily verified using \fullref{thm:12}.
\end{proof}

Corollaries \ref{thm:10}, \ref{thm:23}, \ref{thm:25}, \ref{thm:15} and
\ref{thm:29} have similar versions for bottom knots.

\begin{remark}
  \label{r11}
  One can replace $Z(H)$ in this section with any $\modk $--subalgebra of~$Z(H)$.
\end{remark}

\begin{remark}
  \label{r10}
  For $n\ge 0$, $\tilde{H}^{\otimes n}$ has a natural $Z(H)$--module structure
  induced by multiplication of elements of $Z(H)$ on one of the tensor
  factors in $H^{\otimes n}$.  For each $f\in \modB (m,n)$, the map
  $\tilde{\modJ }(f)$ is a $Z(H)$--module map.  One can show that there is
  a monoidal functor $\tilde{\modJ }'\col \modB \rightarrow \Mod_{Z(H)}$ of $\modB $ into
  the category $\Mod_{Z(H)}$ of $Z(H)$--modules which maps each object
  $\modb ^{\otimes n}$ into $\tilde{H}^{\otimes n}$ and each morphism $f$ into
  $\tilde{\modJ }'(f)=\tilde{\modJ }(f)$.
\end{remark}

\section{Band-reembedding of bottom tangles}
\label{sec10}

\subsection{Refined universal invariants of links}
\label{sec:refin-univ-invar}
As in \fullref{sec:bottom-tangles}, for $n\ge 0$, let $\LL_n$ denote
the set of isotopy classes of $n$--component, framed, oriented, ordered
links for $n\ge 0$.  There is a surjective function
\begin{equation*}
  \clo\col \BT_n\rightarrow \LL_n, \quad T\mapsto \clo(T).
\end{equation*}
We study an algebraic condition for two bottom tangles to yield the
same closure.

\begin{definition}
  \label{def:5}
  Two bottom tangles $T,T'\in \BT_n$ are said to be related
  by a {\em band-reembedding} if there is $W\in \BT_{2n}$ such that
  \begin{gather}
    \label{eq:29}
    T = F_{\modb} ((\epsilon _{\HH}\otimes \HH)^{\otimes n})(W),\quad
    T' = F_{\modb} (\ad_{\HH}^{\otimes n})(W).
  \end{gather}
  See \fullref{fig:reemb} for an example.
\labellist\small
\pinlabel {(a)} [b] at 130 -5
\pinlabel {$W = $} [r] at 40 140
\pinlabel {(b)} [b] at 365 -5
\pinlabel {$T = $} [r] at 265 140
\pinlabel {(c)} [b] at 600 -5
\pinlabel {$T' = $} [r] at 500 140
\endlabellist
  \FIG{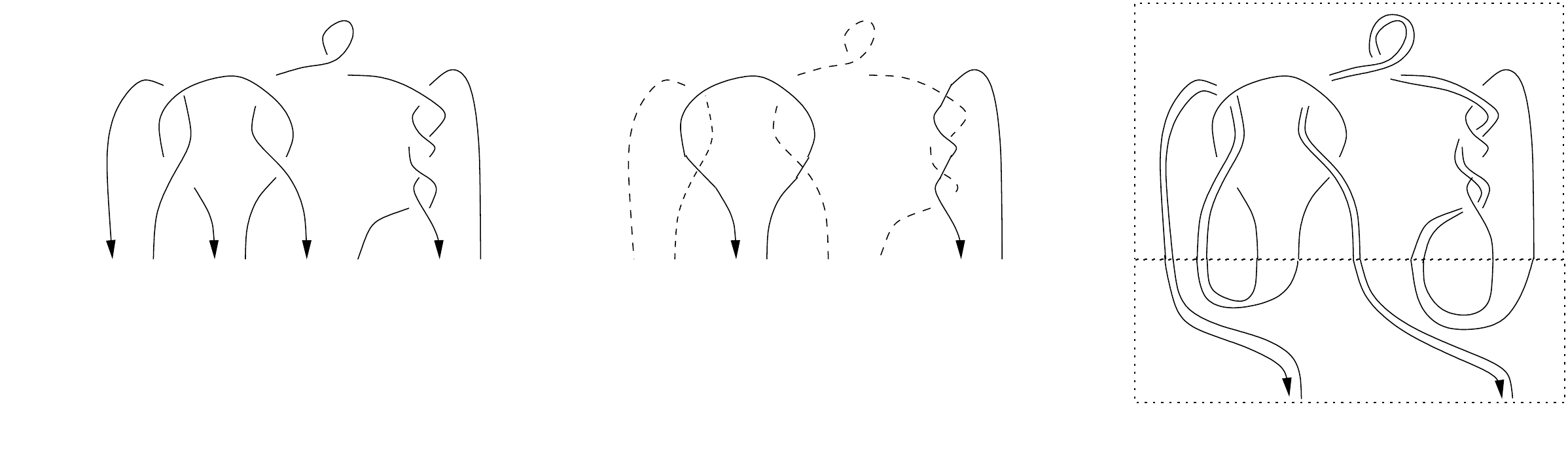}{(a) A bottom tangle $W\in \BT_4$.  (b) The tangle
    $T=F_{\modb} ((\epsilon _{\HH}\otimes \HH)^{\otimes 2})(W)\in \BT_2$.  (c) The tangle
    $T'=F_{\modb} (\ad_{\HH}^{\otimes 2})(W)\in \BT_2$.}{height=34mm}
\end{definition}

If we regard a bottom tangle as a based link in a natural way, then
band-reembedding corresponds to changing the basing.

\begin{proposition}
  \label{thm:22}
  Two bottom tangles $T,T'\in \BT_n$ are related by a
  band-reembedding if and only if $\clo(T)=\clo(T')$.
\end{proposition}

\begin{proof}
  Suppose that $T$ and $T'$ are related by a band-reembedding with
  $W\in \BT_{2n}$ as in \fullref{def:5}.  Note that the tangle
  $T$ is obtained from $W$ by removing the components of $W$ of odd
  indices, and the tangle $T'$ is obtained from the composition tangle
  $1_{\modb} ^{\otimes n}T$ by reembedding the $n$ bands in $1_{\modb ^{\otimes n}}$ along
  the components of $W$ of odd indices.  Hence we easily see that
  $\clo(T)=\clo(T')$.

  Conversely, if $\clo(T)=\clo(T')$, then we can express $T'$ as a
  result from $1_{\modb ^{\otimes n}}T$ by reembedding the $n$ bands in
  $1_{\modb ^{\otimes n}}$, and we can arrange by isotopy that there is
  $W\in \BT_{2n}$ satisfying \eqref{eq:29}.
\end{proof}

Let $H$ be a ribbon Hopf algebra.  \fullref{thm:22} implies
that if two bottom tangles $T,T'\in \BT_n$ satisfies $\clo(T)=\clo(T')$,
then there is $W\in \BT_{2n}$ such that
\begin{gather*}
  J_T = (\epsilon \otimes 1_H)^{\otimes n}(J_W),\quad J_{T'} = \ad^{\otimes n}(J_W).
\end{gather*}
If $K_i\subset H^{\otimes i}$ for $i=0,1,2,\ldots$ are as in \fullref{thm:12},
then we have
\begin{equation}
  \label{eq:5}
  J_{T'}-J_T \in  (\ad^{\otimes n}-(\epsilon \otimes 1_H)^{\otimes n})(K_{2n}).
\end{equation}
Hence we have the following.

\begin{theorem}
  \label{thm:32}
  Let $K_n\subset H^{\otimes n}$ for $n=0,1,2,\ldots$ be $\modZ $--submodules satisfying
  the conditions of \fullref{thm:12}.  Set
  \begin{equation}
    \label{eq:40}
    K'_n = (\ad^{\otimes n}-(\epsilon \otimes 1_H)^{\otimes n})(K_{2n}).
  \end{equation}
 Then, for each $n\ge 0$, the
  function
  \begin{equation*}
    J\col \BT_n \rightarrow  K_n, \quad T\mapsto J_T,
  \end{equation*}
  induces a link invariant
  \begin{equation*}
    \bar{J} \col  \LL_n \rightarrow  K_n/K'_n
  \end{equation*}
\end{theorem}

Note that if we set $K_n=H^{\otimes n}$ in \fullref{thm:32}, then we
get the usual definition of universal invariant of links.

\begin{remark}
  \label{r3}
  The idea of \fullref{thm:32} can also be used to obtain ``more
  refined'' universal invariants for more special classes of links.
  For example, let us consider links and bottom tangles of zero
  linking matrices.  Suppose in \fullref{thm:22} that $T$ and
  $T'$ are of zero linking matrices.  Note that the tangle
  $W\in \BT_{2n}$ satisfying \eqref{eq:29} is not necessarily of zero
  linking matrix, but the $n$--component bottom tangle
  $F_{\modb} ((\epsilon _{\HH}\otimes \HH)^{\otimes n})(W)$, which is equivalent to $T$, is of
  zero linking matrix.  Hence we can replace the conditions for the
  $K_n$ in \fullref{thm:12} with weaker ones.  We hope to give
  details of this idea in future publications.
\end{remark}

\subsection{Ribbon discs}
\label{sec:ribbon-discs}
We close this section with a result which is closely related to
\fullref{thm:22}.

A {\em ribbon disc} for a bottom knot $T$ is a ribbon disc for the
knot $T\cup \gamma $, where $\gamma \subset [0,1]^2\times \{0\}$ is the line segment such that
$\partial \gamma =\partial T$.  Clearly, a bottom knot admits a ribbon disc if and only
if the closure $\clo(T)$ of $T$ is a ribbon knot.

\begin{theorem}
  \label{thm:35}
  For any bottom knot $T\in \BT_1$, the following conditions are
  equivalent.
  \begin{enumerate}
  \item $T$ admits a ribbon disc.
  \item There is an integer $n\ge 0$ and a bottom tangle $W\in \BT_{2n}$
    such that
    \begin{gather}
      \label{e20}
      \eta _{\modb} ^{\otimes n} = F_{\modb} ((\epsilon _{\HH}\otimes \HH)^{\otimes n})(W),\\
      \label{e21}
      T = \mu _{\modb} ^{[n]}F_{\modb} (\ad_{\HH}^{\otimes n})(W).
    \end{gather}
  \end{enumerate}
\end{theorem}

\begin{proof}
  A ribbon disc bounded by $T$ can be decomposed into $n+1$ mutually
  disjoint discs $D_0,D_1,\ldots,D_n$ and $n$ mutually disjoint bands
  $b_1,\ldots,b_n$ for some $n\ge 0$ satisfying the following
  conditions.
  \begin{enumerate}
  \item For each $i=1,\ldots,n$ the band $b_i$ joins $D_0$ and $D_i$,
  \item $D_0$ is a disc attached to the bottom square of the cube
    along a line segment,
  \item The only singularities of the ribbon disc are ribbon
    singularities in $D_i\cap b_j$ for $1\le i,j\le n$.  (We do not allow
    ribbon singularity in $D_0$.)
  \end{enumerate}
  For example, see \fullref{fig:ribbon} (a).

\labellist\tiny
\hair=1pt
\pinlabel {\small(a)} [b] at 90 -5
\pinlabel {$b_1$} [r] at 53 85
\pinlabel {$b_2$} [l] at 147 95
\pinlabel {$D_0$} [tr] at 142 46
\pinlabel {$D_1$} [l] at 78 133
\pinlabel {$D_2$} [l] at 163 140
\pinlabel {\small(b)} [b] at 330 -5
\endlabellist
\FIG{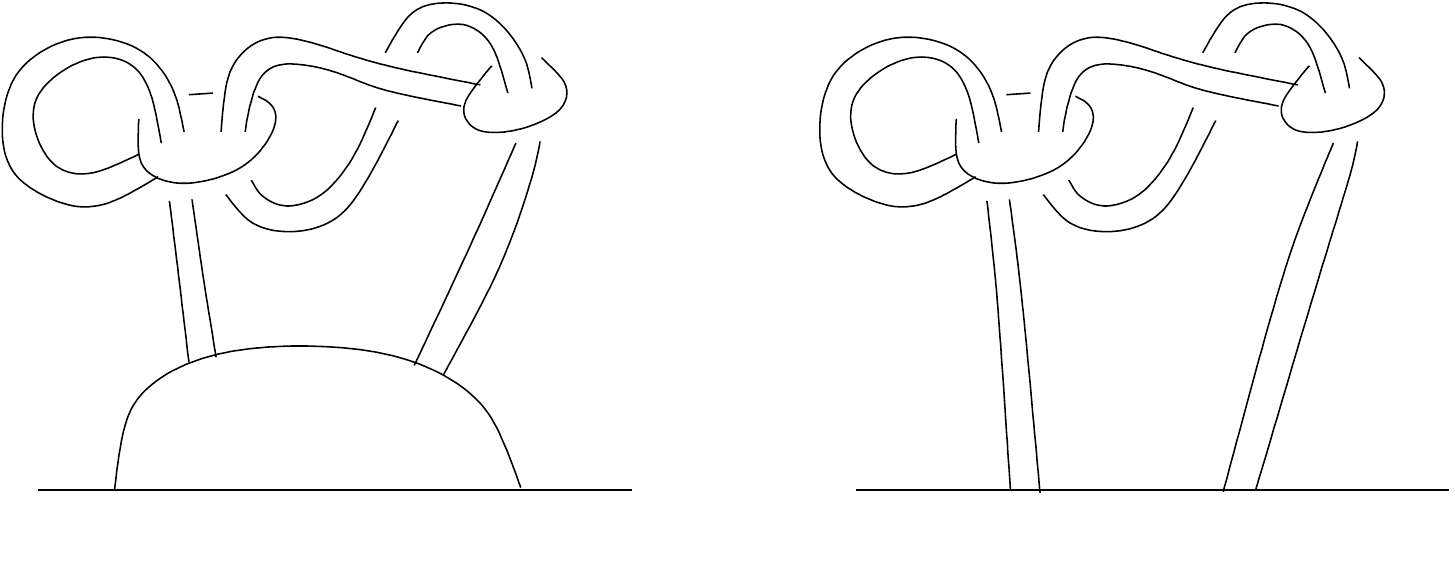}{(a) A ribbon bottom knot $T$ and a ribbon disc decomposed as
  $D_0\cup D_1\cup D_2\cup b_1\cup b_2$.  (b) The bottom tangle $T'$.}{height=30mm}
  Let $T'\in \BT_n$ be the bottom tangle obtained from $T$ by removing
  $D_0$ and regarding the rest as an $n$--component bottom tangle, see
  \fullref{fig:ribbon} (b).  Then $\clo(T')=\clo(\eta _n)$ is an unlink.
  It follows from \fullref{thm:22} that $T'$ and $\eta _n$ are
  related by band-reembedding.  Since $T=\mu _{\modb} ^{[n]}T'$, we have the
  assertion.
\end{proof}

\begin{remark}
  \label{r19}
  In \fullref{thm:35}, one can replace \eqref{e20} with
  \begin{equation}
    \label{e22}
    T = \mu _{\modb} ^{[n]}F_{\modb} (Y_{\HH}^{\otimes n})(W).
  \end{equation}
  We sketch how to prove this claim.  For a ribbon bottom knot $T$,
  there is a Seifert surface $F$ of genus $n\ge 0$ for $T$ and simple
  closed curves $c_1,\ldots ,c_n$ in~$F$ satisfying the following
  conditions.
  \begin{enumerate}
  \item $c_1,\ldots ,c_n$ generates a Lagrangian subgroup of
  $H_1(F;\modZ )\simeq \modZ ^{2n}$.
  \item As a framed link in the cube, $c_1\cup \cdots\cup c_n$ is a $0$--framed
  unlink.  Here the framings of $c_i$ are induced by the surface $F$.
  \end{enumerate}
  (The Seifert surface obtained from a ribbon disc by smoothing the
  singularities in a canonical way has the above property.)  By
  isotopy one can arrange the curves $c_1,\ldots ,c_n$ as in \fullref{fig:ribbonseifert}.
\labellist\small
\pinlabel {$T = $} [r] at 25 50
\pinlabel {$D(`W`)$} at 115 70
\pinlabel {$\cdots$} at 119 55
\pinlabel {$\cdots$} at 119 26
\pinlabel {$c_1$} [t] at 86 20
\pinlabel {$F$} [b] at 145 0
\pinlabel {$c_n$} [t] at 183 20
\endlabellist
\FIGn{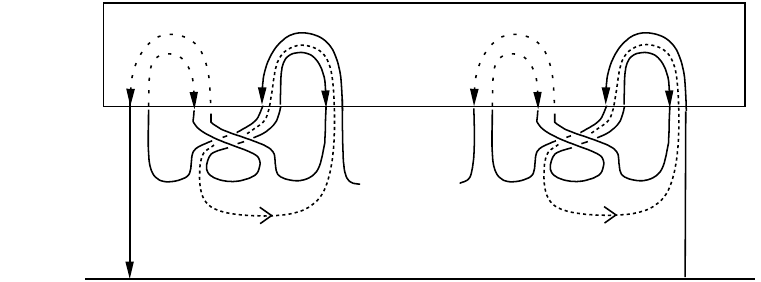}{}{height=25mm} The
  part bounded by a rectangle is a double of a bottom tangle
  $W\in \BT_{2n}$ satisfying \eqref{e20} and \eqref{e22}.
\end{remark}

\section{Algebraic versions of Kirby moves and Hennings
  $3$--man\-i\-fold invariants}
\label{sec11}
An important application of the universal invariants is to the
Hennings invariant of $3$--manifolds and its generalizations.
Hennings \cite{Hennings:96} introduced a class of invariants of
$3$--manifolds associated to quantum groups, which use right integrals
and no finite-dimensional representations.  The Hennings invariants
are studied further by Kauffman and Radford
\cite{Kauffman-Radford:95}, Ohtsuki \cite{Ohtsuki:95}, Lyubashenko
\cite{Lyubashenko}, Kerler \cite{Kerler:97,Kerler:03-2}, Sawin \cite{Sawin},
Virelizier \cite{Virelizier}, etc.  As mentioned in, or at least
obvious from, these papers, the Hennings invariants can be formulated
using universal link invariants.

In this section, we reformulate the Hennings $3$--manifold invariants
using universal invariants of bottom tangles.  For closely related
constructions, see Kerler \cite{Kerler:97} and Virelizier
\cite{Virelizier}.

For a Hopf algebra $A$ in a braided category, we set
\begin{equation*}
  h_A = (\mu _A\otimes A)(A\otimes \Delta _A)\col A^{\otimes 2}\rightarrow A^{\otimes 2}.
\end{equation*}
The morphism $h_A$ is invertible with the inverse
\begin{equation*}
  h_A^{-1}=  (\mu _A\otimes A)(A\otimes S_A\otimes A)(A\otimes \Delta _A).
\end{equation*}
For the transmutation $\bH$ of a ribbon Hopf algebra $H$, we have
\begin{gather*}
  h_{\underline H} = (\mu _H\otimes 1_H)(1_H\otimes \bD),
\end{gather*}
which should not be confused with
\begin{equation*}
  h_H=(\mu _H\otimes 1_H)(1_H\otimes \Delta _H)
\end{equation*}
defined for the Hopf algebra $H$ in the symmetric monoidal category
$\Mod_{\modk} $.

The following is a version of Kirby's theorem \cite{Kirby:78}.  Note
that each move is formulated in an algebraic way.  Therefore we may
regard the following as an {\em algebraic version of Kirby's theorem}.

\begin{theorem}
  \label{thm:34}
  For two bottom tangles $T$ and $T'$, the two $3$--manifolds
  $M_T=(S^3)_{\clo(T)}$ and $M_{T'}=(S^3)_{\clo(T')}$ obtained from
  $S^3$ by surgery along $\clo(T)$ and $\clo(T')$, respectively, are
  orientation-preserving homeomorphic if and only if $T$ and $T'$ are
  related by a sequence of the following moves.
  \begin{enumerate}
  \item Band-reembedding.
  \item Stabilization: Replacing $U\in \BT_n$ with
  $U\otimes v_{\pm} \in \BT_{n+1}$, or its inverse operation.
  \item Handle slide: Replacing $U\in \BT_n$ ($n\ge 2$) with
    $F_{\modb} (h_{\HH}\otimes \HH^{\otimes (n-2)})(U)$, or its inverse
  operation.
  \item Braiding: Replacing $U\in \BT_n$ with $\beta U$, where $\beta \in \modB (n,n)$
    is a doubled braid.
  \end{enumerate}
\end{theorem}

\begin{proof}
  First we see the effects of the moves listed above.

  (1)\qua A band-reembedding does not change the closure, hence the result
  of surgery.

  (2)\qua The effect of stabilization move $U\leftrightarrow U\otimes v_+$ is
  depicted in \fullref{fig:stabilization}.
\labellist\tiny
\pinlabel {$1$} [b] at 28 0
\pinlabel {$2$} [b] at 63 0
\pinlabel {$\cdots$} at 83 35
\pinlabel {$n$} [b] at 99 0
\pinlabel {$1$} [b] at 225 0
\pinlabel {$2$} [b] at 260 0
\pinlabel {$\cdots$} at 281 35
\pinlabel {$n$} [b] at 297 0
\pinlabel {$n{+}1$} [b] at 353 0
\pinlabel {$U$} at 63 90
\pinlabel {$U$} at 265 90
\pinlabel {stabilization} [b] at 160 68
\endlabellist
\FIG{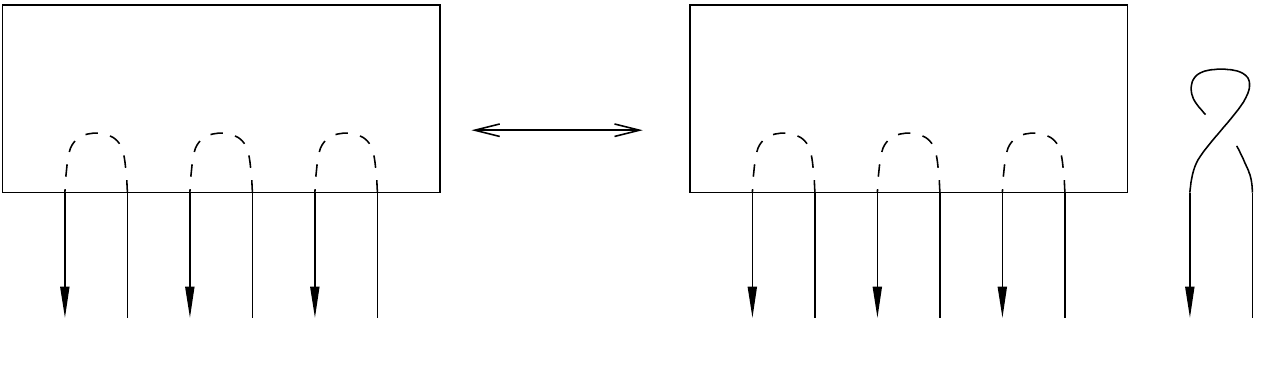}{A stabilization move $U\leftrightarrow U\otimes
  v_+$}{height=25mm} The case of $v_-$ is similar.  The closures
  $\clo(U)$ and $\clo(U\otimes v_{\pm} )=\clo(U)\sqcup\clo(v_{\pm} )$ are
  related by Kirby's stabilization move.  Hence they have the same result
  of surgery.

  (3)\qua The effect of handle slide move $U\leftrightarrow
  U':=F_{\modb} (h_{\HH}\otimes \HH^{\otimes (n-2)})(U)$ is
  depicted in \fullref{fig:handle-slide}.
\labellist\tiny
\pinlabel {$1$} [b] at 28 0
\pinlabel {$2$} [b] at 63 0
\pinlabel {$\cdots$} at 83 35
\pinlabel {$n$} [b] at 99 0
\pinlabel {$1$} [b] at 225 0
\pinlabel {$2$} [b] at 260 0
\pinlabel {$\cdots$} at 281 35
\pinlabel {$n$} [b] at 297 0
\pinlabel {handle-slide} [b] at 160 68
\endlabellist
  \FIG{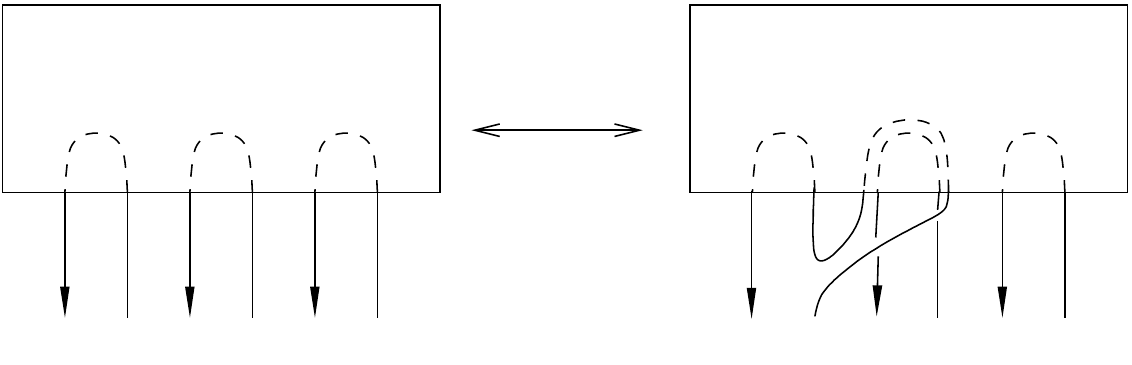}{A
  handle slide move}{height=25mm}  It is easy to see that the
  closures $\clo(U)$ and $\clo(U')$ are related by a Kirby handle slide
  move of the first component over the second.   Hence they have the same
  result of surgery.

  (4)\qua A braiding move (see \fullref{fig:braiding-move}) just changes
  the order of the components at the closure level, and hence does not
  change the result of surgery.
\labellist\tiny
\pinlabel {$1$} [b] at 28 0
\pinlabel {$2$} [b] at 63 0
\pinlabel {$\cdots$} at 83 5
\pinlabel {$n$} [b] at 99 0
\pinlabel {$1$} [b] at 225 0
\pinlabel {$2$} [b] at 260 0
\pinlabel {$\cdots$} at 281 5
\pinlabel {$n$} [b] at 297 0
\pinlabel {braiding} [b] at 160 68
\endlabellist
\FIG{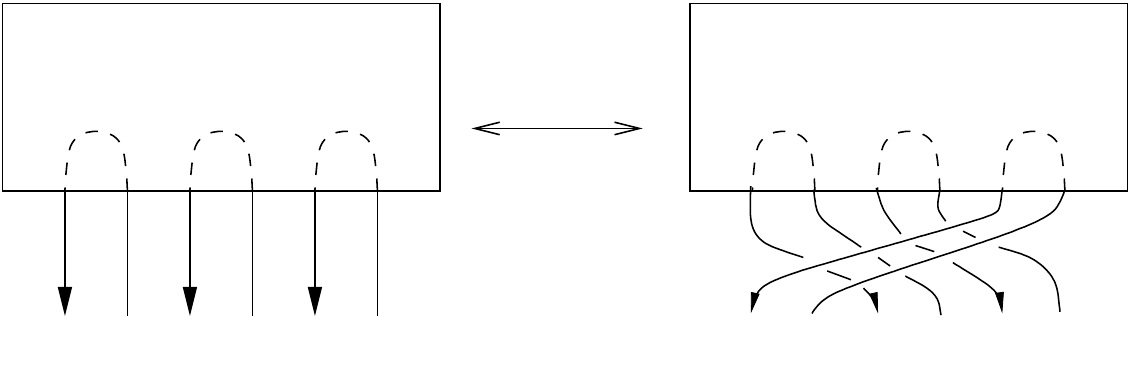}{An example of
  braiding move}{height=25mm}

  The ``if'' part of the theorem follows from the above observations.
  To prove the ``only if'' part, we assume that $M_T$ and $M_{T'}$ are
  orientation-preserving homeomorphic to each other.  By Kirby's
  theorem, there is a sequence from $\clo(T)$ to $\clo(T')$ of
  stabilizations, handle slides, orientation changes of components,
  and changes of ordering.  It suffices to prove that if $\clo(T)$ and
  $\clo(T')$ are related by one of these moves, then $T$ and $T'$ are
  related by the moves listed in the theorem.

  If  $\clo(T)$ and $\clo(T')$ are related by
  change of ordering, then it is easy to see that  $T$ and $T'$ are
  related by a braiding move and a band-reembedding.

  If $\clo(T)$ and $\clo(T')$ are related by stabilization, ie,
  $\clo(T')=\clo(T)\sqcup O_{\pm} $, where $O_{\pm} $ is an unknot of framing
  $\pm 1$, then $T'$ and $T\otimes v_{\mp}$ are related by a band-reembedding.

  Suppose $\clo(T)$ and $\clo(T')$ are related by handle slide of a
  component of $\clo(T)$ over another component of $\clo(T)$.  By
  conjugating with change of ordering, we may assume that the first
  component of $\clo(T)$ is slid over the second component of $\clo(T)$.
  Then there is $T''\in \BT_n$ such that $T''$ is obtained from $T$ by a
  band-reembedding, and $T'$ is obtained from $T''$ by a handle slide
  move or its inverse.  For example, see \fullref{fig:handle-slide-example}.
\labellist\small
\pinlabel {$T = $} [r] at 34 255
\pinlabel {$1$} [t] at 50 198
\pinlabel {$2$} [t] at 170 198
\pinlabel {\parbox{34pt}{\tiny handle slide for closure}} [t] at 222 272
\pinlabel {$1$} [t] at 265 198
\pinlabel {$2$} [t] at 392 198
\pinlabel {$ = T'$} [l] at 415 255
\pinlabel {\parbox{38pt}{\tiny band reembedding}} [r] at 121 159
\pinlabel {\tiny isotopy} [l] at 318 159
\pinlabel {$T'' = $} [r] at 34 94
\pinlabel {$1$} [t] at 50 15
\pinlabel {$2$} [t] at 99 15
\pinlabel {\tiny handle slide} [t] at 220 91
\pinlabel {$1$} [t] at 266 15
\pinlabel {$2$} [t] at 314 15
\endlabellist
\FIG{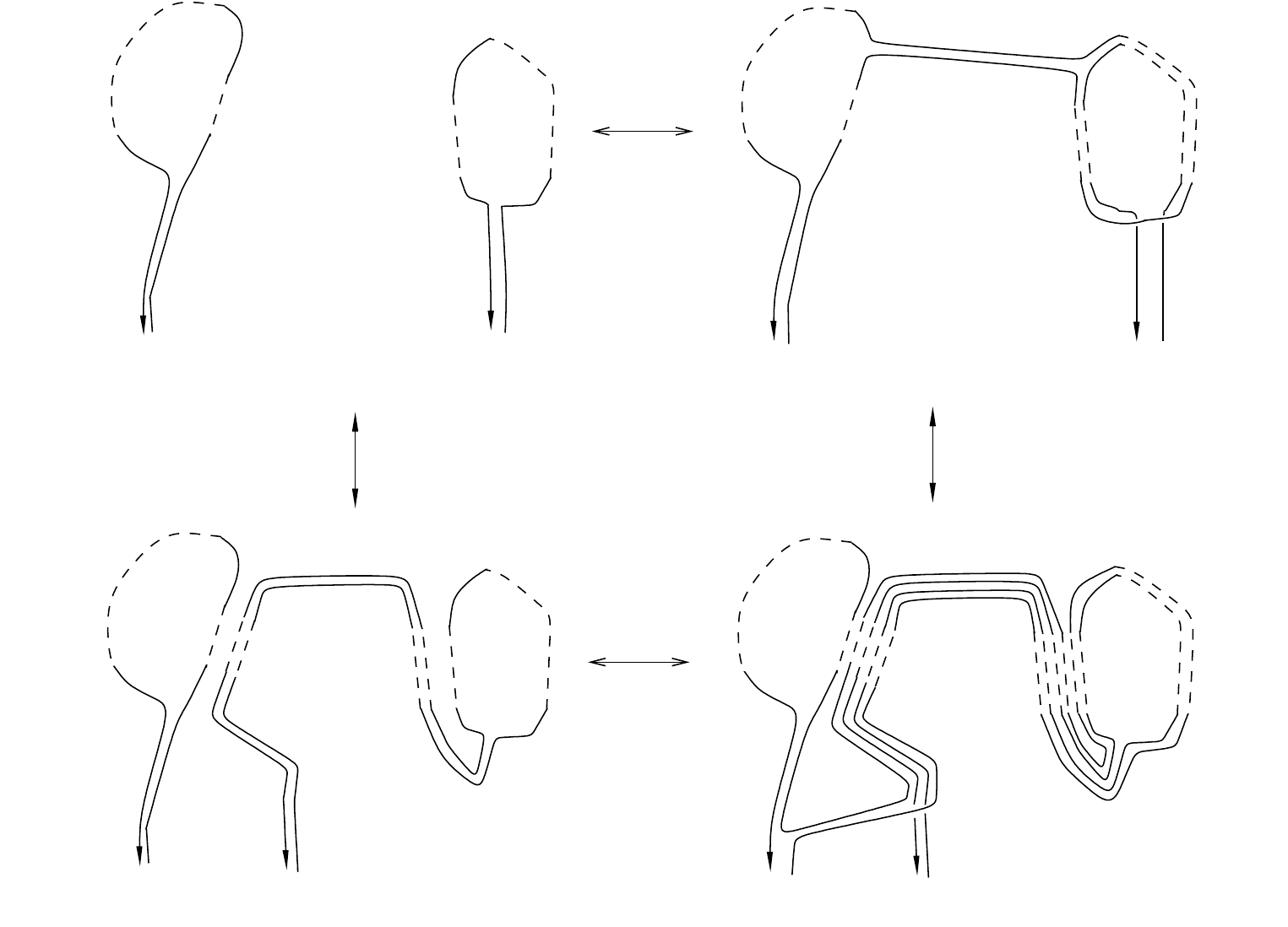}{Here
  only the first and the second components are depicted in each
  figure}{height=65mm}

  Suppose $\clo(T)$ and $\clo(T')$ are related by orientation change of
  $i$th components with $1\le i\le n$.  It suffices to show that $\clo(T)$
  and $\clo(T')$ are related by a sequence of handle slides, changes of
  orientation, changes of ordering and stabilizations.  We may assume
  $i=1$ by change of ordering.  Using stabilization, we may safely
  assume $n\ge 2$.  Now we see that change of orientation of the first
  component can be achieved by a sequence of handle slides and change
  of ordering.  Suppose that $L=L_1\cup L_2$ is an $2$--component link.
  There is a sequence of moves $L=L^0\rightarrow L^1\rightarrow \cdots \rightarrow L^4=-L_1\cup L_2$, with
  $L^i=L^i_1\cup L^i_2$ as depicted in \fullref{fig:ori-change-seq}.
\labellist\tiny
\hair=2pt
\pinlabel {$L_1^0$} [tl] at 32 100
\pinlabel {$L_2^0$} [tr] at 72 100
\pinlabel {\parbox{25pt}{slide $L_1^0$ over $L_2^0$}} [t] at 134 112
\pinlabel {$L_1^1$} [tl] at 195 100
\pinlabel {$L_2^1$} [l] at 265 109
\pinlabel {\parbox{25pt}{slide $L^1_2$ over $L^1_1$}} [t] at 329 112
\pinlabel {$L_2^2$} [br] at 393 97
\pinlabel {$L_1^2$} [tl] at 397 90
\pinlabel {isotopy} [t] at 15 42
\pinlabel {$L_2^2$} [tr] at 86 42
\pinlabel {$L_1^2$} [b] at 128 32
\pinlabel {\parbox{25pt}{slide $L_1^2$ over $L_2^2$}} [t] at 214 32
\pinlabel {$L_2^3$} [tl] at 277 17
\pinlabel {$L_1^3$} [tr] at 318 17
\pinlabel {\parbox{30pt}{change the ordering}} [t] at 385 32
\pinlabel {$L_1^4$} [tl] at 449 17
\pinlabel {$L_2^4$} [tr] at 490 17
\endlabellist
  \FIGn{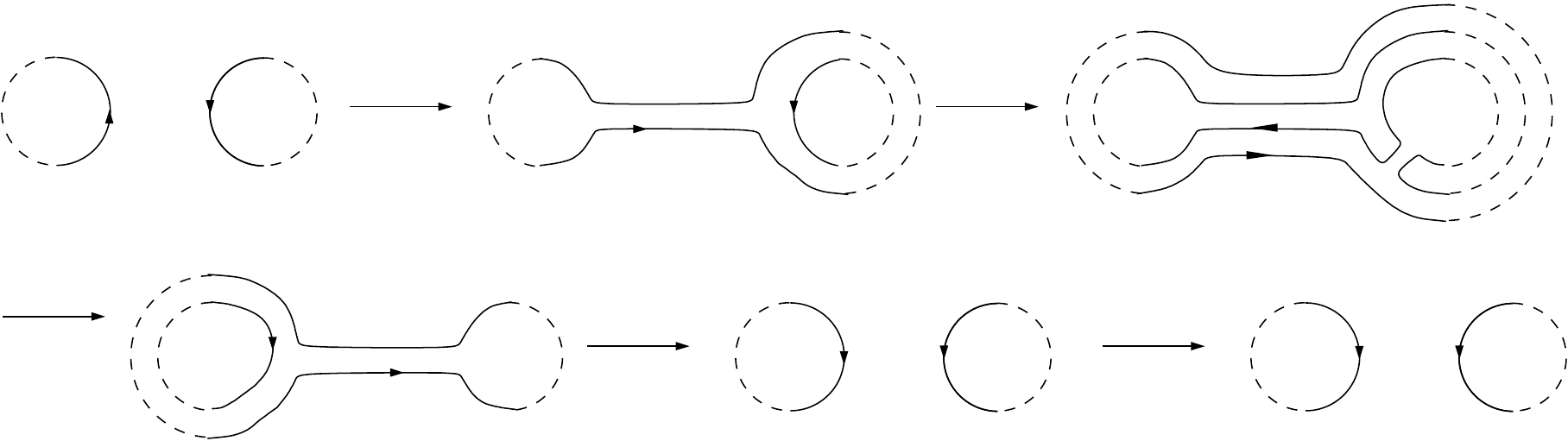}{}{height=35mm} (Note here that, in the move
  $L^1\rightarrow L^2$, we can slide the second component over the first by
  conjugating the handle slide move by changing of the ordering.)
  Hence we have the assertion.
\end{proof}

\begin{lemma}
  \label{lem:3}
  If $f\col H\rightarrow \modk $ is a left $H$--module homomorphism, then we have
  \begin{equation*}
    (1_H\otimes f)\Delta = (1_H\otimes f)\bD = (f\otimes 1_H)\bD.
  \end{equation*}
\end{lemma}

\begin{proof}
  The lemma is probably well known.  We prove it for completeness.

  The first identity is proved using \ref{eq:8} as follows.
  \begin{equation*}
    \begin{split}
      (1\otimes f)\bD(x)
      &=(1\otimes f)(\sum x_{(1)}S(\beta )\otimes \alpha \trr x_{(2)})\\
      &= \sum x_{(1)}S(\beta )\otimes f(\alpha \trr x_{(2)})
      = \sum x_{(1)}S(\beta )\otimes \epsilon (\alpha )f(x_{(2)})\\
      &= \sum x_{(1)}\otimes f(x_{(2)})
      = (1\otimes f)\Delta (x).
    \end{split}
  \end{equation*}
  The identity $(f\otimes 1)\bD(x)=(1\otimes f)\Delta (x)$ follows similarly by using
 the identity
$$\bD(x)=\sum(\beta \trr x_{(2)})\otimes \alpha x_{(1)}\quad
  \text{for } x\in H.\proved$$
\end{proof}

Now we formulate the Hennings invariant in our setting.  If
$\chi \col H\rightarrow \modk $ is a left $H$--module homomorphism and a left integral on
$H$, then we can define a $3$--manifold invariant.  This invariant is
essentially the same as the Hennings invariant defined using the right
integral, since left and right integrals interchange under application
of the antipode.

\begin{proposition}
  \label{thm:37}
  Let $\chi \col H\rightarrow \modk $ be a left $H$--module homomorphism.  Then the
  following conditions are equivalent.
  \begin{enumerate}
  \item $\chi $ is a left integral on $H$, ie,
    \begin{equation}
    \label{eq:37}
    (1\otimes \chi )\Delta  = \eta \chi \col H\rightarrow H.
    \end{equation}
  \item $\chi $ is a two-sided integral on $\bH$ in $\Mod_H$, ie,
    \begin{equation}
      \label{eq:39}
      (1\otimes \chi )\bD=(\chi \otimes 1)\bD = \eta \chi \col H\rightarrow H.
    \end{equation}
  \end{enumerate}
  Suppose either (hence both) of the above holds, and also suppose
  that $\chi (\modr ^{\pm 1})\in \modk $ is invertible.  Then there is a unique
  invariant $\tau _{H,\chi }(M)\in \modk $ of connected, oriented, closed
  $3$--manifolds $M$ such that for each bottom tangle $T\in \BT_n$ we
  have
  \begin{equation}
    \label{eq:23}
    \tau _{H,\chi }(M_T)
    =\frac{\chi ^{\otimes n}(J_T)}{\chi (\modr ^{-1})^{\sigma _+(T)}\chi (\modr )^{\sigma _-(T)}},
  \end{equation}
  where $\sigma _+(T)$ (resp. $\sigma _-(T)$) is the number (with multiplicity)
  of the positive (resp. negative) eigenvalues of the linking matrix
  of $T$, and $M_T=(S^3)_{\clo(T)}$ denote the result from $S^3$ of
  surgery along $\clo(T)$.
\end{proposition}

\begin{proof}
  The first assertion follows from \fullref{lem:3}.

  In the following, we show that the right hand side of \eqref{eq:23}
  is invariant under the moves described in \fullref{thm:34}.

  First we consider the stabilization move.  Suppose $T\in \BT_n$ and
  $T'=T\otimes v_{\pm} \in \BT_{n+1}$.  Then one can easily verify
  $\tau _{H,\chi }(M_T)=\tau _{H,\chi }(M_{T'})$ using
  \begin{equation*}
    \chi ^{\otimes (n+1)}(J_{T\otimes v_{\pm} })
    =\chi ^{\otimes n}(J_T)\cdot\chi (\modr ^{\pm 1}).
  \end{equation*}
  Since the other moves does not change the number of components and
  the number of positive (resp. negative) eigenvalues of the linking
  matrix, it suffices to verify that $\chi ^{\otimes n}(J_T)=\chi ^{\otimes n}(J_{T'})$
  for $T,T'\in \BT_n$ related by each of the other moves.

  Suppose $T$ and $T'$ are related by a band-reembedding.  Then there is
  $W\in \BT_{2n}$ satisfying \eqref{eq:29}.  Since $\chi $ is a left
  $H$--module homomorphism, we have
  \begin{equation*}
    \begin{split}
      \chi ^{\otimes n}(J_T)
      &= \chi ^{\otimes n}(\epsilon \otimes 1_H)^{\otimes n}(J_W)
      = \chi ^{\otimes n}\ad^{\otimes n}(J_W)
      = \chi ^{\otimes n}(J_{T'}).
    \end{split}
  \end{equation*}
  Suppose that $T$ and $T'$ are related by a handle slide, ie,
  $T'=F_{\modb} (h_{\HH}\otimes \HH^{\otimes (n-2)})(T)$.  Since $\chi $ is a two-sided
  integral on $\bH$, we have
  \begin{equation*}
    \begin{split}
      \chi ^{\otimes n}(J_{T'})
      &=\chi ^{\otimes n}(h_{\bH}\otimes 1_H^{\otimes (n-2)})(J_T)\\
      &=\chi ^{\otimes n}((\mu _H\otimes 1_H)(1_H\otimes \bD)\otimes 1_H^{\otimes (n-2)})(J_T)\\
      &=(\chi \mu _H\otimes \chi ^{\otimes (n-2)})(1_H\otimes (1_H\otimes \chi )\bD\otimes 1_H^{\otimes (n-2)})(J_T)\\
      &=(\chi \mu _H\otimes \chi ^{\otimes (n-2)})(1_H\otimes \eta _H\chi \otimes 1_H^{\otimes (n-2)})(J_T)\\
      &=\chi ^{\otimes n}(J_T).
    \end{split}
  \end{equation*}
  Suppose that $T$ and $T'$ are related by a braiding move.  We may
  assume that $T'=(\psi _{\modb ,\modb }^{\pm 1})_{(i-1,n-i-1)}T$ with $1\le i\le n-1$.
  Since $(\chi \otimes \chi )\psi _{H,H}=\chi \otimes \chi $, it follows that
  \begin{equation*}
    \chi ^{\otimes n}(J_{T'}) =\chi ^{\otimes n}(\psi _{H,H}^{\pm 1})_{(i-1,n-i-1)}(J_T)
    =\chi ^{\otimes n}(J_T).
  \end{equation*}
  This completes the proof.
\end{proof}

\begin{remark}
  \label{thm:38}
  One can verify that the invariant $\tau _{H,\chi }(M)$ is equal (up to a
  factor determined only by the first Betti number of $M$) to the
  Hennings invariant of $M$ associated to the right integral
  $\chi S\col H\rightarrow \modk $.
\end{remark}

\begin{remark}
  \label{thm:40}
  Some results in \fullref{sec9} can be used together with
  \fullref{thm:37} to obtain results on the range of values of
  the Hennings invariants for various class of $3$--manifolds.  Recall
  from Hennings \cite{Hennings:96} and Ohtsuki \cite{Ohtsuki:93} (see
  also Virelizier \cite{Virelizier}) that
  the $sl_2$ Reshetikhin--Turaev invariants can be defined using a
  universal link invariant associated to a finite-dimensional quantum
  group $U_q(sl_2)'$ at a root of unity, and a certain trace function
  on $U_q(sl_2)'$.  Hence results in \fullref{sec9} can also be
  used to study the range of values of the Reshetikhin--Turaev
  invariants.
\end{remark}

\section{String links and bottom tangles}
\label{sec:string-links-bottom}
An $n$--component string link $T=T_1\cup \cdots \cup  T_n$ is a tangle consisting
$n$ arcs $T_1,\ldots,T_n$, such that for $i=1,\ldots,n$ the $i$th component
$T_i$ runs from the $i$th upper endpoint to the $i$th lower endpoint.
In other words, $T$ is a morphism $T\in \modT (\downarrow ^{\otimes n},\downarrow ^{\otimes n})$ homotopic
to $\downarrow ^{\otimes n}$.
(Here and in what follows, the endpoints are counted from the left.)
The closure operation is as depicted in \fullref{fig:string-link-new}
(a), (b).
\labellist\small
\pinlabel {(a)} [b] at 46 0
\pinlabel {(b)} [b] at 218 0
\endlabellist
\FIG{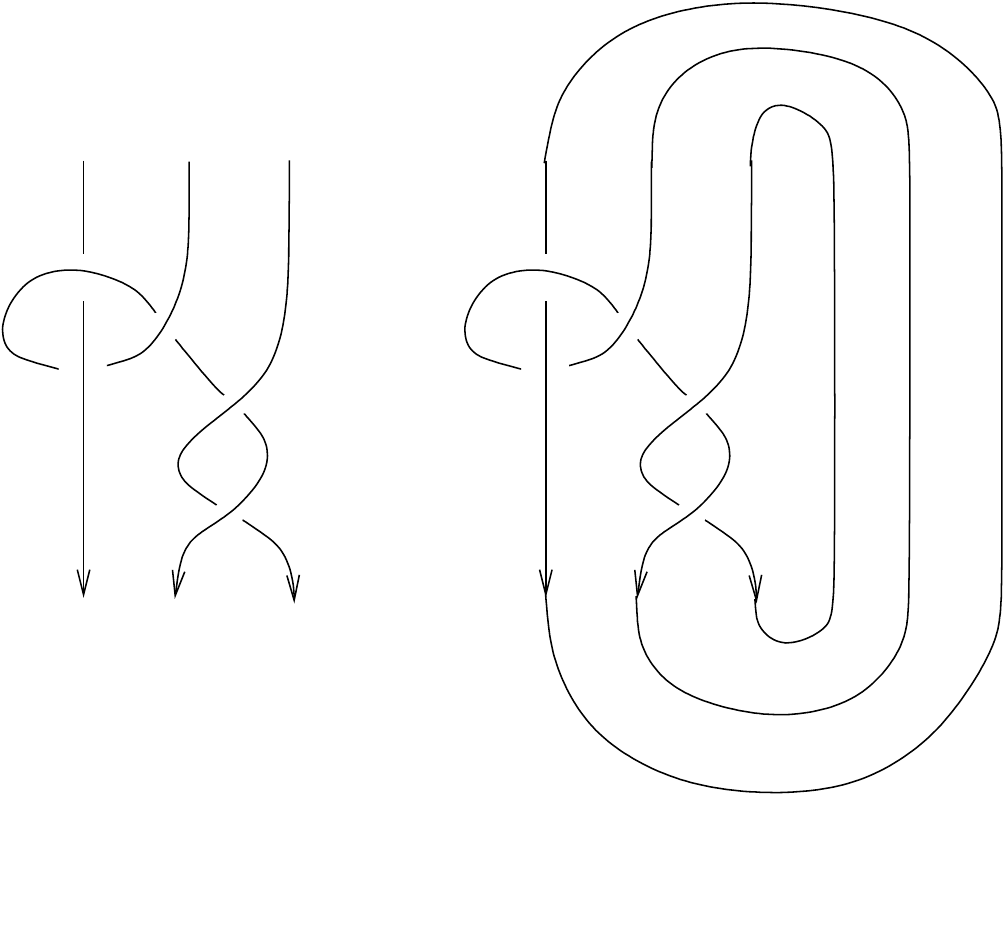}{(a) A string link.  (b) Its
closure.}{height=35mm} One can use string links to study
links via the closure operation.

As in \fullref{sec:pres-string-links}, we denote by $\SL_n$ the
submonoid of $\modT (\downarrow ^{\otimes n},\downarrow ^{\otimes n})$ consisting of the isotopy classes
of the $n$--component framed string links.

Of course, there are many orientation-preserving self-homeomorphisms
of a cube $[0,1]^3$, which transform $n$--component string links into
$n$--component bottom tangles and induces a bijection
$\SL_n\cong\BT_n$.  In this sense, one can think of the notion of
string links and the notion of bottom tangles are equivalent.
However, $\SL_n$ and $\BT_n$ are {\em not} equally convenient.  For
example, the monoid structure in the $\SL_n$ can not be defined in
each $\BT_n$ as conveniently as in $\SL_n$, and also that the external
Hopf algebra structure in the $\BT_n$ can not be defined in the
$\SL_n$ as conveniently as in the $\BT_n$.  It depends on the contexts
which is more useful.

In the following, we define a {\em preferred bijection}
\begin{equation*}
  \tau _n\col \BT_n\rightarrow \SL_n,
\end{equation*}
which enables one to translate results about the bottom tangles into
results about the string links and vice versa.  We define a monoid
structure of each $\BT_n$ such that $\tau _n$ is a monoid homomorphism.
We also study several other structures on $\BT_n$ and $\SL_n$ and
consider the algebraic counterparts for a ribbon Hopf algebra.  The
proofs are straightforward and left to the reader.

For $n\ge 0$, we give $\modb ^{\otimes n}\in \Ob(\modT )$ the standard tensor product
algebra structure (see Majid \cite[Section 2]{Majid:algebras})
\begin{gather*}
  \mu _{\modb ^{\otimes n}}\col \modb ^{\otimes n}\otimes \modb ^{\otimes n}\rightarrow \modb ^{\otimes n},\quad
  \eta _{\modb ^{\otimes n}}\col \one \rightarrow \modb ^{\otimes n},
\end{gather*}
induced by the algebra structure $(\modb ,\mu _{\modb} ,\eta _{\modb} )$, ie,
$\mu _{\modb ^{\otimes 0}}=1_{\one} $, $\mu _{\modb ^{\otimes 1}}=\mu $,
\begin{gather*}
  \mu _{\modb ^{\otimes n}}=
  (\mu _{\modb} \otimes \mu _{\modb ^{\otimes (n-1)}})(\modb \otimes \psi _{\modb ^{\otimes (n-1)},\modb }\otimes \modb ^{\otimes (n-1)})\quad
  \text{for $n\ge 2$},
\end{gather*}
and $\eta _{\modb ^{\otimes n}}=\eta _n$ for $n\ge 0$.
See \fullref{fig:mubn} for example.
\labellist\small
\pinlabel {(a)} [b] at 53 0
\pinlabel {(b)} [b] at 143 0
\pinlabel {(c)} [b] at 297 0
\pinlabel {$=$} at 297 70
\pinlabel {(d)} [b] at 533 0
\pinlabel {$=$} at 510 65
\pinlabel {$=$} at 558 65
\endlabellist
\FIG{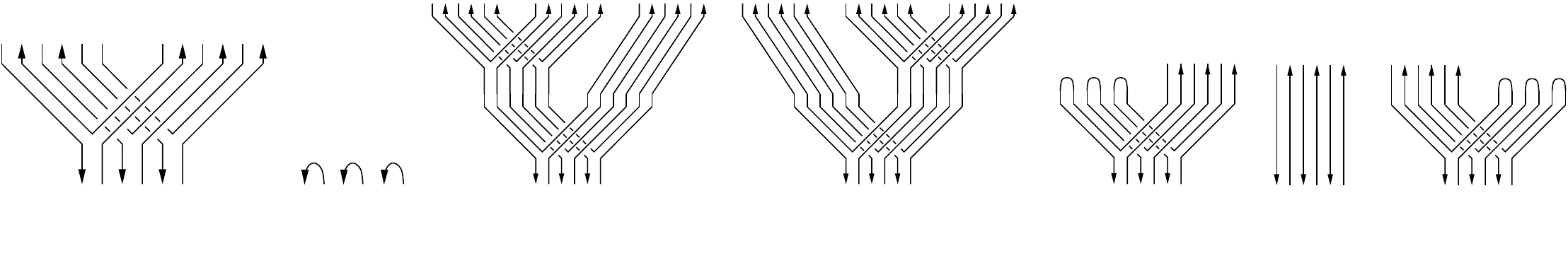}{(a) Multiplication $\mu _{\modb ^{\otimes 3}}$. (b) Unit $\eta _{\modb ^{\otimes 3}}$. (c)
  Associativity. (d) Unitality.}{width=125mm}

We define a monoid structure for $\BT_n$ with multiplication
\begin{equation*}
  \tilde{\mu }_n=\ast\col \BT_n\times \BT_n \rightarrow \BT_n,\quad (T,T')\mapsto T\ast T'
\end{equation*}
defined by
\begin{equation}
  \label{eq:22}
  T\ast T' = \mu _{\modb ^{\otimes n}}(T\otimes T')=
\labellist\tiny
\pinlabel {$T$} at 25 72
\pinlabel {$T'$} at 94 72
\endlabellist
  \raisebox{-7mm}{\incl{15mm}{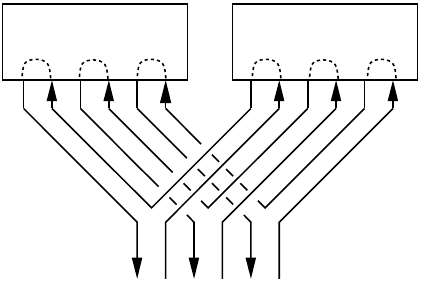}}
\end{equation}
for $T,T'\in \BT_n$, where the figure in the right hand side is for
$n=3$.  Then the set $\BT_n$ has a monoid structure with
multiplication $\tilde{\mu }_n$ and with unit $\eta _n$.

We give $\downarrow \in \Ob(\modT )$ a left $\modb $--module structure defined by the left
action
\begin{equation*}
  \alpha _{\downarrow} =\downarrow \otimes \ev_{\downarrow}
=\raisebox{-3mm}{\incl{7mm}{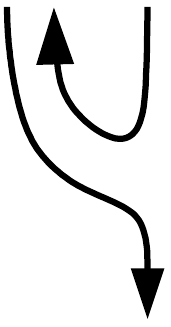}}\col \modb \otimes \downarrow \rightarrow \downarrow .
\end{equation*}
For $n\ge 0$, this left $\modb $--module structure induces in the canonical
way a left $\modb ^{\otimes n}$--module structure for $\downarrow ^{\otimes n}$
\begin{equation*}
  \alpha_{\downarrow^{\otimes n}}\col \modb ^{\otimes n}\otimes \downarrow ^{\otimes n}\rightarrow \downarrow ^{\otimes n},
\end{equation*}
ie, $\alpha_{\downarrow^{\otimes n}}$ is defined inductively by
\begin{gather*}
  \alpha _{\downarrow ^{\otimes 0}} = 1_{\one} ,\quad
  \alpha _{\downarrow ^{\otimes 1}} = \alpha ,\quad
  \\
  \alpha _{\downarrow ^{\otimes n}} =
  (\alpha _{\downarrow} \otimes \alpha _{\downarrow ^{\otimes (n-1)}})(\modb \otimes \psi _{\modb ^{\otimes (n-1)},\downarrow }\otimes \downarrow ^{\otimes (n-1)})
  \quad \text{for $n\ge 2$}.
\end{gather*}
For example, see \fullref{fig:action2}.
\labellist\small
\pinlabel {(a)} [b] at 60 0
\pinlabel {(b)} [b] at 245 0
\pinlabel {$=$} at 240 80
\pinlabel {(c)} [b] at 445 0
\pinlabel {$=$} at 445 80
\endlabellist
\FIG{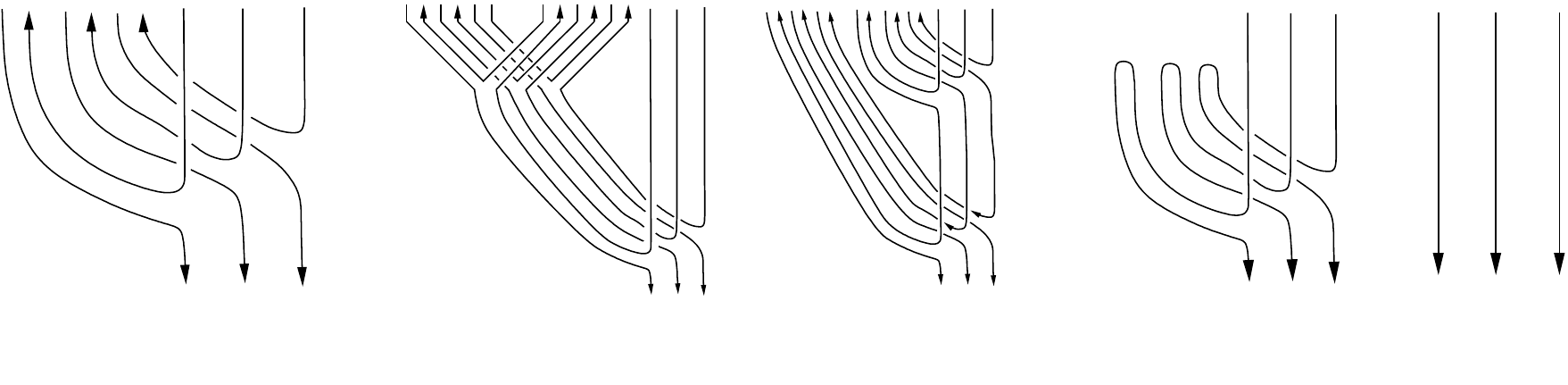}{(a) Left action $\alpha _{\downarrow ^{\otimes 3}}$.  (b) Associativity.
  (c) Unitality.}{height=30mm}

Now we define a function $\tau _n \col  \BT_n \rightarrow  \SL_n$ for $n\ge 0$ by
\begin{equation*}
  \tau _n(T) = \alpha _{\downarrow ^{\otimes n}}(T\otimes \downarrow ^{\otimes n})
\end{equation*}
for $T\in \BT_n$.  In a certain sense, $\tau _n(T)$ is the result of ``letting
$T$ act on $\downarrow ^{\otimes n}$''.  For example, if $T\in \BT_3$, then
\begin{equation*}
  \tau _3(T)
\labellist\tiny
\pinlabel {$T$} at 47 150
\endlabellist
 = \raisebox{-10mm}{\incl{22mm}{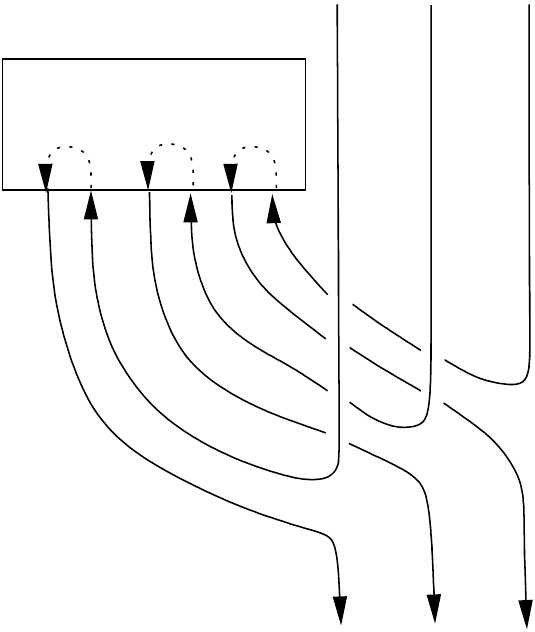}}
\labellist\tiny
\pinlabel {$T$} at 57 46
\endlabellist
 = \raisebox{-7mm}{\incl{14mm}{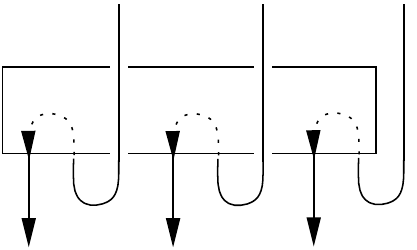}}.
\end{equation*}
The function $\tau _n$ is invertible with the inverse
$\tau _n^{-1}\col \SL_n\rightarrow \BT_n$ given by
\begin{equation*}
  \tau _n^{-1}(L) = \theta _n (T\otimes \uparrow ^{\otimes n})\coev_{\downarrow ^{\otimes n}},
\end{equation*}
where
\begin{equation*}
  \theta _n \col \downarrow ^{\otimes n} \otimes \uparrow ^{\otimes n} \rightarrow \modb ^{\otimes n}
\end{equation*}
is defined inductively by
\begin{gather*}
  \theta _0 = 1_{\one} ,\quad
  \theta _{n+1} = (\downarrow \otimes \psi _{\modb ^{\otimes n},\uparrow })(\downarrow \otimes \theta _n\otimes \uparrow )
\end{gather*}
for $n\ge 1$.  For example, if $L\in \SL_3$, then
\begin{equation*}
  \tau _3^{-1}(L)
\labellist\tiny
\pinlabel {$L$} at 28 143
\endlabellist
= \raisebox{-10mm}{\incl{22mm}{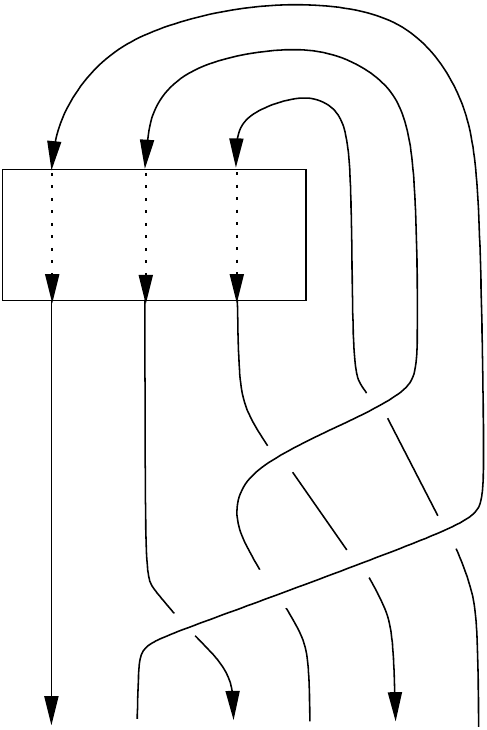}}
\labellist\tiny
\pinlabel {$L$} at 39 42
\endlabellist
= \raisebox{-7mm}{\incl{14mm}{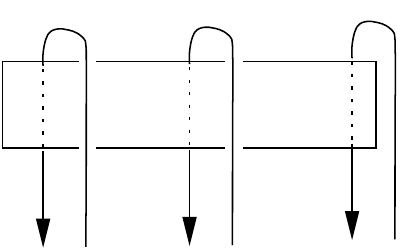}}.
\end{equation*}
The function $\tau _n$ is a monoid isomorphism, ie,
\begin{equation*}
  \tau _n(T\ast T')=\tau _n(T)\tau _n(T'),\quad \tau _n(\eta _n)=\downarrow ^{\otimes n}
\end{equation*}
for $T,T'\in \BT_n$.

As is well known, there is a ``coalgebra-like'' structure on
the $\SL_n$.  For $T\in \SL_n$ and $i=1,\ldots,n$, let $\Delta _i(T)\in \SL_{n+1}$
(resp. $\epsilon _i(T)\in \SL_{n-1}$) be obtained from~$T$ by duplicating
(resp. removing) the $i$th component.  These operations define monoid
homomorphisms
\begin{gather*}
  \Delta _i\col \SL_n \rightarrow  \SL_{n+1},\quad
  \epsilon _i\col \SL_n \rightarrow  \SL_{n-1}.
\end{gather*}
The following diagrams commutes.
\begin{equation}
  \label{eq:107}
  \begin{CD}
    \BT_n @>\check\Delta _{(i-1,n-i)}>> \BT_{n+1}\\
    @V\tau _nVV @VV\tau _{n+1}V\\
    \SL_n @>>\Delta _i> \SL_{n+1}
  \end{CD}
  \qquad
  \begin{CD}
    \BT_n @>\check\epsilon _{(i-1,n-i)}>> \BT_{n-1}\\
    @V\tau _nVV @VV\tau _{n-1}V\\
    \SL_n @>>\epsilon _i> \SL_{n-1} .
  \end{CD}
\end{equation}
Thus the ``coalgebra-like'' structure of the $\SL_n$ corresponds via
the $\tau _n$ to the ``coalgebra-like'' structure in the $\BT_n$.

Now we translate the above observations into the universal invariant
level.  Let $H$ be a ribbon Hopf algebra over a commutative, unital
ring $\modk $.

Define a $\modk$--module homomorphism $\tau '_n\col H^{\otimes n}\rightarrow H^{\otimes n}$, $n\ge 0$,
by $\tau '_0=1_{\modk}  $, $\tau '_1=1_H$, and for $n\ge 2$
\begin{equation*}
  \tau '_n =  \bigl(1_H^{\otimes (n-2)}\otimes \lambda
_2\bigr)\bigl(1_H^{\otimes (n-3)}\otimes \lambda _3\bigr)
  \cdots(1_H\otimes \lambda _{n-1})\lambda _n,
\end{equation*}
where $\lambda _n\col H^{\otimes n}\rightarrow H^{\otimes n}$, $n\ge 2$, is defined by
\begin{equation*}
  \lambda _n\Bigl(\sum x_1\otimes \cdots\otimes x_n\Bigr) =
  \sum x_1\beta  \otimes  (\alpha _{(1)}\trr x_2)\otimes \cdots\otimes (\alpha _{(n-1)}\trr x_n).
\end{equation*}
Then the effects of $\tau _n$ on the universal invariants
is given by
\begin{equation}
  \label{eq:30}
  J_{\tau _n(T)} = \tau '_n(J_T)\quad \text{for $T\in \BT_n$},
\end{equation}
which can be used in translating results about the universal invariant
of bottom tangles into results about the universal invariant of string
links.

We denote by $\bH^{\otimes n}$ the $n$--fold tensor product of $\bH$ in the
braided category $\Mod_H$.  Thus $\bH^{\otimes n}$ is equipped with the
standard algebra structure in $\Mod_H$, with the multiplication
\begin{equation*}
  \ul\mu _n\col \bH^{\otimes n}\otimes \bH^{\otimes n}\rightarrow \bH^{\otimes n}
\end{equation*}
given by $\ul\mu _n = \modJ (\mu _{\modb ^{\otimes n}})$ and with the unit in
$\bH^{\otimes n}$ by $1_H^{\otimes n}$.

The map $\tau '_n$ defines a $\modk $--algebra isomorphism
\begin{equation*}
  \tau '_n\col \bH^{\otimes n}\rightarrow H^{\otimes n},
\end{equation*}
where $H^{\otimes n}$ is equipped with the standard algebra structure.  In
other words, we have
\begin{equation*}
  \tau '_n(\ul\mu _n(x\otimes y))=\tau '_n(x)\tau '_n(y)
\end{equation*}
for $x,y\in H^{\otimes n}$, and $\tau '_n(1_H^{\otimes n})=1_H^{\otimes n}$.

We have algebraic analogues of the diagrams in \eqref{eq:107}
\begin{equation*}
  \begin{CD}
    \bH^{\otimes n} @>\bD_i>>\bH^{\otimes (n+1)}\\
    @V\tau '_nVV @VV\tau '_{n+1}V\\
    H^{\otimes n} @>>\Delta _i> H^{\otimes (n+1)}
  \end{CD}
  \quad \quad
  \begin{CD}
    \bH^{\otimes n} @>\epsilon _i>>\bH^{\otimes (n-1)}\\
    @V\tau '_nVV @VV\tau '_{n-1}V\\
    H^{\otimes n} @>>\epsilon _i> H^{\otimes (n-1)}
  \end{CD}
\end{equation*}
where $\bD_i=1^{\otimes (i-1)}\otimes \bD\otimes 1^{\otimes (n-i)}$,
$\Delta _i=1^{\otimes (i-1)}\otimes \Delta \otimes 1^{\otimes (n-i)}$ and
$\epsilon _i=1^{\otimes (i-1)}\otimes \epsilon \otimes 1^{\otimes (n-i)}$.  If $n=i=1$, then the
commutativity of the diagram on the left, $\bD=\tau _2^{-1}\Delta $, above
coincides \eqref{eq:8}, the definition of the transmuted
comultiplication $\bD$.

\section{Remarks}
\label{sec12}
\subsection{Direct applications of the category $\modB $ for
  representation-colored link invariants}
\label{sec:other-applications}

We have argued that the setting of universal invariants is useful in
the study of representation-colored link invariants.  However, it
would be worth describing how to apply the setting of the category
$\modB $ to the study of representation-colored link invariants without
using universal invariants.

Let $H$ be a ribbon Hopf algebra over a field $\modk $, and let $V$ be a
finite-dimensional left $H$--module.  Let $F^{\modT} _V\col \modT \rightarrow \Mod_H$ denote
the canonical braided functor from the category $\modT $ of framed,
oriented tangles to the category $\Mod_H$ of left $H$--modules, which
maps the object $\downarrow $ to $V$.  Let us denote the restriction of
$F^{\modT} _V$ to $\modB $ by $F_V\col \modB \rightarrow \Mod_H$.  Then $F_V$ maps the object
$\modb $ into $V\otimes V^*$, which we identify with the $\modk $--algebra
$E=\End_{\modk} (V)$ of $\modk $--vector space endomorphisms of~$V$.  As one can
easily verify, the algebra structure of $\modb $ is mapped into that
of~$E$, ie, $F_V(\mu _{\modb} )=\mu _E$, and $F_V(\eta _{\modb} )=\eta _E$, where
$\mu _E\col E\otimes E\rightarrow E$ and $\eta _E\col \modk \rightarrow E$ are the structure morphisms for the
algebra $E$.  Also, the images by $F_V$ of the other generating
morphisms of $\modB $ are determined by
\begin{align*}
  F_V(v_{\pm} )(1_{\modk} )&= \rho _V(\modr ^{\pm 1}),\\
  F_V(c_{\pm} )(1_{\modk} )&=(\rho _V\otimes \rho _V)(c^H_{\pm} ),
\end{align*}
where $\rho _V\col H\rightarrow E$ denotes the left action of $H$ on $V$, ie,
$\rho _V(x)(v)=x\cdot v$ for $x\in H$, $v\in V$.  In fact, we have for each
$T\in \BT_n$
\begin{equation*}
  F_V(T)(1_{\modk} )= \rho _V^{\otimes n}(J_T).
\end{equation*}
Many  results for universal invariants in the previous sections can be
modified into versions for $F_V$.  For example, the following is a
version of \fullref{thm:12}.

\begin{proposition}
  \label{r15}
  Let $K_i\subset E^{\otimes i}$ for $i\ge 0$, be subsets satisfying the following.
  \begin{enumerate}
  \item $1_{\modk} \in K_0$, $1_E,\rho _V(\modr ^{\pm 1})\in K_1$, and
  $(\rho _V\otimes \rho _V)(c^H_{\pm} )\in K_2$.
  \item For $m,n\ge 0$, we have $K_m\otimes K_n\subset K_{m+n}$.
  \item For $p,q\ge 0$ we have
    \begin{align*}
      (1^{\otimes p}\otimes \psi _{E,E}^{\pm 1}\otimes
       1^{\otimes q})(K_{p+q+2}) &\subset K_{p+q+2},\\
      (1^{\otimes p}\otimes \mu _E\otimes 1^{\otimes q})(K_{p+q+2})
      &\subset K_{p+q+1},
    \end{align*}
    where $\psi _{E,E}=F_V(\psi _{\modb ,\modb })\col E\otimes
    E\rightarrow E\otimes E$ is the braiding of
    two copies of $E$ in $\Mod_H$.
  \end{enumerate}
  Then, for any $T\in \BT_n$, $n\ge 0$, we have $F_V(T)\in K_n$.
\end{proposition}

Note that $\rho _V$ can be regarded as a morphism $\rho _V\col \bH\rightarrow E$ in
$\Mod_H$.  In fact, $\rho _V$ is an algebra-morphism, ie,
$\rho _V\mu _H=\mu _E(\rho _V\otimes \rho _V)$, $\rho _V\eta _H=\eta _E$.  One can check that the
morphisms $\rho _V^{\otimes i}\col \bH^{\otimes i}\rightarrow E^{\otimes i}$, $i\ge 0$, form a natural
transformation
\begin{equation*}
  \rho \col  \modJ \Rightarrow F_V\col  \modB \rightarrow \Mod_H,
\end{equation*}
ie, the following diagram is commutative for $i,j\ge 0$, $T\in \modB (i,j)$
\begin{equation*}
  \begin{CD}
    \bH^{\otimes i} @>{\modJ (T)}>> \bH^{\otimes j}\\
    @V{\rho _V^{\otimes i}}VV @VV{\rho _V^{\otimes j}}V\\
    E^{\otimes i} @>>{F_V(T)}> E^{\otimes j}.
  \end{CD}
\end{equation*}

\subsection{Generalizations}

\subsubsection{Non-strict version $\modB ^q$ of $\modB $ and the Kontsevich invariant}
\label{sec:non-strict-version-q}

Recall that there is a non-strict version $\modT ^q$ of $\modT $, ie, the
category of {\em $q$--tangles} (see Le and Murakami \cite{Le-Murakami}), whose objects
are parenthesized tensor words in $\downarrow $ and $\uparrow $ such as
$$(\downarrow \otimes \downarrow )\otimes
  (\uparrow \otimes (\downarrow \otimes \uparrow ))$$
and whose morphisms are isotopy classes of
tangles.  In a natural way, one can define non-strict braided
subcategories $\modB ^q$ (resp. $\modB ^q_0$) of~$\modB $ whose objects are
parenthesized tensor words of $\modb =\downarrow \otimes \uparrow $, such as $(\modb \otimes \modb )\otimes \modb $, and
the morphisms are isotopy classes of tangles in $\modB $ (resp. $\modB _0$).
The results for $\modB $ and $\modB _0$ in Sections \ref{sec4}, \ref{sec5} and
\ref{sec6} can be easily generalized to results for the non-strict
braided categories $\modB ^q$ and $\modB ^q_0$.  Recall \cite{Le-Murakami}
that the Kontsevich invariant can be formulated as a (non-strict)
monoidal functor $Z\col \modT ^q\rightarrow \modA $ of $\modT ^q$ into a certain ``category of
diagrams''.  It is natural to expect that the non-strict versions of
the results for $\modB $ and $\modB _0$ in the present paper can be applied to
the Kontsevich invariant and can give some integrality results for the
Kontsevich invariant.

\subsubsection{Ribbon Hopf algebras in symmetric monoidal category}
\label{sec:ribbon-hopf-algebras}
Universal invariant of links and tangles can be defined for any ribbon
Hopf algebra $H$ in any symmetric monoidal category $\modM $.  If $T$ is a
tangle consisting of $n$ arcs and no circles, then the universal
invariant $J_T$ takes values in $\modM (\one _{\modM} ,H^{\otimes n})$,
where $\one _{\modM} $ is
the unit object in $\modM $.  Most of the results in \fullref{sec8}
can be generalized to this setting.  In particular, there is a braided
functor $\modJ \col \modB \rightarrow \Mod_H$ such that $\modJ (\modb)=H$, where $\Mod_H$ is the
category of left $H$--modules in $\modM $, and $H$ is given the left
$H$--module structure via the adjoint action.

We comment on two interesting special cases below.

\subsubsection{Complete ribbon Hopf algebras and quantized enveloping
  algebras}
\label{sec:complete-ribbon-hopf}
The universal invariant of tangles can also be defined for any ribbon
{\em complete} Hopf algebra $H$ over a linearly topologized,
commutative, unital ring $\modk $.  The construction of universal
invariant can be generalized to ribbon complete Hopf algebras in an
obvious way.  This case may be considered as the special case of
\fullref{sec:ribbon-hopf-algebras}, since $H$ is a ribbon Hopf
algebra in the category of complete $\modk $--modules.

An important example of a complete ribbon Hopf algebra is the
$h$--adic quantized enveloping algebra $U_h(\g)$ of a simple Lie
algebra $\g$.  In future papers
\cite{Habiro:in-preparation,Habiro-Le:in-preparation}, we will
consider this case and prove some integrality results of the universal
invariants.

\subsubsection{Universal invariants and virtual tangles}
\label{sec:univ-invar-virt}

Kauffman \cite{Kauffman:99} introduced {\em virtual knot theory} (see
also Goussarov, Polyak and Viro \cite{Goussarov-Polyak-Viro}, Kamada and
Kamada \cite{Kamada-Kamada}, and Sakai \cite{Sakai1,Sakai2}).  A
virtual link is a diagram in a plane similar to a link diagram but
allowing ``virtual crossings''.  There is a preferred equivalence
relation among virtual links called ``virtual isotopy'', and two
virtually isotopic virtual links are usually regarded as the same.
There is also a weaker notion of equivalence called ``virtual regular
isotopy'', and it is observed by Kauffman \cite{Kauffman:99} that many quantum
link invariants can be extended to invariants of virtual regular
isotopy classes of virtual links.  The notion of virtual links are
naturally generalized to that of virtual tangles.  {\em Virtual framed
isotopy} is generated by virtual regular isotopy and the move
\begin{equation*}
  \def\s{15mm}
  \raisebox{-6mm}{\incl{\s}{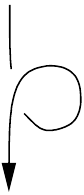}}
  \;\leftrightarrow\;\raisebox{-6mm}{\incl{\s}{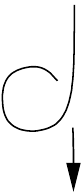}}.
\end{equation*}
An extreme case of \fullref{sec:ribbon-hopf-algebras} is the case
of the symmetric monoidal category $\langle \HH_r\rangle $ freely generated by a
ribbon Hopf algebra $\HH_r$.  The universal tangle invariant $J_T$ in
this case is very closely related to virtual tangles.  We can
construct a canonical bijection between the set
$\langle \HH_r\rangle (\one ,\HH_r^{\otimes n})$ and the set of the virtual framed isotopy
classes of $n$--component ``virtual bottom tangles''.  For $T\in \BT_n$,
the universal invariant $J_T$ takes values in
$\Mod_{\HH_r}(\one ,\HH_r^{\otimes n})\subset \langle \HH_r\rangle (\one ,\HH_r^{\otimes n})$.  Thus, the
universal invariant associated to $\HH_r$ takes values in the virtual
bottom tangles.  For $n\ge 0$, the function
\begin{equation*}
  \modJ _{0,n}=\modJ \col \BT_n\bigl(=\modB (0,n)\bigr)
  \rightarrow  \Mod_{\HH_r}(\one ,\HH_r^{\otimes n}),\quad
  T\mapsto \modJ (T),
\end{equation*}
is injective.  We conjecture that $\modJ _{0,n}$ is surjective.  If this
is true, then we can regard it as an algebraic characterization of
bottom tangles among virtual bottom tangles.  We can formulate similar
conjectures for general tangles and links.  This may be regarded as a
new way to view virtual knot theory in an algebraically natural way.
(We also remark here that there is another (perhaps more natural) way
to formulate virtual knot theory in a category-theoretic setting,
which uses the tangle invariant associated to the symmetric monoidal
category $\langle \HH_r,\modV \rangle $ freely generated by a ribbon Hopf algebra
$\HH_r$ and a left $\HH_r$--module $\modV $ with left dual.)  We plan to
give the details of the above in future publications.

\subsubsection{Quasitriangular Hopf algebras and even bottom tangles}
\label{sec:quasitriangular}

We can generalize our setting to a quasitriangular Hopf algebra, which
may not be ribbon.  For a similar idea of defining quantum invariants
associated to quasitriangular Hopf algebras, see Sawin \cite{Sawin}.

It is convenient to restrict our attention to {\em even-framed}
bottom tangles up to regular isotopy.  Here a tangle $T$ is
even-framed if the closure of each component of $T$ is of even
framing, or, in other words, each component of $T$ has even number of
self crossings.  Let $\modB^{\ev}$ denote the subcategory of $\modB $ such that
$\Ob(\modB ^{\ev})=\Ob(\modB )$ and $\modB^{\ev}(m,n)$ consists of $T\in \modB (m,n)$
even-framed.  Then we have the following.
\begin{enumerate}
\item $\modB ^{\ev}$ is a braided subcategory of $\modB $.
\item $\modB ^{\ev}$ is generated as a braided subcategory of $\modB $ by the
  object $\modb $ and the morphism $\mu _{\modb} ,\eta _{\modb} ,c_+,c_-$.  (This follows
  from \fullref{r2}.)
\item $\modB ^{\ev}$ inherits from $\modB $ an external Hopf algebra structure.
\item There is a braided functor $\modJ ^{\ev}\col \modB ^{\ev}\rightarrow \Mod_H$, defined
  similarly to the ribbon case, and we have an analogue of \fullref{thm:63}.  Hence we have a topological interpretation of
  transmutation of a quasitriangular Hopf algebra.  If $H$ is ribbon,
  then $\modJ ^{\ev}$ is the restriction of $\modJ $ to $\modB ^{\ev}$.
\end{enumerate}

\subsection{The functor $\check\modJ $}
\label{sec:functor-check}

The functor $\modJ \col \modB \rightarrow \Mod_H$ is not faithful for any ribbon Hopf
algebra $H$.  For example, we have $\modJ (t_{\downarrow} \otimes \uparrow
)=\modJ (\downarrow \otimes t_{\uparrow} )$ for any
$H$ but we have $t_{\downarrow} \otimes \uparrow \neq\downarrow \otimes
t_{\uparrow} $.  Using the construction by
Kauffman \cite{Kauffman:93}, one can construct a functor
$\check\modJ \col \modB \rightarrow \operatorname{Cat}(H)$, which distinguishes more 
tangles than $\modJ $.  Here $\operatorname{Cat}(H)$ is the category
defined in \cite{Kauffman:93}, and $\check\modJ $ is just the restriction
to $\modB $ of the functor $F\col \modT \rightarrow \operatorname{Cat}(H)$ defined in
\cite{Kauffman:93}.  Since each $T\in \modB (m,n)$ consists of $m+n$ arc
components, $\check\modJ (T)$ can be defined as an element of
$H^{\otimes (m+n)}$.  If we take $H$ as the Hopf algebra in the braided
category $\langle \HH_r\rangle $ freely generated by a ribbon Hopf algebra as in
\fullref{sec:univ-invar-virt}, then $\check\modJ $ is faithful.

We have not studied the functor $\check\modJ $ in the present paper
because our aim of introducing the category $\modB $ is to provide a
useful tool to study bottom tangles.  The functor $\modJ \col \modB \rightarrow \Mod_H$ is
more suitable than $\check\modJ $ for this purpose.

\subsection{The category $\B$ of bottom tangles in handlebodies}
\label{sec:category-0}

In future papers, we will give details of the following.

Let $\B$ denote the category of {\em bottom tangles in handlebodies},
which is roughly defined as follows.  For $n\ge 0$, let $V_n$ denote a
``standard handlebody of genus~$n$'', which is obtained from the cube
$[0,1]^3$ by adding $n$ handles in a canonical way, see \fullref{fig:handlebody} (a).
\labellist\small
\hair=2pt
\pinlabel {(a)} [b] at 100 0
\pinlabel {$V_2$} at 100 80
\pinlabel {(b)} [b] at 330 0
\pinlabel {$T_1$} [r] at 293 55
\pinlabel {$T_2$} [tr] at 293 100
\pinlabel {$T_3$} [l] at 376 57
\pinlabel {(c)} [b] at 570 0
\pinlabel {$T_1$} [r] at 527 55
\pinlabel {$T_2$} [tr] at 520 105
\pinlabel {$T_3$} [l] at 612 60
\endlabellist
\FIG{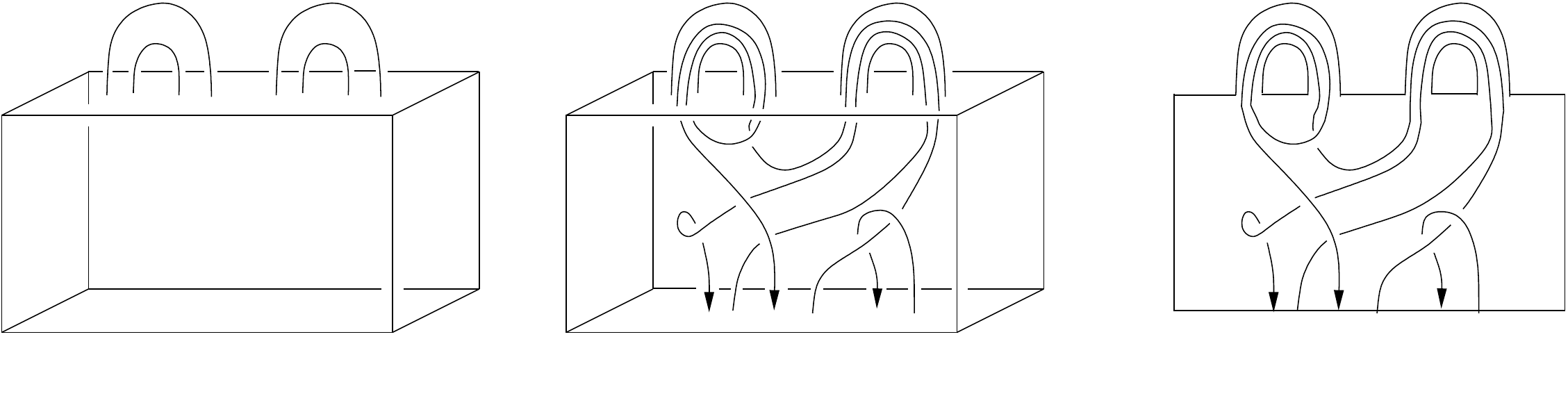}{(a) A standard handlebody
$V_2$ of genus $2$.  (b) A $3$--component bottom tangle
$T=T_1\cup T_2\cup T_3$ in $V_2$.  (c) A diagram for $T$.}{height=33mm} An
$n$--component {\em bottom tangle in $V_n$} is a framed, oriented
tangle $T$ in $V_n$ consisting of $n$ arc components $T_1,\ldots,T_n$,
such that, for $i=1,\ldots,n$, $T_i$ starts at the $2i$th endpoint on
the bottom and end at the $(2i-1)$st endpoint on the bottom.  See
\fullref{fig:handlebody} (b) for example, which we usually draw as
the projected diagram as in (c).

The category $\B$ is defined as follows.  Set
$\Ob(\B)=\{0,1,2,\ldots\}$.  For $m,n\ge 0$, the set $\B(m,n)$ is the set of
isotopy classes of $n$--component bottom tangles in $V_m$.  For
$T\in \B(l,m)$ and $T'\in \B(m,n)$, the composition $T'T\in \B(l,n)$ is
represented by the $l$--component tangle in $V_n$ obtained as follows.
First let $E_T$ denote the ``exterior of $T$ in $V_l$'', ie, the
closure of $V_l\setminus N_T$, with $N_T$ a tubular neighborhood of $T$ in
$V_l$.  Note that $E_T$ may be regarded as a cobordism from the
connected oriented surfaces $F_{l,1}$ of genus $l$ with one boundary
component to $F_{m,1}$.  Thus there is a natural way to identify the
``bottom surface'' of $N_T$ with the ``top surface'' of $V_m$.  By
gluing $E_T$ and $V_m$ along these surfaces, we obtain a $3$--manifold
$E_T\cup V_m$, naturally identified with $V_l$.  The tangle $T'$, viewed
as a tangle in $E_T\cup V_m\cong V_l$, represents the composition $T'T$.
\fullref{fig:composite} shows an example.
\labellist\small
\pinlabel {(a)} [b] at 80 0
\pinlabel {(b)} [b] at 270 0
\pinlabel {(c)} [b] at 450 0
\endlabellist
\FIG{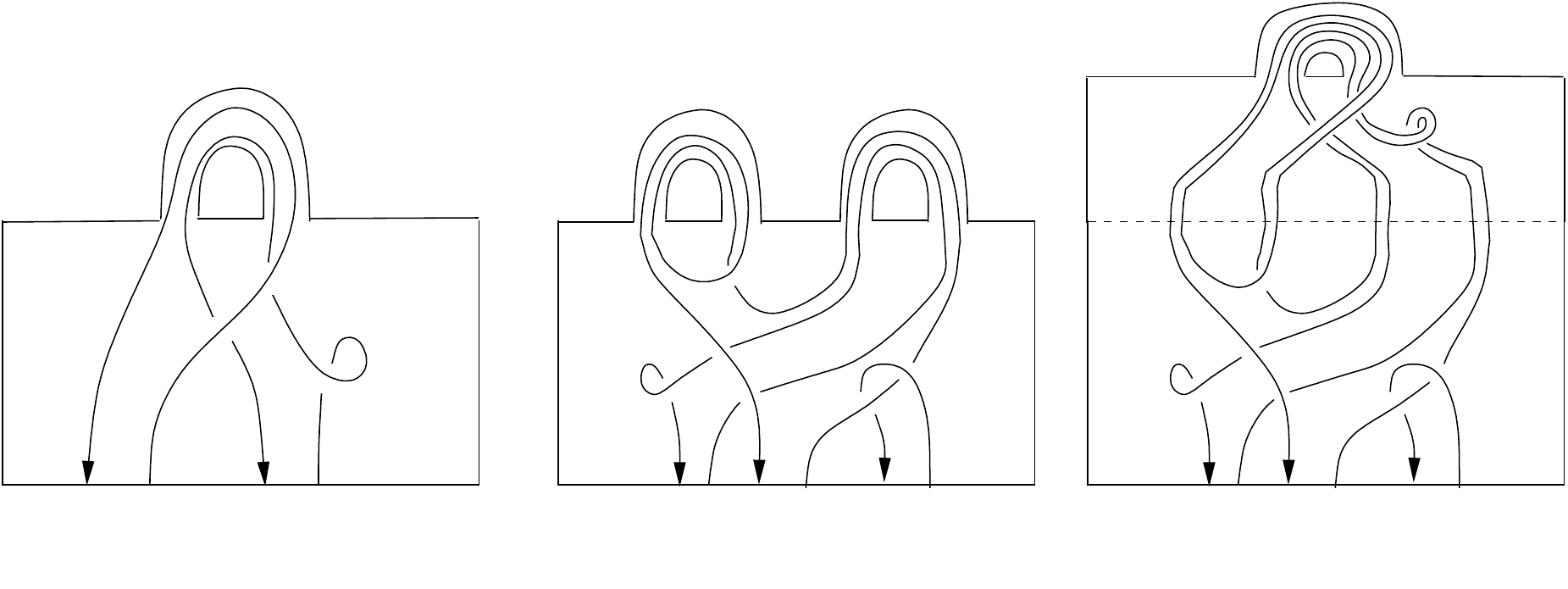}{(a) A
bottom tangle $T\in \B(1,2)$.  (b) A bottom tangle $T'\in \B(2,3)$.
(c) The composition $T'T\in \B(1,3)$.}{height=35mm} For $n\ge 0$, the identity
morphism $1_n\col n\rightarrow n$ is represented by the bottom tangle depicted in
\fullref{fig:cob}.
\labellist\small
\pinlabel {$1_2$} [b] at 60 105
\pinlabel {$\psi_{1,1}$} [b] at 160 105
\pinlabel {$\psi_{1,1}^{-1}$} [b] at 245 105
\pinlabel {$\mu_{\B}$} [b] at 40 0
\pinlabel {$\eta_{\B}$} [b] at 120 0
\pinlabel {$\Delta_{\B}$} [b] at 180 0
\pinlabel {$\epsilon_{\B}$} [b] at 235 0
\pinlabel {$S_{\B}$} [b] at 290 0
\endlabellist
\FIG{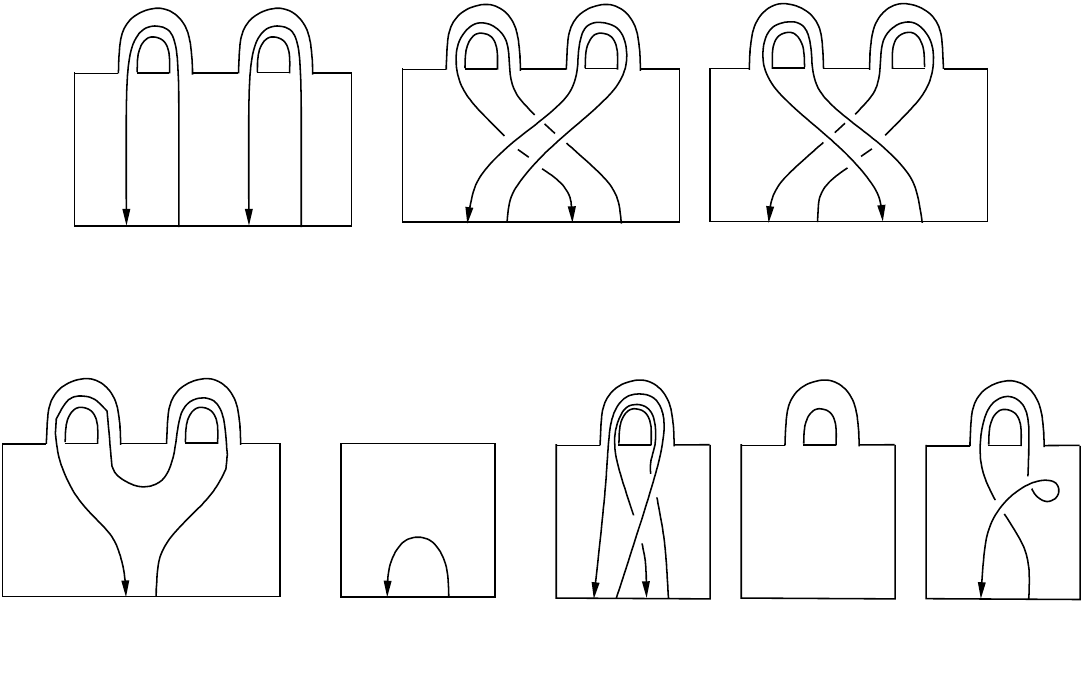}{The morphisms $1_2$, $\psi _{1,1}$,
$\psi _{1,1}^{-1}$, $\mu _{\B}$, $\eta _{\B}$, $\Delta _{\B}$, $\epsilon
_{\B}$ and $S_{\B}$ in
$\B$}{height=50mm} We can prove that the category $\B$ is well
defined.  The category $\B$ has the monoidal structure given by
horizontal pasting.

There is a functor $\xi \col \modB \rightarrow \B$ such that $\xi (\modb ^{\otimes n})=n$, and,
for $T\in \modB (m,n)$, the tangle $\xi (T)\in \B(m,n)$ is obtained from $T$ by
pasting a copy of the identity bottom tangle $1_{\modb ^{\otimes m}}$ on the top
of $T$.  This functor is monoidal, and the braiding structure for $\modB $
induces that for $\B$ via $\xi $, see \fullref{fig:cob}.

Also, there is a Hopf algebra structure
$H_{\B}=(1,\mu_{\B},\eta _{\B},\Delta _{\B},\epsilon _{\B},S_{\B})$ in the usual sense for the
object $1\in \Ob(\B)$.  Graphically, the structure morphisms
$\mu _{\B},\eta _{\B},\Delta _{\B},\epsilon _{\B},S_{\B}$ for $\B$ is as depicted in \fullref{fig:cob}.

The external Hopf algebra structure in $\modB $ is mapped by $\xi $ into the
external Hopf algebra structure in $\B$ associated to $H_{\B}$.

For $n\ge 0$, the function
\[
\xi \col \BT_n\rightarrow \B(0,n)
\]
is bijective.  Hence we can identify the set $\B(0,n)$ with the set of
$n$--component bottom tangles.  (However, $\xi \col \modB (m,n)\rightarrow \B(m,n)$ is
neither injective nor surjective in general.  Hence the functor $\xi $ is
neither full nor faithful.)  Using \fullref{thm:2}, we can prove
that $\B$ is generated as a braided category by the morphisms
\begin{equation}
  \label{e16}
  \mu _{\B},\eta _{\B},\Delta _{\B},\epsilon _{\B},S_{\B},S_{\B}^{-1},v_{\B,+},v_{\B,-},
\end{equation}
where $v_{\B,\pm }=\xi (v_{\pm} )$.

There is a natural, faithful, braided functor $i\col \B\rightarrow \modC $, where
$\modC $ is the category of cobordisms of surfaces with connected
boundary as introduced by Crane and Yetter \cite{Crane-Yetter:99} and
by Kerler \cite{Kerler:99}, independently.  The objects of $\modC $ are
the nonnegative integers $0,1,2,\ldots$, the morphisms from $m$ to $n$ are
(certain equivalence classes of) cobordisms from $F_m$ to $F_n$, where
$F_m$ is a compact, connected, oriented surface of genus $n$ with
$\partial F_m\cong S^1$.  (See also
Habiro \cite{Habiro:claspers}, Kerler \cite{Kerler:97,Kerler:03},
Kerler--Lyubashenko \cite{Kerler-Lyubashenko} and Yetter \cite{Yetter:97}
for descriptions of $\modC $.)  This functor $i$ maps $T\in \B(m,n)$ into
the cobordism $E_T$ defined above.  In the following, we regard
$\B$ as a braided subcategory of $\modC $ via~$i$.
Recall from \cite{Crane-Yetter:99} and \cite{Kerler:99} that $\modC $ is a braided
category, and there is a Hopf algebra
$\modh =(1,\mu _{\modh} ,\eta _{\modh} ,\Delta _{\modh} ,\epsilon _{\modh}
,S_{\modh} )$ in $\modC $ with the underlying
object $1$, which corresponds to the punctured torus $F_1$.

The category $\B$ can be identified with the subcategory of $\modC $ such that
$\Ob(\B)=\Ob(\modC )$ and
\begin{equation*}
  \B(m,n) =
  \bigl\{f\in \modC (m,n)\ver \epsilon _{\modh} ^{\otimes n}f=\epsilon
_{\modh} ^{\otimes m}\bigr\}
\end{equation*}
for $m,n\ge 0$.  Recall from \cite{Kerler:03} that $\modC $ is generated as
a braided category by the generators of $\B$ listed in \eqref{e16} and
an integral $\chi _{\modh} $ of the Hopf algebra $\modh $.  This integral
$\chi _{\modh} $
for $\modh $ is not contained in $\B$.

For each ribbon Hopf algebra $H$, we can define a braided functor
\begin{equation*}
  \modJ ^{\B}\col \B\rightarrow \Mod_H
\end{equation*}
such that $\modJ =\modJ ^{\B} \xi $.  The functor $\modJ ^{\B}$ maps the Hopf algebra
$H_{\B}$ in $\B$ into the transmutation $\bH$ of $H$.  If $H$ is
finite-dimensional ribbon Hopf algebra over a field $\modk $, and is {\em
factorizable} (see Reshetikhin and Semenov-Tian-Shansky
\cite{Reshetikhin-SemenovTianShanskii:89}) -- that is, the
function $\Hom_{\modk} (H,\modk )\rightarrow H$, $f\mapsto (1\otimes f)(c^H_+)$, is an
isomorphism -- then $\modJ ^{\B}$ extends to Kerler's functorial version of
the Hennings invariant \cite{Kerler:97,Kerler:03-2}
\begin{equation}
  \label{e1}
  \modJ ^{\tilde{\modC }}\col \tilde{\modC }\rightarrow \Mod_H.
\end{equation}
An interesting extension of $\B$ is the braided subcategory $\bar{\B}$
of $\modC $ generated by the objects and morphisms of $\B$ and the
morphism $B^*\in \modC (3,0)$ described in \fullref{fig:dualB}.
\labellist\small
\pinlabel {$C$} [tr] at 100 50
\endlabellist
\FIG{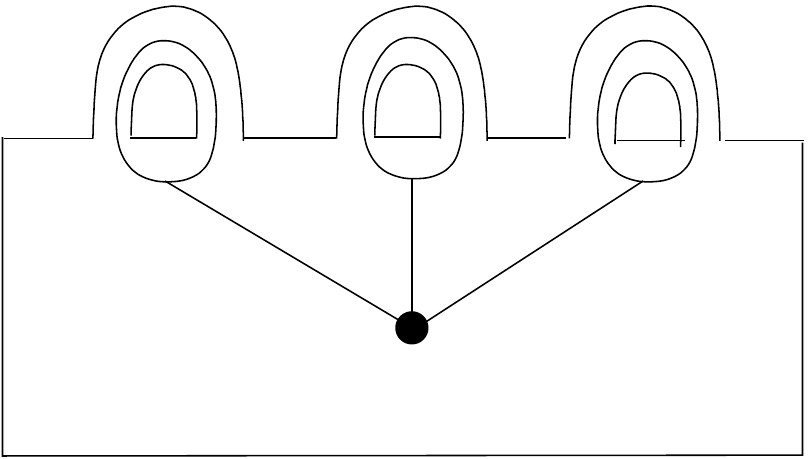}{The cobordism $B^*$ is obtained from the trivial cobordism
$\epsilon _{\B}^{\otimes 3}$ by surgery along the $Y$--graph $C$.}{height=25mm} One
can show that the category $\bar{\B}$ is the same as the category of
{\em bottom tangles in homology handlebodies}, which is defined in the
same way as $\B$ but the bottom tangles are contained in a homology
handlebody.  (Recall that a homology handlebody can be characterized
as a $3$--manifold which is obtainable as the result from a standard
handlebody of surgery along finitely many $Y$--graphs, see
Habegger \cite{Habegger}.)  For each $n\ge 0$, the monoid $\bar{\B}(n,n)$
contains the Lagrangian submonoid $\modL _n$ (see Levine \cite{Levine:01}) of the
monoid of homology cobordisms
(see Goussarov \cite{Goussarov:Y-graphs} and Habiro
\cite{Habiro:claspers}) of a compact, connected,
oriented surface $\Sigma _{n,1}$ of genus $1$ and with one boundary
component, and hence contains the monoid of homology cylinders (or
homologically trivial cobordisms) $\mathcal{H}\mathcal{C}_{n,1}$ over
$\Sigma _{n,1}$ and the Torelli group $\modI _{n,1}$ of $\Sigma _{n,1}$.  Here the
Lagrangian subgroup of $H_1\Sigma _{n,1}\simeq \modZ ^{2n}$, which is necessary
to specify $\modL _n$, is generated by the meridians of the handles in
$V_n$.  The monoid $\bar{\B}(n,n)$ does not contain the whole mapping
class group $\modM _{n,1}$.  In fact, $\bar{\B}(n,n)\cap \modM _{n,1}$ is
precisely the Lagrangian subgroup of $\modM _{n,1}$.  The subgroup
$\B(n,n)\cap \modM _{n,1}$ of $\bar{\B}(n,n)\cap \modM _{n,1}$ corresponds to the
handlebody group $H_{n,1}$, which is the group of isotopy classes of
self-homeomorphisms of handlebody of genus $n$ fixing a disc in the
boundary pointwise.

\subsection{Surgery on $3$--manifolds as monoidal relation}
\label{sec:surg-3-manif}

The idea of identifying the relations on tangles defined by local
moves with monoidal relations in the monoidal category $\modT $ of tangles
can be generalized to $3$--manifolds as follows.  Matveev
\cite{Matveev} defined a class of surgery operations on $3$--manifolds
called {\em $\V$--surgeries}.  A special case of $\V$--surgery removes
a handlebody from a $3$--manifold and reglues it back in a different
way.  Here $\V=(V_1,V_2)$ is a pair of two handlebodies $V_1,V_2$ of
the same genus with boundaries identified, and determines a type of
surgery.  We call such surgery {\em admissible}.  For each such $\V$,
there is a (not unique) pair $(f_1,f_2)$ of morphisms $f_1$ and $f_2$
of the same source and target in the monoidal category $\modC $ of
cobordisms of surfaces with connected boundary
(see Crane--Yetter \cite{Crane-Yetter:99}, Kerler \cite{Kerler:99} and
\fullref{sec:category-0}).  For $3$--manifolds representing morphisms in
$\modC $, the equivalence relations of $3$--manifolds generated by
$\V$--surgeries is the same as the monoidal relation in $\modC $ generated
by the pair $(f_1,f_2)$.  Thus, one can formulate the theory of
admissible surgeries in an algebraic way.

\bibliographystyle{gtart}
\bibliography{link}

\end{document}